\documentclass[twoside]{book}
\usepackage{lectures-ru}

\newcommand{\arxiv}[2]{#1} 

\newcommand{\spell}[2]{#1} 

\arxiv{\geometry{top=0.9in, bottom=0.9in,inner=0.9in, outer=0.7in, paperwidth=6in, paperheight=9in}}{\geometry{top=1.025in, bottom=1.025in,inner=0.9in, outer=0.825in, paperwidth=6.125in, paperheight=9.25in}}
\makeindex

\spell{}{\makeatletter\usepackage{environ}\RenewEnviron{figure}[1][]{}\RenewEnviron{figure*}[1][]{}\RenewEnviron{wrapfigure}[1][]{}\RenewEnviron{Figure}[1][]{}\pagestyle{empty}\let\ps@plain\ps@empty\makeatother\usepackage[none]{hyphenat}\def\emph{\textit}\renewcommand{\footnote}[1]{\textup{\ (#1)}}\renewcommand{\frac}[2]{(#1)/(#2)}\renewcommand{\tfrac}[2]{(#1/#2)}\everydisplay{\textstyle}}

\begin{document}

\hypersetup
{
pdftitle={Что такое дифференциальная геометрия?: кривые и поверхности},
pdfauthor={Антон Петрунин и  Серхио Замора Баррера}
}

\title{\tt Что такое\\ дифференциальная геометрия?:\\
кривые и поверхности}

\author{\tt Антон Петрунин и \tt  Серхио Замора Баррера\\
\\
\\ \tt Перевод и редакция Н. Д. Лебедевой и Г. Ю. Паниной}

\date{}

\maketitle

\thispagestyle{empty}

Художественное оформление Аны Кристины Чавес Калис.
\null
\vfill
\noindent{\includegraphics[scale=0.5]{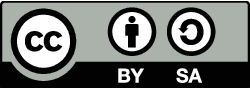}
\vspace*{1mm}
\\
\hbox{\parbox{.8\textwidth}
{Это произведение распространяется на условиях лицензии CC BY-SA 4.0, с ней можно ознакомиться по ссылке
\texttt{https://creativecommons.org/licenses/by-sa/4.0}}}}

\thispagestyle{empty}
\newpage
\addtocontents{toc}{\cftpagenumbersoff{part}}
{
\clearpage
\phantomsection
\pdfbookmark[1]{\contentsname}{bm:toc}
\sloppy
\footnotesize
\tableofcontents

}
\vfill

\begin{figure}[h!]
\centering
\begin{tikzpicture}[->,>=stealth',shorten >=1pt,auto,scale=.22,
  thick,main node/.style={circle,draw,font=\sffamily\bfseries,minimum size=8mm}]

  \node[main node] (1) at (0,0) {\ref{chap:curves-def}};
  \node[main node] (2) at (3,-5) {\ref{chap:length}};
  \node[main node] (3) at (6,-10) {\ref{chap:curve-curvature}};
  \node[main node] (4) at (9,-5) {\ref{chap:poly}};
  \node[main node] (5) at (9,-15) {\ref{chap:torsion}};
  \node[main node] (6) at (12,-10) {\ref{chap:signed-curvature}};
  \node[main node] (7) at (15,-15) {\ref{chap:supporting-curves}};
  \node[main node] (8) at (18,0) {\ref{chap:surfaces-def}};
  \node[main node] (9) at (15,-5) {\ref{chap:first-order}};
  \node[main node] (10) at (18,-10) {\ref{chap:surface-curvature}};
  \node[main node] (11) at (21,-15) {\ref{chap:Curves in a surface}};
  \node[main node] (12) at (18,-20) {\ref{chap:surface-support}};
  \node[main node] (13) at (21,-5) {\ref{chap:shortest}};
  \node[main node] (14) at (24,-10) {\ref{chap:geodesics}};
  \node[main node] (15) at (27,-15) {\ref{chap:parallel-transport}};
  \node[main node] (16) at (33,-15) {\ref{chap:gauss-bonnet}};
  \node[main node] (17) at (30,-10) {\ref{chap:semigeodesic}};
  \node[main node] (18) at (36,-10) {\ref{chap:comparison}};

  \path[every node/.style={font=\sffamily\small}]
   (1) edge node{}(2)
   (2) edge node{}(3)
   (3) edge node{}(5)
   (3) edge node{}(4)
   (3) edge node{}(6)
   (4) edge[dashed] node{}(6)
   (6) edge node{}(7)
   (6) edge node{}(10)
   (7) edge node{}(12)
   (8) edge node{}(9)
   (8) edge node{}(13)
   (9) edge node{}(10)
   (10) edge node{}(11)
   (11) edge node{}(12)
   (10) edge node{}(14)
   (13) edge node{}(14)
   (14) edge node{}(15)
   (14) edge[bend left= 15] node{}(17)
(17) edge[dashed, bend left=15] node{}(14)
   (15) edge node{}(16)
   (16) edge[dashed] node{}(18)
   (17) edge[dashed] node{}(15)
   (17) edge[dashed] node{}(16)
   (17) edge node{}(18);
\end{tikzpicture}
\end{figure}

\vfill

\newpage

\chapter*{Предисловие}
\addcontentsline{toc}{chapter}{Предисловие}
\thispagestyle{myheadings}
\markboth{ПРЕДИСЛОВИЕ}{ПРЕДИСЛОВИЕ}

Этот учебник рассчитан на тех, кто решил заниматься дифференциальной геометрией или же хочет найти вескую причину этого не делать.
Материала хватит на один семестр, и ещё останется.

Дифференциальная геометрия опирается на несколько разделов математики, включая
вещественный анализ,
теорию меры,
вариационное исчисление,
дифференциальные уравнения, топологию, элементарную и выпуклую геометрию.
Кроме того, физическая интуиция помогает разобраться во многих её аспектах.
В эту науку уйма входов, поэтому её и интересно, и трудно и преподавать, и изучать.

Гладкие кривые и поверхности дают важнейший источник примеров и идей дифференциальной геометрии.
Разумно хорошо разобраться в этой области, прежде чем идти дальше --- не стоит спешить.

В книжке делается упор на задачи,
доказательства элементарны, наглядны и почти строги (иногда пропускается кое-что из других разделов, в основном тех, что обсуждаются в приложении).
Мы сосредоточились на нескольких идеях, которые точно пригодятся в дальнейшем.
Поэтому обошли вниманием ряд тем, традиционно включаемых в вводные тексты;
например, мы почти не касаемся минимальных поверхностей и формул Петерсона --- Кодацци.

В то же время включены теоремы, которые обычно не обсуждаются в вводных курсах.

{\sloppy

Первый пример --- теорема о луне в луже Владимира Ионина и Германа Пестова (\ref{thm:moon-orginal}).
Это простейший значимый пример так называемых теорем от локального к глобальному, которые лежат в основе всей дифференциальной геометрии;
он даёт хороший ответ на главный вопрос книги --- «Что такое дифференциальная геометрия?».
Другие примеры включают теорему о седловых графиках Сергея Бернштейна (\ref{thm:bernshtein}) и теорему о бесконечной двусторонней геодезической Стефана Кон-Фоссена (\ref{thm:cohn-vossen}).

}

Учебник основан на наших лекциях, прочитанных осенью 2018 года на MASS-программе Университета штата Пенсильвания.
Многие из этих тем использовались Юрием Бураго в его лекциях, читаемых в Ленинградском университете, когда первый автор был его студентом.
При написании мы подглядывали в учебники
Вильгельма Бляшке~\cite{blaschke},
Виктора Топоногова~\cite{toponogov-book}
и Алексея Чернавского \cite{chernavsky}, а также в лекции Сергея Иванова \cite{ivanov};
многие продвинутые упражнения взяты из \cite{petrunin2020}.
Последняя глава основана на вводном материале из книги Стефани Александер, Виталия Каповича и первого автора~\cite{alexander-kapovitch-petrunin2027}.
Мы хотим поблагодарить
Стефани Александер,
Юрия Бураго,
Берка Джейлана,
Нину Лебедеву,
Александра Лычака,
Бенджамина Маккея
и студентов нашего класса
за помощь.

Работа частично поддержана грантом NSF DMS-2005279 и грантом Фонда Саймонса  № 584781.

\begin{flushright}
Антон Петрунин и
\\
Серхио Замора Баррера.
\end{flushright}

\arxiv{\cleardoublepage
\phantomsection
\AddToShipoutPictureBG*{\includegraphics[width=\paperwidth]{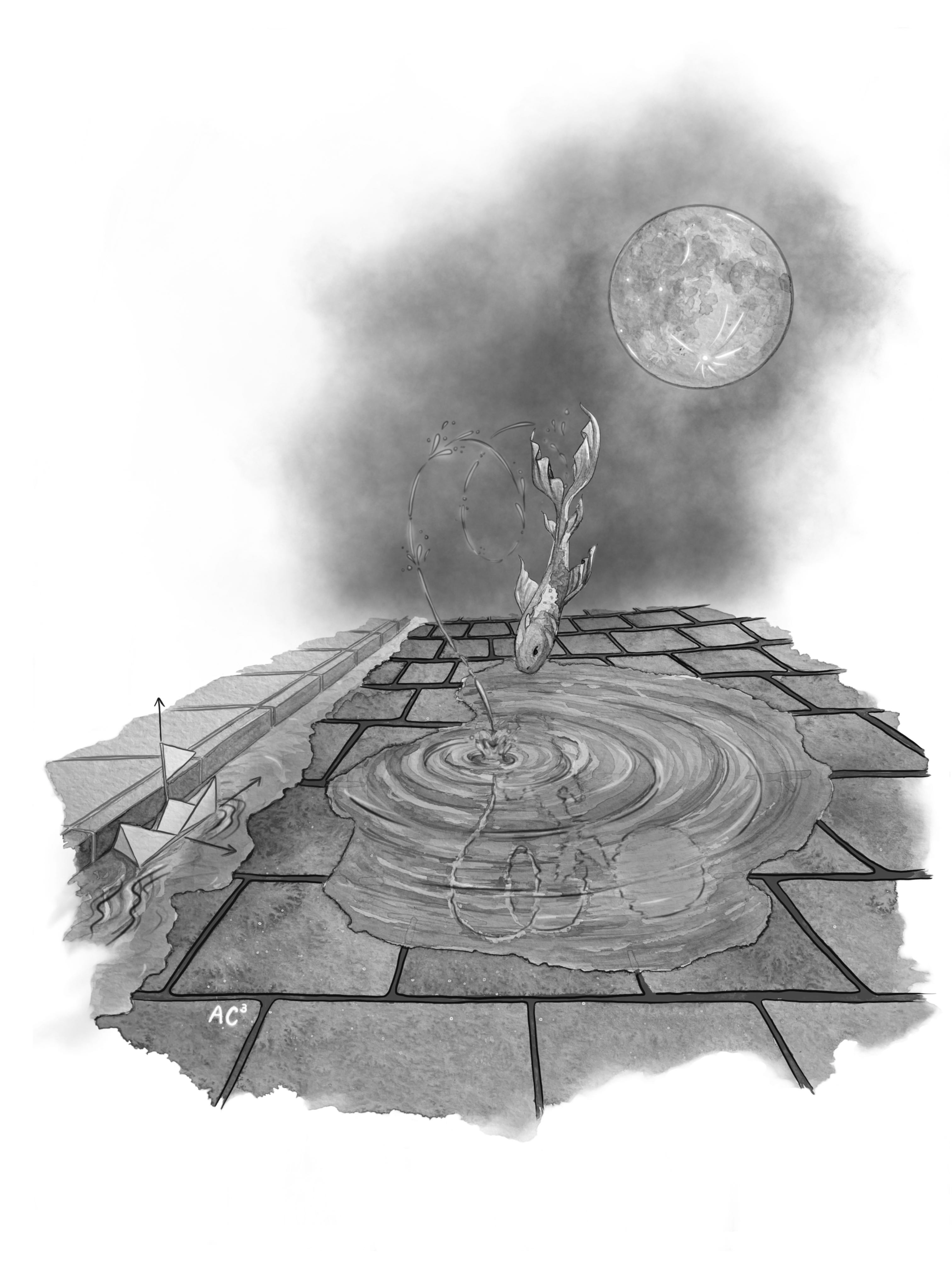}}
\cleardoublepage
\thispagestyle{empty}
\stepcounter{part}
\begin{center}
{\Huge\textbf{Часть \Roman{part}\quad Кривые}}\
\end{center}
\addcontentsline{toc}{part}{\texorpdfstring{Часть \Roman{part}\quad  Кривые}{Часть \Roman{part} Кривые}}
\clearpage}
{\backgroundsetup{
scale=.95,
opacity=1,
angle=0,
contents={
\includegraphics[width=\paperwidth
]{pics/Curves}
}%
}
\cleardoublepage
\phantomsection
\thispagestyle{empty}
\BgThispage
\stepcounter{part}
\addcontentsline{toc}{part}{\texorpdfstring{Часть \Roman{part}\quad  Кривые}{Часть \Roman{part} Кривые}}
\begin{center}
{\Huge\textbf{Часть \Roman{part}\qquad Кривые}}\
\end{center}
}

\chapter{Определения}
\label{chap:curves-def}

По следам велосипеда часто можно определить направление его движения, а также отличить следы переднего колеса от заднего, и найти расстояние между ними.

Подумайте, как всё это можно сделать.
Это заставит вас переоткрыть ощутимую часть дифференциальной геометрии кривых и вам станет легче читать первую часть;
см. также ссылку сразу после упражнения~\ref{ex:bike}.

\section{Прежде чем начать}

Понятие кривой имеет много вариаций.
Некоторые описываются существительными (путь, дуга и так далее), 
другие --- прилагательными (замкнутая, открытая, собственная, простая, гладкая и так далее).
Следующий рисунок даёт представление о четырёх типах кривых.

\vskip-0mm
\begin{figure}[h!]
\begin{minipage}{.48\textwidth}
\centering
\includegraphics{mppics/pic-110}
\end{minipage}\hfill
\begin{minipage}{.48\textwidth}
\centering
\includegraphics{mppics/pic-115}
\end{minipage}
\bigskip
\begin{minipage}{.48\textwidth}
\centering
\includegraphics{mppics/pic-120}
\end{minipage}\hfill
\begin{minipage}{.48\textwidth}
\centering
\includegraphics{mppics/pic-125}
\end{minipage}
\end{figure}
\vskip-0mm

Далее мы обсудим все эти определения.
Эту главу можно пропустить и пользоваться ею в дальнейшем как справочником.

\section{Простые кривые}

Мы полагаем, что читатель знаком с метрическими пространствами (см. приложение~\ref{sec:metric-spcaes});
главными примерами которых будут евклидова плоскость $\mathbb{R}^2$ и евклидово пространство $\mathbb{R}^3$.

Напомним, что \index{интервал}\emph{интервал} --- это связное подмножество вещественной прямой.
Непрерывная биекция $f\:X\to Y$ между подмножествами метрических пространств называется {}\emph{гомеоморфизмом}, если обратное отображение $f^{-1}\:Y\to X$ также непрерывно.  

\begin{thm}{Определение} 
Связное подмножество $\gamma$ метрического пространства называется \index{простая кривая}\emph{простой кривой}, если оно \index{локально гомеоморфно}\emph{локально} гомеоморфно интервалу;
то есть любая точка $p\in \gamma$ имеет окрестность в $\gamma$, гомеоморфную интервалу.
\end{thm}

Известно, что любую простую кривую $\gamma$ можно \index{параметризация}\emph{параметризовать} интервалом или окружностью.
То есть существует гомеоморфизм $\GG\to\gamma$, где $\GG$ --- интервал (открытый, замкнутый или полуоткрытый) или же окружность
\[\mathbb{S}^1=\set{(x,y)\in\mathbb{R}^2}{x^2+y^2=1}.\] 
Это доказано в статье Дэвида Гейла \cite{gale};
доказательство несложное, но оно отвлекло бы нас от основной темы.

Параметризация $\GG\to\gamma$ полностью описывает кривую.
Часто кривая и её параметризация обозначаются той же буквой;
например, можно сказать, что кривая $\gamma$ задана параметризацией $\gamma\:(a,b]\to \mathbb{R}^2$.
Однако любая простая кривая допускает различные параметризации.

\begin{thm}{Упражнение}\label{ex:9}
{\sloppy
\begin{subthm}{ex:9:compact}
Покажите, что образ любого непрерывного инъективного отображения $\gamma\:[0,1]\to\mathbb{R}^2$ является простой кривой.
\end{subthm}

}

\begin{subthm}{ex:9:9}
Постройте непрерывное инъективное отображение $\gamma\:(0,1)\z\to\mathbb{R}^2$, образ которого \textit{не} является простой кривой.
\end{subthm}

\end{thm}

\section{Параметризованные кривые}\label{sec:Parametrized curves}

Пусть $\GG$ --- окружность или интервал (открытый, замкнутый или полуоткрытый), и $\spc{X}$ --- метрическое пространство.
\index{параметризованная кривая}\emph{Параметризованная кривая} определяется как непрерывное отображение $\gamma\:\GG\z\to\spc{X}$.
Для параметризованной кривой мы \textit{не} предполагаем, что отображение инъективно; другими словами, допускаются {}\emph{самопересечения}.

Про простую кривую можно всегда думать как про параметризованную.
При этом термин \index{кривая}\emph{кривая} можно использовать, когда мы не хотим уточнять, параметризованная она или простая.

Если область определения $\GG$ является открытым интервалом или окружностью, то $\gamma$ называется {}\emph{кривой без концов};
в противном случае --- {}\emph{кривой с концами}.
В случае, если $\GG$ --- окружность, $\gamma$ называется {}\index{замкнутая!кривая}\emph{замкнутой кривой}. 
Если $\GG$ --- замкнутый интервал $[a,b]$, то кривая называется \index{дуга}\emph{дугой}.
Если вдобавок $\GG$ --- единичный интервал $[0,1]$, то кривая также называется \index{путь}\emph{путём}.

\begin{wrapfigure}{o}{15 mm}
\vskip-0mm
\centering
\includegraphics{mppics/pic-130}
\vskip-4mm
\end{wrapfigure}

Если для дуги $\gamma \: [a,b] \to \spc{X}$ выполнено условие $p\z=\gamma (a)=\gamma (b)$, то она называется \index{петля}\emph{петлёй}, и $p$ называется её \index{базовая точка}\emph{базовой точкой}.

Пусть $\GG_1$ и $\GG_2$, оба либо интервалы, либо окружности.
Непрерывное сюръективное отображение $\tau\:\GG_1\to\GG_2$ называется \index{монотонное отображение}\emph{монотонным}, если для любого $t\in \GG_2$ множество $\tau^{-1}\{t\}$ связно.
Если $\GG_1$ и $\GG_2$ --- интервалы, то, по теореме о промежуточных значениях, монотонное отображение либо неубывающее, либо невозрастающее;
то есть наше определение совпадает со стандартным, в случае если $\GG_1$ и $\GG_2$ --- интервалы.

\begin{thm}{Упражнение}\label{ex:mono}
Приведите пример монотонного (в частности сюръективного) отображения $(0,1)\to [0,1]$.
\end{thm}

Пусть $\gamma_1\:\GG_1\to \spc{X}$ и $\gamma_2\:\GG_2\to \spc{X}$ --- две параметризованные кривые, такие, что
$\gamma_1=\gamma_2\circ\tau$ для монотонного отображения $\tau\:\GG_1\to\GG_2$.
Тогда мы говорим, что $\gamma_2$ --- \index{репараметризация}\emph{репараметризация}%
\footnote{При этом $\gamma_1$ может \textit{не} быть репараметризацией $\gamma_2$.
Другими словами, с нашим определением, отношение \textit{быть репараметризацией} не является симметричным;
в частности, это \textit{не} отношение эквивалентности.
Дело можно поправить, перейдя к минимальному отношению эквивалентности, как это сделано в \cite[2.5.1]{burago-burago-ivanov},
но мы останемся верны нашему варианту определения.}
$\gamma_1$ с помощью $\tau$.

\begin{thm}{Продвинутое упражнение}\label{aex:simple-curve}
Пусть $X$ --- подмножество плоскости.
Предположим, что две различные точки $p,q\in X$ можно соединить путём в~$X$.
Покажите, что существует простая дуга в~$X$, из $p$ в~$q$.
\end{thm}

Любую петлю (как и замкнутую кривую) можно задать {}\emph{периодической} параметризацией $\gamma\: \mathbb{R}\to \spc{X}$;
то есть такой, что $\gamma(t+\ell)=\gamma(t)$ для некоторого периода $\ell>0$ и всех~$t$.
Например, единичная окружность на плоскости описывается $2{\cdot}\pi$-периодической параметризацией $\gamma(t)\z=(\cos t,\sin t)$.

Верно и обратное: про кривую с периодической параметризацией можно думать как про замкнутую кривую или как про петлю.

\section{Гладкие кривые}\label{sec:Smooth curves}

Кривые в евклидовом пространстве или на плоскости называются соответственно {}\emph{пространственными} кривыми или {}\emph{плоскими}.

Параметризованную пространственную кривую можно описать её координатными функциями 
$\gamma(t)=(x(t),y(t),z(t))$.
Плоские кривые можно рассматривать как частный случай пространственных кривых с $z(t)\equiv 0$.

Напомним, что вещественная функция называется \index{гладкая!функция}\emph{гладкой}, если её производные всех порядков определены всюду в области определения.  
Если каждая из координатных функций $x(t), y(t)$ и $z(t)$ гладкая, то и параметризация называется \index{гладкая!параметризация}\emph{гладкой}.

Если \index{вектор скорости}\emph{вектор скорости} 
$\gamma'(t)=(x'(t),y'(t),z'(t))$
не равен нулю ни в одной точке, то параметризация $\gamma$ называется \index{регулярная!параметризация}\emph{регулярной}.

Параметризованная кривая называется {}\emph{гладкой}, если её параметризация гладкая и регулярная.
Простая пространственная кривая называется \index{гладкая!кривая}\emph{гладкой}, если она допускает регулярную гладкую параметризацию;
для замкнутой кривой предполагается, что параметризация периодическая.
Это главные объекты первой части книги.
Гладкие кривые можно было бы называть {}\emph{регулярными гладкими кривыми};
получилось бы точнее и длиннее.

Гладкая петля может определять негладкую замкнутую кривую;
пример показан на рисунке.

\begin{wrapfigure}{o}{17 mm}
\vskip-4mm
\centering
\includegraphics{mppics/pic-51}
\bigskip
\includegraphics{mppics/pic-140}
\vskip-8mm
\end{wrapfigure}

Согласно следующему упражнению, кривые с гладкими параметризациями могут оказаться негладкими.

\begin{thm}{Упражнение}\label{ex:L-shape}
Из приложения \ref{sec:analysis} видно, что следующая функция является гладкой:
\[f(t)=
\begin{cases}
0&\text{если}\ t\le 0,
\\
\frac{t}{e^{1\!/\!t}}&\text{если}\ t> 0.
\end{cases}
\]

Покажите, что $\alpha(t)=(f(t),f(-t))$ описывает гладкую параметризацию кривой на рисунке;
это простая кривая, образованная объединением двух полуосей на плоскости.

{\sloppy

Покажите, что любая гладкая параметризация этой кривой имеет нулевой вектор скорости в начале координат.
Выведите отсюда, что эта кривая негладкая;
то есть не допускает регулярной гладкой параметризации.

}

\end{thm}

\begin{thm}{Упражнение}\label{ex:cycloid}
Опишите множество вещественных чисел $\ell$, 
при которых параметризация $\gamma_\ell (t)= (t+\ell \cdot \sin t,\ell \cdot \cos t)$, $t\in\mathbb{R}$ является

\begin{minipage}{.30\textwidth}
\begin{subthm}{ex:cycloid:smooth}
гладкой; 
\end{subthm}
\end{minipage}
\hfill
\begin{minipage}{.30\textwidth}
\begin{subthm}{ex:cycloid:regular}
регулярной;
\end{subthm}
\end{minipage}
\hfill
\begin{minipage}{.30\textwidth}
\begin{subthm}{ex:cycloid:simple}
простой.
\end{subthm}
\end{minipage}

\end{thm}

\begin{thm}{Упражнение}\label{ex:nonregular}
Постройте гладкую, но \textit{не}регулярную параметризацию кубической параболы $y=x^3$ на плоскости.
\end{thm}

\section{Неявно заданные кривые}\label{sec:implicit-curves}

Предположим, что $f\:\mathbb{R}^2\to \mathbb{R}$ --- гладкая функция; 
то есть все её частные производные везде определены.
Пусть $\gamma\subset \mathbb{R}^2$ --- её множество уровня, описанное уравнением $f(x,y)=0$.

Предположим, что $0$ является \index{регулярное значение}\emph{регулярным значением}~$f$; то есть градиент $\nabla_p f$ не обращается в ноль ни в одной точке $p\in \gamma$.
Другими словами, если $f(p)=0$, то   
$f_x(p)\ne 0$ или $f_y(p)\ne 0$.%
\footnote{Здесь $f_x$ является сокращённой записью для частной производной
$\tfrac{\partial f}{\partial x}$.\index{10f@$f_x$ (частная производная)}}
Если при этом множество $\gamma$ связно, то по теореме о неявной функции (\ref{thm:imlicit}), $\gamma$ --- гладкая простая кривая.

Описанное условие достаточно, но \textit{не необходимо}.
Например, ноль \textit{не} является регулярным значением функции $f(x,y)\z=y^2$, однако уравнение $f(x,y)=0$ описывает гладкую кривую --- ось $x$.

Аналогично, предположим, что $(f,h)$ --- пара гладких функций на $\mathbb{R}^3$.
Опять же, из теоремы о неявной функции (\ref{thm:imlicit}) следует, что система уравнений
\[\begin{cases}
   f(x,y,z)=0,
   \\
   h(x,y,z)=0
  \end{cases}
\]
описывает гладкую пространственную кривую, если множество решений $\gamma$ связно, и $0$ является регулярным значением отображения $F\:\mathbb{R}^3\to\mathbb{R}^2$, определённого как
\[F\:(x,y,z)\mapsto (f(x,y,z),h(x,y,z)).\]
Это означает, что градиенты $\nabla f$ и $\nabla h$ линейно независимы в любой точке $p\in \gamma$.
Другими словами, якобиан
\[\Jac_pF=
\begin{pmatrix}
f_x&f_y&f_z\\
h_x&h_y&h_z
\end{pmatrix}
\]
отображения $F\:\mathbb{R}^3\to\mathbb{R}^2$ имеет ранг $2$ в любой точке $p \in \gamma$.

Если кривая $\gamma$ описана таким образом,
то мы говорим, что она \index{неявно заданная кривая}\emph{задана неявно}.

Теорема об обратной функции гарантирует существование регулярных гладких параметризаций для любой неявно заданной кривой.
Однако часто удобнее работать непосредственно с неявным представлением.

\begin{thm}{Упражнение}\label{ex:y^2=x^3}
Рассмотрим множество на плоскости, описанное уравнением
$y^2=x^3$.
Является ли оно простой кривой?
Является ли оно гладкой кривой?
\end{thm}

\begin{thm}{Упражнение}\label{ex:viviani}
Опишите множество вещественных чисел $\ell$, 
при которых система уравнений
\[\begin{cases}
x^2+y^2+z^2&=1
\\
x^2+\ell\cdot x+y^2&=0
\end{cases}\]
описывает гладкую кривую.
\end{thm}

\section{Собственные, замкнутые и открытые}\label{sec:proper-curves}

{\sloppy

Параметризованная кривая $\gamma$ в метрическом пространстве $\spc{X}$ называется \index{собственная!кривая}\emph{собственной}, если для любого компактного множества $K \z\subset \spc{X}$ его прообраз $\gamma^{-1}(K)$ компактен.

}

Вот пример несобственной кривой, определённой на всей вещественной прямой: $\gamma(t)=(e^t,0,0)$.
Действительно, множество $(-\infty,0]$ --- это прообраз замкнутого единичного шара с центром в начале координат, и $(-\infty,0]$ некомпактно.

\begin{thm}{Упражнение}\label{ex:open-curve}
Покажите, что кривая $\gamma\:\mathbb{R}\to\mathbb{R}^3$ собственная тогда и только тогда, когда $|\gamma(t)|\z\to\infty$ при $t\to\pm\infty$.
\end{thm}

Напомним, что замкнутый интервал компактен, и замкнутые подмножества компактного множества также компактны
(приложение \ref{sec:topology}).
Поскольку окружности и замкнутые конечные интервалы компактны, замкнутые кривые и дуги являются собственными кривыми.

Простая кривая называется собственной, если она допускает собственную параметризацию.

\begin{thm}{Упражнение}\label{ex:proper-closed}
Покажите, что простая пространственная кривая является собственной тогда и только тогда, когда она образует замкнутое множество.
\end{thm}

Собственная простая кривая называется \index{открытая!кривая}\emph{открытой}, если она не замкнута и не имеет концов.
Таким образом, любая простая собственная кривая без концов является либо замкнутой, либо открытой.
Обратите внимание, что термины \textit{открытая кривая} и \textit{замкнутая кривая} не имеют отношения к открытым и замкнутым множествам.

\begin{thm}{Упражнение}\label{ex:proper-curve}
Используя теорему Жордана (\ref{thm:jordan}), покажите, что любая простая открытая кривая на плоскости делит её на две связные компоненты.
\end{thm}

\chapter{Длина}
\label{chap:length}

\section{Определения}

Напомним, что последовательность вида 
\[a=t_0 < t_1 < \cdots < t_k=b.\]
называется \index{разбиение}\emph{разбиением} интервала $[a,b]$.

\begin{thm}{Определение}\label{def:length}
Пусть $\gamma\:[a,b]\to \spc{X}$ --- кривая в метрическом пространстве.
\index{длина кривой}\emph{Длина} $\gamma$ определяется как
\begin{align*}
\length \gamma
&= 
\sup
\set{\dist{\gamma(t_0)}{\gamma(t_1)}{\spc{X}}
+\dots+
\dist{\gamma(t_{k-1})}{\gamma(t_k)}{\spc{X}}}{},
\end{align*}
где верхняя грань берётся по всем разбиениям $t_0,\dots,t_k$ интервала $[a,b]$.

Длина замкнутой кривой определяется как длина соответствующей петли.
Если кривая параметризована открытым или полуоткрытым интервалом, то её длина определяется как верхняя грань длин всех её сужений на замкнутые интервалы.
\end{thm}

Кривая конечной длины называется \index{спрямляемая кривая}\emph{спрямляемой}.

\begin{thm}{Упражнение}\label{ex:integral-length-0}
Предположим, что $\gamma_1\:[a_1,b_1] \to\mathbb{R}^3$ есть репараметризация кривой $\gamma_2\:[a_2,b_2] \to\mathbb{R}^3$. 
Покажите, что
\[\length \gamma_1 = \length \gamma_2.\]
\end{thm}

\begin{wrapfigure}[4]{r}{33 mm}
\vskip-4mm
\centering
\includegraphics{mppics/pic-224}
\end{wrapfigure}

Пусть $\gamma\:[a,b]\to \mathbb{R}^3$ --- параметризованная пространственная кривая.
Выберем разбиение $a\z=t_0 < t_1 < \cdots < t_k=b$, и пусть $p_i=\gamma(t_i)$.
Тогда ломаная $p_0\dots p_k$ называется \index{вписанная ломаная}\emph{вписанной} в~$\gamma$.
Если $\gamma$ замкнута, то $p_0=p_k$, и, значит, её вписанная ломаная также замкнута.

Отметим, что длину пространственной кривой можно определить как верхнюю грань длин вписанных в неё ломаных.

{\sloppy

\begin{thm}{Упражнение}\label{ex:length-chain}
Рассмотрим путь $\gamma\:[0,1]\to\mathbb{R}^3$.
Предположим, что $\beta_n$ --- последовательность вписанных в $\gamma$ ломаных с вершинами $\gamma(\tfrac in)$ для $i\z\in\{0,\dots,n\}$.
Докажите, что 
\[\length\beta_n\to\length\gamma
\quad\text{при}\quad
n\to \infty.
\]
\end{thm}

}

\begin{thm}{Упражнение}\label{ex:length-image}
Рассмотрим простой путь $\gamma\:[0,1]\to\mathbb{R}^3$.
Предположим, что путь $\beta\:[0,1]\to\mathbb{R}^3$ имеет тот же образ, что и $\gamma$;
то есть $\beta([0,1])=\gamma([0,1])$.
Докажите, что 
\[\length \beta\ge \length \gamma.\]

Попробуйте доказать это неравенство, ослабив предположение до $\beta([0,1])\z\supset\gamma([0,1])$.
\end{thm}

\begin{thm}{Упражнение}\label{ex:integral-length}
Пусть $\gamma\:[a,b]\to\mathbb{R}^3$ --- гладкая кривая.
Докажите, что
\vskip1mm
\begin{minipage}{.45\textwidth}
\begin{subthm}{ex:integral-length>}
$\length \gamma
\le
\int_a^b|\gamma'(t)|\cdot dt$;
\end{subthm}
\end{minipage}
\hfill
\begin{minipage}{.45\textwidth}
\begin{subthm}{ex:integral-length<}
$\length \gamma
\ge
\int_a^b|\gamma'(t)|\cdot dt$.
\end{subthm}
\end{minipage}
\vskip1mm

В частности 
\[\length \gamma
=
\int_a^b|\gamma'(t)|\cdot dt.\eqlbl{eq:length}\]

\end{thm}

\begin{thm}{Продвинутые упражнения}\label{adex:integral-length}

\begin{subthm}{adex:integral-length:a}
Докажите, что формула \ref{eq:length} справедлива для любой липшицевой кривой $\gamma\:[a,b]\z\to\mathbb{R}^3$.
\end{subthm}

\begin{subthm}{adex:integral-length:b}
Постройте такую непостоянную кривую $\gamma\:[a,b]\to\mathbb{R}^3$, что $\gamma'(t)=0$ почти всюду.
(Для такой кривой не выполняется равенство \ref{eq:length}, хотя обе его стороны определены.)
\end{subthm}

\end{thm}

\section{Неспрямляемые кривые}

Опишем так называемую \index{снежинка Коха}\emph{снежинку Коха} ---
классический пример неспрямляемой кривой.

Начнём с равностороннего треугольника.
Разделим каждую сторону на три равные отрезка, и добавим равносторонний треугольник с основанием на среднем отрезке.
К полученному многоугольнику применим ту же операцию и так будем продолжать рекурсивно.
Снежинка Коха --- это граница объединения всех полученных многоугольников.
Две итерации и сама снежинка показаны на рисунке.

\begin{figure}[ht!]
\centering
\includegraphics{mppics/pic-225}
\end{figure}

\begin{thm}{Упражнение}\label{ex:nonrectifiable-curve}

\begin{subthm}{ex:nonrectifiable-curve:a}
Докажите, что снежинка Коха --- простая замкнутая кривая; в частности её можно параметризовать окружностью.
\end{subthm}

\begin{subthm}{ex:nonrectifiable-curve:b}
Докажите, что снежинка Коха не спрямляема. 
\end{subthm}
\end{thm}

\section{Полунепрерывность длины}

Нижний предел последовательности $x_n$ будет обозначаться как
\[\liminf_{n\to\infty} x_n.\] 
Он определяется как наименьший частичный предел; то есть наименьший возможный предел подпоследовательности $x_n$.
Нижний предел определён для любой последовательности вещественных чисел и принимает значения в расширенной числовой прямой $[-\infty,\infty]$.

{\sloppy

\begin{thm}{Теорема}
Предположим, что последовательность кривых $\gamma_n\:[a,b]\to \spc{X}$ в метрическом пространстве $\spc{X}$ сходится поточечно к кривой $\gamma_\infty\:[a,b]\to \spc{X}$;
то есть $\gamma_n(t)\z\to\gamma_\infty(t)$
при любом $t \in [a,b]$ и $n\to\infty$. 
Тогда 
$$\liminf_{n\to\infty} \length\gamma_n \ge \length\gamma_\infty.\eqlbl{eq:semicont-length}$$
\end{thm}

}

\begin{thm}{Следствие}\label{thm:length-semicont}
Длина полунепрерывна снизу по отношению к поточечной сходимости кривых. 
\end{thm}

\parbf{Доказательство.}
Выберем разбиение $a=t_0<t_1<\dots<t_k=b$.
Пусть
\begin{align*}\Sigma_n
&\df
\dist{\gamma_n(t_0)}{\gamma_n(t_1)}{}
+\dots+
\dist{\gamma_n(t_{k-1})}{\gamma_n(t_k)}{},
\\
\Sigma_\infty
&\df
\dist{\gamma_\infty(t_0)}{\gamma_\infty(t_1)}{}
+\dots+
\dist{\gamma_\infty(t_{k-1})}{\gamma_\infty(t_k)}{}.
\end{align*}

Для каждого $i$,
\[\dist{\gamma_n(t_{i-1})}{\gamma_n(t_i)}{}
\to
\dist{\gamma_\infty(t_{i-1})}{\gamma_\infty(t_i)}{},\]
при $n\to\infty$, и, значит, $\Sigma_n\to \Sigma_\infty$.
Поскольку 
$\Sigma_n\le\length\gamma_n$
для каждого~$n$, получаем, что
$$\liminf_{n\to\infty} \length\gamma_n \ge \Sigma_\infty$$
для любого разбиения.
Теперь неравенство \ref{eq:semicont-length} следует из определения длины.
\qeds

\begin{wrapfigure}{o}{20 mm}
\vskip0mm
\centering
\includegraphics{mppics/pic-6}
\end{wrapfigure}

Неравенство \ref{eq:semicont-length} может оказаться строгим.
Например, диагональ $\gamma_\infty$ единичного квадрата 
можно аппроксимировать ступенчатыми ломаными $\gamma_n$,
стороны которых параллельны сторонам квадрата ($\gamma_6$ и $\gamma_\infty$ показаны на рисунке).
В этом случае,
\[\length\gamma_\infty=\sqrt{2}\quad
\text{и}\quad \length\gamma_n=2\quad
\text{для любого}\quad n.\]

\section{Параметризация длиной}

Будем говорить, что кривая $\gamma$ \index{параметризация!длиной}\emph{параметризована длиной} или имеет \index{естественная параметризация}\index{натуральная параметризация}\emph{естественную (натуральную) параметризацию},
если 
\[t_2-t_1=\length \gamma|_{[t_1,t_2]}\]
при любых двух значений параметра $t_1<t_2$;
то есть дуга $\gamma$ от $t_1$ до $t_2$ имеет длину $t_2-t_1$.

\begin{thm}{Упражнение}\label{ex:cont-length}
Пусть $\gamma\:[a,b]\to \spc{X}$ --- спрямляемая кривая в метрическом пространстве.
Для данного $t\in [a,b]$ обозначим через $s(t)$ длину дуги $\gamma|_{[a,t]}$.
Докажите, что функция $t\mapsto s(t)$ непрерывна.

Выведите отсюда, что $\gamma$ допускает параметризацию длиной.
\end{thm}

По упражнению~\ref{ex:integral-length},
гладкая кривая $\gamma(t)=(x(t),y(t),z(t))$ параметризована длиной тогда и только тогда, когда её вектор скорости единичный в любой момент;
то есть
\[|\gamma'(t)|=\sqrt{x'(t)^2+y'(t)^2+z'(t)^2}=1\]
для всех $t$.
Поэтому гладкие кривые, параметризованные длиной, называют также кривыми с \index{кривая с единичной скоростью}\emph{единичной скоростью}.
Отметим, что гладкие параметризации с единичной скоростью всегда регулярны (см.~\ref{sec:Smooth curves}).

\begin{thm}{Предложение}\label{prop:arc-length-smooth}
Если $t\mapsto \gamma(t)$ --- гладкая кривая,
то её параметризация длиной также гладкая и регулярная.
Более того, параметр дуги $s$ для $\gamma$ можно записать в виде интеграла
\[s(t)=\int_{t_0}^t |\gamma'(\tau)|\cdot d\tau.
\eqlbl{s(t)}\]

\end{thm}

Как правило, буквой $s$ будет обозначаться параметр длины.

\parbf{Доказательство.}
Так как $\gamma$ гладкая, $|\gamma'(t)|>0$ для любого~$t$,
а значит, и функция $t\mapsto|\gamma'(t)|$ гладкая.

Согласно формуле Ньютона --- Лейбница, $s'(t)=|\gamma'(t)|$.
В частности, $t\mapsto s(t)$ --- гладкая возрастающая функция с положительной производной.

По теореме об обратной функции (\ref{thm:inverse}), функция $s^{-1}(t)$ также гладкая
и $|(\gamma\circ s^{-1})'|\equiv1$.
Следовательно, репараметризация $\gamma\circ s^{-1}$ имеет единичную скорость.
По построению, $\gamma\circ s^{-1}$ гладкая, и так как $|(\gamma\circ s^{-1})'|\equiv1$, параметризация $\gamma\circ s^{-1}$ регулярна.
\qeds

\begin{thm}{Упражнение}\label{ex:arc-length-helix}
Параметризуйте \index{винтовая линия}\emph{винтовую линию} 
\[\gamma_{a,b}(t)=(a\cdot\cos t,a\cdot \sin t, b\cdot t)\]
длиной.
\end{thm}

Нас будут интересовать свойства кривых, которые не меняются при репараметризациях.
Поэтому всегда можно предполагать, что данная гладкая кривая параметризована длиной.
Выбор этих параметризаций почти каноничен --- они отличаются только знаком и сдвигом, и, значит, с ними проще определить величины, не зависящие от параметризации.
Это наблюдение будет использовано при определении кривизны и кручения.

На практике же проще работать с исходной параметризацией.
Более того, обычно параметр длины в явном виде найти невозможно.

\section{Выпуклые кривые}

Простая кривая на плоскости называется \index{выпуклая!кривая}\emph{выпуклой}, если она ограничивает выпуклую область.
Поскольку граница любой области --- замкнутое множество, любая выпуклая кривая либо замкнута, либо открыта (см. \ref{sec:proper-curves}).

\begin{thm}{Предложение}\label{prop:convex-curve}
Предположим, что замкнутая выпуклая кривая $\alpha$ лежит внутри области, ограниченной простой замкнутой кривой $\beta$ на плоскости.
Тогда
\[\length\alpha\le \length\beta.\]
\end{thm}

Достаточно показать, что периметр любого многоугольника, вписанного в $\alpha$, меньше или равен длине $\beta$.
Поскольку любой многоугольник, вписанный в $\alpha$, выпуклый, достаточно доказать следующее.

\begin{thm}{Лемма}\label{lem:perimeter}
Предположим, что выпуклый многоугольник $P$ лежит в фигуре $F$, ограниченной простой замкнутой кривой.
Тогда
\[\perim P\le \perim F,\]
где $\perim F$ обозначает периметр~$F$.
\end{thm}

\parbf{Доказательство.} 
\index{хорда}\emph{Хорда} фигуры $F$ определяется как отрезок прямой в $F$ с концами на её границе.
Предположим, что $F'$ --- это фигура, полученная из $F$ отрезанием куска по её хорде.
По неравенству треугольника,
\[\perim F'\le \perim F.\]

\begin{wrapfigure}{o}{24 mm}
\vskip-10mm
\centering
\includegraphics{mppics/pic-7}
\vskip3mm
\end{wrapfigure}

Заметим, что существует убывающая последовательность фигур 
\[F=F_0\supset F_1\supset\dots\supset F_n=P,\]
такая что $F_{i+1}$ получается из $F_{i}$ отрезанием куска по хорде.
Следовательно, 
\begin{align*}
\perim F=\perim F_0&\ge\perim F_1\ge\dots\ge\perim F_n=\perim P.
\end{align*}
\qedsf

\parit{Замечание.}
Другие доказательства получаются применением формулы Крофтона (\ref{ex:convex-croftons}) и проецированием $F$ на $P$ как в \ref{lem:nearest-point-projection}.

\begin{thm}{Следствие}\label{cor:convex=>rectifiable}
Любая выпуклая замкнутая кривая на плоскости спрямляема.
\end{thm}

\parbf{Доказательство.}
Любая замкнутая кривая ограничена.
Действительно, кривую можно задать некоторой петлёй $\alpha\:[0,1]\to\mathbb{R}^2$, $\alpha(t)\z=(x(t),y(t))$.
Координатные функции $t\mapsto x(t)$ и $t\mapsto y(t)$  непрерывны, и определены на $[0,1]$.
В частности, обе функции ограничены некоторой константой~$C$.
Следовательно, $\alpha$ лежит в квадрате, определённом неравенствами $|x|\le C$ и $|y|\le C$.

Согласно~\ref{prop:convex-curve}, длина $\alpha$ не превосходит периметр этого квадрата; отсюда результат.
\qeds

Напомним, что \index{выпуклая!оболочка}\emph{выпуклая оболочка} множества $X$ --- это наименьшее выпуклое множество, содержащее $X$; иначе говоря, это пересечение всех выпуклых множеств, содержащих~$X$.

\begin{thm}{Упражнение}\label{ex:convex-hull}
{\sloppy
Пусть $\alpha$ --- простая замкнутая кривая на плоскости.
Обозначим через $K$ выпуклую оболочку $\alpha$; пусть $\beta$ --- кривая, ограничивающая~$K$.
Докажите, что
\[\length \alpha\ge \length \beta.\]

}

Попробуйте доказать, что утверждение остаётся верным для любых замкнутых кривых на плоскости, предполагая только, что $K$ имеет непустую внутренность.
\end{thm}

\section{Формулы Крофтона}
\label{sec:crofton}
\index{формула Крофтона}

Для функции $f\: \mathbb{S}^1 \to \mathbb{R}$ обозначим её среднее значение как $\overline{f(\vec u)}$; то есть
\[\overline{f(\vec u)}=\frac1{2\cdot \pi}\cdot\int_{\vec u \in\mathbb{S}^1} f(\vec u).\]
Для вектора $\vec w$ и единичного вектора $\vec u$ обозначим через $\vec w_{\vec u}$ ортогональную проекцию $\vec w$ на прямую в направлении $\vec u$;
то есть
\[\vec w_{\vec u}=\langle\vec u,\vec w\rangle\cdot\vec u.\] 

\begin{thm}{Теорема}
Для любой  кривой $\gamma$ на плоскости выполняется равенство
\[
\length\gamma
=\tfrac\pi2\cdot \overline{\length\gamma_{\vec u}}, \eqlbl{crofton-formula}
\]
где $\gamma_{\vec u}$ определяется как $\gamma_{\vec u}(t) \df (\gamma (t))_{\vec u}$.
\end{thm}

\parbf{Доказательство.}
Длина вектора ${\vec w}$ пропорциональна средней длине его проекций; то есть
\[|{\vec w}|=k\cdot \overline{|{\vec w}_{\vec u}|}\]
для некоторого $k \in \mathbb{R}$.
(Коэффициент $k$ можно найти интегрированием\footnote{Это среднее значение $|\cos x|$ при $x\in [0,2\cdot\pi]$.}, но мы найдём его другим способом.)
Пусть $\gamma\:[a,b]\to\mathbb{R}^2$ --- гладкая кривая.
Для любого $t \in [a,b]$,
\[\gamma_{\vec u}'(t)=(\gamma'(t))_{\vec u}
\quad\text{и}\quad
|\gamma'_{\vec u}(t)|=|\langle\vec u,\gamma'(t)\rangle|.\]
Из упражнения~\ref{ex:integral-length}, получаем
\begin{align*}
\length\gamma
&=\int_a^b|\gamma'(t)|\cdot dt=
\\
&=\int_a^b  k\cdot \overline{|\gamma_{\vec u}'(t)|}\cdot dt=
\\
&=k\cdot \overline{\length\gamma_{\vec u}}.
\end{align*}

Поскольку $k$ --- константа, её достаточно вычислить, предположив что $\gamma$ --- единичная окружность.
В этом случае,
\[\length \gamma=2\cdot\pi.\]
Кривая $\gamma_{\vec u}$ пробегает отрезок длины~$2$ туда-сюда.
Значит $\length\gamma_{\vec u}=4$ для любого $\vec u$, и 
\[\overline{\length\gamma_{\vec u}} =4.\]
Отсюда $2\cdot \pi =k\cdot 4$.
Итак, мы доказали \ref{crofton-formula} для гладких кривых.

Применяя то же рассуждение вместе с \ref{adex:integral-length}, получаем, что \ref{crofton-formula} выполняется для всех липшицевых кривых.
Далее, \ref{ex:cont-length} влечёт \ref{crofton-formula} для произвольных спрямляемых кривых,
ибо параметризация длиной превращает спрямляемую кривую в липшицеву.

Остаётся рассмотреть случай неспрямляемых кривых;
достаточно показать, что 
\[\length\gamma=\infty
\quad\Longrightarrow\quad
\overline{\length\gamma_{\vec u}}=\infty.
\]

Из определения длины следует, что
\[\length\gamma_{\vec u}+\length\gamma_{\vec v}\ge \length\gamma\]
для любой пары $(\vec u , \vec v )$ ортонормированных векторов в $\mathbb{R}^2$.
Следовательно, если $\gamma$ имеет бесконечную длину, то средняя длина $\gamma_{\vec u}$ также бесконечна.
\qeds

\begin{thm}{Упражнение}\label{ex:convex-croftons}
Рассмотрим простую замкнутую кривую $\gamma$, ограничивающую фигуру~$F$ на плоскости.
Обозначим через $s$ среднюю длину проекций $F$ на прямые.
Покажите, что $\length\gamma\ge \pi \cdot s$.
Более того, равенство выполняется тогда и только тогда, когда $\gamma$ выпукла.

Решите \ref{ex:convex-hull}, используя данное утверждение.
\end{thm}

Следующее упражнение даёт аналогичные формулы для евклидова пространства.

Как и раньше, обозначим через $\vec w_{\vec u}$ ортогональную проекцию вектора $\vec w$ на прямую в направлении $\vec u$.
Далее, обозначим через $\vec w_{\vec u}^\bot$ проекцию $\vec w$ на плоскость, ортогональную к $\vec u$;
то есть
\[\vec w_\vec u^\bot=\vec w - \vec w_{\vec u}.\]
Далее будем пользоваться обозначением 
$\overline{f(\vec u)}$ для среднего значения
функции $f$, определённой на $\mathbb{S}^2$.

\begin{thm}{Продвинутое упражнение}\label{adex:more-croftons}
{\sloppy
Покажите, что длина пространственной кривой пропорциональна 

}
\begin{subthm}{}
средней длине её проекций на все прямые; то есть
\[\length\gamma=k_1\cdot\overline{\length\gamma_{\vec u}}\]
для некоторого $k_1 \in \mathbb{R}$.
\end{subthm}
\begin{subthm}{}средней длине её проекций на все плоскости; то есть
\[\length\gamma=k_2\cdot\overline{\length\gamma_{\vec u}^\bot}\]
для некоторого $k_2 \in \mathbb{R}$.
\end{subthm}
Найдите значения $k_1$ и $k_2$.
\end{thm}

\section{Внутренняя метрика}\label{sec:Length metric}

Пусть $\spc{X}$ --- метрическое пространство.
Для двух точек $x$ и $y$ в $\spc{X}$ обозначим через $\ell(x,y)$ точную нижнюю грань длин всех путей, соединяющих $x$ с $y$;
если же такого пути нет, считаем, что $\ell(x,y)=\infty$.

Легко видеть, что функция $\ell$ удовлетворяет всем аксиомам метрики, за исключением того, что она может равняться бесконечности.
Следовательно, если любые две точки в $\spc{X}$ можно соединить спрямляемой кривой, то $\ell$ определяет новую метрику на $\spc{X}$;
в этом случае $\ell$ называется \index{индуцированная внутренняя метрика}\emph{индуцированной внутренней метрикой} на $\spc{X}$.

Ясно, что $\ell(x,y)\ge \dist{x}{y}{}$ для любой пары точек $x,y\in \spc{X}$.
Если равенство выполняется для всех пар, то метрика $\dist{{*}}{{*}}{}$ называется \index{внутренняя метрика}\emph{внутренней}.

\begin{thm}{Упражнение}\label{ex:induced-is-length}
Пусть $(x,y)\mapsto \ell(x,y)$ --- индуцированная внутренняя метрика на метрическом пространстве $\spc{X}$.
Докажите, что $\ell$ --- внутренняя метрика.
\end{thm}

Большей частью мы рассматриваем пространства с внутренней метрикой.
Например, у евклидова пространства метрика внутренняя.

Метрика на подпространстве $A$ может быть не внутренней, даже если его объемлющее пространство $\spc{X}$ имеет внутреннюю метрику.
Расстояние в индуцированной внутренней метрике на подпространстве $A$ будет обозначаться $\dist{x}{y}A$;
то есть $\dist{x}{y}A$ --- это точная нижняя грань длин путей в $A$ от $x$ до $y$.\index{10aaa@$\lvert x-y\rvert_A$ (внутренняя метрика)}

\begin{thm}{Упражнение}\label{ex:intrinsic-convex}
Пусть $A\subset \mathbb{R}^3$ --- замкнутое подмножество.
Покажите, что $A$ выпукло тогда и только тогда, когда
\[\dist{x}{y}A=\dist{x}{y}{\mathbb{R}^3}\]
при любых $x,y\in A$.
\end{thm}

\section{Сферические кривые}

Обозначим через $\mathbb{S}^2$ единичную сферу в пространстве; то есть
\[\mathbb{S}^2=\set{(x,y,z)\in\mathbb{R}^3}{x^2+y^2+z^2=1}.\]
Пространственная кривая $\gamma$ называется \index{сферическая кривая}\emph{сферической}, если она лежит на $\mathbb{S}^2$;
то есть $|\gamma(t)|=1$ при любом $t$.

Напомним, что $\measuredangle(u,v)$ обозначает угол между векторами $u$ и~$v$.

\begin{thm}{Наблюдение}\label{obs:S2-length}
Для любых $u,v\in \mathbb{S}^2$ выполняется равенство
\[\dist{u}{v}{\mathbb{S}^2}=\measuredangle(u,v).\]

\end{thm}

\parbf{Доказательство.}
Пусть $\gamma$ --- короткая дуга большой окружности%
\footnote{Большая окружность --- это пересечение сферы с плоскостью, проходящей через её центр.}
от $u$ до $v$ в $\mathbb{S}^2$.
Тогда $\length\gamma=\measuredangle(u,v)$.
Следовательно,
\[\dist{u}{v}{\mathbb{S}^2}\le\measuredangle(u,v).\]

Остаётся доказать обратное неравенство.
Другими словами, нам нужно показать, что для любой ломаной $\beta=p_0\dots p_n$, вписанной в~$\gamma$, существует ломаная
$\beta_1=q_0\dots q_n$, вписанная в любую заданную сферическую кривую $\gamma_1$, соединяющую $u$ с $v$, такая, что 
\[\length\beta_1\ge \length \beta.\eqlbl{eq:length beta=<length beta}\]

Определим $q_i$ как первую точку на $\gamma_1$ для которой выполнено равенство $|u-p_i|=|u-q_i|$, но зададим $q_n=v$.
Ясно, что $\beta_1$ вписана в~$\gamma_1$.
По неравенству треугольника для углов (\ref{thm:spherical-triangle-inq}),
\begin{align*}
 \measuredangle(q_{i-1},q_i) &\ge  \measuredangle (u, q_i) - \measuredangle ( u , q_{i-1})  =
\\
&= \measuredangle (u,p_i) - \measuredangle (u,p_{i-1}) =
\\
& = \measuredangle(p_{i-1},p_i).
\end{align*}
По монотонности угла (\ref{lem:angle-monotonicity}),
\[|q_{i-1}-q_i|\ge|p_{i-1}-p_i|\]
и \ref{eq:length beta=<length beta} следует.
\qeds

\begin{thm}{Лемма о полусфере}\label{lem:hemisphere}
Любая замкнутая сферическая кривая длины меньше $2\cdot \pi$ лежит в открытой полусфере.
\end{thm}

Эта лемма проявится в доказательстве теоремы Фенхеля (\ref{thm:fenchel}).
Следующее доказательство принадлежит Стефани Александер.
Оно не так просто, как может показаться.
Прежде чем его читать, стоит попытаться придумать своё.

\parbf{Доказательство.}
Пусть $\gamma$ --- замкнутая кривая в $\mathbb{S}^2$ длины $2\cdot\ell$.
Допустим, что $\ell<\pi$.

{

\begin{wrapfigure}{o}{35mm}
\vskip-3mm
\centering
\includegraphics{mppics/pic-52}
\end{wrapfigure}

Разделим $\gamma$ на две дуги $\gamma_1$ и $\gamma_2$ длины~$\ell$ каждая;
пусть $p$ и $q$ --- их общие концы.
Из \ref{obs:S2-length}, получаем
\begin{align*}
\measuredangle(p,q)&\le \length \gamma_1=
\\
&= \ell<
\\
&<\pi.
\end{align*}

}

Обозначим через $z$ середину между $p$ и $q$ в $\mathbb{S}^2$;
то есть $z$ --- середина короткой дуги большой окружности от $p$ до $q$ в $\mathbb{S}^2$.
Мы утверждаем, что $\gamma$ лежит в открытой полусфере с полюсом в~$z$.
Если это не так, то $\gamma$ пересекает экватор в некоторой точке~$r$.
Не умаляя общности, можно считать, что $r$ принадлежит~$\gamma_1$.

Повернём дугу $\gamma_1$ на угол $\pi$ вокруг прямой, проходящей через $z$ и центр сферы.
Полученная дуга $\gamma_1^{*}$ вместе с $\gamma_1$ образует замкнутую кривую длины $2\cdot \ell$, проходящую через $r$ и её антиподальную точку $r^{*}$.
Снова применив \ref{obs:S2-length}, получим
\[\tfrac12\cdot\length \gamma=\ell\ge \measuredangle(r,r^{*})=\pi\]
--- противоречие.
\qeds

\begin{thm}{Упражнение}\label{ex:antipodal}
Опишите простую замкнутую сферическую кривую, которая не содержит пары антиподальных точек и при этом не лежит ни в какой полусфере.
\end{thm}

\begin{thm}{Упражнение}\label{ex:bisection-of-S2}
Предположим, что простая замкнутая сферическая кривая $\gamma$ делит $\mathbb{S}^2$ на две области равной площади.
Покажите, что 
\[\length\gamma\ge2\cdot\pi.\]
\end{thm}

\begin{thm}{Упражнение}\label{ex:flaw}
Найдите ошибку в решении следующей задачи.
Предложите правильное решение.
\end{thm}

\parbf{Задача.}
Предположим, что замкнутая кривая на плоскости имеет длину не более $4$.
Покажите, что она лежит в единичном круге.

\parbf{Неправильное решение.}
Достаточно показать, что \index{диаметр}\emph{диаметр} нашей кривой, скажем $\gamma$, не превышает $2$;
то есть
\[|p-q|\le 2\eqlbl{eq:|pq|=<2}\]
при любых двух точек $p$ и $q$ на~$\gamma$.

Длина $\gamma$ не меньше длины замкнутой вписанной ломаной, идущей от $p$ до $q$ и обратно до~$p$.
Следовательно,
\[2\cdot |p-q|\le\length \gamma\le 4;\]
откуда вытекает \ref{eq:|pq|=<2}.
\qedsf

\begin{thm}{Продвинутые упражнения} \label{adex:crofton}
Пусть ${\vec u},{\vec w}\in\mathbb{S}^2$.
Обозначим через ${\vec w}^*_{\vec u}$ ближайшую к ${\vec w}$ точку на экваторе с полюсом в ${\vec u}$;
другими словами, если ${\vec w}^\perp_{\vec u}$ --- проекция ${\vec w}$ на плоскость, перпендикулярную ${\vec u}$, то ${\vec w}^*_{\vec u}$ --- единичный вектор в направлении ${\vec w}^\perp_{\vec u}$.
Вектор ${\vec w}^*_{\vec u}$ определён при ${\vec w}\ne\pm {\vec u}$.

\begin{subthm}{adex:crofton:crofton}
Покажите, что для любой спрямляемой 
сферической кривой $\gamma$ выполняется
\[\length\gamma=\overline{\length\gamma^*_{\vec u}},\]
где $\overline{\length\gamma^*_{\vec u}}$ обозначает среднюю длину $\gamma^*_{\vec u}$ при изменении ${\vec u}$ на~$\mathbb{S}^2$.
\end{subthm}

\begin{subthm}{adex:crofton:hemisphere}
Найдите другое доказательство леммы о полусфере (\ref{lem:hemisphere}),
используя \ref{SHORT.adex:crofton:crofton}. 
\end{subthm}
 
\end{thm}

Часть \ref{SHORT.adex:crofton:crofton} упражнения --- это сферическая формула Крофтона.
Её можно переписать следующим образом:
\[\length\gamma=\overline n\cdot \pi,\]
где $\overline n$ --- среднее число точек пересечений $\gamma$ с экваторами.
Эквивалентность вытекает из теоремы Леви о монотонной сходимости.

\chapter{Кривизна}
\label{chap:curve-curvature}

Кривизной полагается называть способ мерить, насколько геометрический объект отклоняется \textit{прямого}, чего бы последнее ни значило.
Соответственно, кривизна кривой мерит, насколько она отклоняется от прямой в данной точке.

\section{Ускорение и скорость}

Напомним, что любую гладкую кривую можно параметризовать длиной (\ref{prop:arc-length-smooth}).
Полученная кривая, пусть это будет $\gamma$, остаётся гладкой и имеет единичную скорость;
то есть $|\gamma'|=1$.
Из курса механики должно быть известно, что при равномерном движении скорость перпендикулярна ускорению;
сейчас мы это переформулируем.

\begin{thm}{Предложение}\label{prop:a'-pertp-a''}
Пусть $\gamma$ --- гладкая пространственная кривая с единичной скоростью.
Тогда $\gamma'(s)\perp \gamma''(s)$ для любого~$s$.
\end{thm}

\index{скалярное произведение}\emph{Скалярное произведение} двух векторов $\vec v$ и $\vec w$ будет обозначаться $\langle \vec v,\vec w\rangle$.
Напомним, что производная скалярного произведения удовлетворяет правилу произведения (также известное как тождество Лейбница);
то есть, если $\vec v=\vec v(t)$ и $\vec w=\vec w(t)$ --- гладкие векторозначные функции вещественной переменной $t$, то
\[\langle \vec v,\vec w\rangle'=\langle \vec v',\vec w\rangle+\langle \vec v,\vec w'\rangle.\]

\parbf{Доказательство.}
Продифференцировав тождество $\langle\gamma',\gamma'\rangle\z=1$, получим
$2\cdot \langle\gamma'',\gamma' \rangle=\langle\gamma',\gamma'\rangle'=0$;
следовательно, $\gamma''\perp\gamma'$.
\qeds

\section{Кривизна}\label{sec:curvature}

Пусть $\gamma$ --- гладкая пространственная кривая с единичной скоростью.
Тогда величина $|\gamma''(s)|$ называется \index{10k@$\kur$ (кривизна)}\index{кривизна}\emph{кривизной} $\gamma$ в момент времени~$s$.
Можно также сказать, что $|\gamma''(s)|$ --- кривизна в точке $p=\gamma(s)$;
если наша кривая простая, то это не приводит к неоднозначности.
Кривизна обозначается через $\kur(s)$ или $\kur(s)_\gamma$;
в случае простых кривых можно также пользоваться обозначениями $\kur(p)$ или $\kur(p)_\gamma$.

{\sloppy

\begin{thm}{Упражнение}\label{ex:zero-curvature-curve}
Покажите, что среди гладких простых пространственных кривых, только отрезки прямой имеют нулевую кривизну в каждой точке.
\end{thm}

}

\begin{thm}{Упражнение}\label{ex:scaled-curvature}
Пусть $\gamma$ --- гладкая простая пространственная кривая и $\gamma_{\lambda}$ --- её гомотетия с коэффициентом $\lambda >0$;
то есть $\gamma_{\lambda}(t)=\lambda \cdot\gamma(t)$ для любого~$t$.
Покажите, что
\[\kur(\lambda \cdot p )_{\gamma_{\lambda}}
= \frac{\kur(p)_{\gamma}}\lambda\]
для любого $p \in \gamma$.
\end{thm}

\begin{thm}{Упражнение}\label{ex:curvature-of-spherical-curve}
Покажите, что кривизна любой гладкой сферической кривой не меньше $1$.
\end{thm}

\section{Касательная индикатриса}\label{sec:Tangent indicatrix}

Пусть $\gamma$ --- гладкая пространственная кривая.
Кривая 
\[\tan(t)=\tfrac{\gamma'(t)}{|\gamma'(t)|};
\eqlbl{eq:tantrix}\] 
называется её \index{касательная!индикатриса}\emph{касательной индикатрисой}.
Заметим, что $|\tan(t)|\z=1$ для любого $t$;
то есть $\tan$ является сферической кривой.

Если $s\mapsto \gamma(s)$ --- параметризация длиной, то $\tan(s)=\gamma'(s)$.
В этом случае у нас есть следующее выражение для кривизны:
\[\kur(s)\z=|\tan'(s)|\z=|\gamma''(s)|.\]

Для общей параметризации $t\mapsto \gamma(t)$,
кривизна выражается как
\[ \kur(t)=\frac{|\tan'(t)|}{|\gamma'(t)|}.\eqlbl{eq:curvature}\]
Действительно, если $s(t)$ --- параметр длины, то $s'(t)=|\gamma'(t)|$, и
\begin{align*}
\kur&=\left|\frac{d\tan}{ ds}\right|=
\left|\frac{d\tan}{ dt}\right|/\left|\frac{ds}{ dt}\right|=
\frac{|\tan'|}{|\gamma'|}.
\end{align*}

Напомним, что гладкая кривая имеет гладкую регулярную параметризацию; см. \ref{sec:Smooth curves}.
Если $\gamma$ гладкая, то $t\mapsto \tan(t)$ задаёт гладкую параметризацию,
и она регулярна, если $|\tan'|=\kur\cdot|\gamma'|$ не обращается в нуль.
Отсюда следует, что \textit{если кривизна гладкой кривой не обнуляется, то её касательная индикатриса --- гладкая кривая}.

\begin{thm}{Упражнение}\label{ex:curvature-formulas}
Пусть $\gamma$ --- гладкая пространственная кривая.
Используя \ref{eq:tantrix} и \ref{eq:curvature}, выведите следующие формулы для её кривизны:
\setlength{\columnseprule}{0.4pt}
\begin{multicols}{2}

\begin{subthm}{ex:curvature-formulas:a} 
\[\kur=\frac{|\vec w|}{|\gamma'|^2},\]
где $\vec w$ --- проекция $\gamma''(t)$ на нормальную плоскость к $\gamma'(t)$;
\end{subthm}

\columnbreak

\begin{subthm}{ex:curvature-formulas:b}
\[\kur=\frac{|\gamma''\times \gamma'|}{|\gamma'|^{3}},\]
где $\times$ обозначает \index{векторное произведение}\emph{векторное произведение}.
\end{subthm}

\vfill\null
\end{multicols}
\end{thm}

{\sloppy

\begin{thm}{Упражнение}\label{ex:curvature-graph}
Пусть $f$ --- гладкая вещественная функция.
Применив формулы предыдущего упражнения, докажите, что 
\[\kur=\frac{|f''(x)|}{(1+f'(x)^2)^{\frac32}}\]
есть кривизна графика $y=f(x)$ в точке $(x,f(x))$. 
\end{thm}

}

\begin{thm}{Продвинутое упражнение}\label{ex:approximation-const-curvature}
Покажите, что любую гладкую кривую $\gamma\:\mathbb{I}\z\to\mathbb{R}^3$ с кривизной не более $1$ можно аппроксимировать гладкими кривыми с постоянной кривизной~$1$.
То есть, найдётся последовательность гладких кривых $\gamma_n\:\mathbb{I}\z\to\mathbb{R}^3$ с постоянной кривизной $1$ такая, что $\gamma_n(t)\to \gamma(t)$ для любого~$t$.
\end{thm}

\section{Касательные}

{

\begin{wrapfigure}[6]{r}{32 mm}
\vskip-14mm
\centering
\includegraphics{mppics/pic-3400}
\vskip0mm
\end{wrapfigure}

Пусть $\gamma$ --- гладкая пространственная кривая, и $\tan$ --- её касательная индикатриса.
Прямая, проходящая через $\gamma(t)$ в направлении $\tan(t)$, называется \index{касательная!прямая}\emph{касательной} к $\gamma$ при~$t$.
Любой вектор, пропорциональный $\tan(t)$, называется \index{касательная!вектор}\emph{касательным} к $\gamma$ при~$t$.

Касательную прямую можно также определить как единственную прямую \index{порядок касания}\emph{первого порядка касания} с $\gamma$ при $t$;
то есть $\rho(\ell)=o(\ell)$, где $\rho(\ell)$ обозначает расстояние от $\gamma(t+\ell)$ до прямой.

}

\begin{thm}{Продвинутое упражнение}\label{ex:no-parallel-tangents}
Постройте гладкую замкнутую пространственную кривую без параллельных касательных прямых.
\end{thm}

Говорят, что гладкие кривые $\gamma_1$ и $\gamma_2$ \index{касательная!кривые}\emph{касаются} при $s_1$ и $s_2$, если $\gamma_1(s_1)=\gamma_2(s_2)$ и касательная к $\gamma_1$ при $s_1$ совпадает с касательной к $\gamma_2$ при $s_2$.
Можно также сказать, что эти кривые касаются друг друга в точке $p=\gamma_1(s_1)\z=\gamma_2(s_2)$;
если обе кривые простые, то это не приводит к неоднозначности.

\section{Полная кривизна}\label{sec:Total curvature}

Пусть $\gamma\:\mathbb{I}\to\mathbb{R}^3$ --- гладкая кривая с единичной скоростью.
Интеграл 
\[\tc\gamma\df\int_{\mathbb{I}}\kur(s)\cdot ds\]
называется \index{полная!кривизна}\emph{полной кривизной}\label{page:total curvature of:smooth-def}
$\gamma$.
Полная кривизна мерит насколько кривая поворачивает;
её также называют \index{вариация поворота}\emph{вариацией поворота}.

Формула замены переменной позволяет выразить полную кривизну кривой с общей параметризацией $t\mapsto \gamma(t)$:
\[
\tc\gamma\df\int_{\mathbb{I}}\kur(t)\cdot|\gamma'(t)| \cdot dt.
\eqlbl{eq:tocurv}
\]

Поскольку у гладкой кривой, параметризованной длиной, кривизна равна скорости касательной индикатрисы,
получаем следующее.

\begin{thm}{Наблюдение}\label{obs:tantrix}
Полная кривизна гладкой кривой равна длине её касательной индикатрисы.
\end{thm}

\begin{thm}{Упражнение}\label{ex:helix-curvature}
Найдите кривизну винтовой линии 
\[\gamma_{a,b}(t)=(a\cdot \cos t,a\cdot \sin t,b\cdot t),\]
её касательную индикатрису и полную кривизну её дуги $\gamma_{a,b}|_{[0,2\cdot\pi]}$.
\end{thm}

{\sloppy

\begin{thm}{Теорема Фенхеля}
\label{thm:fenchel}
\index{теорема Фенхеля}
Полная кривизна любой замкнутой гладкой пространственной кривой не меньше $2\cdot\pi$.
\end{thm}

}

\parbf{Доказательство.}
Пусть $\gamma$ --- замкнутая гладкая пространственная кривая.
Давайте считать, что $\gamma$ описывается петлёй $\gamma\:[a,b]\z\to \mathbb{R}^3$, параметризованной длиной;
в этом случае, $\gamma(a)=\gamma(b)$ и $\gamma'(a)\z=\gamma'(b)$.

Рассмотрим её касательную индикатрису $\tan=\gamma'$.
Напомним, что $|\tan(s)|=1$ для любого $s$; то есть $\tan$ --- замкнутая сферическая кривая.

Покажем, что $\tan$ не может лежать в полусфере.
Рассуждая от противного, можно предположить, что она лежит в полусфере, определяемой неравенством $z>0$ в координатах $(x,y,z)$.
Другими словами, если $\gamma(s)=(x(s), y(s), z(s))$, то $z'(t)>0$ для любого~$t$.
Следовательно,
\[z(b)-z(a)=\int_a^b z'(s)\cdot ds>0.\]
В частности, $\gamma(a)\ne \gamma(b)$ --- противоречие.

Применив наблюдение (\ref{obs:tantrix}) и лемму о полусфере (\ref{lem:hemisphere}), получим  
\[\tc\gamma=\length \tan\ge2\cdot\pi.\]
\qedsf

\begin{thm}{Упражнение}\label{ex:length>=2pi}
Докажите, что замкнутая пространственная кривая $\gamma$ с кривизной не более~$1$ не короче единичной окружности;
то есть
\[\length\gamma\ge 2\cdot \pi.\]

\end{thm}

\begin{thm}{Продвинутое упражнение}\label{ex:gamma/|gamma|}
{\sloppy
Предположим, что гладкая пространственная кривая $\gamma$ не проходит через начало координат.
Рассмотрим сферическую кривую $\sigma(t)\z\df\frac{\gamma(t)}{|\gamma(t)|}$.
Докажите, что
\[\length \sigma< \tc\gamma+\pi.\]
Более того, если $\gamma$ замкнута, то
\[\length \sigma\le \tc\gamma.\]

}
\end{thm}

Последнее неравенство даёт другое доказательство теоремы Фенхеля.
Действительно, можно предположить, что начало координат лежит на хорде~$\gamma$.
В этом случае, замкнутая сферическая кривая $\sigma$ идёт от одной точки к её антиподу и возвращается обратно,
преодолевая расстояние $\pi$ в каждую сторону.
Значит
\[\length\sigma\ge 2\cdot\pi.\]

Напомним, что кривизна сферической кривой хотя бы $1$
(см. \ref{ex:curvature-of-spherical-curve}).
В частности, длина сферической кривой не превосходит её полной кривизны.
Следующая теорема показывает, что то же неравенство выполняется для \textit{замкнутых} кривых в единичном шаре.

\begin{thm}{Теорема о ДНК}\label{thm:DNA}
Пусть $\gamma$ --- гладкая замкнутая кривая, которая лежит в единичном шаре.
Тогда
\[\tc\gamma\ge \length\gamma.\]

\end{thm}

Несколько доказательств этой теоремы собраны в статье Сергея Табачникова~\cite{tabachnikov}.
Двумерный случай был доказан Иштваном Фари \cite{fary1950}.
Дон Чакериан \cite{chakerian1962} обобщил теорему на старшие размерности.
Следующее упражнение основано на другом его доказательстве \cite{chakerian1964},
и ещё одно обсуждается в~\ref{sec:DNA-poly}.

\begin{thm}{Упражнение}\label{ex:DNA}
Рассмотрим гладкую кривую с единичной скоростью $\gamma\:[0,\ell]\to\mathbb{R}^3$, лежащую в единичном шаре; то есть $|\gamma|\le 1$.

\begin{subthm}{ex:DNA:c''c>=k}
Докажите, что
\[\langle\gamma''(s),\gamma(s)\rangle\ge-\kur(s)\]
для любого~$s$.
\end{subthm}

\begin{subthm}{ex:DNA:int>=length-tc}
Воспользовавшись \ref{SHORT.ex:DNA:c''c>=k}, докажите, что
\[\int_0^\ell\langle\gamma(s),\gamma'(s)\rangle'\cdot ds\ge
\ell-\tc\gamma.\]

\end{subthm}

\begin{subthm}{ex:DNA:end}
Пусть $\gamma(0)=\gamma(\ell)$ и $\gamma'(0)=\gamma'(\ell)$.
Докажите, что
\[\int_0^\ell\langle\gamma(s),\gamma'(s)\rangle'\cdot ds=0.\]
Докажите \ref{thm:DNA}, воспользовавшись \ref{SHORT.ex:DNA:int>=length-tc} и этим равенством.
\end{subthm}
\end{thm}

\section{Выпуклые кривые}

В этом разделе мы покажем, что касательная индикатриса выпуклой кривой вращается монотонно.
Следующее упражнение будет использовано в доказательстве.

\begin{thm}{Упражнение}\label{ex:tangent-support}
Пусть $F$ --- выпуклое множество на плоскости ограниченное гладкой кривой $\gamma$.
Докажите, что прямая $\ell$ касается $\gamma$ в точке $p$ тогда и только тогда, когда $\ell$ \index{опорная!прямая}\emph{подпирает} $F$ в точке $p$;
то есть $\ell\ni p$ и $F$ лежит в полуплоскости, ограниченной~$\ell$.
\end{thm}

Напомним, что отображение монотонно, если прообраз любой точки в целевом пространстве связен (и в частности непуст).

\begin{thm}{Предложение}\label{prop:convex-monotone}
Пусть $\gamma$ --- гладкая выпуклая плоская кривая.

\begin{subthm}{prop:convex-monotone:closed}
Если $\gamma$ замкнута и $\gamma\:\mathbb{S}^1\to \mathbb{R}^2$ --- её параметризация, то её касательная индикатриса $\tan\:\mathbb{S}^1\to \mathbb{S}^1$ --- монотонное отображение.
\end{subthm}

\begin{subthm}{prop:convex-monotone:open}
Если $\gamma$ открыта и $\gamma\:\mathbb{R}\to \mathbb{R}^2$ --- её параметризация, то её касательная индикатриса $\tan\:\mathbb{R}\to\mathbb{S}^1$ определяет монотонное отображение на некоторый интервал в замкнутой полуокружности.
\end{subthm}

\end{thm}

Согласно следствию ниже, \textit{для выпуклых кривых выполняется равенство в теореме Фенхеля} (\ref{thm:fenchel}).
Позже, в \ref{prop:fenchel=}, мы покажем, что равенство выполняется \textit{только} для выпуклых кривых.

\begin{thm}{Следствие}\label{cor:fenchel=convex}
Пусть $\gamma$ --- выпуклая кривая на плоскости.

\begin{subthm}{}
Если $\gamma$ замкнута, то $\tc\gamma=2\cdot\pi$.
\end{subthm}

\begin{subthm}{}
Если $\gamma$ открыта, то $\tc\gamma\le\pi$.
\end{subthm}

\end{thm}

\parbf{Доказательство.}
Следует из \ref{prop:convex-monotone:closed}, \ref{obs:tantrix} и \ref{ex:integral-length-0}.
\qeds

\begin{wrapfigure}{r}{32 mm}
\vskip-0mm
\centering
\includegraphics{mppics/pic-3500}
\vskip0mm
\end{wrapfigure}

\parbf{Доказательство \ref{prop:convex-monotone};} \ref{SHORT.prop:convex-monotone:closed}.
Так как $\gamma$ замкнута, она ограничивает компактное выпуклое множество $F$.
Можно предположить, что $F$ лежит слева от $\gamma$.

Выберем единичный вектор $\vec u$ и координаты, с осью $x$ направленной по $\vec u$.

Из упражнения~\ref{ex:tangent-support} вытекает, что $\vec u=\tan(s)$ тогда и только тогда, когда $p=\gamma(s)$ является точкой минимума $y$-координаты на~$F$.
Действительно, пусть $p$ --- точка минимума.
Проведём через $p$ горизонтальную прямую $\ell$; она подпирает $F$ в $p$.
Согласно упражнению, $\ell$ касается $\gamma$ в $p$.
Так как $F$ лежит слева от $\gamma$, получаем, что $\tan(s)=\vec u$.
И наоборот, если $\tan(s)=\vec u$, то касательная в точке $p=\gamma(s)$ горизонтальна,
и, согласно упражнению, она подпирает $F$.
Так как $F$ лежит от $\gamma$ слева, $y$-координата достигает минимума на $F$ в точке $p$.

Поскольку множество $F$ компактно, $y$-координата достигает минимума на $F$ в некоторой точке $p=\gamma(s)$.
Точка $p$ может оказаться единственной, в этом случае, $\tan^{-1}\{\vec u\}=\{s\}$,
или же $y$-координата может иметь отрезок минимальных точек на $\gamma$, в этом случае, $\tan^{-1}\{\vec u\}$ --- дуга в~$\mathbb{S}^1$.
Следовательно, отображение $\tan\:\mathbb{S}^1\to \mathbb{S}^1$ монотонно.

\parit{\ref{SHORT.prop:convex-monotone:open}} 
То же рассуждение, что и в \ref{SHORT.prop:convex-monotone:closed}, показывает, что $\tan$ является монотонным отображением на свой образ.
Ясно, что образ связен в $\mathbb{S}^1$.
Остаётся показать, что образ лежит в полуокружности;
другими словами,  
\[\measuredangle(\vec w,\tan(s))\ge\tfrac\pi2
\quad\text{для некоторого}\quad \vec w\quad\text{и любого}\quad
s.
\eqlbl{eq:<(w,tan).pi/2}
\]

\begin{figure}[ht!]
\centering
\includegraphics{mppics/pic-3502}
\end{figure}

Поскольку $\gamma$ открыта, она ограничивает некомпактное выпуклое замкнутое множество $F$.
Как и прежде, считаем, что $F$ лежит слева от $\gamma$.

Докажем, что существует луч, скажем $h$, лежащий в~$F$.
Можно предположить, что начало координат $o$ лежит в~$F$.
Рассмотрим последовательность точек $q_n\in F$ таких, что $|q_n|\z\to \infty$ при $n\to \infty$.
Обозначим через $\vec v_n$ единичный вектор в направлении $q_n$; то есть~$\vec v_n=\tfrac{q_n}{|q_n|}$.

Так как единичная окружность компактна, перейдя к подпоследовательности $q_n$, можно считать, что $\vec v_n$ сходится к единичному вектору, скажем к $\vec v$.
Проведём луч $h$ из $o$ в направлении $\vec v$.
Любую точку на $h$ можно приблизить точками из отрезков $[o,q_n]$ при $n\to\infty$.
Значит луч $h$ лежит в $F$, ведь $F$ замкнуто.

Пусть $\vec w$ --- поворот $\vec v$ против часовой стрелки на угол $\tfrac\pi 2$,
$\ell$ --- касательная прямая к $\gamma$ в точке $p=\gamma(s)$,
и $H$ --- замкнутая левая полуплоскость ограниченная $\ell$;
то есть $H$ лежит слева от $\ell$ по направлению $\tan(s)$.
То же рассуждение, что и в \ref{SHORT.prop:convex-monotone:closed}, показывает, что $F$, а значит, и $h$,
лежат в $H$.
В частности, $\vec v$ указывает из $p$ в $H$, а это эквивалентно~\ref{eq:<(w,tan).pi/2}.
\qeds

\section{Лемма о луке}

Следующая лемма была доказана Эрхардом Шмидтом \cite{schmidt}; она обобщает результат Акселя Шура \cite{shur}.

\begin{wrapfigure}[9]{r}{39 mm}
\vskip-6mm
\centering
\includegraphics{mppics/pic-3510}
\vskip0mm
\end{wrapfigure}

{\sloppy

Лемма является дифференциально-геометрическим аналогом так называемой {}\emph{леммы о руке} Огюстена-Луи Коши, согласно которой, \textit{если для выпуклого многоугольника
$p_0\dots p_n$ на плоскости и пространственной ломаной $q_0\dots q_n$ выполняется 
\begin{align*}
|p_i-p_{i-1}|&=|q_i-q_{i-1}|,
\\
\measuredangle\hinge{p_i}{p_{i+1}}{p_{i-1}}&\le \measuredangle\hinge{q_i}{q_{i+1}}{q_{i-1}}
\end{align*}
для всех $i$, то $|p_0-p_n|\le |q_0-q_n|$.}
(Можно думать про это так: если распрямлять все суставы руки, то расстояние от плеча до кончика среднего пальца увеличится.)

}

\begin{thm}{Лемма}\label{lem:bow}\index{лемма о луке}
Пусть $\gamma_1\:[a,b]\to\mathbb{R}^2$ и $\gamma_2\:[a,b] \to\mathbb{R}^3$ --- две гладкие кривые с единичной скоростью.
Предположим, что $\kur(s)_{\gamma_1}\ge\kur(s)_{\gamma_2}$ для любого $s$ 
и кривая
$\gamma_1$ --- дуга выпуклой кривой; то есть она идёт по границе выпуклой плоской фигуры.
Тогда расстояние между конечными точками $\gamma_1$ не превосходит расстояние между конечными точками $\gamma_2$; то есть
\[|\gamma_1(b)-\gamma_1(a)|\le |\gamma_2(b)-\gamma_2(a)|.\]

\end{thm}

Следующее упражнение говорит, что выпуклость $\gamma_1$ необходима.
Его стоит решить, прежде чем читать доказательство.

\begin{thm}{Упражнение}\label{ex:anti-bow}
Постройте две простые гладкие кривые с единичной скоростью $\gamma_1,\gamma_2\:[a,b]\to\mathbb{R}^2$, что $\kur(s)_{\gamma_1}>\kur(s)_{\gamma_2}>0$ для любого $s$ и
\[|\gamma_1(b)-\gamma_1(a)|> |\gamma_2(b)-\gamma_2(a)|.\]

\end{thm}

\parbf{Доказательство леммы.}
Можно предположить, что $\gamma_1(a)\ne \gamma_1(b)$;
иначе нечего доказывать.
По выпуклости, кривая $\gamma_1$ лежит с одной стороны от прямой $\ell$, проходящей через $\gamma_1(a)$ и $\gamma_1(b)$;
будем думать, что $\ell$ направлена горизонтально и $\gamma_1$ лежит под ней.

Пусть $\gamma_1(s_0)$ будет самой низкой точкой на $\gamma_1$;
то есть $\gamma_1(s_0)$ имеет минимальную $y$-координату.

Обозначим через $\tan_1$ и $\tan_2$ касательные индикатрисы $\gamma_1$ и $\gamma_2$, соответственно.
Рассмотрим два единичных вектора 
\[
\vec u_1=\tan_1(s_0)=\gamma_1'(s_0),
\quad\text{и}\quad
\vec u_2=\tan_2(s_0)=\gamma_2'(s_0).
\]
Отметим, что $\gamma_1(b)$ лежит в направлении $\vec u_1$ от $\gamma_1(a)$.

\begin{wrapfigure}[10]{r}{39 mm}
\vskip-3mm
\centering
\includegraphics{mppics/pic-57}
\vskip0mm
\end{wrapfigure}

Покажем, что 
\[\measuredangle(\gamma'_1(s),\vec u_1)\ge \measuredangle(\gamma'_2(s),\vec u_2)
\eqlbl{<gamma',u}
\]
для любого $s$.
Давайте считать, что $s\le s_0$; случай $s\ge s_0$ аналогичен.

Заметим, что
\[
\begin{aligned}
\measuredangle(\gamma'_1(s),\vec u_1)&=\measuredangle(\tan_1(s),\vec u_1)=
\\
&=\length (\tan_1|_{[s,s_0]}).
\end{aligned}
\eqlbl{<=length}\]
для любого $s\le s_0$.
Действительно, по \ref{ex:tangent-support}, $y$-координата $\gamma_1$ не убывает на интервале $[a,s_0]$.
Следовательно, дуга $\tan_1|_{[a,s_0]}$ лежит в одном из единичных полукругов с концами $\vec u_1$ и $-\vec u_1$.
Остаётся применить \ref{obs:tantrix}, \ref{prop:convex-monotone} и \ref{ex:integral-length-0}.

Из \ref{obs:S2-length}, также получаем, что 
\[\measuredangle(\gamma'_2(s),\vec u_2)=\measuredangle(\tan_2(s),\vec u_2)\le \length (\tan_2|_{[s,s_0]}).
\eqlbl{<=<length}\]

{

Далее,
\begin{align*}
\length (\tan_1|_{[s,s_0]})
&=\int_s^{s_0}|\tan_1'(t)|\cdot d t=
\\
&=\int_s^{s_0}\kur_1(t)\cdot d t\ge
\int_s^{s_0}\kur_2(t)\cdot d t=
\\
&=\int_s^{s_0}|\tan_2'(t)|\cdot d t= 
\length (\tan_2|_{[s,s_0]}).
\end{align*}
Это неравенство, вместе с \ref{<=length} и \ref{<=<length}, влечёт \ref{<gamma',u}.
}

Поскольку $1=|\gamma_1'(s)|=|\gamma_2'(s)|=|\vec u_1|=|\vec u_2|$,
\[\langle\gamma'_1(s),\vec u_1\rangle=\cos \measuredangle(\gamma'_1(s),\vec u_1)
\quad\text{и}\quad
\langle\gamma'_2(s),\vec u_2\rangle=\cos \measuredangle(\gamma'_2(s),\vec u_2).
\]
Косинус убывает на интервале $[0,\pi]$; следовательно, \ref{<gamma',u} влечёт 
\[\langle\gamma'_1(s),\vec u_1\rangle\le \langle\gamma'_2(s),\vec u_2\rangle\eqlbl{<gamma',u>}\]
для любого~$s$.
Далее, 
\[|\gamma_1(b)-\gamma_1(a)|=\langle \vec u_1,\gamma_1(b)-\gamma_1(a)\rangle,\]
ибо $\gamma_1(b)$ лежит в направлении $\vec u_1$ от $\gamma_1(a)$.
Поскольку $\vec u_2$ --- единичный вектор,
\[|\gamma_2(b)-\gamma_2(a)|\ge\langle \vec u_2,\gamma_2(b)-\gamma_2(a)\rangle.\]

Проинтегрировав \ref{<gamma',u>}, получим 
\begin{align*}
|\gamma_1(b)-\gamma_1(a)|&=\langle \vec u_1,\gamma_1(b)-\gamma_1(a)\rangle=
\\
&=
\int_a^b\langle \vec u_1,\gamma'_1(s)\rangle\cdot ds \le 
\int_a^b\langle \vec u_2,\gamma'_2(s)\rangle\cdot ds 
=
\\
&=\langle \vec u_2,\gamma_2(b)-\gamma_2(a)\rangle
\le |\gamma_2(b)-\gamma_2(a)|.
\end{align*}
\qedsf

\begin{thm}{Продвинутое упражнение}\label{ex:bow'}
Предположим, что $\gamma_1$ и $\gamma_2$ как в лемме о луке (\ref{lem:bow}),
а $\tan_1$ и $\tan_2$ --- их касательные индикатрисы.

{

\begin{wrapfigure}{r}{50 mm}
\vskip-0mm
\centering
\includegraphics{mppics/pic-251}
\vskip-4mm
\end{wrapfigure}

Пусть 
\begin{align*}
\vec w_i&=\gamma_i(b)-\gamma_i(a),
\\
\alpha_i&=\measuredangle(\tan_i(a),\vec w_i),
\\
\beta_i&=\measuredangle(\tan_i(b),\vec w_i).
\end{align*}

}

\begin{subthm}{ex:bow'+}
Предположим, что $\beta_1\le\tfrac\pi2$.
Покажите, что $\alpha_1\ge \alpha_2$.
\end{subthm}

\begin{subthm}{ex:bow'-}
Постройте пример, показывающий, что неравенство $\alpha_1\ge \alpha_2$ не выполняется в общем случае.
\end{subthm}

\end{thm}

\begin{thm}{Упражнение}\label{ex:length-dist}
Пусть $\gamma\:[a,b]\to \mathbb{R}^3$ --- гладкая кривая и $0\z<\theta\z\le\tfrac\pi2$.
Предположим, что
\[\tc\gamma\le 2\cdot\theta.\]

\begin{subthm}{ex:length-dist:>}
Покажите, что
\[|\gamma(b)-\gamma(a)|> \cos\theta\cdot\length\gamma.\]
\end{subthm}

\begin{subthm}{ex:length-dist:self-intersection:>pi}
Покажите, что $\tc\gamma>\pi$, если $\gamma$ самопересекается.
Нарисуйте гладкую кривую $\gamma$ с $\tc\gamma<2\cdot\pi$ и самопересечением.
\end{subthm}


\begin{subthm}{ex:length-dist:=}
Покажите, что неравенство в \ref{SHORT.ex:length-dist:>} оптимально; то есть для заданного 
$\theta$ существует гладкая кривая $\gamma$, такая что $\tc\gamma\z\le 
2\cdot\theta$, и $\frac{|\gamma(b)-\gamma(a)|}{\length\gamma}$ произвольно 
близко к $\cos\theta$.
\end{subthm}

\end{thm}

\begin{thm}{Упражнение}\label{ex:schwartz}
Пусть $p$ и $q$ --- точки на единичной окружности, делящие её на две дуги с длинами $\ell_1<\ell_2$.
Предположим, что пространственная кривая $\gamma$ соединяет $p$ с $q$ и имеет кривизну не более $1$.
Покажите, что 
\[\length \gamma\le \ell_1
\quad\text{либо}\quad
\length \gamma\ge \ell_2.
\]
\end{thm}

Следующее упражнение обобщает \ref{ex:length>=2pi}.

\begin{thm}{Упражнение}\label{ex:loop}
Предположим, что $\gamma\:[a,b]\to \mathbb{R}^3$ --- гладкая петля с кривизной не более $1$.
Покажите, что 
\[\length\gamma\ge2\cdot\pi.\]

\end{thm}

\begin{thm}{Упражнение}\label{ex:bow-upper}
Пусть $\kur$ --- гладкая неотрицательная функция, определённая на интервале $[0,\ell]$.
Постройте гладкую кривую $\gamma\:[0,\ell]\to\mathbb{R}^3$ с единичной скоростью и кривизной $\kur(s)$ при любом $s$ так, чтобы расстояние $|\gamma(\ell)-\gamma(0)|$ было произвольно близко к $\ell$.
\end{thm}

\begin{wrapfigure}{r}{41 mm}
\vskip-0mm
\centering
\includegraphics{mppics/pic-283}
\vskip0mm
\end{wrapfigure}

\begin{thm}{Продвинутое упражнение}\label{ex:gromov-twist}\\
Пусть $\gamma$ --- замкнутая гладкая пространственная кривая с кривизной не более $2$.
Предположим, что $|\gamma(t)|\le 1$ для любого~$t$.
Покажите, что если $\gamma(t)\ne 0$, то 
\[|\gamma(t)| \le \sin (\alpha(t)),\]
где $\alpha(t)$ обозначает угол между $\gamma(t)$ и~$\gamma'(t)$.
\end{thm}

\chapter{Ломаные}
\label{chap:poly}

Эта глава связывает кривизну кривой и углы вписанных в неё ломаных;
она поможет обзавестись правильной геометрической интуицией кривизны.

\section{Кусочно-гладкие кривые}

{

\begin{wrapfigure}{o}{25 mm}
\vskip-2mm
\centering
\includegraphics{mppics/pic-54}
\end{wrapfigure}

Предположим, что $\alpha\:[a,b]\to \mathbb{R}^3$ и $\beta\:[b,c]\z\to \mathbb{R}^3$ --- такие две кривые, что $\alpha(b)\z=\beta(b)$.
Тогда их можно объединить в одну кривую $\gamma\:[a,c]\z\to \mathbb{R}^3$, определяемую как 
\[\gamma(t)=
\begin{cases}
\alpha(t)&\text{если}\quad t\le b,
\\
\beta(t)&\text{если}\quad t> b.
\end{cases}
\]
Кривая $\gamma$ называется 
\index{произведение кривых}\emph{произведением} $\alpha$ и $\beta$.
(Условие $\alpha(b)=\beta(b)$ обеспечивает непрерывность отображения $t\mapsto\gamma(t)$.)

}

То же определение можно применить если $\alpha$ и/или $\beta$ определены на полуоткрытых интервалах
$(a,b]$ и/или $[b,c)$.

Предположение, что интервалы $[a,b]$ и $[b,c]$ продолжают друг друга, не существенно.
Достаточно, чтобы конечная точка $\alpha$ совпадала с начальной точкой $\beta$.
Тогда интервалы можно подвинуть так, чтоб один продолжал другой. 

Если дополнительно $\beta(c)=\alpha(a)$, то можно рассмотреть циклическое произведение этих кривых;
так мы получим замкнутую кривую.

Если $\alpha'(b)$ и $\beta'(b)$ определены, то угол $\theta\z=\measuredangle(\alpha'(b),\beta'(b))$ называется \index{внешний угол}\emph{внешним углом} $\gamma$ при~$b$.
Если $\theta=\pi$, то говорят, что $\gamma$ имеет \index{точка!возврата}\emph{точку возврата} при~$b$.

Пространственная кривая $\gamma$ называется \index{кусочно-гладкая кривая}\emph{кусочно-гладкой}, если её можно представить как произведение конечного числа гладких кривых; если $\gamma$ замкнута, то предполагается, что произведение циклическое.

Если $\gamma$ --- произведение гладких дуг $\gamma_1,\dots,\gamma_n$, то полная кривизна $\gamma$ определяется как сумма полных кривизн $\gamma_i$ и внешних углов;
то есть 
\[\tc\gamma=\tc{\gamma_1}+\dots+\tc{\gamma_n}+\theta_1+\dots+\theta_{n-1},\]
где $\theta_i$ --- внешний угол на стыке между $\gamma_i$ и $\gamma_{i+1}$.

Если же $\gamma$ замкнута, то её полная кривизна определяется как
\[\tc\gamma=\tc{\gamma_1}+\dots+\tc{\gamma_n}+\theta_1+\dots+\theta_{n},\]
где $\theta_n$ --- внешний угол на стыке между $\gamma_n$ и $\gamma_1$.

{

\begin{wrapfigure}{r}{23 mm}
\vskip-3mm
\centering
\includegraphics{mppics/pic-354}
\end{wrapfigure}

В частности, если $\gamma\:[a,b] \z\to \mathbb{R}^3$ --- гладкая петля, и $\hat\gamma$ соответствующая ей замкнутая кривая, то
\[\tc{\hat\gamma}\df\tc\gamma + \theta,\]
где $\theta=\measuredangle(\gamma'(a),\gamma'(b))$.

}

\section{Обобщённая теорема Фенхеля}

\begin{thm}{Теорема}\label{thm:gen-fenchel}
Пусть $\gamma$ --- замкнутая кусочно-гладкая пространственная кривая.
Тогда
\[\tc\gamma\ge2\cdot\pi.\]

\end{thm}

{\sloppy

\parbf{Доказательство.}
Пусть $\gamma$ --- циклическое произведение дуг $\gamma_1,\dots,\gamma_n$,
и $\theta_1,\dots,\theta_n$ --- её внешние углы.
Нам нужно показать, что \index{10phi@$\tc{\gamma}$ (полная кривизна)}
\[\tc{\gamma_1}+\dots+\tc{\gamma_n}+\theta_1+\dots+\theta_n\ge2\cdot\pi.\eqlbl{eq:gen-fenchel}\]

}

Рассмотрим касательную индикатрису $\tan_i$ для каждой дуги $\gamma_i$;
все они сферические дуги.

Рассуждение в доказательстве теоремы Фенхеля (\ref{thm:fenchel}), показывает, что все $\tan_1,\dots,\tan_n$ не могут лежать в открытой полусфере.

Сферическое расстояние от конечной точки $\tan_i$ до начальной точки $\tan_{i+1}$ равно внешнему углу $\theta_i$ (мы нумеруем дуги циклически, поэтому $\gamma_{n+1}=\gamma_1$).
Соединив конечную точку $\tan_i$ с начальной точкой $\tan_{i+1}$ короткой дугой большой окружности на сфере,
получаем замкнутую сферическую кривую, которая на $\theta_1+\dots+\theta_n$ длиннее общей длины $\tan_1,\dots,\tan_n$.

Применив к ней лемму о полусфере (\ref{lem:hemisphere}), получим
\[\length\tan_1+\dots+\length\tan_n+\theta_1+\dots+\theta_n\ge 2\cdot\pi.\]
Остаётся применить \ref{obs:tantrix}.
\qedsf

\begin{thm}{Лемма о хорде}\label{lem:chord}
Пусть $\gamma\:[a,b]\z\to\mathbb{R}^3$
--- гладкая дуга, и
$\ell$ --- её хорда.
Предположим, что $\gamma$ подходит к $\ell$ под углами $\alpha$ и $\beta$ в точках $\gamma(a)$ и $\gamma(b)$ соответственно;
то есть
\[\alpha=\measuredangle(\vec w,\gamma'(a))\quad\text{и}\quad \beta=\measuredangle(\vec w,\gamma'(b)),\]
где $\vec w=\gamma(b)-\gamma(a)$.
Тогда 
\[\tc\gamma\ge \alpha+\beta.\eqlbl{tc>a+b}\] 

\end{thm}

\parbf{Доказательство.}
Пропараметризуем хорду $\ell$ от $\gamma(b)$ до $\gamma(a)$, и
 пусть $\hat\gamma$ --- циклическое произведение $\gamma$ и $\ell$.
Замкнутая кривая $\hat\gamma$ имеет два внешних угла $\pi-\alpha$ и $\pi-\beta$.

\begin{wrapfigure}{r}{45 mm}
\vskip-5mm
\centering
\includegraphics{mppics/pic-53}
\vskip0mm
\end{wrapfigure}

Поскольку кривизна $\ell$ равна нулю, получаем 
\[\tc{\hat\gamma}=\tc\gamma+(\pi-\alpha)+(\pi-\beta).\]
По обобщённой теореме Фенхеля (\ref{thm:gen-fenchel}),
$\tc{\hat\gamma}\ge 2\cdot\pi$; отсюда \ref{tc>a+b}.
\qeds

\begin{thm}{Упражнение}\label{ex:chord-lemma-optimal}
Покажите, что оценка в лемме оптимальна.

Точнее, по данной паре различных точек $p, q$ и паре единичных векторов $\vec u,\vec v$ в $\mathbb{R}^3$,
постройте гладкую кривую $\gamma$, которая начинается в $p$ в направлении $\vec u$, заканчивается в $q$ в направлении $\vec v$, и при этом полная кривизна 
$\tc\gamma$ сколь угодно близка к $\measuredangle(\vec w,\vec u)+\measuredangle(\vec w,\vec v)$, где $\vec w$ --- вектор в направлении из $p$ в $q$.

\end{thm}

\section{Ломаные} 

Ломаные --- это частный случай кусочно-гладких кривых;
каждая дуга в их произведении является отрезком прямой.
Поскольку кривизна отрезка равна нулю, полная кривизна ломаной есть сумма её внешних углов.

\begin{thm}{Упражнение}\label{ex:monotonic-tc}
Пусть $a$, $b$, $c$, $d$ и $x$ --- различные точки в $\mathbb{R}^3$.
Покажите, что полная кривизна ломаной $abcd$ не может превысить полную кривизну $abxcd$; то есть 
\[\tc {abcd} \le \tc {abxcd}.\]

Воспользуйтесь этим, чтобы показать, что кривизна любой замкнутой ломаной не меньше $2\cdot\pi$.
\end{thm}

\begin{thm}{Предложение}\label{prop:inscribed-total-curvature}
Пусть ломаная $p_1\dots p_n$ вписана в гладкую кривую~$\gamma$.
Тогда 
\[\tc\gamma\ge \tc{p_1\dots p_n}.\]
Более того, если $\gamma$ замкнута, то можно считать, что ломаная $p_1\dots p_n$ также замкнута.

\end{thm}

\begin{wrapfigure}[7]{o}{40 mm}
\vskip-4mm
\centering
\includegraphics{mppics/pic-55}
\vskip0mm
\end{wrapfigure}

\parbf{Доказательство.}
Будем считать, что $\gamma$ замкнута. 
Введём обозначения:
\begin{align*}
p_i&=\gamma(t_i),
&
\alpha_i&=\measuredangle(\vec w_i,\vec v_i),
\\
\vec w_i&=p_{i+1}-p_i,
& 
\beta_i&=\measuredangle(\vec w_{i-1},\vec v_i),
\\
\vec v_i&=\gamma'(t_i),
&
\theta_i&=\measuredangle(\vec w_{i-1},\vec w_i).
\end{align*}
Индексы считаем по модулю $n$;
так что $p_{n+1}\z=p_1$.

Поскольку кривизна отрезков равна нулю, 
полная кривизна ломаной $p_1\dots p_n$ равна сумме внешних углов $\theta_i$.

По неравенству треугольника для углов \ref{thm:spherical-triangle-inq}, 
\[\theta_i\le \alpha_i+\beta_i.\]
По лемме о хорде (\ref{lem:chord}), полная кривизна дуги $\gamma$ от $p_i$ до $p_{i+1}$ не меньше $\alpha_i+\beta_{i+1}$.
Таким образом, 
\begin{align*}
\tc{p_1\dots p_n}&=\theta_1+\dots+\theta_n\le
\\
&\le\beta_1+\alpha_1+\dots+\beta_n+\alpha_n = 
\\
&=(\alpha_1+\beta_2)+\dots+(\alpha_n+\beta_1) \le 
\\
&\le \tc\gamma.
\end{align*}

Если $\gamma$ незамкнута, то вычисления аналогичны
\begin{align*}
\tc{p_1\dots p_n}&=\theta_2+\dots+\theta_{n-1}\le
\\
&\le\beta_2+\alpha_2+\dots+\beta_{n-1}+\alpha_{n-1} \le
\\
&\le (\alpha_1+\beta_2)+\dots+(\alpha_{n-1}+\beta_n) \le
\\
&\le \tc\gamma.
\end{align*}
\qedsf

\begin{thm}{Упражнение}\label{ex:sef-intersection}\label{ex:sef-intersection:>pi}  
Решите упражнение \ref{ex:length-dist:self-intersection:>pi}, используя \ref{prop:inscribed-total-curvature}.
\end{thm}

\begin{wrapfigure}{r}{30 mm}
\vskip-0mm
\centering
\includegraphics{mppics/pic-20}
\vskip0mm
\end{wrapfigure}

\begin{thm}{Упражнение}\label{ex:quadrisecant}
Замкнутая кривая $\gamma$ пересекает прямую в четырёх точках $a$, $b$, $c$ и~$d$.
Предположим, что эти точки появляются на прямой в порядке $a$, $b$, $c$, $d$, а на кривой $\gamma$ в порядке $a$, $c$, $b$,~$d$.
Покажите, что 
\[\tc\gamma\ge 4\cdot\pi.\]

\end{thm}


\section[\texorpdfstring{А что если $\Phi(\gamma)=2\cdot \pi$?}{А что если Φ(γ)=2·π?}]{А что если $\bm{\Phi(\gamma)=2\cdot \pi}$?}

\begin{thm}{Предложение}\label{prop:fenchel=}
Случай равенства в теореме Фенхеля выполняется только для плоских выпуклых кривых;
то есть полная кривизна гладкой пространственной кривой $\gamma$ равна $2\cdot\pi$ тогда и только тогда, когда $\gamma$ является выпуклой плоской кривой.
\end{thm}

\parbf{Доказательство \ref{prop:fenchel=}.}
Достаточность следует из \ref{cor:fenchel=convex};
остаётся доказать необходимость.

Пусть $abcd$ --- четырёхугольник, вписанный в~$\gamma$.
По определению полной кривизны,
\begin{align*}
\tc{abcd}&=(\pi-
\measuredangle\hinge adb)+(\pi-
\measuredangle\hinge bac)+(\pi-
\measuredangle\hinge cbd)+(\pi-
\measuredangle\hinge dca)=
\\
&=4\cdot\pi -(
\measuredangle\hinge adb
+
\measuredangle\hinge bac
+
\measuredangle\hinge cbd
+
\measuredangle\hinge dca)
\end{align*}

По неравенству треугольника для углов (\ref{thm:spherical-triangle-inq}),
\[
\measuredangle\hinge bac
\le
\measuredangle\hinge bad
+ 
\measuredangle\hinge bdc
\quad\text{и}\quad
\measuredangle\hinge dca\le
\measuredangle\hinge dcb
+ 
\measuredangle\hinge dba.
\eqlbl{eq:spheric-triangle}
\]

\begin{wrapfigure}{r}{30 mm}
\vskip-5mm
\centering
\includegraphics{mppics/pic-56}
\vskip0mm
\end{wrapfigure}

Сумма углов в любом треугольнике равна $\pi$.
Значит из вышесказанного получаем
\begin{align*}
\tc{abcd}\ge 4\cdot \pi 
&- (\measuredangle\hinge adb+\measuredangle\hinge bad+ 
\measuredangle\hinge dba)-
\\
&-(\measuredangle\hinge cbd+\measuredangle\hinge dcb 
+\measuredangle\hinge  bdc)=
\\
=2\cdot\pi.&
\end{align*}

По \ref{prop:inscribed-total-curvature},
\[\tc{abcd}\le \tc\gamma\le 2\cdot\pi.\]
Следовательно, в \ref{eq:spheric-triangle} достигаются равенства.
Это означает, что точка $d$ принадлежит углу $abc$, 
а точка $b$ --- углу $cda$.
Последнее влечёт, что $abcd$ является выпуклым четырёхугольником на плоскости.

{\sloppy

То есть, любой четырёхугольник, вписанный в $\gamma$, плоский и выпуклый.
Таким образом, все точки $\gamma$ лежат в одной плоскости, определяемой тремя точками на~$\gamma$.
А поскольку любой четырёхугольник, вписанный в $\gamma$, выпуклый,
кривая $\gamma$ и сама выпукла. 
\qeds

}

\section{Обобщённая теорема о ДНК}\label{sec:DNA-poly}

\begin{thm}{Теорема}\label{thm:DNA-poly}
Пусть $p_1\dots p_n$ --- замкнутая ломаная в единичном шаре.
Тогда 
\[\tc{p_1\dots p_n}>\length(p_1\dots p_n).\]
\end{thm}

Согласно упражнению \ref{ex:total-curvature=}, эта теорема влечёт гладкий вариант теоремы (\ref{thm:DNA}).
То есть \ref{thm:DNA-poly} обобщает \ref{thm:DNA}.

\parbf{Доказательство.}
Будем предполагать, что $p_n=p_0$, $p_{n+1}=p_1$ и так далее.
Обозначим через $\theta_i$ внешний угол при $p_i$.

\begin{figure}[ht!]
\vskip-0mm
\centering
\includegraphics{mppics/pic-16}
\vskip0mm
\end{figure}

Пусть $o$ --- центр шара.
Рассмотрим последовательность треугольников 
\[\triangle q_0q_1s_0\cong \triangle p_0p_1o,\ \ \triangle q_1q_2s_1\cong \triangle p_1p_2o,\ \dots\]
таких, что точки $q_0,q_1,\dots$ лежат на одной прямой в том же порядке, и все точки $s_i$ лежат с одной стороны от этой прямой.

Заметим, что $s_0s_nq_nq_0$ --- параллелограмм, и, значит,
\[|s_n-s_0|=|q_n-q_0|=\length (p_1\dots p_n).\]
Следовательно, 
\[|s_0-s_1|+\dots+|s_{n-1}-s_n|\ge \length (p_1\dots p_n).\]

Далее, 
\[|q_i-s_{i-1}|=|q_i-s_i|=|p_i-o|\le 1\]
и
\[\theta_i\ge\measuredangle \hinge{q_i}{s_{i-1}}{s_i}\]
для каждого $i$.
Следовательно,
\[\theta_i>|s_{i-1}-s_i|\]
для каждого $i$.
Таким образом,
\begin{align*}
\tc {p_1\dots p_n}
&=\theta_1+\dots+\theta_n>
\\
&> |s_{0}-s_1|+\dots +|s_{n-1}-s_n|\ge 
\\
&\ge\length (p_1\dots p_n).
\end{align*}
\qedsf

Упомянем ещё одно обобщение теоремы о ДНК;
оно было получено Джеффри Лагариасом и Томасом Ричардсоном \cite{lagarias-richardso}; другое доказательство нашли Александр Назаров и Фёдор Петров \cite{nazarov-petrov}.
Оба доказательства до обидного сложны.

\begin{thm}{Теорема}
Предположим, что $\alpha$ --- замкнутая кривая, лежащая в выпуклой фигуре на плоскости, ограниченной кривой $\gamma$.
Тогда средняя кривизна $\alpha$ не меньше средней кривизны $\gamma$.

\end{thm}

\section{Обобщённая кривизна}

Следующее упражнение говорит, что неравенство в \ref{prop:inscribed-total-curvature} оптимально.

\begin{thm}{Упражнение}\label{ex:total-curvature=}
Покажите, что для любой гладкой пространственной кривой $\gamma$ выполнено равенство
\[\tc\gamma=\sup\{\tc\beta\},\]
где точная верхняя грань берётся по всем ломаным $\beta$, вписанным в $\gamma$
(если $\gamma$ замкнута, то предполагаем, что и $\beta$ замкнута).
\end{thm}

Равенство в упражнении используется в определении полной кривизны произвольной кривой~$\gamma$ --- её определяют как \textit{точную верхнюю грань полной кривизны невырожденных ломаных, вписанных в~$\gamma$.}

Эта теория была переоткрыта несколько раз; см. \cite[III §~1]{pogorelov}, \cite{aleksandrov-reshetnyak} и \cite{sullivan-curves}.
Большинство утверждений этой главы можно распространить на кривые конечной полной кривизны в этом обобщённом смысле.

\begin{thm}{Упражнение}\label{ex:tc-rectifiable}
Предположим, что полная кривизна кривой $\gamma\:[0,1]\to\mathbb{R}^3$ ограничена в обобщённом смысле;
то есть существует верхняя грань на полные кривизны ломаных, вписанных в~$\gamma$.

Покажите, что $\gamma$ спрямляема.
Постройте пример, показывающий, что обратное неверно.
\end{thm}

\chapter{Кручение}
\label{chap:torsion}

Эта глава в основном предоставляет практику в вычислениях.
За исключением определений в разделе~\ref{sec:frenet-frame}, она не будет использоваться в дальнейшем.

Кручение во многом напоминает кривизну.
Кривизна показывает, насколько кривая отклоняется от прямой линии, а кручение показывает, насколько пространственная кривая отклоняется от плоской (см. \ref{ex:lancret}).

\section{Базис Френе}\label{sec:frenet-frame}

Пусть $\gamma$ --- гладкая пространственная кривая.
Будем считать, что $\gamma$ параметризована длиной,
так что вектор скорости $\tan(s)\z=\gamma'(s)$ единичный.

Предположим, что $\gamma$ имеет ненулевую кривизну в момент времени~$s$;
другими словами, $\gamma''(s)\z\ne 0$.
Тогда \index{нормальный вектор}\emph{нормальный вектор} при $s$ определяется как
\[\norm(s)=\frac{\gamma''(s)}{|\gamma''(s)|}.\]
Заметим, что \index{10tnb@$\tan$, $\norm$, $\bi$ (базис Френе)}
\[\tan'(s)=\gamma''(s)=\kur(s)\cdot\norm(s).\]

Согласно \ref{prop:a'-pertp-a''}, $\norm(s)\perp \tan(s)$.
Поэтому векторное произведение 
\[\bi(s)=\tan(s)\times \norm(s)\]
является единичным вектором.
Более того, тройка векторов $\tan(s)$, $\norm(s)$, $\bi(s)$ образует ориентированный ортонормированный базис в $\mathbb{R}^3$;
он называется \index{базис Френе}\emph{базис Френе} кривой $\gamma$ при~$s$.
В частности, 
\[\begin{aligned}
\langle\tan,\tan\rangle&=1,
&
\langle\norm,\norm\rangle&=1,
&\langle\bi,\bi\rangle&=1,
\\
\langle\tan,\norm\rangle&=0,
&
\langle\norm,\bi\rangle&=0,
&
\langle\bi,\tan\rangle&=0.
\end{aligned}
\eqlbl{eq:orthornomal}
\]
Векторы $\tan(s)$, $\norm(s)$ и $\bi(s)$ называются \index{касательный вектор}\emph{касательным}, \index{нормальный вектор}\emph{нормальным} и \index{бинормальный вектор}\emph{бинормальным} базиса Френе, соответственно.
Отметим, что базис Френе определён только если $\kur(s)\z\ne 0$.

Плоскость $\Pi_s$ через $\gamma(s)$, натянутая на $\tan(s)$ и $\norm(s)$, называется \index{соприкасающаяся!плоскость}\emph{соприкасающейся плоскостью} при $s$;
она может определяться и как плоскость через $\gamma(s)$, перпендикулярная бинормальному вектору $\bi(s)$.
Это единственная плоскость, имеющая \index{порядок касания}\emph{второй порядок касания} с $\gamma$ при $s$;
то есть $\rho(\ell)=o(\ell^2)$, где $\rho(\ell)$ --- расстояние от $\gamma(s+\ell)$ до~$\Pi_s$.

\section{Кручение}

Пусть $\gamma$ --- гладкая пространственная кривая с единичной скоростью,
и $\tan,\norm,\bi$ --- её базис Френе.
Величина \index{10tau@$\tor$ (кручение)}
\[\tor(s)=\langle \norm'(s),\bi(s)\rangle\]
называется \index{кручение}\emph{кручением} $\gamma$ при~$s$.

Как и базис Френе, кручение $\tor(s_0)$ определено, только если $\kur(s_0)\z\ne0$.
Действительно, поскольку функция $s\mapsto \kur(s)$ является непрерывной, 
$\kur(s_0)\z\ne 0$ означает, что $\kur(s)\z\ne 0$ для всех $s\approx s_0$.
Следовательно, базис Френе определён на открытом интервале, содержащем~$s_0$.
Ясно, что $\tan(s)$, $\norm(s)$ и $\bi(s)$ гладко зависят от $s$ в своих областях определения.
Поэтому производная $\norm'(s_0)$ определена, а, значит, определено и кручение.

\begin{thm}{Упражнение}\label{ex:helix-torsion}
При данных $a$ и $b$,
вычислите кривизну и кручение винтовой линии
$\gamma_{a,b}(t)=(a\cdot \cos t,a\cdot\sin t, b\cdot t)$.

Выведите отсюда, что для данных $\kur>0$ и $\tor$
существует винтовая линия с постоянной кривизной $\kur$ и постоянным кручением~$\tor$.
\end{thm}

\section{Формулы Френе}

Предположим, что базис Френе $\tan(s),\norm(s),\bi(s)$ кривой $\gamma$ определён в точке~$s$.
Напомним, что 
\[\tan'=\kur\cdot \norm.
\eqlbl{eq:frenet-tau}\]
Давайте выразим оставшиеся производные $\norm'$ и $\bi'$ в базисе $\tan,\norm,\bi$.

Сначала покажем, что
\[\norm'=-\kur\cdot\tan+\tor\cdot\bi.\eqlbl{eq:frenet-nu}\]

Действительно, поскольку базис $\tan,\norm,\bi$ ортонормирован, приведённая выше формула эквивалентна следующим трём тождествам:
\[\begin{aligned}
\langle \norm',\tan\rangle&=-\kur,
&
\langle \norm',\norm\rangle&=0,
&
\langle \norm',\bi\rangle&=\tor,
\end{aligned}\eqlbl{eq:<N',?>}\]
Последнее тождество следует из определения кручения.
Второе следует из тождества $\langle \norm,\norm\rangle\z=1$ в \ref{eq:orthornomal}. 
Продифференцировав тождество $\langle\tan,\norm\rangle\z=0$ в \ref{eq:orthornomal}, получим 
\[\langle\tan',\norm\rangle+\langle\tan,\norm'\rangle=0.\]
Применив \ref{eq:frenet-tau}, получим первое тождество в \ref{eq:<N',?>}.

Дифференцируя третье тождество в \ref{eq:orthornomal}, получаем, что $\bi'\perp\bi$.
Продифференцировав остальные два тождества с $\bi$ в \ref{eq:orthornomal}, получим
\begin{align*}
\langle\bi',\tan\rangle&=-\langle\bi,\tan'\rangle=-\kur\cdot\langle\bi,\norm\rangle=0,
\\
\langle\bi,\norm'\rangle&=-\langle\bi',\norm\rangle=\tor.
\end{align*}
И поскольку базис $\tan,\norm,\bi$ ортонормирован, получаем, что
\[\bi'=-\tor\cdot\norm.\eqlbl{eq:frenet-beta}\]

Уравнения \ref{eq:frenet-tau}, \ref{eq:frenet-nu} и \ref{eq:frenet-beta} называются \index{формулы Френе}\emph{формулами Френе}.
Все три можно записать в виде одного матричного тождества:
\[
\begin{pmatrix}
\tan'
\\
\norm'
\\
\bi'
\end{pmatrix}
=
\begin{pmatrix}
0&\kur&0
\\
-\kur&0&\tor
\\
0&-\tor&0
\end{pmatrix}
\cdot
\begin{pmatrix}
\tan
\\
\norm
\\
\bi
\end{pmatrix}.
\]

Напомним, что бинормаль $\bi$ ортогональна соприкасающейся плоскости.
Поэтому уравнение \ref{eq:frenet-beta} означает, что кручение мерит, как быстро крутится соприкасающаяся плоскость при движении вдоль~$\gamma$.

\begin{thm}{Упражнение}\label{ex:beta-from-tau+nu}
Выведите формулу \ref{eq:frenet-beta} из \ref{eq:frenet-tau} и \ref{eq:frenet-nu}, дифференцируя тождество
$\bi=\tan\times \norm$.
\end{thm}

\begin{thm}{Упражнение}\label{ex:torsion=0}
Пусть $\gamma$ --- гладкая пространственная кривая с ненулевой кривизной.
Докажите, что $\gamma$ лежит в плоскости тогда и только тогда, когда её кручение тождественно равно нулю.
\end{thm}

\begin{thm}{Упражнение}\label{ex:+B}
Пусть $\tan,\norm,\bi$ --- базис Френе гладкой кривой
$\gamma_0\:[a,b]\to \mathbb{R}^3$.
Рассмотрим кривую $\gamma_1(t)\z\df\gamma_0(t)\z+\bi(t)$.
Докажите, что
\[\length\gamma_1\ge\length\gamma_0.\]
\end{thm}

\begin{thm}{Упражнение}\label{ex:frenet}
Пусть $\gamma$ --- гладкая пространственная кривая,
$\tan,\norm,\bi$ --- её базис Френе, и $\tor$ --- её кручение.
Выведите равенства
\[\bi=\frac{\gamma'\times\gamma''}{|\gamma'\times\gamma''|}
\quad\text{и}\quad
\tor=\frac{\langle\gamma'\times\gamma'',\gamma'''\rangle}{|\gamma'\times\gamma''|^2}.
\]

\end{thm}

\begin{thm}{Упражнение}\label{ex:moment-curve}
Найдите кривизну $\kur(t)$ и кручение $\tor(t)$ \index{кривая моментов}\emph{кривой моментов} $\gamma\:t\z\mapsto (t,t^2,t^3)$ в точке $\gamma(t)$.
\end{thm}

Следующее упражнение тесно связано с леммой о луке (\ref{lem:bow}).

\begin{thm}{Упражнение}\label{ex:bow-converse}
Пусть $\gamma_1,\gamma_2\:[a,b]\to\mathbb{R}^3$ --- две гладкие кривые с единичной скоростью.
Предположим, что 
\[\dist{\gamma_1(t_1)}{\gamma_1(t_2)}{}\ge \dist{\gamma_2(t_1)}{\gamma_2(t_2)}{}\]
при любых $t_1$ и $t_2$.
Покажите, что $\kur(t_0)_{\gamma_1}\le \kur(t_0)_{\gamma_2}$ для любого $t_0$.
\end{thm}

\begin{thm}{Продвинутое упражнение}\label{ex:torsion-indicatrix}
Пусть $\gamma$ --- замкнутая гладкая пространственная кривая с положительным кручением.
Докажите, что её касательная индикатриса имеет самопересечение.
\end{thm}

\section{Линии откоса}

Гладкая пространственная кривая $\gamma$ называется \index{линия откоса}\emph{линией откоса}, если она идёт под постоянным углом к фиксированному направлению.
Следующую теорему доказал Мишель Анж Ланкре~\cite{lancret}.

\begin{thm}{Теорема}\label{thm:const-slope}
Пусть $\gamma$ --- гладкая кривая с кривизной и кручением $\kur$ и $\tor$.
Предположим, что $\kur(s)>0$ для всех~$s$.
Тогда $\gamma$ --- линия откоса тогда и только тогда, когда отношение $\tfrac\tor\kur$ постоянно.
\end{thm}

\begin{thm}{Доказательство и упражнение}\label{ex:lancret}
{\sloppy
Пусть $\gamma$ --- гладкая пространственная кривая с ненулевой кривизной, $\tan,\norm,\bi$ --- её базис Френе, а $\kur$ и $\tor$ --- её кривизна и кручение.

}

\begin{subthm}{ex:lancret:a}
Предположим, что для некоторого фиксированного ненулевого вектора $\vec w$ значение $\langle \vec w,\tan\rangle$ постоянно.
Докажите, что $\langle \vec w, \norm\rangle \z=0$.
Выведите, что $\langle \vec w,\bi\rangle$ постоянно.
Докажите, что \[\tor\cdot\langle \vec w,\bi\rangle -\kur\cdot\langle \vec w,\tan\rangle =0.\]
Выведите отсюда, что отношение $\tfrac\tor\kur$ постоянно.
\end{subthm}

{\sloppy

\begin{subthm}{ex:lancret:b}
Предположим, что отношение $\tfrac\tor\kur$ не меняется.
Докажите, что вектор $\vec w\z=\tfrac\tor\kur\cdot \tan+\bi$ также не меняется.
Выведите отсюда, что $\gamma$ --- линия откоса.
\end{subthm}

}

\end{thm}

Пусть $\gamma$ --- гладкая кривая с единичной скоростью и $s_0$ --- фиксированное значение её параметра. 
Тогда кривая 
\[\alpha(s)=\gamma(s)+(s_0-s)\cdot \gamma'(s)\]
называется \index{эвольвента}\emph{эвольвентой}~$\gamma$.
Заметим, что если $\ell(s)$ обозначает касательную к $\gamma$ в точке $s$,
то $\alpha(s)\in \ell(s)$ и $\alpha'(s)\perp \ell(s)$ для всех~$s$.

\begin{thm}{Упражнение}\label{ex:evolvent-constant-slope}
Докажите, что эвольвента линии откоса лежит на плоскости.
\end{thm}

\section{Сферические кривые}

\begin{thm}{Теорема}
Пусть $\gamma$ --- гладкая пространственная кривая с ненулевым кручением $\tor$ (и, следовательно, с ненулевой кривизной~$\kur$).
Тогда $\gamma$ лежит на единичной сфере тогда и только тогда, когда выполняется тождество
\[\left|\frac{\kur'}{\tor}\right|=\kur\cdot\sqrt{\kur^2-1}.\]
\end{thm}

\begin{thm}{Доказательство и упражнение}\label{ex:spherical-frenet}
Предположим, что $\gamma$ --- гладкая пространственная кривая с единичной скоростью,
$\tan,\norm,\bi$ --- её базис Френе,
и $\kur$, $\tor$ --- её кривизна и кручение.

\smallskip

Предположим, что $\gamma$ сферическая; то есть $|\gamma(s)|=1$ для любого~$s$.
Покажите, что

\begin{subthm}{ex:spherical-frenet:tau} $\langle\tan,\gamma\rangle=0$ и $\langle\norm,\gamma\rangle^2+\langle\bi,\gamma\rangle^2=1$.
\end{subthm}

\begin{subthm}{ex:spherical-frenet:nu} $\langle\norm,\gamma\rangle=-\tfrac1\kur$.
\end{subthm}

\begin{subthm}{ex:spherical-frenet:beta} $\langle\bi,\gamma\rangle'=\tfrac\tor\kur$.
\end{subthm}

\begin{subthm}{ex:spherical-frenet:beta+}
Воспользовавшись \ref{SHORT.ex:spherical-frenet:beta}, докажите, что если $\gamma$ замкнута, то $\tor(s)=0$ для некоторого~$s$.
\end{subthm}

\begin{subthm}{ex:spherical-frenet:kur-tor}
Предположим, что кручение $\gamma$ не равно нулю.
Используйте \ref{SHORT.ex:spherical-frenet:tau}--\ref{SHORT.ex:spherical-frenet:beta}, чтобы показать, что
\[\left|\frac{\kur'}{\tor}\right|=\kur\cdot\sqrt{\kur^2-1}.\]
\end{subthm}
Пусть теперь $\gamma$ удовлетворяет тождеству в \ref{SHORT.ex:spherical-frenet:kur-tor}.
\begin{subthm}{ex:spherical-frenet:f}
Покажите, что точка $\gamma+\tfrac1\kur\cdot \norm+\tfrac{\kur'}{\kur^2\cdot\tor}\cdot\bi$ не меняется.
Выведите, что $\gamma$ лежит на единичной сфере с центром в этой точке.
\end{subthm}

\end{thm}

Для кривой $\gamma$ с единичной скоростью и ненулевой кривизной и кручением при~$s$,
сфера $\Sigma_s$, проходящая через $\gamma(s)$ с центром в
\[p(s)=\gamma(s)+\tfrac1{\kur(s)}\cdot \norm(s)+\tfrac{\kur'(s)}{\kur^2(s)\cdot\tor(s)}\cdot\bi(s)\]
называется \index{соприкасающаяся!сфера}\emph{соприкасающейся сферой} кривой $\gamma$ при~$s$.
Это единственная сфера, которая имеет \index{порядок касания}\emph{касание третьего порядка} с $\gamma$ при~$s$;
то есть $\rho(\ell)=o(\ell^3)$, где $\rho(\ell)$ обозначает расстояние от $\gamma(s+\ell)$ до $\Sigma_s$.
 
\section{Основная теорема}

\begin{thm}{Теорема}\label{thm:fund-curves}
Пусть $s\mapsto \kur(s)$ и $s\mapsto \tor(s)$ --- две гладкие функции, определённые на интервале $\mathbb{I}$.
Предположим, что $\kur(s)>0$ для всех~$s$.
Тогда существует гладкая кривая $\gamma\:\mathbb{I}\to\mathbb{R}^3$ с единичной скоростью, кривизной $\kur(s)$ и кручением $\tor(s)$ при любом~$s\in \mathbb{I}$.
Более того, $\gamma$ однозначно определена с точностью до изометрии пространства, сохраняющей ориентацию.
\end{thm}

Доказательство основывается на теореме о существовании и единственности решений дифференциальных уравнений (\ref{thm:ODE}).

\parbf{Доказательство.}
Выберем значение параметра $s_0$, точку $\gamma(s_0)$ и ориентированный ортонормированный базис $\tan(s_0)$, $\norm(s_0)$, $\bi(s_0)$.
Рассмотрим следующую систему дифференциальных уравнений
\[
\begin{cases}
\gamma'=\tan,
\\
\tan'=\kur\cdot\norm,
\\
\norm'=-\kur\cdot\tan+\tor\cdot\bi,
\\
\bi'=-\tor\cdot\norm
\end{cases}
\eqlbl{eq:gamma'tan'norm'bi'}
\]
с начальными данными $\gamma(s_0)$, $\tan(s_0)$, $\norm(s_0)$, $\bi(s_0)$.
(Наша система уравнений состоит из четырёх векторных уравнений, так что её можно переписать как систему из 12-и скалярных.)

Согласно \ref{thm:ODE}, эта система имеет единственное решение, которое определено на максимальном подинтервале $\mathbb{J}\subset \mathbb{I}$, содержащем $s_0$.
Покажем, что $\mathbb{J}=\mathbb{I}$.

Сначала заметим, что 
\[\begin{aligned}
\langle\tan,\tan\rangle&=1,
&
\langle\norm,\norm\rangle&=1,
&
\langle\bi,\bi\rangle&=1,
\\
\langle\tan,\norm\rangle&=0,&
\langle\norm,\bi\rangle&=0,&
\langle\bi,\tan\rangle&=0
\end{aligned}
\eqlbl{eq:111000}
\]
для любого значения параметра $s$.

Действительно, из \ref{eq:gamma'tan'norm'bi'} получаем следующую систему:
\[
\begin{cases}
\langle\tan,\tan\rangle'
&=
2\cdot\langle\tan,\tan'\rangle
=
2\cdot\kur\cdot \langle\tan,\norm\rangle,
\\
\langle\norm,\norm\rangle'
&=
2\cdot\langle\norm,\norm'\rangle
=
-
2\cdot\kur\cdot\langle\norm,\tan\rangle
+
2\cdot\tor\cdot\langle\norm,\bi\rangle,
\\
\langle\bi,\bi\rangle'
&=
2\cdot\langle\bi,\bi'\rangle
=
-2\cdot\tor\langle\bi,\norm\rangle,
\\
\langle\tan,\norm\rangle'
&=
\langle\tan',\norm\rangle
+
\langle\tan,\norm'\rangle
=
\kur\cdot\langle\norm,\norm\rangle
-
\kur\cdot\langle\tan,\tan\rangle
+
\tor\cdot\langle\tan,\bi\rangle,
\\
\langle\norm,\bi\rangle'
&=
\langle\norm',\bi\rangle+\langle\norm,\bi'\rangle
=\kur\cdot\langle\tan,\bi\rangle+\tor\cdot\langle\bi,\bi\rangle-\tor\cdot\langle\norm,\norm\rangle,
\\
\langle\bi,\tan\rangle'
&=
\langle\bi',\tan\rangle+\langle\bi,\tan'\rangle
=
-\tor\cdot \langle\norm,\tan\rangle
+\kur\cdot\langle\bi,\norm\rangle.
\end{cases}
\eqlbl{eq:<gamma'tan'norm'bi'>}
\]

Константы в \ref{eq:111000} удовлетворяют этой системе.
Более того, поскольку базис $\tan(s_0)$, $\norm(s_0)$, $\bi(s_0)$ ориентирован и ортонормирован,
\ref{eq:111000} решает нашу задачу Коши для системы \ref{eq:<gamma'tan'norm'bi'>}.

Допустим, что $\mathbb{J} \varsubsetneq \mathbb{I}$.
Тогда один из концов $\mathbb{J}$, скажем $b$, лежит во внутренней части $\mathbb{I}$.
Теорема \ref{thm:ODE} применима для $\Omega=\mathbb{R}^{12}\times \mathbb{I}$.
Значит, одно из значений $\gamma(s)$, $\tan(s)$, $\norm(s)$, $\bi(s)$
устремляется к бесконечности при $s\to b$.
Но это невозможно --- векторы $\tan(s)$, $\norm(s)$, $\bi(s)$ остаются единичными, и $|\gamma'(s)|=|\tan(s)|=1$;
так что $\gamma$ проходит лишь конечное расстояние при $s\to b$.
Таким образом, $\mathbb{J}= \mathbb{I}$, и первая часть теоремы доказана.

Теперь предположим, что есть две кривые $\gamma_1$ и $\gamma_2$ с заданными функциями кривизны и кручения.
Применив изометрию пространства, можно считать, что $\gamma_1(s_0)=\gamma_2(s_0)$ и базисы Френе кривых совпадают при $s_0$.
Тогда $\gamma_1=\gamma_2$ по единственности решений системы (\ref{thm:ODE}), что завершает доказательство.
\qeds

\begin{thm}{Упражнение}\label{ex:cur+tor=helix}
Предположим, что кривая $\gamma\:\mathbb{R}\to\mathbb{R}^3$ с постоянной скоростью, кривизной и кручением.
Покажите, что $\gamma$ является винтовой линией, возможно, вырождающейся в окружность;
то есть в подходящей системе координат 
$\gamma(t)=(a\cdot \cos t,a\cdot\sin t, b\cdot t)$
для некоторых констант $a$ и~$b$.
\end{thm}

\begin{thm}{Продвинутое упражнение}\label{ex:const-dist}
Пусть $\gamma$ --- гладкая пространственная кривая такая, что расстояние $|\gamma(t)-\gamma(t+\ell)|$ зависит только от $\ell$.
Покажите, что $\gamma$ является винтовой линией, возможно, вырождающейся в прямую или окружность.
\end{thm}

\chapter{Кривизна со знаком}\label{chap:signed-curvature}

Кривая на плоскости может поворачивать влево или вправо.
Это позволяет ввести знак кривизны плоских кривых.
Если вести машину по плоскости вдоль кривой, то ориентированная кривизна будет задавать положение руля.

\section{Определения}\label{sec:def(skur)}

Пусть $\gamma$ --- гладкая плоская кривая с единичной скоростью.
Как обычно, $\tan(s)=\gamma'(s)$ --- её единичный касательный вектор при~$s$.

Повернём $\tan(s)$ на угол $\tfrac\pi 2$ против часовой стрелки; 
обозначим полученный вектор через $\norm(s)$.
Пара $\tan(s),\norm(s)$ образует ориентированный ортонормированный базис плоскости.
Он аналогичен базису Френе (см. \ref{sec:frenet-frame}) и также будет назваться \index{базис Френе}\emph{базисом Френе}.

Напомним, что $\gamma''(s)\perp \gamma'(s)$ (см. \ref{prop:a'-pertp-a''}).
Следовательно, \index{10k@$\skur$ (ориентированная кривизна)}
\[\tan'(s)=\skur(s)\cdot \norm(s).\eqlbl{eq:tau'}\]
для некоторого числа $\skur(s)$.
Величина $\skur(s)$ называется \index{ориентированная кривизна}\index{кривизна! ориентированная}\emph{ориентированной кривизной} $\gamma$ при~$s$;
её также можно обозначать через $\skur(s)_\gamma$.

Заметим, что 
\[\kur(s)=|\skur(s)|;\]
то есть, с точностью до знака, ориентированная кривизна $\skur(s)$ совпадает с обычной кривизной $\kur(s)$ (см. \ref{sec:curvature}).
Знак ориентированной кривизны указывает направление поворота: $\skur (s)>0$, если кривая уходит влево.
Если захочется подчеркнуть, что мы работаем с обычной кривизной, 
то можно сказать \index{кривизна!без знака}\emph{кривизна без знака}.

Ориентированная кривизна меняет знак как при обращении параметризации кривой, так и при обращении ориентации плоскости.

Поскольку $\tan(s),\norm(s)$ --- ортонормированный базис, 
\begin{align*}
\langle\tan,\tan\rangle&=1,
&
\langle\norm,\norm\rangle&=1, 
&
\langle\tan,\norm\rangle&=0.
\end{align*}
Продифференцировав эти тождества, получим 
\begin{align*}
\langle\tan',\tan\rangle&=0,
&
\langle\norm',\norm\rangle&=0,
&
\langle\tan',\norm\rangle+\langle\tan,\norm'\rangle&=0.
\end{align*}
Согласно \ref{eq:tau'}, $\langle\tan',\norm\rangle=\skur$. 
Следовательно, $\langle\tan,\norm'\rangle=-\skur$, и 
\[\norm'(s)=-\skur(s)\cdot \tan(s).\eqlbl{eq:nu'}\]
Уравнения \ref{eq:tau'} и \ref{eq:nu'} называются \index{формулы Френе}\emph{формулами Френе} для плоских кривых. 
Их можно записать одним матричным уравнением:
\[
\begin{pmatrix}
\tan'
\\
\norm'
\end{pmatrix}
=
\begin{pmatrix}
0&\skur
\\
-\skur&0
\end{pmatrix}
\cdot
\begin{pmatrix}
\tan
\\
\norm
\end{pmatrix}.
\]

\begin{thm}{Упражнение}\label{ex:bike}
Пусть $\tan$ --- касательная индикатриса гладкой кривой $\gamma_0\:[a,b]\to\mathbb{R}^2$.
Рассмотрим кривую $\gamma_1\:[a,b]\to\mathbb{R}^2$, определённую как $\gamma_1(t)\z\df\gamma_0(t)+\tan(t)$.
Докажите, что

\begin{minipage}{.47\textwidth}
\begin{subthm}{ex:bike:length}$\length\gamma_0\le \length\gamma_1$;
\end{subthm}
\end{minipage}
\hfill
\begin{minipage}{.47\textwidth}
\begin{subthm}{ex:bike:tc}$\tc{\gamma_0}\le \length\gamma_1$.
\end{subthm}
\end{minipage}

\end{thm}

Кривые $\gamma_0$ и $\gamma_1$ выше описывают следы шин идеализированного велосипеда с единичным расстоянием между колёсами.
Из упражнения следует, что след  переднего колеса обычно длинней.
Больше по теме можно найти в обзоре Роберта Фута, Марка Леви и Сергея Табачникова \cite{foote-levi-tabachnikov}.

\section{Основная теорема}

\begin{thm}{Теорема}\label{thm:fund-curves-2D}
Для любой гладкой функции $s\mapsto \skur(s)$, определённой на интервале $\mathbb{I}$,
существует гладкая кривая $\gamma\:\mathbb{I}\to\mathbb{R}^2$ с единичной скоростью и ориентированной кривизной $\skur(s)$ при любом $s$.
Более того, $\gamma$ определяется однозначно с точностью до изометрии плоскости, сохраняющей ориентацию.
\end{thm}

Теорема является частным случаем своего трёхмерного аналога (\ref{thm:fund-curves}), но мы проведём прямое доказательство.

\parbf{Доказательство.}
Выберем $s_0\in\mathbb{I}$.
Рассмотрим функцию
\[\theta(s)
=
\int_{s_0}^s\skur(t)\cdot dt.\]
Согласно формуле Ньютона --- Лейбница, $\theta'(s)\z=\skur(s)$ для всех~$s$.
Пусть 
$\tan(s)=(\,\cos(\theta(s)),\,\sin(\theta(s))\,)$,
и $\norm(s)=(\,-\sin(\theta(s)),\,\cos(\theta(s))\,)$;
то есть $\norm(s)$ --- поворот $\tan(s)$ на угол $\tfrac\pi2$ против часовой стрелки.
Рассмотрим кривую 
\[\gamma(s)=\int_{s_0}^s\tan(s)\cdot ds.\]
Так как $|\gamma'|=|\tan|=1$, кривая $\gamma$ имеет единичную скорость, и $\tan,\norm$ --- её базис Френе. 

Заметим, что
\begin{align*}
\gamma''(s)&=\tan'(s)=
\\
&=(\,\cos(\theta(s))',\,\sin(\theta(s))'\,)=
\\
&=\theta'(s)\cdot (\,-\sin(\theta(s)),\,\cos(\theta(s))\,)=
\\
&=\skur(s)\cdot \norm(s).
\end{align*}
То есть, ориентированная кривизна кривой $\gamma$ при~$s$ равна $\skur(s)$. 

Существование доказано; осталась единственность.

Предположим, что $\gamma_1$ и $\gamma_2$ --- две кривые, которые удовлетворяют условиям теоремы.
Применив изометрию плоскости, можно добиться, чтобы совпали точки $\gamma_1(s_0)$ и $\gamma_2(s_0)$, а также базисы Френе обеих кривых при $s_0$.
Пусть $\tan_1,\norm_1$ и $\tan_2,\norm_2$ --- базисы Френе кривых $\gamma_1$ и $\gamma_2$, соответственно.
Обе тройки $\gamma_i,\tan_i,\norm_i$ удовлетворяют следующей системе обыкновенных дифференциальных уравнений 
\[
\begin{cases}
\gamma_i'=\tan_i,
\\
\tan_i'=\skur\cdot\norm_i,
\\
\norm_i'=-\skur\cdot\tan_i.
\end{cases}
\]
Более того, у обеих троек те же начальные значения в $s_0$,
а значит, $\gamma_1=\gamma_2$ по единственности решений задачи Коши (\ref{thm:ODE}).
\qeds

Пусть $\gamma\:\mathbb{I}\z\to\mathbb{R}^2$ --- кривая с единичной скоростью.
Функция $\theta\:\mathbb{I}\z\to\mathbb{R}$ называется \index{непрерывный аргумент}\emph{непрерывным аргументом} $\gamma$, если она непрерывна,~и
\[\gamma'(s)=(\,\cos (\theta(s)),\,\sin(\theta(s))\,)\]
для любого $s$.
Из доказательства теоремы видно следующее.

\begin{thm}{Следствие}\label{cor:2D-angle}
Для любой гладкой кривой с единичной скоростью $\gamma\:\mathbb{I}\to\mathbb{R}^2$ существует непрерывный аргумент $\theta\:\mathbb{I}\to\mathbb{R}$.
Более того,
\[\theta'(s)=\skur(s),\]
где $\skur$ --- ориентированная кривизна~$\gamma$.
\end{thm}

\section{Полная ориентированная кривизна}\label{sec:Total signed curvature}

Пусть $\gamma\:\mathbb{I}\to\mathbb{R}^2$ --- гладкая плоская кривая с единичной скоростью.
\index{кривизна}\index{полная!ориентированная кривизна}\emph{Полная ориентированная кривизна} $\gamma$ будет обозначаться $\tgc\gamma$; она определяется как интеграл \index{10psi@$\tgc\gamma$ (полная ориентированная кривизна)}
\[\tgc\gamma
=
\int_\mathbb{I} \skur(s)\cdot ds,\eqlbl{eq:tsc-k}\]
где $\skur$ --- ориентированная кривизна~$\gamma$.

Если $\mathbb{I}=[a,b]$, то 
\[\tgc\gamma=\theta(b)-\theta(a),\eqlbl{eq:tsc-theta}\]
где $\theta$ --- непрерывный аргумент $\gamma$ (см. \ref{cor:2D-angle}).

Если $\gamma$ --- кусочно-гладкая, то её полная ориентированная кривизна определяется как сумма полных ориентированных кривизн её дуг плюс сумма \textit{ориентированных} внешних углов на стыках;
эти углы положительны, когда $\gamma$ поворачивает влево, отрицательны, когда $\gamma$ поворачивает вправо, и равны $0$, когда $\gamma$ идёт прямо.
Ориентированный угол не определён в точке возврата.

Другими словами, если $\gamma$ --- произведение гладких дуг $\gamma_1,\dots,\gamma_n$, то 
\[\tgc\gamma=\tgc{\gamma_1}+\dots+\tgc{\gamma_n}+\theta_1+\dots+\theta_{n-1},\]
где $\theta_i$ --- внешний ориентированный угол на стыке между $\gamma_i$ и $\gamma_{i+1}$.
Если $\gamma$ замкнута, то произведение циклическое, и
\[\tgc\gamma=\tgc{\gamma_1}+\dots+\tgc{\gamma_n}+\theta_1+\dots+\theta_{n},\]
где $\theta_n$ --- внешний ориентированный угол на стыке между $\gamma_n$ и $\gamma_1$.

Поскольку $\left|\int \skur(s)\cdot ds\right|\le \int|\skur(s)|\cdot ds$, неравенство
\[|\tgc\gamma|\le \tc\gamma\eqlbl{eq:tsc-tc}\] 
выполняется для любой гладкой кривой $\gamma$ на плоскости;
то есть полная ориентированная кривизна $\tgc{}$ не превышает полную кривизну $\tc{}$ по абсолютной величине.
При этом равенство достигается тогда и только тогда, когда ориентированная кривизна не меняет знак.

\begin{thm}{Упражнение}\label{ex:trochoids}
Трохоида --- это кривая, описываемая точкой, закреплённой на колесе, катящемся по прямой линии.
\begin{figure}[!ht]
\centering
\begin{lpic}[t(-0mm),b(0mm),r(0mm),l(0mm)]{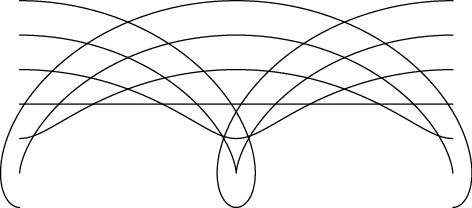}
\lbl[l]{4,0;{\footnotesize $-\tfrac32$}}
\lbl[l]{4,6;{\footnotesize $-1$}}
\lbl[tl]{10,15;{\footnotesize $-\tfrac12$}}
\lbl[t]{22,17;{\footnotesize $0$}}
\lbl[r]{3,23.6;{\footnotesize $\tfrac12$}}
\lbl[r]{3,29.4;{\footnotesize $1$}}
\lbl[r]{3,35.2;{\footnotesize $\tfrac32$}}
\end{lpic}
\end{figure}
Семейство \index{трохоида}\emph{трохоид} $\gamma_a\:[0,2\cdot\pi]\to \mathbb{R}^2$ (см. рисунок) можно параметризовать как
\[\gamma_a(t)=(t+a\cdot \sin t, a\cdot \cos t).\]
\begin{enumerate}[(a)]
\item Для данного $a\in \mathbb{R}$ найдите $\tgc{\gamma_a}$, если она определена.
\item Для данного $a\in \mathbb{R}$ найдите $\tc{\gamma_a}$.
\end{enumerate}
\end{thm}

\begin{thm}{Предложение}\label{prop:total-signed-curvature}
Полная ориентированная кривизна любой простой замкнутой гладкой плоской кривой $\gamma$ равна $\pm2\cdot\pi$; это значение равно $+2\cdot\pi$, если область, ограниченная $\gamma$, находится слева от неё, и $-2\cdot\pi$ в противном случае.

Более того, то же выполняется и для любой простой замкнутой кусочно-гладкой кривой $\gamma$ на плоскости, если её полная ориентированная кривизна определена.
\end{thm}

Предложение даёт дифференциально-геометрический аналог теоремы о сумме внутренних углов многоугольника (\ref{thm:sum=(n-2)pi}), которую мы будем использовать в доказательстве.
Более концептуальное доказательство нашёл Хайнц Хопф \cite{hopf1935}, \cite[с. 42]{hopf1989}.

\parbf{Доказательство.}
Не умаляя общности, можно предположить, что $\gamma$ ориентирована так, что область, ограниченная $\gamma$, находится слева от неё.
Также можно предположить, что $\gamma$ имеет единичную скорость.

Рассмотрим замкнутую ломаную $p_1\dots p_n$, вписанную в~$\gamma$.
Можно считать, что дуги между вершинами настолько малы, что ломаная проста, и каждая дуга $\gamma_i$ от $p_i$ до $p_{i+1}$ имеет малую полную кривизну; скажем, $\tc{\gamma_i}<\pi$ при любом~$i$.

{

\begin{wrapfigure}[13]{o}{41 mm}
\vskip-4mm
\centering
\includegraphics{mppics/pic-59}
\vskip0mm
\end{wrapfigure}

Пусть $p_i=\gamma(t_i)$;
как обычно, будем считать индексы по модулю $n$; в частности, $p_{n+1}\z=p_1$.
Тогда
\begin{align*}
\vec w_i&=p_{i+1}-p_i,& \vec v_i&=\gamma'(t_i),
\\
\alpha_i&=\measuredangle(\vec v_i,\vec w_i),&\beta_i&=\measuredangle(\vec w_{i-1},\vec v_i),
\end{align*}
здесь $\alpha_i,\beta_i\in(-\pi,\pi)$ --- ориентированные углы: $\alpha_i$ положительный, если $\vec w_i$ направлен левее~$\vec v_i$.

}

Согласно \ref{eq:tsc-theta}, величина
\[\tgc{\gamma_i}-\alpha_i-\beta_{i+1}\eqlbl{eq:Psi-alpha-beta}\]
 кратна $2\cdot\pi$.
Поскольку $\tc{\gamma_i}<\pi$, из леммы о хорде (\ref{lem:chord}) получаем, что $|\alpha_i|\z+|\beta_{i+1}|<\pi$.
Из \ref{eq:tsc-tc}, получаем $|\tgc{\gamma_i}|\z\le\tc{\gamma_i}$;
следовательно, величина в \ref{eq:Psi-alpha-beta} равна нулю.
То есть
\[\tgc{\gamma_i}=\alpha_i+\beta_{i+1}\]
для каждого $i$ имеем

Заметим, что 
\[\delta_i=\pi-\alpha_i-\beta_i\eqlbl{eq:delta=pi-alpha-beta}\] 
--- внутренний угол ломаной при $p_i$;
$\delta_i\in (0,2\cdot\pi)$ для каждого~$i$.
Напомним, что сумма внутренних углов $n$-угольника равна $(n-2)\cdot \pi$ (\ref{thm:sum=(n-2)pi}); то есть
\[\delta_1+\dots+\delta_n=(n-2)\cdot \pi.\]
Следовательно, 
\[
\begin{aligned}
\tgc\gamma&=\tgc{\gamma_1}+\dots+\tgc{\gamma_n}=
\\
&=(\alpha_1+\beta_2)+\dots+(\alpha_n+\beta_1)=
\\
&=(\beta_1+\alpha_1)+\dots+(\beta_n+\alpha_n)=
\\
&=(\pi-\delta_1)+\dots+(\pi-\delta_n)=
\\
&=n\cdot\pi-(n-2)\cdot \pi=
\\
&=2\cdot\pi.
\end{aligned}\eqlbl{eq:delta=pi-alpha-beta-sum}\]

{\sloppy

Случай кусочно-гладких кривых аналогичен;
потребуется дополнительно подразбить дуги циклического произведения чтобы выполнялось вышеуказанное условие.
Вместо \ref{eq:delta=pi-alpha-beta}, придётся пользоваться уравнением
\[\delta_i=\pi-\alpha_i-\beta_i-\theta_i,\]
где $\theta_i$ --- ориентированный внешний угол кривой $\gamma$ при $p_i$;
он равен нулю, если кривая $\gamma$ гладкая в $p_i$.
Вычисления в \ref{eq:delta=pi-alpha-beta-sum} превращаются в следующее:
\begin{align*}
\tgc\gamma&=\tgc{\gamma_1}+\dots+\tgc{\gamma_n}+\theta_1+\dots+\theta_n=
\\
&=(\alpha_1+\beta_2)+\dots+(\alpha_n+\beta_1)+\theta_1+\dots+\theta_n=
\\
&=(\beta_1+\alpha_1+\theta_1)+\dots+(\beta_n+\alpha_n+\theta_n)=
\\
&=(\pi-\delta_1)+\dots+(\pi-\delta_n)=
\\
&=n\cdot\pi-(n-2)\cdot \pi=
\\
&=2\cdot\pi.
\end{align*}
\qedsf

}

\begin{thm}{Упражнение}\label{ex:zero-tsc}
Нарисуйте такую гладкую замкнутую кривую $\gamma$ на плоскости, что

\begin{subthm}{ex:zero-tsc:0}
$\tgc\gamma=0$;
\end{subthm}
 
\begin{subthm}{ex:zero-tsc:5}
$\tgc\gamma=\tc\gamma=10\cdot\pi$;
\end{subthm}

\begin{subthm}{ex:zero-tsc:2-4}
$\tgc\gamma=2\cdot\pi$ и $\tc\gamma=4\cdot\pi$.
\end{subthm}

\end{thm}

\begin{thm}{Упражнение}\label{ex:length'}
Пусть $\gamma\:[a,b]\to\mathbb{R}^2$ --- гладкая кривая на плоскости, и $\tan,\norm$ --- её базис Френе.
Для данного параметра $\ell$, рассмотрим
кривую $\gamma_\ell(t)\z=\gamma(t)\z+\ell\cdot\norm(t)$; она называется кривой, \index{параллельная кривая}\emph{параллельной} к $\gamma$ с параметром $\ell$.

\begin{subthm}{ex:length':reg}
Покажите, что параметризация $\gamma_\ell$ регулярна, если $\ell\cdot \skur(t)_\gamma\ne 1$ для всех $t$.
\end{subthm}
 
\begin{subthm}{ex:length':formula}
Рассмотрим функцию $L(\ell)=\length\gamma_\ell$.
Покажите, что 
\[L(\ell)=L(0)-\ell\cdot\tgc\gamma\eqlbl{eq:length(parallel-curve)}\]
при всех $\ell$, достаточно близких к $0$. 
\end{subthm}

\begin{subthm}{ex:length':antiformula}
Приведите пример, показывающий, что формула \ref{eq:length(parallel-curve)} может не выполняться при некоторых значениях~$\ell$.
\end{subthm}

\end{thm}

\section{Соприкасающаяся окружность}

\begin{thm}{Предложение}\label{prop:circle}
Для данной точки $p\in\mathbb{R}^2$,
единичного вектора $\tan$ 
и вещественного числа $\skur$ существует единственная гладкая кривая $\sigma\:\mathbb{R}\to\mathbb{R}^2$ с единичной скоростью, которая начинается в $p$ в направлении $\tan$ и имеет постоянную ориентированную кривизну $\skur$.

Более того, если $\skur=0$, то это прямая $\sigma(s)=p+s\cdot \tan$;
если $\skur\ne 0$, то $\sigma$ описывает окружность радиуса $\tfrac1{|\skur|}$ с центром в $p+\tfrac1\skur\cdot \norm$, где $\tan,\norm$ --- это ориентированный ортонормальный базис.
\end{thm}

\parbf{Доказательство.}
Выберем систему координат с началом отсчёта в $p$, чтобы $\tan$ шёл по оси $x$.
В этом случае $\norm$ идёт по оси $y$.

Заметим, что
\begin{align*}\theta(s)&=\int_{0}^s\skur\cdot dt=\skur\cdot s
\end{align*}
является непрерывным аргументом $\sigma$, см. \ref{cor:2D-angle}.
Следовательно,
\[\sigma'(s)=(\,\cos(\skur\cdot s),\,\sin(\skur\cdot s)\,).\]
Остаётся проинтегрировать последнее равенство.
Если $\skur=0$, то получаем $\sigma(s)=(s,0)$,
что описывает прямую $\sigma(s)=p+s\cdot \tan$.

Если $\skur\ne 0$, то получаем
\[\sigma(s)=(\,\tfrac1\skur\cdot\sin(\skur\cdot s),\, \tfrac1\skur\cdot(1-\cos(\skur\cdot s))\,);\]
это окружность радиуса $r=\tfrac1{|\skur|}$ и центром в $(0,\tfrac1\skur)=p+\tfrac1\skur\cdot\norm$.
\qeds

\begin{thm}{Определение}
Пусть $\gamma$ --- гладкая плоская кривая с единичной скоростью,
и $\skur(s)$ --- её ориентированная кривизна при~$s$.

Кривая $\sigma_s$ с единичной скоростью постоянной ориентированной кривизны $\skur(s)$, которая начинается в $\gamma(s)$ в направлении $\gamma'(s)$, называется \index{соприкасающаяся!окружность}\emph{соприкасающейся окружностью} кривой $\gamma$ при~$s$.

Центр и радиус соприкасающейся окружности называются \index{центр кривизны}\emph{центром кривизны} и \index{радиус кривизны}\emph{радиусом кривизны} кривой в соответственной точке.
\end{thm}

{

\begin{wrapfigure}{o}{31 mm}
\vskip-0mm
\centering
\includegraphics{mppics/pic-21}
\vskip0mm
\end{wrapfigure}

\textit{Соприкасающаяся окружность} может оказаться окружностью или прямой;
в последнем случае центр кривизны не определён, а радиус кривизны можно считать бесконечным.  

Соприкасающаяся окружность $\sigma_s$ также может быть определена как единственная окружность (или прямая), имеющая \index{порядок касания}\emph{второй порядок касания} с $\gamma$ при $s$;
то есть $\rho(\ell)\z=o(\ell^2)$, где $\rho(\ell)$ обозначает расстояние от $\gamma(s+\ell)$ до $\sigma_s$.

Следующее упражнение рекомендуется читателю, знакомому с понятием \index{инверсия}\emph{инверсии}.

}

\begin{thm}{Продвинутое упражнение}\label{ex:inverse}
Пусть $\gamma$ --- гладкая кривая на плоскости, которая не проходит через начало координат и $\hat \gamma$ её инверсия относительно окружности с центром в начале координат.
Покажите, что инверсия переводит соприкасающуюся окружность к $\gamma$ при~$s$ в соприкасающуюся окружность к $\hat \gamma$ при~$s$.
\end{thm}

\section{Лемма о спирали}
\label{spiral}
\index{лемма о спирали}

{

\begin{wrapfigure}{r}{31 mm}
\vskip-14mm
\includegraphics{mppics/pic-61}
\end{wrapfigure}

Следующая лемма была доказана Питером Тейтом \cite{tait}
и переоткрыта Адольфом Кнезером \cite{kneser}.

\begin{thm}{Лемма}\label{lem:spiral}
Пусть $\gamma$ --- гладкая плоская кривая со строго убывающей положительной кривизной.
Тогда соприкасающиеся к $\gamma$ окружности образуют монотонное семейство;
то есть, если $\sigma_s$ --- соприкасающаяся к $\gamma$ окружность при $s$,
то для любых $s_0<s_1$ окружность $\sigma_{s_0}$ лежит внутри открытого круга, ограниченного $\sigma_{s_1}$.\index{10sigma@$\sigma_s$ (соприкасающаяся окружность)}
\end{thm}

}

На рисунке изображена кривая с её соприкасающимися окружностями.
Они заполняют кольцо и делают это довольно причудливым образом: если гладкая функция постоянна на каждой окружности, то она постоянна во всём кольце \cite[Лекция 10]{fuchs-tabachnikov}.
Также отметим, что кривая $\gamma$ касается окружности из этого семейства в каждой из своих точек.
Однако, она не идёт ни по одной из этих окружностей.

\parbf{Доказательство.}
Пусть $\tan(s),\norm(s)$ --- базис Френе.
Обозначим через $\omega(s)$ и $r(s)$
центр и радиус кривизны~$\gamma$;
согласно \ref{prop:circle},
\[\omega(s)=\gamma(s)+r(s)\cdot \norm(s).\]

Поскольку кривизна положительна, $r(s)\cdot\skur(s)=1$.
Применив формулу Френе \ref{eq:nu'}, получим
\begin{align*}
\omega'(s)&=\gamma'(s)+r'(s)\cdot \norm(s)+r(s)\cdot \norm'(s)=
\\
&=\tan(s)+r'(s)\cdot \norm(s)-r(s)\cdot \skur(s)\cdot \tan(s)=
\\
&=r'(s)\cdot \norm(s).
\end{align*}
Так как $\skur(s)$ убывает, $r(s)$ возрастает;
следовательно, $r'\ge 0$.
Из этого вытекает, что $|\omega'(s)|= r'(s)$ и $\omega'(s)$ направлена вдоль $\norm(s)$.

Поскольку $\norm'(s)=-\skur(s)\cdot\tan(s)$, направление $\omega'(s)$ не может оставаться постоянным на любом интервале;
другими словами, кривая $s\mapsto \omega(s)$ не содержит отрезков прямых.

\begin{wrapfigure}[7]{o}{35 mm}
\vskip-0mm
\centering
\includegraphics{mppics/pic-84}
\end{wrapfigure}

В частности, 
\[|\omega(s_1)-\omega(s_0)|\z<\length(\omega|_{[s_0,s_1]})\]
при любых $s_0<s_1$.
Следовательно, 
\begin{align*}
|\omega(s_1)-\omega(s_0)|&<\length(\omega|_{[s_0,s_1]})=
\\
&=\int_{s_0}^{s_1}|\omega'(s)|\cdot ds=
\\
&=\int_{s_0}^{s_1}r'(s)\cdot ds=
\\
&=r(s_1)-r(s_0).
\end{align*}
Другими словами, расстояние между центрами $\sigma_{s_1}$ и $\sigma_{s_0}$
строго меньше разницы между их радиусами, отсюда лемма.
\qeds

Кривая $s\mapsto \omega(s)$ называется \index{эволюта}\emph{эволютой} $\gamma$;
она описывает центры кривизны данной кривой и записывается как
\[\omega(t)=\gamma(t)+\tfrac1{\skur(t)}\cdot \norm(t).\]
Из доказательства следует, что $(\tfrac1{\skur})'\cdot\norm$ --- её вектор скорости.

\begin{thm}{Упражнение}\label{ex:evolute}
Пусть $\omega$ --- эволюта гладкой кривой~$\gamma$ на плоскости.
Предположим, что у $\gamma$ положительная кривизна $\skur$ и $\skur ^{\prime} \neq 0$ во всех точках.
Найдите базис Френе и ориентированную кривизну эволюты $\omega$ через $\skur$ и базис Френе $(\tan,\norm)$ кривой $\gamma$.
\end{thm}

Про следующую теорему можно думать так:
\textit{если ехать по плоскости на машине и всё время поворачивать руль влево,
то невозможно вернуться на исходное место.}

\begin{thm}{Теорема}\label{thm:spiral}
Любая гладкая плоская кривая $\gamma$ с положительной строго монотонной кривизной не имеет самопересечений.
\end{thm}

Теорема остаётся верной и без предположения о положительности кривизны; доказательство требует лишь незначительных изменений.

\parbf{Доказательство.}
Заметим, что $\gamma(s)\in \sigma_s$, где  $\sigma_s$ --- окружность, соприкасающаяся к $\gamma$ при~$s$.
Из \ref{lem:spiral} следует, что $\sigma_{s_0}$ не пересекает $\sigma_{s_1}$ при $s_1\ne s_0$.
В частности, $\gamma(s_1)\ne \gamma(s_0)$, и теорема следует.
\qeds

{\sloppy

\begin{thm}{Продвинутое упражнение}\label{ex:3D-spiral}
Убедитесь, что трёхмерный аналог теоремы неверен.
То есть постройте гладкую пространственную кривую с пересечением и строго монотонной кривизной.
\end{thm}

}

\begin{thm}{Упражнение}\label{ex:double-tangent}
Предположим, что $\gamma$ --- гладкая плоская кривая с положительной строго монотонной кривизной.

\begin{subthm}{ex:double-tangent:a}
Покажите, что ни одна прямая не касается $\gamma$ в двух различных точках.
\end{subthm}

\begin{subthm}{ex:double-tangent:b}
Покажите, что ни одна окружность не касается $\gamma$ в трёх различных точках.
\end{subthm}

\end{thm}

{

\begin{wrapfigure}{o}{25 mm}
\vskip-4mm
\centering
\includegraphics{mppics/pic-25}
\vskip0mm
\end{wrapfigure}

Часть \ref{SHORT.ex:double-tangent:a} перестаёт быть верной, если разрешить отрицательную кривизну; пример показан на рисунке.

}

\begin{thm}{Продвинутое упражнение}\label{ex:spherical-spiral}
Покажите, что гладкая сферическая кривая с ненулевым кручением не имеет самопересечений.
\end{thm}

\chapter{Опорные кривые}
\label{chap:supporting-curves}

Если одна кривая касается другой, оставаясь на той же её стороне, то кривизну одной можно оценить через кривизну другой.
Мы докажем это утверждение и покажем, как его можно использовать при изучении кривых на плоскости.

\section{Сонаправленность}

Пусть $\gamma_1$ и $\gamma_2$ --- гладкие кривые на плоскости.
Напомним, что кривые $\gamma_1$ и $\gamma_2$ касаются при $t_1$ и $t_2$,
если $\gamma_1(t_1)=\gamma_2(t_2)$
и у них общая касательная при $t_1$ и $t_2$ соответственно.
В этом случае точка $p\z=\gamma_1(t_1)=\gamma_2(t_2)$ называется \index{точка!касания}\emph{точкой касания} кривых.
Заметим, что векторы скорости $\gamma_1'(t_1)$ и $\gamma_2'(t_2)$ параллельны.
\begin{figure}[!ht]
\vskip-0mm
\centering
\includegraphics{mppics/pic-85}
\vskip-0mm
\end{figure}
Касание называется \index{сонаправленное касание}\emph{сонаправленым} или \index{противонаправленное касание}\emph{противонаправленым}, если вектора $\gamma_1'(t_1)$ и $\gamma_2'(t_2)$ сонаправлены или противоположны соответственно.

Обращение параметризации у одной из кривых превращает сонаправленное касание в противонаправленное и наоборот.
Так что всегда можно считать данное касание сонаправленным.

\pagebreak

\section{Опорные кривые}

\begin{wrapfigure}[8]{o}{43 mm}
\vskip-4mm
\centering
\includegraphics{mppics/pic-86}
\vskip0mm
\end{wrapfigure}

Пусть $\gamma_1$ и $\gamma_2$ --- две гладкие плоские кривые с общей точкой 
\[p=\gamma_1(t_1)=\gamma_2(t_2),\] 
и при этом $p$ не является концом ни $\gamma_1$, ни $\gamma_2$.
Предположим, что для некоторого $\epsilon>0$, дуга $\gamma_2|_{[t_2-\epsilon, t_2+\epsilon]}$ лежит в замкнутой области $R$ с дугой $\gamma_1|_{[t_1-\epsilon, t_1+\epsilon]}$ на её границе.
Тогда мы говорим, что $\gamma_1$ --- \index{опорная!кривая}\emph{локальная опорная} к $\gamma_2$ при параметрах $t_1$ и $t_2$.
При этом, если кривые на рисунке параметризованы по стрелкам, то $\gamma_1$ подпирает $\gamma_2$ \emph{справа} в точке $p$ (а также $\gamma_2$ подпирает $\gamma_1$ \emph{слева} в~$p$).

Предположим, что $\gamma_1$ простая, и она разрезает плоскость на две замкнутые области: слева и справа от~$\gamma_1$.
Мы говорим, что $\gamma_1$ \index{опорная!кривая}\emph{глобально подпирает} $\gamma_2$ в точке $p=\gamma_2(t_2)$,
если $\gamma_2$ лежит в одной из этих замкнутых областей, и при этом $p$ лежит на~$\gamma_1$.

Далее, предположим, что $\gamma_2$ --- простая замкнутая кривая на плоскости.
По теореме Жордана (\ref{thm:jordan}), $\gamma_2$ разрезает плоскость на две замкнутые области: одна из них ограниченна, другая --- нет.
Мы говорим, что точка $p$ лежит {}\emph{внутри} (соответственно, {}\emph{снаружи}) $\gamma_2$, если $p$ лежит в ограниченной области (соответственно, неограниченной области).
Если $\gamma_1$ подпирает $\gamma_2$ и лежит внутри $\gamma_2$ (снаружи $\gamma_2$), то мы говорим, что $\gamma_1$ подпирает $\gamma_2$ \index{опорная кривая}\emph{изнутри} (соответственно, {}\emph{снаружи})

Если $p=\gamma_1(t_1)=\gamma_2(t_2)$ и кривые $\gamma_1$ и $\gamma_2$ не касаются друг друга при $t_1$ и $t_2$, то в момент времени $t_2$ кривая $\gamma_2$ переходит с одной стороны $\gamma_1$ на другую.
В этом случае $\gamma_1$ не может подпирать $\gamma_2$ при $t_1$ и $t_2$, то есть верно следующее.

\begin{thm}{Наблюдение и определение}
Пусть $\gamma_1$ и $\gamma_2$ --- две гладкие кривые на плоскости.
Предположим, что $\gamma_1$ локально подпирает $\gamma_2$ при $t_1$ и $t_2$.
Тогда $\gamma_1$ и $\gamma_2$ касаются при $t_1$ и $t_2$.

В этом случае, если кривые сонаправлены, и область $R$ в определении опорных кривых лежит справа (слева) от дуги $\gamma_1$, то мы говорим, что 
$\gamma_1$ подпирает $\gamma_2$ слева (соответственно справа).
\end{thm}

Мы говорим, что гладкая плоская кривая $\gamma$ имеет \index{вершина кривой}\emph{вершину} при $s$, 
если $\skur'(s)_\gamma=0$, то есть функция ориентированной кривизны имеет критическую точку в $s$.
Можно также сказать, что точка $p=\gamma(s)$ есть вершина $\gamma$ (если $\gamma$ простая, то это не приводит к неоднозначности).

\begin{thm}{Упражнение}\label{ex:vertex-support}
Докажите, что точка $p\z=\gamma(s)$ гладкой плоской кривой $\gamma$ является её вершиной, если
соприкасающаяся окружность $\sigma_s$ при $s$ локально подпирает $\gamma$ при $\gamma(s)$.
\end{thm}

\section{Признак опорной}

\begin{thm}{Предложение}\label{prop:supporting-circline}
Пусть $\gamma_1$ и $\gamma_2$ --- гладкие кривые на плоскости.

Предположим, что $\gamma_1$ локально подпирает $\gamma_2$ слева (справа) при $t_1$ и $t_2$.
Тогда
\[\skur_1(t_1)\ge \skur_2(t_2)\quad(\text{соответственно}\quad \skur_1(t_1)\le \skur_2(t_2)),\]
где $\skur_1$ и $\skur_2$ --- ориентированные кривизны $\gamma_1$ и $\gamma_2$, соответственно.

При этом строгое неравенство влечёт обратное.
А именно, если $\gamma_1$ и $\gamma_2$ касаются и сонаправлены при $t_1$ и $t_2$, то
\[\skur_1(t_1)> \skur_2(t_2)\quad(\text{соответственно}\quad \skur_1(t_1)< \skur_2(t_2)),\]
влечёт, что $\gamma_1$ локально подпирает $\gamma_2$ слева (справа) при $t_1$ и $t_2$.

\end{thm}

\parbf{Доказательство.}
Не умаляя общности, можно считать, что $t_1\z=t_2\z=0$, общая точка $\gamma_1(0)=\gamma_2(0)$ --- начало координат, а векторы скорости $\gamma'_1(0)$, $\gamma'_2(0)$ направлены по оси $x$.
Тогда малые дуги $\gamma_1|_{[-\epsilon,+\epsilon]}$ и $\gamma_2|_{[-\epsilon,+\epsilon]}$ задаются графиками 
$y=f_1(x)$ и $y=f_2(x)$ гладких функций $f_1$ и $f_2$, таких что $f_i(0)=0$ и $f_i'(0)=0$.
Из \ref{ex:curvature-graph}, $f_1''(0)=\skur_1(0)$ и $f_2''(0)=\skur_2(0)$.

Очевидно, что $\gamma_1$ подпирает $\gamma_2$ слева (справа), если 
\[f_1(x)\ge f_2(x)\quad(\text{соответственно}\quad f_1(x)\le f_2(x))\]
для всех достаточно малых~$x$.
Остаётся применить  к функции $f_1-f_2$ признак экстремума по второй производной.

Доказательство обратного аналогично.
\qeds

\begin{thm}{Продвинутое упражнение}\label{ex:support}
Предположим, что две гладкие простые плоские кривые с единичной скоростью $\gamma_0$ и $\gamma_1$ касаются и сонаправлены в точке $p\z=\gamma_0(0)\z=\gamma_1(0)$ и при этом $\skur_0(s)\le\skur_1(s)$ для любого~$s$.
Покажите, что $\gamma_0$ локально подпирает $\gamma_1$ справа в точке~$p$.

Приведите пример двух простых кривых $\gamma_0$ и $\gamma_1$, удовлетворяющих указанному условию, таких что $\gamma_0$ замкнута, но
не является глобальной опорной к $\gamma_1$ в точке $p$.
\end{thm}

Напомним (см. \ref{thm:DNA}), что среднее значение кривизны замкнутой гладкой кривой в единичном круге не меньше~$1$.
В частности, на ней должна быть точка с кривизной хотя бы~$1$.
Следующее упражнение говорит, что то же верно и для петель.

\begin{thm}{Упражнение}\label{ex:in-circle}
Предположим, что гладкая петля $\gamma$ лежит в единичном круге на плоскости.
Покажите, что $\gamma$ имеет кривизну не меньше~$1$ в какой-то точке.
\end{thm}

\begin{thm}{Упражнение}\label{ex:between-parallels-1}
{\sloppy
Пусть замкнутая гладкая кривая $\gamma$ лежит между двумя параллельными прямыми на плоскости, расстояние между которыми равно $2$.
Покажите, что существует точка на $\gamma$ с кривизной хотя бы~$1$.

}

Попробуйте доказать то же самое для гладкой плоской петли.
\end{thm}

\begin{thm}{Упражнение}\label{ex:in-triangle}
Пусть замкнутая гладкая плоская кривая $\gamma$ лежит внутри треугольника $\triangle$ с единичным радиусом вписанной окружности.
Покажите, что существует точка на $\gamma$ с кривизной хотя бы~$1$.
\end{thm}

Отметим, что три упражнения выше являются частными случаями упражнения \ref{ex:moon-rad},
однако постарайтесь найти прямые решения.

{

\begin{wrapfigure}{r}{32 mm}
\vskip-4mm
\centering
\includegraphics{mppics/pic-70}
\vskip0mm
\end{wrapfigure}

\begin{thm}{Упражнение}\label{ex:lens}
Пусть $F$ --- фигура на плоскости, ограниченная двумя дугами окружностей $\sigma_1$ и $\sigma_2$ с ориентированной кривизной $1$, которые проходят от $x$ до~$y$.
Предположим, что $\sigma_1$ короче, чем~$\sigma_2$,
а гладкая дуга $\gamma$ лежит в $F$ и имеет обе конечные точки на $\sigma_1$.
Покажите, что кривизна $\gamma$ не меньше $1$ при каком-то значении параметра.

\end{thm}

}

\section{Выпуклые кривые}

Напомним, что плоская кривая называется \index{выпуклая!кривая}\emph{выпуклой}, если она ограничивает выпуклую область.

\begin{thm}{Предложение}\label{prop:convex}
Пусть простая замкнутая гладкая кривая $\gamma$ ограничивает компактное множество~$F$ на плоскости.
Для выпуклости $F$ необходимо и достаточно чтобы ориентированная кривизна $\gamma$ не меняла знак.
\end{thm}

\begin{thm}{Лемма о линзе}\label{lem:lens}
Пусть $\gamma$ --- гладкая простая плоская кривая, идущая от $x$ до~$y$.
Пусть $\gamma$ проходит строго с правой стороны (левой стороны) от ориентированной прямой $xy$, и только её конечные точки $x$ и $y$ лежат на прямой.
Тогда на $\gamma$ найдётся точка с положительной (соответственно, отрицательной) ориентированной кривизной.
\end{thm}

{

\begin{wrapfigure}{o}{35 mm}
\vskip-0mm
\centering
\includegraphics{mppics/pic-22}
\vskip0mm
\end{wrapfigure}

Лемма не выполняется для кривых с самопересечениями.
Например, кривая $\gamma$ на рисунке всегда поворачивает направо, 
поэтому её ориентированная кривизна везде отрицательна, однако она находится с правой стороны от прямой $xy$.

}

\begin{wrapfigure}[6]{i}{50 mm}
\vskip-0mm
\centering
\includegraphics{mppics/pic-24}
\end{wrapfigure}

\parbf{Доказательство.}
Выберем точки $p$ и $q$ на прямой $xy$, 
таким образом, чтобы точки $p, x, y, q$ появлялись на прямой в этом порядке.
Можно предположить, что $p$ и $q$ лежат достаточно далеко от $x$ и $y$, так что полукруг с диаметром $pq$ содержит~$\gamma$.

Рассмотрим наименьший круговой сегмент с хордой $[p,q]$, который содержит~$\gamma$.
Его дуга $\sigma$ подпирает $\gamma$ в некоторой точке $w=\gamma(t_0)$.

Давайте параметризуем $\sigma$ от $p$ до~$q$.
Заметим, что $\gamma$ и $\sigma$ касаются и сонаправлены в~$w$.
В противном случае, дуга $\gamma$ от $w$ до $y$ была бы заперта в криволинейном треугольнике $xwp$, ограниченным отрезком $[p,x]$ и дугами кривых $\sigma$ и $\gamma$.
Но это невозможно, поскольку $y$ не принадлежит этому треугольнику.

\begin{wrapfigure}{o}{50 mm}
\vskip-4mm
\centering
\includegraphics{mppics/pic-23}
\bigskip
\includegraphics{mppics/pic-230}
\end{wrapfigure}

Из этого следует, что $\sigma$ подпирает $\gamma$ справа при $t_0$.
Согласно \ref{prop:supporting-circline},
\[\skur(w)_\gamma\ge \skur_\sigma >0.\]
\qedsf

\parit{Замечание.}
Можно было бы выбрать точку $w$ на $\gamma$, лежащую на максимальном расстоянии от прямой $xy$.
То же рассуждение показывает, что кривизна в $w$ неотрицательна, что немного слабее требуемой положительной кривизны.

\parbf{Доказательство \ref{prop:convex};} \textit{необходимость.}
Если $F$ выпукла, то каждая касательная прямая к $\gamma$ подпирает~$\gamma$.
При движении по $\gamma$, фигура $F$ должна оставаться с одной стороны от её касательной.
Можно предположить, что каждая касательная прямая подпирает $\gamma$ с одной стороны, скажем, справа.
Поскольку прямая имеет нулевую кривизну, по признаку опорной (\ref{prop:supporting-circline}) $\skur\ge 0$ в каждой точке.

\begin{wrapfigure}{r}{35 mm}
\vskip-3mm
\centering
\includegraphics{mppics/pic-68}
\vskip0mm
\end{wrapfigure}

\parit{Достаточность.}
Обозначим через $K$ выпуклую оболочку~$F$.
Если $F$ не выпукла, то $F$ является собственным подмножеством~$K$.
Следовательно, граница $\partial K$ содержит отрезок прямой, который не является частью $\partial F$.
Другими словами, существует прямая, опорная к $\gamma$ в двух точках, скажем $x$ и $y$.
Эти точки делят $\gamma$ на две дуги $\gamma_1$ и $\gamma_2$, обе отличные от отрезка прямой $[x,y]$.

Одна из дуг $\gamma_1$ или $\gamma_2$ параметризована от $x$ до $y$, а другая от $y$ до~$x$.
Применив лемму о линзе, получаем, что на $\gamma$ есть точки с ориентированными кривизнами противоположных знаков.
Возможно, при этом нужно будет перейти к меньшим дугам так, чтобы только их конечные точки лежали на прямой.
\qeds

\begin{thm}{Упражнение}\label{ex:convex small}
Пусть $\gamma$ --- гладкая простая замкнутая плоская кривая диаметра больше~$2$.
Покажите, что на $\gamma$ найдётся точка с кривизной меньше~$1$.
\end{thm}

\begin{wrapfigure}{r}{45 mm}
\vskip-6mm
\centering
\includegraphics{mppics/pic-713}
\vskip0mm
\end{wrapfigure}

\begin{thm}{Упражнение}\label{ex:convex-lens}
Пусть $\gamma$ --- простая гладкая плоская дуга с концами $p$ и~$q$.
Предположим, что $\gamma$ имеет неотрицательную ориентированную кривизну, и $|p-q|$ равно диаметру~$\gamma$.
Покажите, что дуга $\gamma$ и её хорда $[p,q]$ ограничивают выпуклую фигуру на плоскости.
\end{thm}

\begin{thm}{Упражнение}\label{ex:diameter-of-simple-curve}
Покажите, что любая простая гладкая плоская кривая $\gamma$ с кривизной не менее $1$ имеет диаметр не более $2$.

Попробуйте доказать, что $\gamma$ лежит в единичном круге.
\end{thm}

\section{О луне в луже}

Следующая теорема уточняет результат Владимира Ионина и Германа Пестова \cite{ionin-pestov},\index{теорема Ионина --- Пестова}
а для выпуклых кривых он был известен ранее \cite[§ 24]{blaschke-1916}.

\begin{wrapfigure}{r}{18 mm}
\vskip-6mm
\centering
\includegraphics{mppics/pic-67}
\vskip-2mm
\end{wrapfigure}

\begin{thm}{Теорема}\label{thm:moon-orginal}
Предположим, что простая гладкая петля с кривизной не больше $1$ ограничивает фигуру $F$ на плоскости.
Тогда $F$ содержит единичный круг.
\end{thm}

Это простой, но содержательный пример теорем типа {}\emph{от локального к глобальному}.
То есть, исходя из некоторых локальных свойств (в данном случае оценки на кривизну) мы доказываем какое-то глобальное свойство (в данном случае существование единичного круга, окружённого кривой).

{

\begin{wrapfigure}{r}{33 mm}
\vskip-0mm
\centering
\includegraphics{mppics/pic-62}
\vskip0mm
\end{wrapfigure}

Можно было бы попробовать начать с какого-то круга в $F$ раздувать его в надежде достичь единичного радиуса.
Однако, пример на рисунке показывает, что это не всегда приводит к решению --- 
иногда, чтобы достичь единичного радиуса придётся сперва сдуть круг.

\begin{thm}{Основная лемма}\label{thm:moon}
Пусть $\gamma$ --- простая гладкая петля на плоскости.
Тогда соприкасающаяся к ней окружность в некоторой точке (отличной от базовой) глобально подпирает $\gamma$ изнутри.
\end{thm}

}

Сначала покажем, что теорема следует из леммы.

\parbf{Доказательство \ref{thm:moon-orginal}, используя \ref{thm:moon}.}
Поскольку кривизна $\gamma$ не превышает $1$, каждая соприкасающаяся окружность имеет радиус не менее~$1$.
По основной лемме, одна из соприкасающихся окружностей, скажем $\sigma$, подпирает $\gamma$ изнутри.
В частности, $\sigma$ лежит внутри $\gamma$, отсюда результат.
\qeds

\parbf{Доказательство \ref{thm:moon}.}
Обозначим через $F$ замкнутую область, окружённую $\gamma$.
Можно считать, что $F$ находится слева от~$\gamma$.
Рассуждая от противного, допустим, что соприкасающаяся окружность в каждой точке $p\in \gamma$ не лежит в~$F$.

\begin{figure}[!ht]
\vskip-0mm
\centering
\includegraphics{mppics/pic-32}
\vskip-2mm
\end{figure}

Для данной точки $p\in\gamma$ рассмотрим максимальную окружность, которая полностью лежит в $F$ и касается $\gamma$ в~$p$.
Эту окружность, скажем $\sigma$, будем называть {}\emph{вписанной} в $F$ при~$p$.
Её кривизна $\skur_\sigma$ обязана быть больше, чем $\skur(p)_\gamma$.
Действительно, из \ref{prop:supporting-circline}, $\skur_\sigma\ge \skur(p)_\gamma$, ведь $\sigma$ подпирает $\gamma$ слева.
В случае равенства $\sigma$ является соприкасающейся окружностью при~$p$,
а это невозможно в силу нашего предположения.

Следовательно, $\sigma$ обязана коснуться $\gamma$ в другой точке.
Иначе её можно было бы подраздуть, оставляя внутри~$F$.
Действительно, поскольку $\skur_\sigma> \skur(p)_\gamma$, 
по \ref{prop:supporting-circline} найдётся окрестность $U$ точки $p$ так, что после небольшого раздутия $\sigma$, пересечение $U\cap \sigma$ всё ещё в~$F$.
С другой стороны, если $\sigma$ не касается $\gamma$ в другой точке, то после некоторого (возможно меньшего) раздутия $\sigma$ дополнение $\sigma\setminus U$ всё ещё в~$F$.
То есть, немного увеличенная $\sigma$ всё ещё лежит в $F$ --- противоречие.

Выберем точку $p_1$ на петле $\gamma$, отличную от базовой точки. 
Пусть $\sigma_1$ --- вписанная окружность при $p_1$,
и $\gamma_1$ --- дуга $\gamma$ от $p_1$ до первой точки $q_1$ на $\sigma_1$.
Обозначим через $\hat\sigma_1$ и $\check\sigma_1$ две дуги $\sigma_1$ от $p_1$ до $q_1$ так, что циклическое произведение $\hat\sigma_1$ и $\gamma_1$ окружает~$\check\sigma_1$. 

Пусть $p_2$ --- середина $\gamma_1$.
Обозначим через $\sigma_2$ вписанную окружность в точке $p_2$.

Окружность $\sigma_2$ не может пересечь $\hat\sigma_1$,
ибо если $\sigma_2$ пересекает $\hat\sigma_1$ в какой-то точке $s$, то она должна иметь ещё две общие точки с $\check\sigma_1$, скажем $x$ и $y$ --- по одной для каждой дуги $\sigma_2$ от $p_2$ до~$s$.
Следовательно, $\sigma_1=\sigma_2$, ведь у них три общие точки: $s$, $x$ и~$y$. 
С другой стороны, по построению, $p_2\in \sigma_2$ и $p_2\notin \sigma_1$ --- противоречие.

\begin{wrapfigure}{r}{32 mm}
\vskip-0mm
\centering
\includegraphics{mppics/pic-64}
\caption*{Два овала изображают окружности.}
\vskip-2mm
\end{wrapfigure}

Напомним, что $\sigma_2$ должна касаться $\gamma$ в другой точке.
Из вышесказанного вытекает, что $\sigma_2$ может касаться только~$\gamma_1$. 
Следовательно, можно выбрать дугу $\gamma_2\subset \gamma_1$, которая идёт от $p_2$ до первой точки $q_2$ на $\sigma_2$.
Поскольку $p_2$ --- середина $\gamma_1$, 
\[\length \gamma_2< \tfrac12\cdot\length\gamma_1.\eqlbl{eq:length<length/2}\]

Повторяя это построение рекурсивно,
получаем бесконечную последовательность дуг $\gamma_1\z\supset \gamma_2\z\supset\dots$;
по \ref{eq:length<length/2} также получим, что 
\[\length\gamma_n\to0\quad\text{когда}\quad n\to\infty.\] 
Следовательно, пересечение $\gamma_1\cap\gamma_2\cap\dots$
содержит единственную точку, скажем $p_\infty$.

Пусть $\sigma_\infty$ --- вписанная окружность при $p_\infty$; она должна касаться $\gamma$ в другой точке, скажем $q_\infty$.
То же рассуждение, что и выше, показывает, что $q_\infty\in\gamma_n$ для любого~$n$.
Следовательно, $q_\infty =p_\infty$ --- противоречие.
\qeds

\begin{thm}{Упражнение}\label{ex:moon-rad}
Пусть замкнутая гладкая кривая $\gamma$ лежит в фигуре $F$, ограниченной простой замкнутой плоской кривой.
Пусть $R$ --- это максимальный радиус кругов, которые лежат в~$F$.
Покажите, что кривизна $\gamma$ не менее $\tfrac1R$ при некотором значении параметра.
\end{thm}

\section{Теорема о четырёх вершинах}
\index{теорема о четырёх вершинах}

{

\begin{wrapfigure}{r}{20 mm}
\vskip-8mm
\centering
\includegraphics{mppics/pic-26}
\vskip0mm
\end{wrapfigure}

Напомним, что вершина гладкой кривой определяется как критическая точка её ориентированной кривизны.
В частности, точка локального минимума (или максимума) кривизны является вершиной,
а у окружности все точки --- вершины.

\begin{thm}{Теорема}\label{thm:4-vert}
Любая гладкая простая замкнутая плоская кривая имеет хотя бы четыре вершины.
\end{thm}

}

Ясно, что любая замкнутая гладкая кривая имеет как минимум две вершины --- точки, в которых достигаются минимум и максимум кривизны.
На рисунке отмечены вершины двух кривых;
первая имеет одно самопересечение и ровно две вершины;
у второй четыре вершины и нет самопересечений.

Для выпуклых кривых теорему доказал
Шьямадас Мукхопадхьяй \cite{mukhopadhyaya}.
Одно из самых красивых доказательств общего случая предложил Роберт Оссерман \cite{osserman}.
Мы докажем следующее более сильное утверждение, используя основную лемму предыдущего раздела.
Больше по теме можно узнать в нашей статье \cite{petrunin-zamora:moon} и в указанных там ссылках.

{

\begin{wrapfigure}[7]{r}{33 mm}
\vskip-6mm
\centering
\includegraphics{mppics/pic-63}
\vskip0mm
\end{wrapfigure}

\begin{thm}{Теорема}\label{thm:4-vert-supporting}
На любой гладкой простой замкнутой плоской кривой $\gamma$ можно найти четыре точки, соприкасающиеся окружности в которых подпирают $\gamma$; две изнутри и две --- снаружи.
\end{thm}

\parbf{Доказательство \ref{thm:4-vert}, используя \ref{thm:4-vert-supporting}.}
Достаточно доказать следующее: \textit{если соприкасающаяся окружность $\sigma$ в точке $p$ подпирает $\gamma$ в $p$, то $p$ --- вершина~$\gamma$}.

}

Если это не так, то кривизна на малой дуге вокруг $p$ монотонна.
По лемме о спирали (\ref{lem:spiral}), соприкасающиеся окружности на этой дуге образуют монотонное семейство.
В частности, кривая $\gamma$ проходит сквозь $\sigma$ в точке~$p$, и, значит, $\sigma$ не подпирает $\gamma$.
\qeds

\parbf{Доказательство \ref{thm:4-vert-supporting}.}
По основной лемме (\ref{thm:moon}), соприкасающаяся окружность при некоторой точке $p\in\gamma$ подпирает $\gamma$ изнутри.
Кривую $\gamma$ можно рассмотреть как петлю с базовой точкой~$p$.
Значит, по основной лемме, найдётся ещё одна точка $q\in\gamma$ с тем же свойством.

Это даёт пару соприкасающихся окружностей, опорных к $\gamma$ изнутри;
остаётся раздобыть ещё две.

Для получения соприкасающихся окружностей, подпирающих $\gamma$ снаружи, можно повторить доказательство основной леммы, но вместо вписанных окружностей брать окружность (или прямую) максимальной ориентированной кривизны, опорную к кривой снаружи, предполагая, что $\gamma$ ориентирована так, что область слева от неё ограничена.
\qeds

\parbf{Другая концовка доказательства.}
Применим к $\gamma$ инверсию относительно окружности с центром внутри~$\gamma$. Тогда у полученной кривой, скажем $\gamma_1$, также можно найти пару соприкасающихся окружностей, опорных к $\gamma_1$ изнутри.
Согласно \ref{ex:inverse}, их инверсии являются соприкасающимися окружностями к~$\gamma$.
При этом область, лежащая внутри $\gamma$, отображается в область снаружи $\gamma_1$, и наоборот.
Следовательно, эти две новые окружности подпирают $\gamma$ снаружи.\qeds

{

\begin{wrapfigure}{r}{20 mm}
\vskip-0mm
\centering
\includegraphics{mppics/pic-725}
\vskip0mm
\end{wrapfigure}

\begin{thm}{Упражнение}\label{ex:2-squares}
Предположим, что гладкая простая замкнутая плоская кривая \(\gamma\) лежит в квадрате со стороной $2$ и окружает квадрат с диагональю~$2$.
Докажите, что \(\gamma\) содержит точку с кривизной равной~$1$.
\end{thm}

}

\begin{thm}{Упражнение}\label{ex:moon-area}
Предположим, что гладкая простая замкнутая плоская кривая ограничивает область площади $a$.
Докажите, что на ней найдётся точка с кривизной $\sqrt{\pi/a}$.
\end{thm}

\begin{wrapfigure}[5]{r}{25 mm}
\vskip-7mm
\centering
\includegraphics{mppics/pic-65}
\vskip0mm
\end{wrapfigure}

\begin{thm}{Продвинутое упражнение}\label{ex:curve-crosses-circle}
{\sloppy
Предположим, что простая замкнутая гладкая плоская кривая $\gamma$  пересекает окружность $\sigma$ в точках $p_1,\dots,p_{2{\cdot} n}$, и эти точки появляются в том же порядке и на $\gamma$ и на $\sigma$.
Докажите, что у $\gamma$ есть как минимум $2\cdot n$ вершин.

}

Постройте пример простой замкнутой гладкой плоской кривой $\gamma$ у которой только $4$ вершины, и при этом она пересекает данную окружность в произвольно большом числе точек.
\end{thm}

{\sloppy

\begin{thm}{Продвинутое упражнение}\label{ex:berk}
Пусть $\gamma$ --- простая гладкая плоская кривая с кривизной, ограниченной~$1$.
Докажите, что $\gamma$ окружает два непересекающихся открытых единичных диска тогда и только тогда, когда её диаметр не менее $4$;
то есть $\dist{p}{q}{}\ge 4$ для некоторых точек $p,q\in\gamma$.
\end{thm}

}


Следующее упражнение является версией теоремы о четырёх вершинах для пространственных кривых без параллельных касательных; такие кривые существуют по \ref{ex:no-parallel-tangents}.

\begin{thm}{Продвинутое упражнение}\label{ex:4x0-torsion}
Пусть $\gamma$ --- замкнутая гладкая пространственная кривая, у которой нет пары точек с параллельными касательными.
Допустим, что кривизна $\gamma$ не обращается в ноль ни в одной точке.
Докажите, что у $\gamma$ есть как минимум четыре точки с нулевым кручением.
\end{thm}

\arxiv{\cleardoublepage
\phantomsection
\AddToShipoutPictureBG*{\includegraphics[width=\paperwidth]{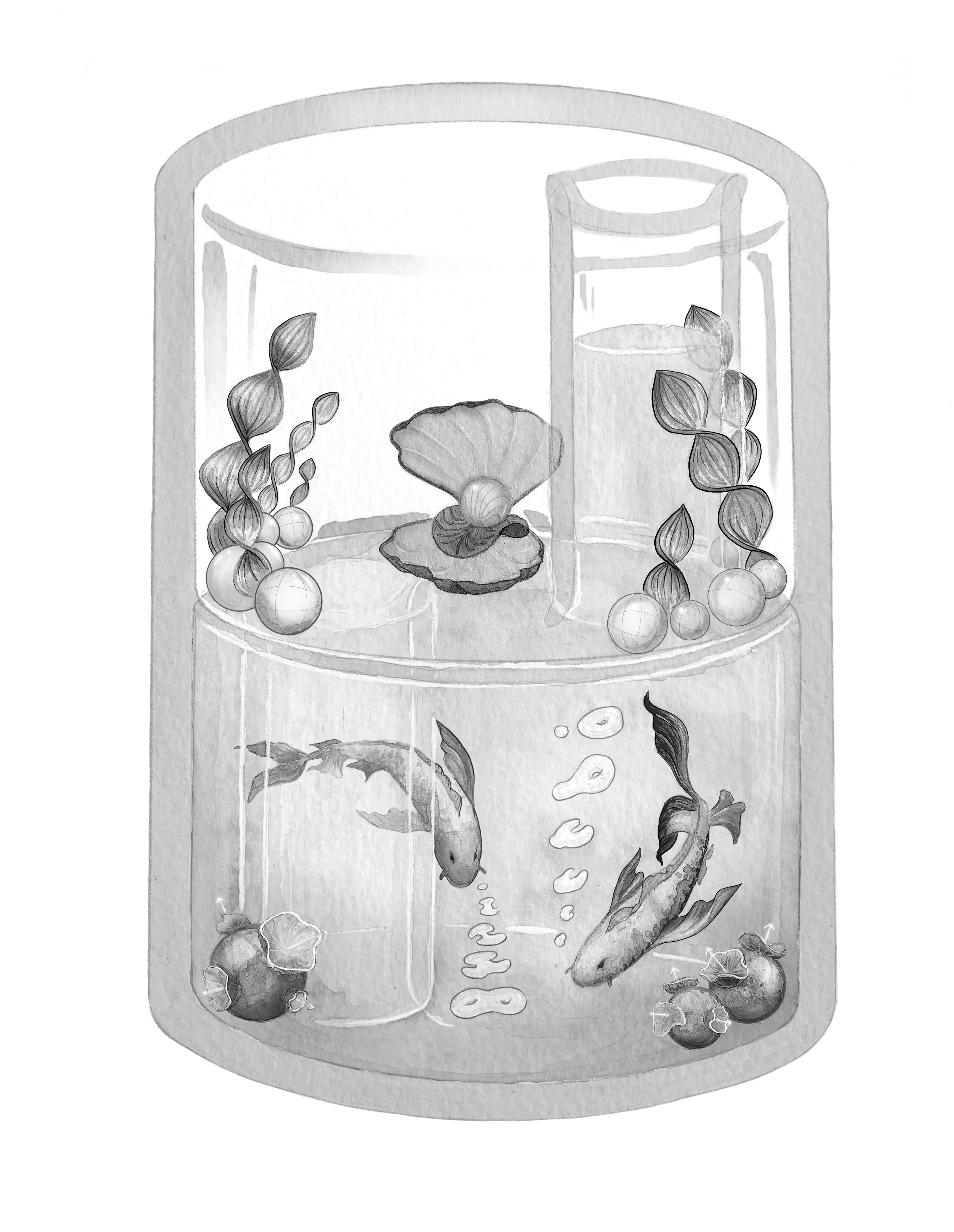}}
\cleardoublepage
\thispagestyle{empty}
\stepcounter{part}
\begin{center}
{\Huge\textbf{Часть \Roman{part}\qquad Поверхности}}\
\end{center}
\addcontentsline{toc}{part}{\texorpdfstring{Часть \Roman{part}\quad Поверхности}{Часть \Roman{part} Поверхности}}
\clearpage
}
{\backgroundsetup{scale=.95,opacity=1,angle=0,vshift=10mm,contents={%
\includegraphics[width=\paperwidth]{pics/Surfaces}
}%
}
\cleardoublepage
\phantomsection
\stepcounter{part}
\addcontentsline{toc}{part}{\texorpdfstring{Часть \Roman{part}\quad Поверхности}{Часть \Roman{part} Поверхности}}
\thispagestyle{empty}
\begin{center}
{\Huge\textbf{Часть \Roman{part}\qquad Поверхности}}\
\end{center}
\BgThispage
}

\chapter{Определения}
\label{chap:surfaces-def}

\section{Топологические поверхности}

В основном нас будут интересовать гладкие поверхности, определённые в следующем разделе.
Сейчас мы дадим более общее определение, которое будет использоваться лишь пару раз.

Связное подмножество $\Sigma$ евклидова пространства $\mathbb{R}^3$
называется \index{поверхность}\index{топологическая поверхность}\emph{топологической поверхностью} (точнее {}\emph{вложенной топологической поверхностью без края}),
если любая точка $p\in \Sigma$ имеет окрестность $W$ в $\Sigma$,
которую можно параметризовать открытым подмножеством евклидовой плоскости;
то есть существует гомеоморфизм $V\to W$ из открытого множества $V\subset \mathbb{R}^2$; см. приложение~\ref{sec:topology}.

\section{Гладкие поверхности}\label{sec:def-smooth-surface}

Напомним, что функция $f$ двух переменных $x$ и $y$ называется \index{гладкая!функция}\emph{гладкой}, если все её частные производные $\frac{\partial^{m+n}}{\partial x^m\partial y^n}f$ определены (а значит, и непрерывны) в области определения~$f$.

{\sloppy

Связное множество $\Sigma \subset \mathbb{R}^3$ называется \index{поверхность}\index{гладкая!поверхность}\emph{гладкой поверхностью}%
\footnote{Это сокращение термина {}\emph{гладкая регулярная вложенная поверхность}.}, если оно локально описывается как график гладкой функции в подходящей системе координат.

Точнее, для любой точки $p\in \Sigma$ найдётся система координат $(x,y,z)$ и окрестность $U\ni p$ так, что
пересечение $W=U\cap \Sigma$ задаётся графиком $z=f(x,y)$ гладкой функции $f$, определённой в некоторой открытой области $V$ плоскости $(x,y)$.

}

Заметим, что отображение $(x,y)\mapsto (x,y,f(x,y))$ задаёт гомеоморфизм $V\to W$.
Таким образом, гладкая поверхность является примером топологической поверхности.

\parbf{Примеры.}
Координатная лоскость $z=0$ даёт пример гладкой поверхности, ведь она является графиком функции $f(x,y)=0$.

Любая другая плоскость также является гладкой поверхностью, ведь в подходящей системе координат  она является координатной.
Можно также представить плоскость как график линейной функции 
$f(x,y)\z=a\cdot x+b\cdot y+c$ для некоторых констант $a$, $b$ и $c$
(предполагая, что плоскость не вертикальна, иначе придётся взять другие координаты).

Более хитрый пример даёт единичная сфера 
\[\mathbb{S}^2=\set{(x,y,z)\in\mathbb{R}^3}{x^2+y^2+z^2=1}.\]
Вся сфера не является графиком,
однако её можно покрыть $6$ графиками:
\begin{align*}
z&=f_\pm(x,y)=\pm \sqrt{1-x^2-y^2},
\\
y&=g_\pm(x,z)=\pm \sqrt{1-x^2-z^2},
\\
x&=h_\pm(y,z)=\pm \sqrt{1-y^2-z^2},
\end{align*}
где каждая функция $f_+$, $f_-$, $g_+$, $g_-$, $h_+$ и $h_-$ определена в открытом единичном круге.
Любая точка $p\in\mathbb{S}^2$ лежит на одном из этих графиков, следовательно, $\mathbb{S}^2$ --- гладкая  поверхность.

\section{Поверхности с краем}

{\sloppy

Связное подмножество поверхности, ограниченное одной или несколькими кусочно-гладкими кривыми, называется \index{поверхность!с краем}\emph{поверхностью с краем}; эти кривые и образуют \index{край}\emph{край} поверхности.

}

Как правило, мы говорим {}\emph{поверхность}, имея в виду {}\emph{гладкую поверхность без края}.
При необходимости, можно использовать термин {}\emph{поверхность с возможно непустым краем}.

\section{Типы поверхностей}

Если поверхность $\Sigma$ образует замкнутое множество в $\mathbb{R}^3$, то она называется \index{собственная!поверхность}\emph{собственной}.
Например, график $z=f(x,y)$ любой гладкой функции $f$, определённой на всей плоскости, является собственной поверхностью.
Сфера $\mathbb{S}^2$ --- другой пример собственной поверхности.

С другой стороны, открытый круг 
\[\set{(x,y,z)\in\mathbb{R}^3}{x^2+y^2<1,\  z=0}\]
--- несобственная поверхность, ведь он не открыт, и не замкнут в $\mathbb{R}^3$.

{\sloppy

Компактная поверхность без края называется \index{замкнутая!поверхность}\emph{замкнутой}
(по аналогии с замкнутой кривой, но без связи с замкнутым множеством).

}

\index{открытая!поверхность}\emph{Открытая} поверхность --- это
собственная некомпактная поверхность без края 
(опять же, термин связан с открытой кривой, но не с открытым множеством).

Для примера, параболоид $z=x^2+y^2$ --- открытая поверхность,
а сфера $\mathbb{S}^2$ --- замкнутая.
Заметим, что \textit{любая собственная поверхность без края либо замкнута либо открыта}.

{\sloppy

Следующий факт является трёхмерным аналогом теоремы Жордана (\ref{ex:proper-curve}).
Он выводится из так называемой {}\emph{двойственности Александера} \cite{hatcher}.
Мы опускаем доказательство, оно увело бы нас далеко в сторону.

}

\begin{thm}{Факт}\label{clm:proper-divides}
Дополнение любой собственной топологической поверхности без края (эквивалентно, любой открытой или замкнутой топологической поверхности) имеет ровно две связные компоненты.
\end{thm}

\section{Неявно заданные поверхности}

\begin{thm}{Предложение}\label{prop:implicit-surface}
Пусть ноль --- регулярное значение гладкой функции $f\:\mathbb{R}^3\to \mathbb{R}$;
то есть $\nabla_p f\ne 0$, если $f(p)=0$.
Тогда любая компонента связности $\Sigma$ множества уровня $f(x,y,z)=0$ является гладкой поверхностью.
\end{thm}

\parbf{Доказательство.}
Выберем точку $p\in\Sigma$.
Так как $\nabla_p f\ne 0$, выполнятся какое-то из условий 
$f_x(p)\ne 0$,
$f_y(p)\ne 0$ или
$f_z(p)\ne 0$.
Можно предположить, что $f_z(p)\ne 0$;
иначе, переименуем координаты $x,y,z$.

По теореме о неявной функции (\ref{thm:imlicit}), окрестность $p$ в $\Sigma$ является графиком $z=h(x,y)$ гладкой функции $h$, определённой на открытом множестве в $\mathbb{R}^2$.
Остаётся применить определение гладкой поверхности (см.~\ref{sec:def-smooth-surface}).
\qeds

\begin{thm}{Упражнение}\label{ex:hyperboloids}
При каких $\ell$ множество уровня $x^2+y^2-z^2=\ell$ является гладкой поверхностью?
\end{thm}

\section{Локальные параметризации}
\index{параметризация}

Пусть $U$ --- открытая область в $\mathbb{R}^2$, и $s\:U\to \mathbb{R}^3$ --- гладкое отображение;
оно называется \index{регулярная!параметризация}\emph{регулярным}, если его матрица Якоби имеет максимальный ранг в любой точке.
В данном случае это означает, что производные $s_u$ и $s_v$ линейно независимы в любой точке $(u,v)\in U$;
эквивалентно $s_u\times s_v\ne 0$, где $\times$ обозначает векторное произведение.

\begin{thm}{Предложение}\label{prop:graph-chart}
Образ гладкого регулярного вложения $s$ открытого связного множества $U\subset \mathbb{R}^2$ является гладкой поверхностью.
\end{thm}

\parbf{Доказательство \ref{prop:graph-chart}.}
Пусть $s(u,v)=(x(u,v),y(u,v),z(u,v))$.
Поскольку $s$ регулярно, матрица Якоби
\[\Jac s=
\renewcommand\arraystretch{1.3}
\begin{pmatrix}
x_u&x_v\\
y_u&y_v\\
z_u&z_v
\end{pmatrix}
\]
имеет ранг два в любой точке $(u,v)\in U$.
Пусть $\Sigma=s(U)$.

Выберем точку $p\in \Sigma$; сдвинув координаты $(x,y,z)$ и $(u,v)$, можно считать, что $p = (0,0,0) =s(0,0)$.
Переименовав координаты $x,y,z$, можно добиться, чтобы
матрица $\left(\begin{smallmatrix}
x_u&x_v\\
y_u&y_v
\end{smallmatrix}\right)$
стала обратима.
Отметим, что это матрица Якоби отображения $(u,v)\mapsto (x(u,v),y(u,v))$.

По теореме об обратной функции (\ref{thm:inverse}), есть гладкое регулярное отображение
$w\:(x,y)\mapsto (u,v)$, определённое на открытом множестве $W\ni 0$ плоскости $(x,y)$ такое, что $w(0,0)=(0,0)$ и $s\z\circ w(x,y)\z=(x,y,f(x,y))$, где $f=z\circ w$.
То есть, подмножество $s\z\circ w(W)\subset \Sigma$ является графиком $f$.

И снова по теореме об обратной функции, образ $w(W)$ открыт в $U$.
Поскольку $s$ --- вложение, наш график открыт в~$\Sigma$;
то есть существует такое открытое множество $V\subset \mathbb{R}^3$, что $s\circ w(W)=V\cap \Sigma$ --- график гладкой функции.
А это уже означает, что $\Sigma$ гладкая поверхность, ведь точка $p\in \Sigma$ выбиралась произвольно.
\qeds

\begin{thm}{Упражнение}\label{ex:9-surf}
Найдите гладкое регулярное инъективное отображение $s\:\mathbb{R}^2\to\mathbb{R}^3$ образ которого \textit{не} является поверхностью.
\end{thm}

Если $s$ и $\Sigma$, как в предложении, то $s$ называется \index{гладкая!параметризация}\emph{гладкой параметризацией} поверхности~$\Sigma$. 

Не все гладкие поверхности можно так запараметризовать;
например, сферу $\mathbb{S}^2$ нельзя.
Однако, \textit{любая гладкая поверхность $\Sigma$ допускает локальную параметризацию в любой точке $p\in\Sigma$}; то есть $p$ имеет открытую окрестность $W\subset \Sigma$ с гладкой регулярной параметризацией~$s$.
В этом случае любая точка в $W$ может быть описана двумя параметрами, обычно обозначаемыми как $u$ и $v$;
они называются \index{локальные координаты}\emph{локальными координатами} при~$p$.
Отображение $s$ называется \index{карта}\emph{картой} поверхности~$\Sigma$.

Если $W$ --- график $z=h(x,y)$ гладкой функции $h$, то отображение 
\[s\:(u,v)\mapsto (u,v,h(u,v))\] является картой.
Действительно, отображение $(u,v,h(u,v))\mapsto (u,v)$ является обратным к $s$, и оно непрерывно;
то есть $s$ --- вложение.
Кроме того,
$s_u\z=(1,0,h_u)$ и $s_v\z=(0,1,h_v)$. 
В частности, $s_u$ и $s_v$ линейно независимы;
то есть $s$ регулярно.

\begin{thm}{Следствие}\label{cor:reg-parmeterization}
Связное множество $\Sigma\subset \mathbb{R}^3$ является гладкой поверхностью тогда и только тогда, когда у любой точки в $\Sigma$ есть 
окрестность в $\Sigma$, которую можно покрыть картой.
\end{thm}

Функция $g\: \Sigma \to \mathbb{R}$, определённая на гладкой поверхности $\Sigma$, называется \index{гладкая!функция}\emph{гладкой}, если для любой карты $s \: U\to \Sigma$,
композиция $g\circ s$ гладкая; то есть все частные производные $\frac{\partial^{m+n}}{\partial u^m\partial v^n}(g\circ s)$ определены и непрерывны в области определения.

\begin{thm}{Упражнение}\label{ex:smooth-fun(surf)}
Пусть $\Sigma\subset \mathbb{R}^3$ --- гладкая поверхность.
Покажите, что функция $g\:\Sigma\to\mathbb{R}$ --- гладкая тогда и только тогда, когда для любой точки $p\in \Sigma$ найдётся окрестность $N\subset \mathbb{R}^3$ и такая гладкая функция $h\:N\to\mathbb{R}$, что равенство $g(q)=h(q)$ выполняется для точки $q\in \Sigma\cap N$.

Постройте гладкую поверхность $\Sigma$ с гладкой функцией $g\:\Sigma\to\mathbb{R}$, которую невозможно продолжить до гладкой функции $h\:\mathbb{R}^3\to\mathbb{R}$.
\end{thm}

\begin{thm}{Упражнение}\label{ex:inversion-chart}
Покажите, что
\[s(u,v)=(\tfrac{2\cdot u}{1+u^2+v^2},\tfrac{2\cdot v}{1+u^2+v^2},\tfrac{2}{1+u^2+v^2}).\]
является картой единичной сферы с центром в точке $(0,0,1)$; опишите образ~$s$.
\end{thm}

\begin{wrapfigure}{r}{31 mm}
\vskip-6mm
\centering
\includegraphics{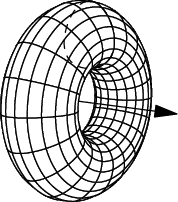}
\vskip0mm
\end{wrapfigure}

Пусть $\gamma(t)=(x(t),y(t))$ --- кривая на плоскости.
Напомним, что \index{поверхность!вращения}\emph{поверхность вращения} кривой $\gamma$ вокруг оси $x$ описывается как образ отображения 
\[(t, s)\mapsto (x(t), y(t)\cdot\cos s,y(t)\cdot\sin s).\]
При этом параметры $t$ и $s$ называются соответственно \index{широта}\emph{широтой} и \index{долгота}\emph{долготой};
для фиксированных $t$ или $s$ полученные кривые называются \index{параллель}\emph{параллелями} или
\index{меридиан}\emph{меридианами} соответственно. 
Отметим, что параллели образованы окружностями в плоскости, перпендикулярной оси вращения.
Кривая $\gamma$ называется \index{образующая}\emph{образующей} этой поверхности.

\begin{thm}{Упражнение}\label{ex:revolution}
Предположим, что $\gamma$ --- простая замкнутая гладкая плоская кривая, которая не пересекает ось $x$.
Покажите, что поверхность вращения вокруг оси $x$ с образующей $\gamma$ является гладкой поверхностью.
\end{thm}

\section{Глобальные параметризации}\label{sec:global-parametrizations}
\index{параметризация}

Поверхность можно описать как вложение другой поверхности.

Для примера, рассмотрим эллипсоид
\[\Theta=\set{(x,y,z)\in\mathbb{R}^3}{\tfrac{x^2}{a^2}+\tfrac{y^2}{b^2}+\tfrac{z^2}{c^2}=1},\]
где $a$, $b$ и~$c$ --- положительные константы.
По \ref{prop:implicit-surface}, $\Theta$ --- гладкая поверхность.
Действительно, пусть $h(x,y,z)\df\tfrac{x^2}{a^2}+\tfrac{y^2}{b^2}+\tfrac{z^2}{c^2}$,
тогда
\[\nabla h(x,y,z)=(\tfrac{2}{a^2}\cdot x,\tfrac{2}{b^2}\cdot y,\tfrac{2}{c^2}\cdot z).\]
Следовательно, $\nabla h\ne0$, если $h=1$; то есть $1$ является регулярным значением~$h$.
Остаётся заметить, что множество $\Theta$ связно.

Поверхность $\Theta$ задаётся как образ сферы $\mathbb{S}^2$ при отображении
\[(x,y,z)\z\mapsto (a\cdot x, b\cdot y,c\cdot z).\]

Отображение $f\:\Sigma \to \mathbb{R}^3$ (или $f\:\Sigma \to \mathbb{R}^2$) считается гладким, если гладкие все его координатные функции.
Кроме того, гладкое отображение $f \: \Sigma \to \mathbb{R}^3$ называется 
\emph{гладкой параметризованной поверхностью}, если оно является вложением, и для любой карты $s \:U\to \Sigma$,
композиция $f\circ s$ регулярна;
то есть два вектора 
$\frac{\partial}{\partial u}(f\circ s)$ и $\frac{\partial}{\partial v}(f\circ s)$ являются линейно независимыми.
В этом случае, образ $\Sigma^{*}=f(\Sigma)$ --- гладкая поверхность, ведь для любой карты $s\:U\to \Sigma$ композиция $f\circ s\:U\to \Sigma^{*}$ задаёт карту на $\Sigma^{*}$. 

Такое отображение $f$ называется \index{диффеоморфизм}\emph{диффеоморфизмом} из $\Sigma$ в $\Sigma^{*}$;
говорят, что поверхности $\Sigma$ и $\Sigma^{*}$ {}\emph{диффеоморфны}, если существует диффеоморфизм $f\:\Sigma\to\Sigma^{*}$.
Следующее упражнение говорит, что \textit{быть диффеоморфными} --- это отношением эквивалентности.

\begin{thm}{Упражнение}\label{ex:inv-diffeomorphism}
Докажите, что обратное отображение к диффеоморфизму также диффеоморфизм.
\end{thm}

\begin{thm}{Продвинутое упражнение}\label{ex:star-shaped-disc}
Докажите, следующее:

\begin{subthm}{ex:plane-n}
Дополнения $n$-точечных подмножеств плоскости диффеоморфны между собой.
\end{subthm}

\begin{subthm}{ex:star-shaped-disc:smooth}
Открытые выпуклые множества плоскости, ограниченные замкнутыми гладкими кривыми диффеоморфны между собой.
\end{subthm}

{\sloppy

\begin{subthm}{ex:star-shaped-disc:nonsmooth}
Любые открытые выпуклые множества плоскости диффеоморфны между собой.
\end{subthm}

\begin{subthm}{ex:star-shaped-disc:star-shaped}
Любые открытые звёздные множества плоскости диффеоморфны между собой.
\end{subthm}

}

\end{thm}

\chapter{Первая производная}
\label{chap:first-order}
\section{Касательная плоскость}

\begin{thm}{Определение}\label{def:tangent-vector}
Пусть $\Sigma$ --- гладкая поверхность.
Вектор $\vec w$ называется \index{касательный вектор}\emph{касательным} к $\Sigma$ в точке $p$, если найдётся кривая $\gamma$ на $\Sigma$ с вектором скорости $\vec w$ в точке $p$;
то есть $p=\gamma(t)$ и $\vec w=\gamma'(t)$ для некоторого~$t$.
\end{thm}

\begin{thm}{Предложение и определение}\label{def:tangent-plane}
Пусть $p$ --- точка гладкой поверхности $\Sigma$.
Тогда множество касательных векторов к $\Sigma$ в точке $p$ образует плоскость;
эта плоскость называется \index{касательная!плоскость}\emph{касательной плоскостью} к $\Sigma$ в точке~$p$.

Более того, если $s\:U\to \Sigma$ --- локальная карта и $p=s(u_0,v_0)$, то 
касательная плоскость к $\Sigma$ в точке $p$ натянута на вектора $s_u(u_0,v_0)$ и $s_v(u_0,v_0)$.
\end{thm}

Касательная плоскость к $\Sigma$ в точке $p$ обозначается как $\T_p$ или $\T_p\Sigma$.
Эту плоскость можно рассматривать как линейное подпространство $\mathbb{R}^3$ или как параллельную плоскость, проходящую через~$p$;
последнюю можно называть \index{касательная плоскость}\emph{аффинной касательной плоскостью}.
Она даёт наилучшее приближение плоскостью окрестности $p$ в $\Sigma$.
Точнее, 
она имеет \index{порядок касания}\emph{первый порядок касания} с $\Sigma$ в точке $p$;
то есть $\rho(q)\z=o(|p\z-q|)$, где $q\in \Sigma$, а $\rho(q)$ обозначает расстояние от $q$ до $\T_p$.

Выбор касательного вектора $\vec v_p$ при каждой точке $p$ поверхности $\Sigma$ называется \emph{полем касательных векторов}, если $\vec v_p$ гладко зависит от $p$;
то есть в любой карте $s$ поле записывается как
\[\vec v_{s(u,v)}=a(u,v)\cdot s_u+b(u,v)\cdot s_v,\]
где $a$ и $b$ --- гладкие функции, определённые в области $s$.

\parbf{Доказательство.}
Выберем карту $s$ при~$p$.
Предположим, что гладкая кривая $\gamma$ начинается в~$p$.
Можно считать, что $\gamma$ накрыта картой;
в частности, 
\[\gamma(t)=s(u(t),v(t))\]
для гладких функций $t\mapsto u(t)$ и $t\mapsto v(t)$.
По правилу дифференцирования сложной функции, 
\[\gamma'=s_u\cdot u'+ s_v\cdot v'.\]
В частности, $\gamma'$ является линейной комбинацией $s_u$ и $s_v$.

Поскольку гладкие функции $t\mapsto u(t)$ и $t\mapsto v(t)$ можно было выбрать произвольно, любая линейная комбинация $s_u$ и $s_v$ является касательным вектором в~$p$. 
\qeds

\begin{thm}{Упражнение}\label{ex:tangent-normal}
Поверхность $\Sigma$ задана как компонента связности множества уровня $f(x,y,z)=0$ гладкой функции $f:\mathbb{R}^3\to\mathbb{R}$ с регулярным значением $0$.
Покажите, что касательная плоскость $\T_p\Sigma$ перпендикулярна градиенту $\nabla_pf$ в любой точке $p\in\Sigma$.
\end{thm}

{\sloppy

\begin{thm}{Упражнение}\label{ex:vertical-tangent}
Пусть $\Sigma$ --- гладкая поверхность в $\mathbb{R}^3$, и $p\in\Sigma$.
Покажите, что окрестность $p$ в $\Sigma$ является графиком $z=f(x,y)$ гладкой функции $f$, определённой на открытом подмножестве в плоскости $(x,y)$, тогда и только тогда, когда касательная плоскость $\T_p$ не является {}\emph{вертикальной}; то есть если $\T_p$ не перпендикулярна горизонтальной плоскости.
\end{thm}

}

\begin{thm}{Упражнение}\label{ex:tangent-single-point}
Покажите, что если гладкая поверхность $\Sigma$ пересекает плоскость $\Pi$ в единственной точке $p$, то $\Pi$ является касательной к $\Sigma$ в~$p$.
\end{thm}

\section{Производная по направлению}\label{sec:dirder}

В этом разделе мы расширим определение производной по направлению на гладкие функции, определённые на гладких поверхностях.
Сначала вспомним стандартное определение для плоскости.

Даны точка $p\in \mathbb{R}^2$, вектор $\vec w\in \mathbb{R}^2$ и функция $f\:\mathbb{R}^2\to\mathbb{R}$.
Рассмотрим функцию
$h(t)=f(p+t\cdot\vec w)$.
Тогда производная по направлению функции $f$ в точке $p$ вдоль $\vec w$ определяется как \index{10d@$D_{\vec{w}}f$ (производная по направлению)}
\[D_{\vec w}f(p)\df h'(0).\]

Напомним, что функция $f\: \Sigma \to \mathbb{R}$ называется гладкой, если для любой карты $s\: U \to \Sigma$ композиция $f \circ s$ --- гладкая функция.

\begin{thm}{Предложение и определение}\label{def:directional-derivative}
Пусть $f$ --- гладкая функция, определённая на гладкой поверхности $\Sigma$.
Предположим, что $\gamma$ --- гладкая кривая в $\Sigma$, которая выходит из $p$ с вектором скорости $\vec{w}\in \T_p$;
то есть $\gamma(0)=p$ и $\gamma'(0)=\vec{w}$.
Тогда производная $(f\circ\gamma)'(0)$
зависит только от $f$, $p$ и $\vec{w}$;
она называется \index{производная по направлению}\emph{производной по направлению} функции $f$ вдоль $\vec{w}$ в точке $p$
и обозначается как
\[D_{\vec{w}}f,\quad D_{\vec{w}}f(p), \quad\text{или}\quad D_{\vec{w}}f(p)_\Sigma\] 
--- позволяется опускать $p$ и $\Sigma$, если всё ясно из контекста.

Более того, если $(u,v)\mapsto s(u,v)$ --- локальная карта, и $\vec{w}=a\cdot s_u \z+b\cdot s_v$ в точке $p$, то 
\[D_{\vec{w}}f=a\cdot (f\circ s)_u+b\cdot (f\circ s)_v.\]

\end{thm}

Если поверхность $\Sigma$ является плоскостью, то определение согласуется со стандартным.
Действительно, в этом случае $\gamma(t)=p+\vec w\cdot t$ является кривой в $\Sigma$, которая начинается в $p$ с вектором скорости~$\vec{w}$.
Однако, в общем случае точка $p+\vec w\cdot t$ не обязана лежать на поверхности,
и тогда значение $f(p+\vec w\cdot t)$ не определено.
В этом случае стандартное определение не работает.

\parbf{Доказательство.}
Не умаляя общности, можно предположить, что $p\z=s(0,0)$ и кривая $\gamma$ покрыта картой $s$.
В этом случае 
\[\gamma(t)=s(u(t),v(t))\]
для гладких функций $u,v$, определённых в окрестности точки $0$, таких что 
$u(0)\z=v(0)\z=0$.

По правилу дифференцирования сложной функции,
\begin{align*}
\gamma'(0)&=u'(0)\cdot s_u+v'(0)\cdot s_v
\end{align*}
в точке $(0,0)$.
Поскольку $\vec{w}=\gamma'(0)$ и векторы $s_u$, $s_v$ линейно независимы, мы получаем, что $a=u'(0)$ и $b=v'(0)$.

Применяя то же правило, получаем 
\[
(f\circ\gamma)'(0)=a\cdot (f\circ s)_u+b\cdot (f\circ s)_v.
\eqlbl{eq:f-gamma}
\]
в точке $(0,0)$.

Левая часть в \ref{eq:f-gamma} не зависит от выбора карты $s$, а правая зависит только от $p$, $\vec w$, $f$ и~$s$. 
Следовательно, $(f\circ\gamma)'(0)$ зависит только от $p$, $\vec w$ и~$f$.

Последнее утверждение следует из \ref{eq:f-gamma}.
\qeds

\begin{thm}{Продвинутое упражнение}\label{ex:lin-ind-chart}
Пусть $\vec x$ и $\vec y$ --- это векторные поля на гладкой поверхности $\Sigma$.
Предположим, что $\vec x_p$ и $\vec y_p$ линейно независимы в некоторой точке $p\in \Sigma$.
Постройте такие две функции $u$ и $v$ в окрестности точки $p$, что 
\begin{align*}
D_{\vec x} u&>0,
&
D_{\vec y} u&=0,
&
D_{\vec x} v&=0,
&
D_{\vec y} v&>0.
\end{align*}

Постройте на $\Sigma$ карту при $p$, в которой поля $\vec x$ и $\vec y$ касательны к координатным линиям.
\end{thm}

\section{Касательные векторы как функционалы}

В этом разделе даётся более концептуальное определение касательных векторов.
Мы не будем им пользоваться, но лучше о нём знать.

Касательный вектор $\vec w\in \T_p$ к гладкой поверхности $\Sigma$ 
определяет линейный функционал%
\footnote{Термин \emph{функционал} используется для функций, которые принимают функцию в качестве аргумента и возвращают число.} $D_{\vec w}$;
он заглатывает гладкую функцию $\phi$ на $\Sigma$, и выплёвывает  производную по направлению $D_{\vec w}\phi$.
Обратите внимание, что $D_{\vec w}$ подчиняется правилу произведения, то есть
\[D_{\vec w}(\phi\cdot\psi)=(D_{\vec w}\phi)\cdot \psi(p)+\phi(p)\cdot(D_{\vec w}\psi).
\eqlbl{eq:tangent-functional}\]

Нетрудно показать, что касательный вектор $\vec w$ полностью определяется функционалом $D_{\vec w}$.
Более того, касательные векторы в точке $p$ могут \textit{определяться} как линейные функционалы на пространстве гладких функций, которые удовлетворяют тождеству \ref{eq:tangent-functional}.

Иначе говоря, тождество \ref{eq:tangent-functional} улавливает всё о касательных векторах.
Хотя этот подход не очень наглядный, он удобен в доказательствах;
например, \ref{def:directional-derivative} превращается в тавтологию.

\section{Дифференциал}\label{sec:differential}

Напомним, что гладкое отображение из области на плоскости в пространство определено в \ref{sec:Multivariable calculus}.
Пусть $\Sigma_0$ --- гладкая поверхность в~$\mathbb{R}^3$.
Отображение $s\:\Sigma_0\to \mathbb{R}^3$ называется \index{гладкое отображение}\emph{гладим}, если для любой карты $w\:U\to \Sigma_0$, композиция $s\circ w$ задаёт гладкое отображение из $U\to\mathbb{R}^3$.

Любое гладкое отображение $s$ из поверхности $\Sigma_0$ в $\mathbb{R}^3$ можно описать его координатными функциями 
$s(p)=(x(p),y(p),z(p))$. 
Его производная по направлению определяется покоординатно, то есть
\[D_{\vec{w}} s\df(D_{\vec{w}}x,D_{\vec{w}}y,D_{\vec{w}}z).\]

Предположим, что $s$ отображает одну гладкую поверхность $\Sigma_0$ в другую $\Sigma_1$.
Пусть $p_0\in \Sigma_0$ и $p_1=s(p_0)$.
Выберем такую кривую $\gamma_0$ на $\Sigma_0$, что $\gamma_0(0)=p_0$ и $\gamma_0'(0)=\vec w$.
Тогда $\gamma_1= s\circ \gamma_0$ --- гладкая кривая на~$\Sigma_1$. 
По определению производной по направлению, $D_{\vec w} s=\gamma_1'(0)$, и, значит, $D_{\vec w} s\in \T_{p_1}\Sigma_1$ для любого $\vec w\in \T_{p_0}$.

Из \ref{def:directional-derivative} вытекает, что
$\vec w \mapsto D_{\vec w} s$ определяет линейное отображение $\T_{p_0}\Sigma_0\z\to \T_{ p_1}\Sigma_1$;
то есть
\[D_{c\cdot \vec w} s=c\cdot D_{\vec w} s
\quad\text{и}\quad D_{\vec v+ \vec w} s=D_{\vec v} s+ D_{\vec w} s\]
при любых $\vec v, \vec w\in\T_{p_0}$ и $c\in\mathbb{R}$.
Отображение $d_{p_0} s\:\T_{p_0}\Sigma_0\z\to \T_{ p_1}\Sigma_1$, определяемое как
\[d_{p_0} s\:\vec w \mapsto D_{\vec w} s\]
называется \index{дифференциал}\emph{дифференциалом} $s$ при~$p_0$.

Дифференциал $d_{p_0} s$ описывается $2{\times}2$-матрицей $M$ в ортонормированных базисах $\T_{p_0}\Sigma_0$ и $\T_{p_1}\Sigma_1$.
Определим \index{якобиан}\emph{якобиан} $s$ при $p_0$ как $\jac_{p_0} s=|\det M|$; он  
не зависит от выбора ортонормированных базисов в $\T_{p_0}\Sigma_0$ и $\T_{p_1}\Sigma_1$.%
\label{page:|L|}\index{10d@$d_p f$ (дифференциал)}%
\footnote{Про $\jac_{p_0} s$ можно думать более геометрически:
если $R_0$ --- область в $\T_{p_0}$ и $R_1=(d_{p_0} s)(R_0)$ её образ, то
\[\area R_1=\jac_{p_0} s \cdot \area R_0.\]
Это равенство пригодится в определении площади поверхности.}

Если $r\:\Sigma_1\to\Sigma_2$ --- ещё одно гладкое отображение между гладкими поверхностями $\Sigma_1$ и $\Sigma_2$, то
\[d_{p_0}( r\circ s)=d_{p_1} r \circ d_{p_0} s;\]
отсюда
\[\jac_{p_0}( r\circ s)
=
\jac_{p_1} r\cdot\jac_{p_0} s .\eqlbl{eq:jac-composition}\]

Приведённые выше построения применимы к карте $s$;
в этом случае поверхность $\Sigma_0$ является открытой областью в $\mathbb{R}^2$.
Тогда значение $\jac_{p_0} s$ можно найти с помощью формул
\begin{align*}
\jac s
&=|s_u\times s_v|=
\\
&=\sqrt{\langle s_u, s_u\rangle\cdot\langle s_v, s_v\rangle -\langle s_u, s_v\rangle^2}=
\\
&=\sqrt{\det[(\Jac s)^\top\cdot \Jac s ]},
\end{align*}
где $\Jac s$ обозначает матрицу Якоби $s$; см. приложение~\ref{sec:Multivariable calculus}.

\section{Интеграл и площадь}

Пусть $\Sigma$ --- гладкая поверхность, и $h\:\Sigma\to\mathbb{R}$ --- гладкая функция.
Определим интеграл $h$ по борелеву множеству $R\subset \Sigma$;
чаще всего это определение будет применяться к поверхностям с краем.

Напомним, что якобиан $\jac_ps$ определён в предыдущем разделе.
Допустим, что существует карта $(u,v)\mapsto s(u,v)$ для $\Sigma$, определённая на открытом множестве $U\subset\mathbb{R}^2$, такая что $R\subset s(U)$.
В этом случае, определим
\[\iint_R h
\df
\iint_{s^{-1}(R)} h\circ s(u,v)\cdot \jac_{(u,v)}s \cdot du\cdot dv.
\eqlbl{eq:area-def}\]

По правилу замены переменной (\ref{thm:mult-substitution}), правая часть в \ref{eq:area-def} не зависит от выбора~$s$.
То есть, если $s_1\:U_1\to \Sigma$ --- другая карта, такая что $s_1(U_1)\supset R$, то 
\[\iint_{s^{-1}(R)} h\circ s(u,v)\cdot \jac_{(u,v)}s \cdot du\cdot dv=\iint_{s_1^{-1}(R)} h\circ s_1(u,v)\cdot \jac_{(u,v)}s_1 \cdot du\cdot dv.\]
Другими словами, равенство \ref{eq:area-def} можно использовать как определение левой части.

В общем случае, область $R$ можно разбить на счётное число областей $R_1, R_2, \dots$, так что каждая $R_i$ лежит в образе карты
и определить интеграл по $R$ как сумму
\[\iint_Rh
\df
\iint_{R_1}h+\iint_{R_2}h+\dots\]
В случае, если $R$ компактно или же $h \geq 0$, легко проверить, что значение $\iint_Rh$ не зависит от выбора такого разбиения.

Площадь области $R$ на гладкой поверхности $\Sigma$ определяется как интеграл
\[\area R=\iint_R1.\]

Следующее предложение даёт правило замены переменной для интеграла по поверхности.

\begin{thm}{Формула площади}\label{prop:surface-integral}
Пусть $s\:\Sigma_0\to \Sigma_1$ --- диффеоморфизм между гладкими поверхностями.
Рассмотрим область $R\subset \Sigma_0$ и гладкую функцию $f\:\Sigma_1\to\mathbb{R}$.
Предположим, что область $R$ компактна или же функция $f$ неотрицательна.
Тогда 
\[\iint_R (f\circ s)\cdot \jac s=\int_{s(R)}f.\]
В частности, при $f\equiv 1$, получаем
\[\iint_R \jac s=\area (s(R)).\]
\end{thm}

\parbf{Доказательство.}
Следует из \ref{eq:jac-composition} и определения интеграла по поверхности.
\qeds

Пусть $\Sigma_1$ и $\Sigma_2$ --- две гладкие поверхности.
Отображение $f\:\Sigma_1\z\to \Sigma_2$ называется \index{нерастягивающее отображение}\emph{нерастягивающим}, если для любой кривой $\gamma$ на $\Sigma_1$ выполняется $\length\gamma\ge \length (f\circ\gamma)$.
Следующая теорема предоставляет более естественное определение площади.
Несмотря на интуитивную формулировку, доказательство выходит далеко за рамки нашей книжки;
оно основано на обобщении формулы площади для липшицевых отображений \cite[3.2.3]{federer}.

\begin{thm}{Теорема}\label{thm:area-axioms}
{\sloppy
Функционал площади удовлетворяет следующим свойствам:

}

\begin{subthm}{thm:area-axioms:aditivity}
Сигма-аддитивность: 
Пусть $R_1,R_2,\dots$ --- последовательность непересекающихся борелевских множеств на гладкой поверхности.
Тогда 
\[\area (R_1\cup R_2\cup \dots)=\area R_1+\area R_2+{}\dots\]
\end{subthm}

{\sloppy

\begin{subthm}{thm:area-axioms:monotonicity}
Монотонность:
Пусть $f\:\Sigma_1\to \Sigma_2$ --- нерастягивающее отображение между двумя гладкими поверхностями.
Предположим, что $R_1\subset \Sigma_1$ и $R_2\subset \Sigma_2$ --- такие борелевские множества, что $f(R_1)\supset R_2$.
Тогда 
\[\area R_1\ge \area R_2.\]
\end{subthm}

}

\begin{subthm}{thm:area-axioms:unit}
Единичный квадрат имеет единичную площадь. 
\end{subthm}

Более того, функционал площади однозначно определяется этими свойствами на всех борелевских множествах.
\end{thm}

\parit{Замечание.}
Понятия площади и длины схожи.
Однако, длина кривой определялась, используя другую идею (через верхнюю грань длин вписанных ломаных).
Известно, что аналогичное определение не работает даже для очень простых поверхностей.
Это видно из классического примера --- так называемого \textit{сапога Шварца}.
Этот пример и различные подходы к понятию площади обсуждаются в популярной статье Владимира Дубровского~\cite{dubrovsky}.

\section{Вектор нормали и ориентация}

Поверхность $\Sigma$ называется \index{ориентированная поверхность}\emph{ориентированной}, если она снабжена полем  \index{10nu@$\Norm$ (поле нормалей)}\index{нормаль}\index{поле нормалей}\emph{нормалей} $\Norm$;
то есть таким непрерывным отображением $p\mapsto \Norm(p)$, что $\Norm(p)\perp\T_p$ и $|\Norm(p)|=1$ для любого~$p$.
Выбор поля $\Norm$ называется {}\emph{ориентацией}~$\Sigma$.
Поверхность $\Sigma$ называется {}\emph{ориентируемой}, если её можно ориентировать.
Каждая ориентируемая поверхность допускает две ориентации: $\Norm$ и $-\Norm$.

Пусть $\Sigma$ --- гладкая ориентированная поверхность и $\Norm$ --- её поле нормалей.
Отображение $\Norm\:\Sigma\to \mathbb{S}^2$, определяемое как $p\mapsto \Norm(p)$, называется \index{сферическое отображение}\emph{сферическим} или \index{гауссово отображение}\emph{гауссовым}.

\begin{wrapfigure}{r}{42 mm}
\vskip-7mm
\centering
\includegraphics{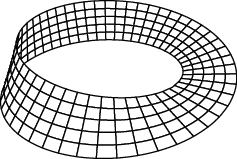}
\vskip-1mm
\end{wrapfigure}

Для поверхностей сферическое отображение играет ту же роль, что и касательная индикатриса для кривых.

Лента Мёбиуса, показанная на рисунке, даёт пример неориентируемой поверхности --- поле нормалей нельзя продолжить непрерывно вдоль средней линии ленты (при попытке обойти ленту вокруг нормаль меняет знак).

Каждая поверхность локально ориентируема.
Более того, каждая карта $s(u,v)$ допускает ориентацию 
\[\Norm=
\frac{s_u\times s_v}
{\left|s_u\times s_v\right|}.\]
Действительно, векторы $s_u$ и $s_v$ являются касательными векторами в точке $p$.
Значит, их векторное произведение перпендикулярно касательной плоскости,
и оно не обращается в ноль, так как $s_u$ и $s_v$ линейно независимы.
Ясно, что отображение $(u,v)\mapsto \Norm(u,v)$ непрерывно.
Следовательно, $\Norm$ --- это поле единичных нормалей.

{\sloppy

\begin{thm}{Упражнение}\label{ex:const-normal}
Предположим, что гладкая поверхность $\Sigma$ имеет постоянное поле нормалей $\nu_0$.
Покажите, что $\Sigma$ лежит в плоскости, перпендикулярной к $\nu_0$.
\end{thm}

}

\begin{thm}{Упражнение}\label{ex:implicit-orientable}
Пусть $0$ --- регулярное значение гладкой функции
$h:\mathbb{R}^3\to\mathbb{R}$.
Покажите, что компонента связности множества уровня $h(x,y,z)=0$ есть ориентируемая поверхность.
\end{thm}

Напомним, что любая собственная поверхность $\Sigma$ без края в евклидовом пространстве делит его на две связные компоненты (\ref{clm:proper-divides}).
Следовательно, мы можем выбрать поле единичных нормалей на $\Sigma$, которое указывает в одну из компонент дополнения; таким образом, мы получаем следующее.

\begin{thm}{Замечание}
Любая гладкая собственная поверхность без края в евклидовом пространстве является ориентируемой.
\end{thm}

В частности, ленту Мёбиуса невозможно продолжить до собственной гладкой поверхности без края.

\section{Сечения}

\begin{thm}{Лемма}\label{lem:reg-section}
Пусть $\Sigma$ --- гладкая поверхность, $f\:\Sigma\to\mathbb{R}$ --- гладкая функция, и $r_0\in\mathbb{R}$.
Тогда найдётся $r\in\mathbb{R}$ произвольно близкое к~$r_0$ такое, что каждая компонента связности линии уровня $L_r\z=\set{x\in\Sigma}{f(x)=r}$ является гладкой кривой.
\end{thm}

\parbf{Доказательство.}
Поверхность $\Sigma$ можно покрыть счётным набором карт $s_i\:U_i\to \Sigma$.
При этом композиция $f\circ s_i$ является гладкой функцией для любого $i$.
По лемме Сарда (\ref{lem:sard}), почти все $r\in \mathbb{R}$ являются регулярным значением каждой композиции $f\circ s_i$.

Зафиксируем такое значение $r$, достаточно близкое к $r_0$, и рассмотрим его линию уровня $L_r\subset\Sigma$.
Любая точка в $L_r$ лежит в образе одной из карт.
Следовательно, эта точка имеет окрестность, которая является гладкой кривой;
отсюда результат.
\qeds

\begin{thm}{Продвинутое упражнение}\label{ex:plane-section}
Пусть $\Pi$ --- координатная плоскость $z=0$, и $A \subset \Pi$ --- любое замкнутое подмножество.
Постройте такую открытую гладкую поверхность $\Sigma$, что $\Sigma \cap \Pi = A$.
\end{thm}

Упражнение выше показывает, что сечения гладких поверхностей бывают довольно противными.
Следствие ниже позволяет подсдвинуть плоскость, чтобы сечение стало приемлемым.

\begin{thm}{Следствие}
Пусть $\Sigma$ --- гладкая поверхность.
Тогда для любой плоскости $\Pi$ существует параллельная плоскость $\Pi^{*}$, лежащая сколь угодно близко к $\Pi$ и такая, что пересечение $\Sigma\cap\Pi^{*}$ есть объединение непересекающихся гладких кривых.
\end{thm}

\chapter{Кривизны}
\label{chap:surface-curvature}

Из предыдущей главы мы знаем, что касательная плоскость к поверхности имеет с ней касание первого порядка.
В частности, вплоть до первого порядка все поверхности выглядят одинаково вблизи одной точки.
Для приближений второго порядка появляются различия, о которых и пойдёт речь.

\section{Касательно-нормальные координаты}
\label{sec:lmn}
\index{касательно-нормальные координаты}

{\sloppy

Выберем точку $p$ на гладкой ориентированной поверхности~$\Sigma$;
пусть $\Norm$ --- её поле нормалей.
Рассмотрим систему координат $(x,y,z)$ с началом в точке $p$, чтобы касательная плоскость $\T_p$ стала горизонтальной, а ось $z$ была направлена по $\Norm(p)$.
Такие координаты будут называться {}\emph{касательно-нормальными}. 
Согласно \ref{ex:vertical-tangent},
окрестность $p$ в $\Sigma$ можно задать графиком $z=f(x,y)$ гладкой функции,
при этом
\begin{align*}
f(0,0)&=0,
&
f_x(0,0)&=0,
&
f_y(0,0)&=0.
\end{align*}
Первое равенство означает, что $p=(0,0,0)$ лежит на графике, а два других --- то, что касательная плоскость в $p$ горизонтальна.

}

\begin{wrapfigure}[7]{o}{42 mm}
\vskip-3mm
\centering
\includegraphics{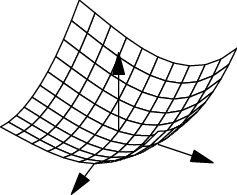}
\vskip-3mm
\end{wrapfigure}

Пусть
\begin{align*}
\ell&=f_{xx}(0,0),
\\
m&=f_{xy}(0,0)=f_{yx}(0,0),
\\
n&=f_{yy}(0,0).
\end{align*}
\textit{Ряд Тейлора} $f$ в $(0,0)$ до второго порядка записывается как
\[f(x,y)=\tfrac12(\ell\cdot x^2+2\cdot m\cdot x\cdot y+n\cdot y^2)+o(x^2+y^2).\phantom{f(x,y)=\tfrac12(\ell\cdot x^2+}\]
Значения $\ell$, $m$ и $n$ полностью определяются этим уравнением.\index{10lmn@$\ell$, $m$, $n$ (компоненты матрицы Гессе)}
График многочлена Тейлора 
\[z=\tfrac12\cdot(\ell\cdot x^2+2\cdot m\cdot x\cdot y+n\cdot y^2)\]
называется \index{соприкасающийся параболоид}\emph{соприкасающимся параболоидом}.
Он имеет \index{порядок касания}\emph{второй порядок касания} с $\Sigma$ в~$p$.

Далее заметим, что 
\[\ell\cdot x^2+2\cdot m\cdot x\cdot y+n\cdot y^2=\langle M_p\cdot (\begin{smallmatrix}
x\\y
\end{smallmatrix}), (\begin{smallmatrix}
x\\y
\end{smallmatrix})\rangle,\]
где 
\[M_p=\begin{pmatrix}
 \ell
 &m
 \\
 m
 &n
 \end{pmatrix}
=\begin{pmatrix}
 f_{xx}(0,0)
 &f_{xy}(0,0)
 \\
 f_{yx}(0,0)
 &f_{yy}(0,0)
 \end{pmatrix}.
\eqlbl{eq:hessian}
\]
--- так называемая \index{матрица Гессе}\emph{матрица Гессе} $f$ в точке $(0,0)$,\index{10m@$M_p$ (матрица Гессе)}

\section{Главные кривизны}\label{sec:Principal curvatures}

Касательно-нормальные координаты почти однозначно определяются точкой на поверхности;
они единственны с точностью до вращения плоскости $(x,y)$.
При вращении этой плоскости  матрица Гессе $M_p$ переписывается в новом базисе.

Поскольку матрица $M_p=\left(\begin{smallmatrix}
 \ell
 &m
 \\
 m
 &n
 \end{smallmatrix}\right)$ 
симметрична, по спектральной теореме (\ref{thm:spectral}), её можно диагонализовать ортогональной матрицей.
То есть, вращая плоскость $(x,y)$, можно добиться того, что $m=0$.
В этом случае,
\[M_p=\begin{pmatrix}
 k_1
 &0
 \\
 0
 &k_2
 \end{pmatrix},
\]
числа $k_1$ и $k_2$  на диагонали называются \index{главные кривизны и направления}\emph{главными кривизнами} $\Sigma$ в точке $p$;\index{10k@$k_1$, $k_2$ (главные кривизны)}
их можно обозначать $k_1(p)$ и $k_2(p)$, или $k_1(p)_\Sigma$ и $k_2(p)_\Sigma$, если нужно подчеркнуть, что они вычисляются в точке $p$ для поверхности~$\Sigma$.
Если не указано обратное, мы будем считать, что 
\[k_1\le k_2.\]

{\sloppy

Если в этих координатах $\Sigma$ локально задаётся графиком $z\z=f(x,y)$, то 
\[f(x,y)=\tfrac12\cdot(k_1\cdot x^2+k_2\cdot y^2)+o(x^2+y^2).\]

}

Главные кривизны можно также определить как собственные значения матрицы $M_p$.
Собственные направления $M_p$ называются {}\emph{главными направлениями} $\Sigma$ в точке~$p$.
Если $k_1(p)\ne k_2(p)$, то в точке $p$ есть ровно два главных направления, которые перпендикулярны друг другу;
если же $k_1(p) = k_2(p)$, то все касательные направления в точке $p$ главные.

Если обратить ориентацию $\Sigma$, то главные кривизны сменят свои знаки и индексы:
$k_1$ превратится в $-k_2$, а $k_2$ в $-k_1$.

Гладкая кривая на поверхности $\Sigma$, всё время идущая в главных направлениях, называется \index{линия кривизны}\emph{линией кривизны}~$\Sigma$.
Если $k_1(p)\z\ne k_2(p)$, то в точке $p$ существует карта $(u,v)\mapsto s(u,v)$ с координатными линиями, образованными линиями кривизны (это вытекает из \ref{ex:lin-ind-chart}), при этом $s_u\perp s_v$, ведь главные направления перпендикулярны.

\begin{thm}{Упражнение}\label{ex:line-of-curvature}
Пусть $\Sigma$ --- гладкая поверхность, зеркально-симметричная относительно плоскости $\Pi$.
Допустим, что $\Sigma$ и $\Pi$ пересекаются по гладкой кривой~$\gamma$.
Покажите, что $\gamma$ --- линия кривизны.
\end{thm}

\section{Ещё кривизны}\label{sec:More curvatures}

Пусть $p$ --- точка ориентированной гладкой поверхности~$\Sigma$.
Напомним, что $k_1(p)$ и $k_2(p)$ обозначают главные кривизны в точке~$p$.

Произведение 
\[K(p)=k_1(p)\cdot k_2(p)\]
называется \index{10k@$K$ (гауссова кривизна)}\index{гауссова кривизна}\emph{гауссовой кривизной} в точке~$p$.
Её можно обозначать $K$, $K(p)$ и даже $K(p)_\Sigma$, если хочется указать и точку и поверхность.

Поскольку определитель равен произведению собственных значений, получаем, что
\[K=\ell\cdot n-m^2,\]
где 
$M_p=
(\begin{smallmatrix}
\ell&m
\\
m&n
\end{smallmatrix}
)$
--- матрица Гессе.

Отметим, что обращение ориентации $\Sigma$ не меняет гауссову кривизну,
и, значит, она определена и для неориентированных поверхностей.

\begin{thm}{Упражнение}\label{ex:gauss+orientable}
Покажите, что любая поверхность с положительной гауссовой кривизной ориентируема. 
\end{thm}

Сумма \index{10h@$H$ (средняя кривизна)}
\[H(p)=k_1(p)+ k_2(p)\] 
называется \index{кривизна}\index{средняя кривизна}\emph{средней кривизной}%
\footnote{Некоторые авторы определяют её как $\tfrac12\cdot(k_1(p)+ k_2(p))$ --- среднее значение главных кривизн.
Это лучше подходит к названию, но оказывается неудобно.}
в точке~$p$.
Если нужно, её можно обозначать как $H(p)_\Sigma$.
Можно думать, что средняя кривизна --- это след матрицы Гессе $M_p=
(\begin{smallmatrix}
\ell&m
\\
m&n
\end{smallmatrix}
)$;
то есть
\[H=\ell+n.\] 
Обращение ориентации $\Sigma$ меняет знак её средней кривизны.

Поверхность с нулевой средней кривизной называется \index{минимальная поверхность}\emph{минимальной}.

\begin{thm}{Упражнение}\label{ex:re-scale-surface-curvature}
Пусть $\Sigma$ --- ориентированная поверхность, а $\Sigma_{\lambda}$ --- её гомотетия с коэффициентом $\lambda > 0$; то есть $\Sigma_{\lambda}$ состоит из точек $\lambda \cdot x$, где $x \in \Sigma$.
Покажите, что
\[K(\lambda\cdot p)_{\Sigma_{\lambda}}
= \tfrac{1}{\lambda^2}\cdot K(p)_{\Sigma}
\quad\text{и}\quad
H(\lambda \cdot p)_{\Sigma_{\lambda}} = \tfrac1\lambda\cdot H(p)_{\Sigma}\]
для любой точки $p\in \Sigma$.  
\end{thm}

\section{Оператор формы}\label{sec:shape}

В следующих определениях используется понятие производной по направлению и дифференциала, см.~\ref{sec:dirder} и~\ref{sec:differential}.

Пусть $p$ --- точка гладкой поверхности $\Sigma$ с ориентацией, заданной полем нормалей $\Norm$,
и пусть $\vec w$ --- касательный вектор при $p$.
\index{оператор формы}\emph{Оператор формы}%
\footnote{
Следующие билинейные формы на касательной плоскости  
\begin{align*}
\mathrm{I}(\vec v,\vec w)&=\langle\vec v,\vec w\rangle,
&
\mathrm{II}(\vec v,\vec w)&=\langle\Shape\vec v,\vec w\rangle,
&
\mathrm{III}(\vec v,\vec w)&=\langle\Shape\vec v,\Shape\vec w\rangle
\end{align*}
называются \index{фундаментальная форма}\emph{первой, второй и третьей фундаментальными формами}, соответственно.
Исторически, эти формы были введены до оператора формы, но мы вовсе не станем их рассматривать.
} от $\vec w$ определяется как
\[\Shape_p\vec w=-D_{\vec w}\Norm.\]
Его также можно определить как
\[\Shape=-d\Norm,\eqlbl{eq:shape=-L}\] 
где $d\Norm$ --- дифференциал сферического отображения $\Norm\:\Sigma\to\mathbb{S}^2$; то есть $d_p\Norm(\vec v)=(D_{\vec v}\Norm)(p)$.

Напомним, что $d_p\Norm$ --- это линейное отображение $\T_p\Sigma\to \T_{\Norm(p)}\mathbb{S}^2$.
Заметим, что $\T_p\Sigma$ совпадает с $\T_{\Norm(p)}\mathbb{S}^2$, ведь обе плоскости перпендикулярны $\Norm(p)$.
Поэтому $\Shape_p$ действительно является линейным оператором $\T_p\Sigma\to \T_p\Sigma$ (это также будет следовать из \ref{thm:shape-chart}).

Оператор формы от касательного вектора $\vec w\z\in \T_p\Sigma$ может обозначаться как $\Shape\vec w$, если из контекста ясно, о какой базовой точке $p$ и о какой поверхности $\Sigma$ идет речь;
в противном случае можно писать 
\[\Shape_p(\vec w)\quad\text{или}\quad \Shape_p(\vec w)_\Sigma.\]

\begin{thm}{Теорема}\label{thm:shape-chart}
Пусть $(u,v)\mapsto s(u,v)$ --- гладкое отображение в гладкую поверхность $\Sigma$ с полем нормалей $\Norm$.
Тогда 
\begin{align*}
\langle \Shape(s_u), s_u\rangle 
&=\langle s_{uu},\Norm\rangle,
&
\langle \Shape(s_v), s_u\rangle 
&=\langle s_{uv},\Norm\rangle,
\\
\langle \Shape(s_u), s_v\rangle 
&=\langle s_{uv},\Norm\rangle,
&
\langle \Shape(s_v), s_v\rangle 
&=\langle s_{vv},\Norm\rangle,
\\
\langle \Shape(s_u), \Norm\rangle 
&=0,
&
\langle \Shape(s_v), \Norm\rangle 
&=0
\end{align*}
при любых $(u,v)$.

\end{thm}

\enlargethispage{2\baselineskip}

\parbf{Доказательство.}
Будем использовать сокращение $\Norm=\Norm(u,v)$ для $\Norm(s(u,v))$,
так что 
\[
\begin{aligned}
\Shape(s_u)&=-D_{s_u}\Norm=-\Norm_u,
&
\Shape(s_v)&=-D_{s_v}\Norm=-\Norm_v.
\end{aligned}
\eqlbl{eq:shape=norm_u}
\]

Поскольку $\Norm$ --- единичный вектор, ортогональный к $s_u$ и $s_v$,
\begin{align*}
\langle \Norm,s_u\rangle&\equiv0,
&
\langle \Norm,s_v\rangle&\equiv0,
&
\langle \Norm,\Norm\rangle&\equiv1.
\end{align*}
Взяв частные производные этих тождеств, получим
\begin{align*}
\langle \Norm_u,s_u\rangle+\langle \Norm,s_{uu}\rangle&=0,
&
\langle \Norm_v,s_u\rangle+\langle \Norm,s_{uv}\rangle&=0,
\\
\langle \Norm_u,s_v\rangle+\langle \Norm,s_{uv}\rangle&=0,
&
\langle \Norm_v,s_v\rangle+\langle \Norm,s_{vv}\rangle&=0,
\\
2\cdot\langle \Norm_u,\Norm\rangle&=0,
&
2\cdot\langle \Norm_v,\Norm\rangle&=0.
\end{align*}
Остаётся подставить сюда выражения из \ref{eq:shape=norm_u}.
\qeds

\begin{thm}{Упражнение}\label{ex:self-adjoint}
Покажите, что оператор формы \index{самосопряжённый оператор}\emph{самосопряжён}; то есть
\[\langle \Shape\vec u,\vec v\rangle=\langle \vec u,\Shape\vec v\rangle\]
при любых $\vec u,\vec v\in\T_p$.
\end{thm}

Напомним, что компоненты $\ell$, $m$ и $n$ матрицы Гессе определены в~\ref{sec:lmn}.

\begin{thm}{Следствие}\label{cor:Shape(ij)}
Пусть $z=f(x,y)$ --- локальное представление гладкой поверхности $\Sigma$ в касательно-нормальных координатах при точке~$p$,
и пусть $(\begin{smallmatrix}
\ell&m\\ m&n
\end{smallmatrix})$ --- её матрица Гессе.
Тогда 
\begin{align*}
\Shape\vec i&=\ell\cdot \vec i+m\cdot \vec j,
&
\Shape\vec j&=m\cdot \vec i+n\cdot\vec j,
\end{align*}
где $\vec i,\vec j,\vec k$\index{10i@$\vec i$, $\vec j$, $\vec k$ (стандартный базис)} --- стандартный базис в $\mathbb{R}^3$.
Иначе говоря, оператор формы задаётся умножением на матрицу Гессе.
\end{thm}

Это следствие обнажает связь между кривизнами и оператором формы.
Главные кривизны $\Sigma$ в точке $p$ --- это собственные значения $\Shape_p$, главные направления --- это собственные направления $\Shape_p$, гауссова кривизна --- это определитель $\Shape_p$, а средняя кривизна --- это след $\Shape_p$.

Поскольку матрица Гессе симметрична, из следствия получаем, что оператор формы самосопряжён;
это решает упражнение \ref{ex:self-adjoint}.

\parbf{Доказательство.}
Заметим, что $s\:(u,v)\mapsto (u,v,f(u,v))$ --- карта на $\Sigma$ при~$p$.
Кроме того, 
\begin{align*}
s_u(0,0)&=\vec i,
&
s_v(0,0)&=\vec j,
&
\Norm(0,0)&=\vec k,
\\
s_{uu}(0,0)&=\ell\cdot \vec k,
&
s_{uv}(0,0)&=m\cdot \vec k,
&
s_{vv}(0,0)&=n\cdot \vec k.
\end{align*}
Остаётся применить \ref{thm:shape-chart}.
\qeds

\begin{thm}{Следствие}\label{cor:intK}
Пусть $\Sigma$ --- гладкая поверхность гауссова кривизна которой не обращается в ноль.
Допустим, что её сферическое отображение $\Norm\:\Sigma\to\mathbb{S}^2$ инъективно.
Тогда 
\[\iint_\Sigma|K|=\area(\Norm(\Sigma)).\]
\end{thm}

{\sloppy

\parbf{Доказательство.}
Выберем ортонормированный базис в $\T_p$, в главных направлениях.
Тогда оператор формы выражается матрицей 
$(\begin{smallmatrix}
 k_1
 &0
 \\
 0
 &k_2
\end{smallmatrix})$.

}

Поскольку $\Shape_p=-d_p\Norm$, из \ref{cor:Shape(ij)} получаем, что
\[\jac_p\Norm
=
|\det(\begin{smallmatrix}
k_1
&0
\\
0
&k_2
\end{smallmatrix})|
=|K(p)|\ne 0\]
для любой точки $p\in\Sigma$.
В частности отображение $\Norm$ регулярно
и, поскольку оно инъективно, задаёт диффеоморфизм  $\Sigma\z\to\Norm(\Sigma)$.
Остаётся применить формулу площади (\ref{prop:surface-integral}).
\qeds

\enlargethispage{2\baselineskip}

\begin{thm}{Упражнение}\label{ex:normal-curvature=const}
Пусть $\Sigma$ --- гладкая поверхность с полем нормалей $\Norm$, и её главные кривизны  равны $1$ во всех точках.

\begin{subthm}{ex:normal-curvature=const:a}
Докажите, что $\Shape_p(\vec w)=\vec w$ при любых $p\in\Sigma$ и $\vec w\in \T_p\Sigma$.
\end{subthm}

\begin{subthm}{ex:normal-curvature=const:b}
Докажите, что точка $q\z= p+\Norm(p)$ не зависит от выбора $p\in\Sigma$.
Выведите отсюда, что $\Sigma$ является подмножеством единичной сферы с центром в~$q$.
\end{subthm}

\end{thm}

{\sloppy

\begin{thm}{Продвинутое упражнение}\label{ex:normal-curvature=0}
Пусть $\Sigma$ --- гладкая поверхность, $\Norm$ --- её поле нормалей, и $Z_0\subset \Sigma$ --- связное множество с нулевым оператором формы. 
Докажите, что $Z_0$ лежит в плоскости.
\end{thm}

}

{}\emph{Угол} между двумя ориентированными поверхностями в точке их пересечения $p$ определяется как угол между их нормалями в~$p$.

\begin{thm}{Упражнение}\label{ex:shape-curvature-line}
Предположим, что две гладкие ориентированные поверхности $\Sigma_1$ и $\Sigma_2$ пересекаются по гладкой кривой~$\gamma$ под постоянным углом.
Пусть $\gamma$ --- линия кривизны на $\Sigma_1$.
Докажите, что $\gamma$ также является линией кривизны на $\Sigma_2$.

Выведите отсюда, что если гладкая поверхность $\Sigma$ пересекает плоскость или сферу вдоль гладкой кривой $\gamma$ под постоянным углом, то $\gamma$ является линией кривизны на~$\Sigma$.
\end{thm}

\begin{thm}{Упражнение}\label{ex:equidistant}
Пусть $\Sigma$ --- замкнутая гладкая поверхность с ориентацией, определённой полем нормалей $\Norm$.

\begin{subthm}{ex:equidistant:smooth}
Докажите, что множество 
\[\Sigma_t=\set{p+t\cdot \Norm(p)}{p\in\Sigma}\] 
является гладкой поверхностью при всех $t$ достаточно близко к нулю.
\end{subthm}

\begin{subthm}{ex:equidistant:area}
Докажите, что для всех $t$, достаточно близких к нулю, справедливо равенство
\[\area\Sigma_t=\area\Sigma-t\cdot \iint_\Sigma H+t^2\cdot \iint_\Sigma K,\]
где $H$ и $K$ обозначают среднюю и гауссову кривизну~$\Sigma$.
\end{subthm}

\end{thm}

\begin{thm}{Продвинутое упражнение}\label{ex:flat-plane}{\sloppy
Пусть $\Sigma$ --- гладкая ориентированная поверхность, параметризованная координатным прямоугольником $(u,v)\mapsto s(u,v)$,
и пусть $\vec u=\tfrac{s_u}{|s_u|}$, $\vec v=\tfrac{s_v}{|s_v|}$ и $\Norm(u,v)$ --- нормаль в точке $s(u,v)$.

}

Предположим, что $\vec u$ и $\vec v$ --- поля главных направлений;
пусть $0$ и $k\z=k(u,v)$ будут их главными кривизнами,
при этом $k$ не обращается в ноль и $|s_u|=1$ на одной из $u$-координатных линий.

\begin{subthm}{ex:flat-plane:orthonormal}
Докажите, что $\Norm(u,v)$, $\vec u(u,v)$ и $\vec v(u,v)$ образуют ортонормированный базис при любых $(u,v)$.
\end{subthm}

\begin{subthm}{ex:flat-plane:depend}
Докажите, что $\Norm(u,v)$, $\vec u(u,v)$ и $\vec v(u,v)$ зависят только от $v$.
Выведите отсюда, что $u$-координатные линии являются отрезками прямых.
\end{subthm}

\begin{subthm}{ex:flat-plane:depend-u}
Докажите, что $s_{uu}$ пропорционален $s_u$ во всех точках.
Выведите отсюда, что $|s_u|=1$ во всех точках.
\end{subthm}

\begin{subthm}{ex:flat-plane:linear}
Докажите, что при фиксированном $v_0$ значение $\tfrac1{k(u,v_0)}$ линейно зависит от $u$.
\end{subthm}

\end{thm}

Упражнение содержит основную часть доказательства следующей теоремы:
\textit{Любая открытая поверхность с нулевой гауссовой кривизной является \index{цилиндрическая поверхность}\emph{цилиндрической}}\,;
то есть она заметается семейством параллельных прямых.
Упражнения \ref{ex:lin-ind-chart} и \ref{ex:line-cylinder} иллюстрируют другие части этого доказательства.
Упомянутая теорема была доказана Алексеем Погореловым \cite[II §3 теорема 2]{pogorelov1956} в гораздо более общих предположениях и позже была переоткрыта пару раз \cite{hartman-nirenberg,massey1962}.

\chapter{Кривые на поверхности}\label{chap:Curves in a surface}

\section{Базис Дарбу}\label{sec:Darboux}

\begin{wrapfigure}{r}{42 mm}
\vskip-20mm
\centering
\begin{lpic}[t(-0mm),b(0mm),r(0mm),l(0mm)]{asy/paraboloid+curve(1)}
\lbl[ul]{34,14;$\tan$}
\lbl[b]{20,43;$\Norm$}
\lbl[bl]{38,35;$\mu$}
\end{lpic}
\vskip-0mm
\end{wrapfigure}

Пусть $\gamma$ --- гладкая кривая на гладкой ориентированной поверхности~$\Sigma$,
с полем нормалей $\Norm$.
Будем считать, что $\gamma$ параметризована длиной, а значит, $\tan(s)\z=\gamma'(s)$ --- её касательная индикатриса.
Далее будем пользоваться сокращением $\Norm(s)\z=\Norm(\gamma(s))$.

Единичные векторы $\tan(s)$ и $\Norm(s)$ ортогональны.
Поэтому существует единственный единичный вектор $\mu(s)$\index{10tmn@$\tan$, $\mu$, $\Norm$ (базис Дарбу)}, такой что 
$\tan(s),\mu(s),\Norm(s)$ образуют ориентированный ортонормированный базис.
Он называется \index{базис Дарбу}\emph{базис Дарбу} кривой $\gamma$ на~$\Sigma$.

Поскольку $\T_{\gamma(s)}\z\perp\Norm(s)$, вектор $\mu(s)$ касается $\Sigma$ в точке $\gamma(s)$.
На самом деле, $\mu(s)$ можно получить повернув $\tan(s)$ в $\T_{\gamma(s)}$ против часовой стрелки на угол $\tfrac\pi2$.
Этот вектор также можно определить через векторное произведение $\mu(s)\z\df\Norm(s)\times \tan(s)$.

Поскольку у $\gamma$ единичная скорость, $\gamma''\perp \gamma'$ (см. \ref{prop:a'-pertp-a''}).
В частности, $\gamma''$ можно записать как линейную комбинацию $\mu$ и $\Norm$;
то есть \index{10k@$k_g$ (геодезическая кривизна)}\index{10k@$k_n$ (нормальная кривизна)}
\[\gamma''(s)=k_g(s)\cdot \mu(s)+k_n(s)\cdot\Norm(s).\]
Величины $k_g(s)$ и $k_n(s)$ называются соответственно \index{кривизна}\index{геодезическая!кривизна}\emph{геодезической} и \index{нормальная кривизна}\emph{нормальной кривизной} кривой $\gamma$ при $s$.
Поскольку базис $\tan(s)$, $\mu(s)$, $\Norm(s)$ ортонормирован, эти кривизны переписываются через скалярные произведения
\begin{align*}
k_g(s)&=\langle \gamma''(s),\mu(s)\rangle= 
&
k_n(s)&=\langle \gamma''(s),\Norm(s)\rangle=
\\
&=\langle \tan'(s),\mu(s)\rangle,
&
&=\langle \tan'(s),\Norm(s)\rangle.
\end{align*}

Продифференцировав $0=\langle \tan(s),\Norm(s)\rangle$, получим
\begin{align*}
0&=\langle \tan(s),\Norm(s)\rangle'=
\\
&=\langle \tan'(s),\Norm(s)\rangle+\langle \tan(s),\Norm'(s)\rangle=
\\
&=k_n(s)+\langle \tan(s),D_{\tan(s)}\Norm\rangle.
\end{align*}
И, применив определение оператора формы,
получим следующее.

\begin{thm}{Предложение}\label{prop:normal-shape}
Пусть $\gamma$ --- гладкая кривая с единичной скоростью на гладкой поверхности~$\Sigma$, $p=\gamma(s_0)$ и $\vec v=\gamma'(s_0)$.
Тогда 
\[k_n(s_0)=\langle \Shape_p(\vec v),\vec v\rangle,\]
где $k_n$ --- нормальная кривизна $\gamma$ при $s_0$, а $\Shape_p$ --- оператор формы в~точке~$p$.
\end{thm}

Согласно предложению, нормальная кривизна гладкой кривой на $\Sigma$ полностью определяется её вектором скорости~$\vec v$.
Поэтому нормальную кривизну в направлении $\vec v$ допустимо обозначать как $k_{\vec v}$\index{10k@$k_{\vec v}$ (нормальная кривизна)} и
\[k_{\vec v}=\langle \Shape_p(\vec v),\vec v\rangle\]
для любого единичного вектора $\vec v$ в $\T_p$.

\section{Формула Эйлера}

Пусть $p$ --- точка гладкой поверхности~$\Sigma$.
Выберем касательно-нормальные координаты при $p$ с диагональной матрицей Гессе
\[M_p=\begin{pmatrix}
 k_1(p)
 &0
 \\
 0
 &k_2(p)
 \end{pmatrix}.
\]
Рассмотрим вектор ${\vec w}=a\cdot\vec i+b\cdot\vec j$ касательный в $p$.
Согласно \ref{cor:Shape(ij)},
\[
\langle \Shape\vec w,\vec w\rangle
=a^2\cdot k_1(p) +b^2\cdot k_2(p). 
\]
Поскольку $|{\vec w}|^2=a^2+b^2$, получаем следующее.

\begin{thm}{Наблюдение}\label{obs:k1-k2}
Для любой точки $p$ на ориентированной гладкой поверхности $\Sigma$,
главные кривизны $k_1(p)$ и $k_2(p)$ являются соответственно минимумом и максимумом нормальных кривизн в~точке $p$.
Более того, если $\theta$ --- угол между единичным вектором ${\vec w}\in\T_p$ и первым главным направлением, то 
\[k_{\vec w}(p)=k_1(p)\cdot(\cos\theta)^2+k_2(p)\cdot(\sin\theta)^2.\]

\end{thm}

Последнее тождество называется \index{формула Эйлера}\emph{формулой Эйлера}.

\begin{thm}{Упражнение}\label{ex:mean-curvature}
Покажите, что средняя кривизна $H(p)$ в точке~$p$ на гладкой поверхности равна сумме нормальных кривизн для любой пары ортогональных касательных направлений в $p$. 
\end{thm}

\begin{thm}{Упражнение}\label{ex:average}
Покажите, что $\tfrac38\cdot H(p)^2-\tfrac12\cdot K(p)$ равно среднему значению квадрата нормальных кривизн при $p$;
то есть среднему значению $k_{\vec w}^2$ для всех единичных векторов ${\vec w}\in\T_p$.
\end{thm}

{\sloppy

\begin{thm}{Теорема Мёнь\'{е}}
\label{thm:meusnier}
\index{теорема Мёнье}
Пусть $\gamma$ --- гладкая кривая на гладкой ориентированной поверхности~$\Sigma$, $p=\gamma(t_0)$, ${\vec v}\z=\gamma'(t_0)$ и $\alpha\z=\measuredangle(\Norm(p),\norm(t_0))$;
то есть $\alpha$ --- угол между нормалью к $\Sigma$ в точке $p$ и нормальным вектором в базисе Френе кривой $\gamma$ при~$t_0$.
Тогда 
\[\kur(t_0)\cdot\cos\alpha=k_{n}(t_0);\]
здесь $\kur(t_0)$ и $k_n(t_0)$ --- кривизна и нормальная кривизна $\gamma$ при $t_0$, соответственно. 
\end{thm}

}

\parbf{Доказательство.}
Так как $\gamma''=\tan'=\kur\cdot \norm$, 
\begin{align*}
k_{n}(t_0)&=\langle\gamma''(t_0),\Norm(p)\rangle=
\\
&=\kur(t_0)\cdot\langle\norm(t_0),\Norm(p)\rangle=
\\
&=\kur(t_0)\cdot\cos\alpha.
\end{align*}
\qedsf

\begin{thm}{Упражнение}\label{ex:meusnier}
Пусть ${\vec v}\z\in \T_p\Sigma$ --- единичный касательный вектор к гладкой поверхности $\Sigma$ в точке $p\in\Sigma$.
Допустим, что нормальная кривизна $k_{\vec v}(p)$ не равна нулю.

Покажите, что соприкасающиеся окружности в точке $p$ гладких кривых в $\Sigma$, которые идут в направлении ${\vec v}$, заметают сферу $S$ с центром в точке $p+\tfrac1{k_{\vec v}}\cdot\Norm$ и радиусом $r=\tfrac1{|k_{\vec v}|}$.
\end{thm}

\begin{thm}{Упражнение}\label{ex:principal-revolution}
Пусть $\Sigma$ --- поверхность вращения вокруг оси $x$ гладкой простой кривой с единичной скоростью $\gamma(s)=(x(s),y(s))$, лежащей в верхней полуплоскости.
Допустим, что $x'\ne 0$.

\begin{subthm}{ex:principal-revolution:a}
Покажите, что параллели и меридианы поверхности являются её линиями кривизны.
\end{subthm}

{\sloppy

\begin{subthm}{ex:principal-revolution:formula}
Покажите, что 
\[\frac{|x'(s)|}{y(s)}
\quad
\text{и}
\quad
\frac{-y''(s)}{|x'(s)|}
\]
являются главными кривизнами $\Sigma$ в точке $(x(s),y(s),0)$
в направлениях соответствующих параллели и меридиана.
\end{subthm}

}

\begin{subthm}{ex:principal-revolution:pseudosphere}
Покажите, что $\Sigma$ имеет гауссову кривизну $-1$ во всех точках тогда и только тогда, когда $y$ удовлетворяет дифференциальному уравнению $y''=y$. 

\end{subthm}

\end{thm}

{

\begin{wrapfigure}{r}{31 mm}
\vskip-0mm
\includegraphics{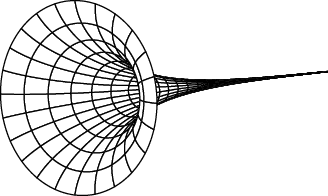}
\vskip-3mm
\end{wrapfigure}

В случае $y=e^{-s}$ кривая \( \gamma \) называется \index{трактриса}\emph{трактрисой}, а поверхность — \index{псевдосфера}\emph{псевдосферой};
согласно \ref{SHORT.ex:principal-revolution:pseudosphere}, её гауссова кривизна равна \( -1 \).

\begin{thm}{Упражнение}\label{ex:catenoid-is-minimal}
Покажите, что \index{катеноид}\emph{катеноид}, заданный уравнением
\[(\cosh z)^2=x^2+y^2\]
является минимальной поверхностью.
\end{thm}

}

\begin{thm}{Упражнение}\label{ex:helicoid-is-minimal}
Покажите, что \index{геликоид}\emph{геликоид}, определяемый следующей параметризацией
\[s(u,v)=(u\cdot \sin v,u\cdot \cos v,v)\]
является минимальной поверхностью.
\end{thm}

\begin{wrapfigure}{r}{51 mm}
\vskip-6mm
\centering
\includegraphics{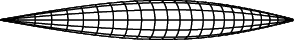}
\vskip0mm
\end{wrapfigure}

\begin{thm}{Упражнение}\label{ex:rev(sin)}
Пусть $\Sigma$ --- поверхность вращения вокруг оси $x$
с образующей $y=a\cdot \sin x$ для постоянной $a>0$ и $x\in (0,\pi)$.
Покажите, что гауссова кривизна $\Sigma$ не превышает $1$.
\end{thm}

\begin{thm}{Упражнение}\label{ex:rev(lin)}
Пусть $f\:(a,b)\to\mathbb{R}$ --- гладкая положительная функция, а $\Sigma$ --- поверхность вращения графика $y=f(x)$ вокруг оси $x$.
Предположим, что $\Sigma$ имеет нулевую гауссову кривизну.
Покажите, что $f$ является линейной функцией; то есть $f(x)=c\cdot x+d$ для некоторых констант $c$ и $d$.
\end{thm}

\section{Аквариум Лагунова}
\index{аквариум Лагунова}

Начнём с трёхмерного варианта задачи о луне в луже (\ref{thm:moon-orginal}).

\begin{thm}{Вопрос}\label{quest:lagunov}
Предположим, что тело $V\subset \mathbb{R}^3$ ограничено замкнутой поверхностью $\Sigma$ с 
главными кривизнами не больше $1$ по абсолютной величине.
Верно ли, что $V$ содержит шар радиуса~$1$?
\end{thm}

\begin{thm}{Упражнение}\label{ex:moon-revolution}
Покажите, что для тел вращения ответ на вопрос \ref{quest:lagunov} положителен.
\end{thm}

Позже (см. \ref{ex:convex-lagunov} и \ref{ex:rad=2})
мы увидим, что для выпуклых тел и поверхностей в открытом шаре радиуса $2$ ответ также положительный.
Сейчас же мы построим пример Владимира Лагунова \cite{lagunov-1961}, показывающий, что в общем случае ответ отрицательный.

\begin{figure}[t!]
\centering
\includegraphics{mppics/pic-33}
\vskip0mm
\end{figure}

\parbf{Построение.}
Начнём с тела вращения $V_1$, поперечное сечение которого показано на следующем рисунке.
Граница поперечного сечения состоит из $6$ длинных горизонтальных отрезков, включённых в $3$ простые замкнутые гладкие кривые.
(При их построении надо использовать процедуру сглаживания из приложения \ref{sec:analysis}.)
Граница $V_1$ состоит из трёх гладких сфер.

\begin{wrapfigure}{o}{21 mm}
\vskip-0mm
\centering
\includegraphics{mppics/pic-910}
\vskip0mm
\end{wrapfigure}

Мы предполагаем, что кривизна кривых не превышает~$1$.
Более того, за исключением почти горизонтальных частей, эта кривизна близка к $1$.
Самое толстое место в теле $V_1$ там, где все три граничных сферы подходят близко друг к другу;
остальная часть $V_1$ предполагается очень тонкой.
Можно добиться того, чтобы радиус $r$ максимального шара в $V_1$ только чуть превышал $r_2=\tfrac2{\sqrt{3}}-1$.
(Маленький чёрный круг на рисунке имеет радиус $r_2$, если три больших круга единичные.)
В частности, можно считать, что $r<\tfrac16$.

Упражнение \ref{ex:principal-revolution} приводит формулы для главных кривизн, применимые к граничным сферам $V_1$;
из них следует, что обе главные кривизны не превышают $1$ по абсолютной величине.

Остаётся сделать границу $V_1$ связной, не испортив при этом оценки на главные кривизны и не изменив максимальный радиус шаров внутри тела.

Каждая граничная сфера содержит два плоских круга; они образуют три пары, лежащие близко друг к другу.
Просверлим дыру через две такие пары и соединим отверстия трубкой --- телом вращения, ось которого смещена, но остаётся параллельной оси $V_1$.
Пусть $V_2$ --- полученное тело; его поперечное сечение показано на рисунке.

Теперь повторим операцию для других двух пар.
Пусть $V_3$ --- полученное тело.

Заметим, что граница $V_3$ связна.
Предполагая, что отверстия большие, можно добиться, чтобы её главные кривизны по-прежнему не превышали $1$; последнее доказывается так же, как для~$V_1$.
\qeds

Оценка $r_2=\tfrac2{\sqrt{3}}-1$ оптимальна;
это и прочее было доказано Владимиром Лагуновым и Абрамом Фетом \cite{lagunov-1960, lagunov-1961, lagunov-fet-1963, lagunov-fet-1965}.

\begin{wrapfigure}{o}{55 mm}
\centering
\vskip-0mm
\includegraphics{mppics/pic-920}
\vskip0mm
\end{wrapfigure}

Поверхность $V_3$ в аквариуме Лагунова имеет род $2$; то есть её можно параметризовать сферой с двумя ручками.

Действительно, граница $V_1$ состоит из трёх гладких сфер.

При первом сверлении, мы получаем четыре дырки в поверхноти тела --- две сферы получают по одной дырке, и одна сфера две дырки.
Соединив две сферы трубкой, получаем одну сферу.
При соединении двух дырок оставшейся сферы, получаем тор.
Он изображён слева на рисунке.
Таким образом, $V_2$ ограничено сферой и тором.

Чтобы построить $V_3$ из $V_2$, мы делаем тор из оставшейся сферы и соединяем его трубкой с другим тором, и получается сфера с двумя ручками.

\begin{thm}{Упражнение}\label{ex:lagunov-genus4}
Измените построение Лагунова так, чтобы граничная поверхность была сферой с четырьмя ручками.

Попробуйте изменить построение, чтобы граничная поверхность была сферой с тремя ручками.
\end{thm}



\chapter{Опорные поверхности}
\label{chap:surface-support}

\section{Определения}

Две гладкие поверхности $\Sigma_1$ и $\Sigma_2$ \index{касательные поверхности}\emph{касаются} в точке $p \z\in \Sigma_1 \cap \Sigma_2$, если у них в этой точке общая касательная плоскость.
В этом случае их нормали в $p$ либо равны, либо противонаправлены.
То есть, если поверхности ориентированы, то можно говорить о \index{сонаправленные и противонаправленные!поверхности}\emph{сонаправленном} или {}\emph{противонаправленном} касании в $p$.

Поверхность $\Sigma_1$ \index{локальная опорная}\emph{локально подпирает} поверхность $\Sigma_2$ в точке~$p$, если $p$ --- общая точка поверхностей и нашлась такая окрестность $U$ точки $p$, что $\Sigma_2\cap U$ целиком лежит по одну сторону от $\Sigma_1$ в $U$.
Согласно следующему упражнению, такие поверхности обязаны касаться в~$p$; выбрав их ориентации, можно считать, что они сонаправлены в $p$.

\begin{thm}{Упражнение}\label{ex:supp>tan}
Пусть гладкая поверхность $\Sigma_1$ локально подпирает гладкую поверхность $\Sigma_2$ в точке~$p$.
Докажите, что $\Sigma_1$ касается $\Sigma_2$ в точке $p$.
\end{thm}

Допустим, что поверхности $\Sigma_1$ и $\Sigma_2$ касаются и сонаправлены в~$p$.
Представим их локально графиками функций, скажем $f_1$ и $f_2$, в общих касательно-нормальных координатах при $p$.
Отметим, что $\Sigma_2$ локально подпирает $\Sigma_1$ в $p$ тогда и только тогда, когда неравенство
\[ f_1(x,y)\ge f_2(x,y)
\quad\text{(\,или}\quad
f_1(x,y)\le f_2(x,y)\,)\]
выполнено при всех $(x,y)$, достаточно близких к нулю;
в этом случае $\Sigma_2$ локально подпирает $\Sigma_1$ \index{опорная!поверхность}\emph{снаружи} (соответственно, \emph{изнутри}).

\begin{thm}{Предложение}\label{prop:surf-support}
Пусть $\Sigma_1$ и $\Sigma_2$ --- гладкие ориентированные поверхности.
Предположим, что $\Sigma_1$ локально подпирает $\Sigma_2$ изнутри в точке $p$ (эквивалентно, $\Sigma_2$ локально подпирает $\Sigma_1$ снаружи).
Тогда 
\[k_1(p)_{\Sigma_1}\ge k_1(p)_{\Sigma_2}\quad\text{и}\quad k_2(p)_{\Sigma_1}\z\ge k_2(p)_{\Sigma_2}.\]
\end{thm}

\begin{thm}{Упражнение}\label{ex:surf-support}
Постройте две поверхности $\Sigma_1$ и $\Sigma_2$ с общей точкой $p$ и общей нормалью при $p$, которые локально не подпирают друг друга в~$p$, но
$k_1(p)_{\Sigma_1}\z> k_1(p)_{\Sigma_2}$ и $k_2(p)_{\Sigma_1}\z> k_2(p)_{\Sigma_2}$.
\end{thm}

\parbf{Доказательство \ref{prop:surf-support}.}
Можно считать, что $\Sigma_1$ и $\Sigma_2$ --- графики $z\z=f_1(x,y)$ и $z=f_2(x,y)$ в общих касательно-нормальных координатах при $p$ и $f_1\ge f_2$.

Пусть $\vec w\z\in \T_p\Sigma_1=\T_p\Sigma_2$ --- общий касательный единичный вектор.
Построим плоскость $\Pi$, через $p$, в направлении нормали $\Norm_p$ и ${\vec w}$.
Пусть $\gamma_1$ и $\gamma_2$ --- кривые пересечения $\Sigma_1$ и $\Sigma_2$ с $\Pi$.

Выберем ориентацию на $\Pi$, чтобы $\Norm_p$ указывал влево от ${\vec w}$.
Далее, параметризуем обе кривые в направлении ${\vec w}$ при $p$;
в этом случае они сонаправленны и $\gamma_1$ подпирает $\gamma_2$ слева в $p$.

\begin{wrapfigure}{o}{35 mm}
\vskip-6mm
\centering
\includegraphics{mppics/pic-80}
\vskip-0mm
\end{wrapfigure}

Из \ref{prop:supporting-circline}, получаем следующее неравенство на нормальные кривизны $\Sigma_1$ и $\Sigma_2$ в $p$ в направлении ${\vec w}$:
\[k_{\vec w}(p)_{\Sigma_1}\ge k_{\vec w}(p)_{\Sigma_2}.\eqlbl{kw>=kw}\]

\setlength{\columnseprule}{0.4pt}
\begin{multicols}{2}
Далее, из \ref{obs:k1-k2}, 
\[k_1(p)_{\Sigma_i}=\min_{\vec w}\set{k_{\vec w}(p)_{\Sigma_i}}{},\]
где $\vec w\in\T_p$ и $|\vec w|=1$.
Выберем ${\vec w}$ так, что 
\[k_1(p)_{\Sigma_1}\z=k_{\vec w}(p)_{\Sigma_1}.\]
Тогда из \ref{kw>=kw}, получим
\begin{align*}
k_1(p)_{\Sigma_1}&=k_{\vec w}(p)_{\Sigma_1}\ge
\\
&\ge k_{\vec w}(p)_{\Sigma_2}\ge
\\
&\ge\min_{\vec v}\set{k_{\vec v}(p)_{\Sigma_2}}{}=
\\
&=k_1(p)_{\Sigma_2}.
\end{align*}

\columnbreak

Аналогично, из \ref{obs:k1-k2},
\[k_2(p)_{\Sigma_i}=\max_{\vec w}\set{k_{\vec w}(p)_{\Sigma_i}}{},\]
где $\vec w\in\T_p$ и $|\vec w|=1$.
Выберем ${\vec w}$ так, что 
\[k_2(p)_{\Sigma_2}=k_{\vec w}(p)_{\Sigma_2},\]
и из \ref{kw>=kw},
\begin{align*}
k_2(p)_{\Sigma_2}&=k_{\vec w}(p)_{\Sigma_2}\le
\\
&\le k_{\vec w}(p)_{\Sigma_1}\le
\\
&\le\max_{\vec v}\set{k_{\vec v}(p)_{\Sigma_1}}{}=
\\
&=k_2(p)_{\Sigma_1}.
\end{align*}
\end{multicols}

В обоих случаях предполагается, что $\vec v\in\T_p$ и $|\vec v|=1$.\qeds

\enlargethispage{2\baselineskip}

\begin{thm}{Следствие}\label{cor:surf-support}
Пусть $\Sigma_1$ и $\Sigma_2$ --- ориентированные гладкие поверхности.
Если $\Sigma_1$ локально подпирает $\Sigma_2$ изнутри в точке~$p$, то

\begin{subthm}{cor:surf-support:mean}
$H(p)_{\Sigma_1}\ge H(p)_{\Sigma_2}$;
\end{subthm}

\begin{subthm}{cor:surf-support:gauss}
Если $k_1(p)_{\Sigma_2}\ge 0$, то $K(p)_{\Sigma_1}\ge K(p)_{\Sigma_2}$.
\end{subthm}
 
\end{thm}

\parbf{Доказательство;} \ref{SHORT.cor:surf-support:mean}.
Следует из \ref{prop:surf-support}, поскольку $H=k_1+k_2$. 

\parit{\ref{SHORT.cor:surf-support:gauss}.}
Поскольку $k_2(p)_{\Sigma_i}\ge k_1(p)_{\Sigma_i}$, и $k_1(p)_{\Sigma_2}\ge 0$, все главные кривизны 
$k_1(p)_{\Sigma_1}$, 
$k_1(p)_{\Sigma_2}$, 
$k_2(p)_{\Sigma_1}$ и 
$k_2(p)_{\Sigma_2}$ неотрицательны.
Из \ref{prop:surf-support}, получаем
\begin{align*}
K(p)_{\Sigma_1}&=k_1(p)_{\Sigma_1}\cdot k_2(p)_{\Sigma_1}\ge 
\\
&\ge k_1(p)_{\Sigma_2}\cdot k_2(p)_{\Sigma_2}=K(p)_{\Sigma_2}.
\end{align*}
\qedsf

\begin{thm}{Упражнение}\label{ex:positive-gauss-0}
Покажите, что на любой замкнутой поверхности, лежащей в единичном шаре, найдётся точка с гауссовой кривизной хотя бы~$1$.
Выведите отсюда, что любая замкнутая поверхность имеет точку с положительной гауссовой кривизной.
\end{thm}

{\sloppy

\begin{thm}{Упражнение}\label{ex:positive-gauss}
Покажите, что любая замкнутая поверхность, лежащая на расстоянии не более $1$ от некоторой прямой, имеет точку с гауссовой кривизной не менее~$1$.
\end{thm}

}

\section{Выпуклые поверхности}

Поверхность называется \index{выпуклая!поверхность}\emph{выпуклой}, если она ограничивает выпуклую область.

\begin{thm}{Упражнение}\label{ex:convex-surf}
Покажите, что гауссова кривизна гладких выпуклых поверхностей неотрицательна.
\end{thm}

\begin{thm}{Продвинутое упражнение}\label{ex:convex-lagunov}
Пусть $R$ --- выпуклое тело в $\mathbb{R}^3$, ограниченное поверхностью с главными кривизнами, не превышающими~$1$.
Покажите, что $R$ содержит единичный шар.
\end{thm}

{\sloppy

Напомним, что область $R$ евклидова пространства называется {}\emph{строго выпуклой}, если для любых двух точек $x,y\in R$ любая точка $z$, лежащая между $x$ и $y$, находится во внутренности~$R$.

Любое открытое выпуклое множество строго выпукло. 
Куб, как и любой выпуклый многогранник, нестрого выпуклы.
Ясно, что \textit{замкнутая выпуклая область является строго выпуклой тогда и только тогда, когда её граница не содержит отрезков}.

}

\enlargethispage{2\baselineskip}

\begin{thm}{Лемма}\label{lem:gauss+=>convexity}
Пусть $z=f(x,y)$ --- локальное представление гладкой поверхности $\Sigma$ в касательно-нормальных координатах при $p\in\Sigma$.
Если обе главные кривизны $\Sigma$ в~$p$ положительны,
то функция $f$ строго выпукла в окрестности нуля и имеет в нём локальный минимум.

В частности, касательная плоскость $\T_p$ локально подпирает $\Sigma$ снаружи в~$p$.
\end{thm}

\parbf{Доказательство.}
Поскольку обе главные кривизны положительны, \ref{cor:Shape(ij)} влечёт, что 
\[D^2_{\vec w}f(0,0)=\langle \Shape_p({\vec w}),{\vec w}\rangle\ge k_1(p)>0\] 
для любого единичного касательного вектора $\vec w\in\T_p\Sigma$ (в наших координатах $\vec w$ --- горизонтальный вектор).

Так как функция $(x,y,{\vec w})\mapsto D^2_{\vec w}f(x,y)$ непрерывна, $D^2_{\vec w}f(x,y)>0$, при любом $\vec w\ne 0$ и в любой точке $(x,y)$ в малой окрестности нуля.
Значит, $f$ строго выпукла в некоторой окрестности нуля (\ref{thm:Jensen}).

{\sloppy

И наконец, у $f$ есть строгий локальный минимум в нуле,
ибо $\nabla f(0,0)\z=0$ и $f$ строго выпукла в окрестности нуля.
\qeds

}

\begin{thm}{Упражнение}\label{ex:section-of-convex}
Пусть $\Sigma$ --- гладкая поверхность (без края) с положительной гауссовой кривизной.
Покажите, что любая компонента связности пересечения $\Sigma$ с плоскостью либо одна точка, либо гладкая кривая с ориентированной кривизной одного знака.
\end{thm}

\section{Признак выпуклости}

\begin{thm}{Теорема}\label{thm:convex-embedded}
Пусть $\Sigma$ --- открытая или замкнутая гладкая поверхность с положительной гауссовой кривизной.
Тогда $\Sigma$ ограничивает строго выпуклую область в $\mathbb{R}^3$.
\end{thm}

В доказательстве потребуется, что поверхность связна (это верно по определению);
иначе пара сфер дала бы контрпример.
Верно также, что \textit{любая замкнутая гладкая поверхность с неотрицательной гауссовой кривизной ограничивает выпуклую область},
но доказательство ощутимо сложней \cite{hadamard,gomes,sacksteder}.
Теорема остаётся верной для поверхностей с самопересечениями. 
Это доказано Джеймсом Стокером \cite{stoker}, который приписал результат Жаку Адамару, доказавшему близкое утверждение \cite[§ 23]{hadamard}.
Из доказательства ниже можно извлечь, что \textit{любая замкнутая связная локально выпуклая область евклидова пространства выпукла}.

\enlargethispage{2\baselineskip}

\parbf{Доказательство.}
Поскольку гауссова кривизна положительна, можно считать, что главные кривизны положительны во всех точках.
Обозначим через $R$ область, ограниченную $\Sigma$, в которую указывает поле нормалей $\Norm$ на $\Sigma$.
(Область $R$ существует по \ref{clm:proper-divides}.)

Сначала покажем, что $R$ {}\emph{локально строго выпукла};
то есть для любой точки $p\in \Sigma$, пересечение $R$ и малого шара с центром в~$p$ строго выпукло.
Действительно, пусть $\Sigma$ локально задаётся графиком $z=f(x,y)$ в касательно-нормальных координатах при~$p$.
По \ref{lem:gauss+=>convexity}, функция $f$ строго выпукла в окрестности нуля.
В частности, пересечение малого шара с центром в $p$ и эпиграфа $z\ge f(x,y)$ строго выпукло.

Поскольку $\Sigma$ связна, связна и область $R$.
Более того, любые две точки внутри $R$ можно соединить ломаной во внутренности~$R$.

Допустим, что внутренность $R$ не выпукла;
то есть найдутся точки $x,y\in R$, и точка $z$ между $x$ и $y$, которая не лежит во внутренности~$R$.
Рассмотрим ломаную $\beta$ от $x$ до $y$ во внутренности~$R$.
Пусть $y_0$ --- первая такая точка на $\beta$, что хорда $[x,y_0]$ касается $\Sigma$, скажем в~$z_0$.

{

\begin{wrapfigure}{r}{43 mm}
\vskip-4mm
\centering
\includegraphics{mppics/pic-37}
\vskip-0mm
\end{wrapfigure}

Поскольку $R$ локально строго выпукла, $R\z\cap B(z_0,\epsilon)$ является строго выпуклой при достаточно малых $\epsilon>0$.
С другой стороны, $z_0$ лежит между двумя точками в пересечении $[x,y_0]\cap B(z_0,\epsilon)$.
Поскольку $[x,y_0]\subset R$, мы приходим к противоречию.

}

\enlargethispage{2\baselineskip}

Следовательно, внутренность $R$ выпукла.
Отметим, что область $R$ есть замыкание своей внутренней части, поэтому $R$ также выпукла.

Поскольку $R$ является локально строго выпуклой, её граница $\Sigma$ не содержит отрезков прямых.
Следовательно, $R$  строго выпукла.
\qeds

\begin{thm}{Упражнение}\label{ex:surrounds-disc}
Предположим, что замкнутая поверхность $\Sigma$ ограничивает область, содержащую единичную окружность.
Покажите, что $K(p)\le 1$ в некоторой точки $p \in \Sigma$.
\end{thm}

\begin{thm}{Упражнение}\label{ex:small-gauss}
Пусть $\Sigma$ --- замкнутая гладкая выпуклая поверхность диаметра не меньше $\pi$;
то есть существует пара точек $p,q\in\Sigma$ таких, что $|p-q|\ge \pi$.
Покажите, что $\Sigma$ имеет точку с гауссовой кривизной не более~$1$.
\end{thm}

\section{Ещё о выпуклости}

\begin{thm}{Теорема}\label{thm:convex-closed}
Любая замкнутая гладкая выпуклая поверхность $\Sigma$ является гладкой сферой;
то есть $\Sigma$ допускает гладкую регулярную параметризацию стандартной сферой $\mathbb{S}^2$.
\end{thm}

{

\begin{wrapfigure}[7]{r}{33 mm}
\vskip-10mm
\centering
\includegraphics{mppics/pic-78}
\end{wrapfigure}

\begin{thm}{Доказательство и упражнение}\label{ex:convex-proper-sphere}\\
Пусть $R$ --- выпуклая компактная область ограниченая гладкой поверхностью~$\Sigma$.
Предположим, что $R$ содержит начало координат.

\begin{subthm}{ex:convex-proper-sphere:single}
Покажите, что любой луч, исходящий из начала координат, пересекает $\Sigma$ в единственной точке;
иначе говоря, существует такая положительная функция $\rho\:\mathbb{S}^2\z\to\mathbb{R}$, что $\Sigma$ состоит из точек $q\z=\rho(\xi)\cdot \xi$ для $\xi\in \mathbb{S}^2$.
\end{subthm}

\begin{subthm}{ex:convex-proper-sphere:smooth}
Покажите, что функция $\rho\:\mathbb{S}^2\to\mathbb{R}$ гладкая.
Выведите, что $\xi\z\mapsto \rho(\xi)\cdot \xi$ --- гладкая регулярная параметризация $\mathbb{S}^2\z\to \Sigma$.
\end{subthm}

\end{thm}

}

\begin{thm}{Теорема}\label{thm:convex-open}
Пусть $\Sigma$ --- открытая гладкая поверхность, ограничивающая строго выпуклую замкнутую область.
Тогда, в некоторой системе координат, $\Sigma$ задаётся графиком $z\z=f(x,y)$ выпуклой функции $f$, определённой на выпуклой открытой области $\Omega$ в $(x,y)$-плоскости.
Более того, $f(x_n,y_n)\to\infty$ при $(x_n,y_n)\to(x_\infty,y_\infty)\z\in \partial\Omega$.

\end{thm}

{

\begin{wrapfigure}{r}{28 mm}
\vskip-0mm
\centering
\includegraphics{mppics/pic-1181}
\vskip2mm
\end{wrapfigure}

\begin{thm}{Доказательство и упражнение}\label{ex:convex-proper-plane}
Предположим, что строго выпуклая замкнутая некомпактная область $R$ содержит начало координат во внутренней части и ограничена гладкой поверхностью~$\Sigma$.
Докажите следующее.

\begin{subthm}{ex:convex-proper-plane:a}
$R$ содержит луч, скажем $\ell$.
\end{subthm}

\begin{subthm}{ex:convex-proper-plane:b}
Любая прямая $m$, параллельная $\ell$, пересекает $\Sigma$ не более чем в одной точке.
\end{subthm}

\begin{subthm}{ex:convex-proper-plane:c}
В системе координат, с осью $z$ направленной по $\ell$, проекция $\Sigma$ на $(x,y)$-плоскость есть открытое выпуклое множество, скажем $\Omega$.
\end{subthm}

\begin{subthm}{ex:convex-proper-plane:d}
$\Sigma$ является графиком $z=f(x,y)$ гладкой выпуклой функции $f$, определённой на $\Omega$.
\end{subthm}

\begin{subthm}{ex:convex-proper-plane:e}
Докажите последнее утверждение в \ref{thm:convex-open}.
\end{subthm}

\end{thm}

}

\begin{thm}{Упражнение}\label{ex:open+convex=plane}
{\sloppy
Покажите, что любая открытая поверхность $\Sigma$ с положительной гауссовой кривизной является топологической плоскостью;
то есть существует вложение $\mathbb{R}^2\hookrightarrow\mathbb{R}^3$ с образом~$\Sigma$.

}

Попробуйте доказать, что $\Sigma$ является гладкой плоскостью;
то есть вложение $f$ можно сделать гладким и регулярным.
\end{thm}

{\sloppy

\begin{thm}{Упражнение}\label{ex:circular-cone}
Покажите, что любая открытая гладкая поверхность $\Sigma$ с положительной гауссовой кривизной
помещается внутри бесконечного кругового конуса.
Другими словами, в некоторой системе координат, $\Sigma$ лежит в области $z \z\ge m \cdot\sqrt{x^2 + y^2}$ для положительной константы $m$.
\end{thm}

}

\begin{thm}{Упражнение}\label{ex:intK}
Пусть $\Sigma$ --- гладкая выпуклая поверхность положительной гауссовой кривизны и $\Norm\:\Sigma\to \mathbb{S}^2$ --- её сферическое отображение.

\begin{multicols}{2}

\begin{subthm}{ex:intK:4pi}
Предполагая, что $\Sigma$ замкнута,
покажите, что $\Norm\:\Sigma\to \mathbb{S}^2$ является биекцией.
Выведите отсюда, что 
\[\iint_\Sigma K=4\cdot\pi.\]
\end{subthm}

\columnbreak

\begin{subthm}{ex:intK:2pi}
Предполагая, что $\Sigma$ открыта,
покажите, что $\Norm\:\Sigma\to \mathbb{S}^2$
задаёт биекцию на подмножество полусферы,
и выведите, что 
\[\iint_\Sigma K\le 2\cdot\pi.\]
\end{subthm}

\end{multicols}

\end{thm}

\chapter{Седловые поверхности}

\section{Определения}\label{sec:saddle}

\begin{wrapfigure}{r}{43 mm}
\vskip-16mm
\centering
\includegraphics{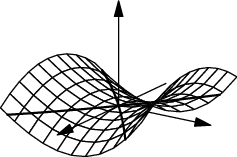}
\vskip0mm
\end{wrapfigure}

Поверхность называется \index{седловая поверхность}\emph{седловой}, если её гауссова кривизна в каждой точке неположительна;
другими словами, главные кривизны в каждой точке имеют противоположные знаки или хотя бы одна из них равна нулю.

Если гауссова кривизна отрицательна в каждой точке,
то поверхность называется {}\emph{строго седловой};
это означает, что главные кривизны имеют противоположные знаки в каждой точке.
В этом случае касательная плоскость не может подпирать поверхность даже локально --- при движении по поверхности в главных направлениях из точки, можно очутиться и выше и ниже касательной плоскости.

\begin{thm}{Упражнение}\label{ex:convex-revolution}
Пусть $f\:\mathbb{R}\to\mathbb{R}$ --- гладкая положительная функция.
Покажите, что поверхность вращения графика $y\z=f(x)$ вокруг оси $x$
является седловой тогда и только тогда, когда $f$ выпукла; то есть если $f''(x)\ge0$ для любого~$x$.
\end{thm}

Поверхность $\Sigma$ называется \index{линейчатая поверхность}\emph{линейчатой}, если для каждой точки $p\in \Sigma$ существует прямая $\ell_p\subset \Sigma$, проходящая через $p$.

\begin{thm}{Упражнение}\label{ex:ruled=>saddle}
Покажите, что любая линейчатая поверхность седловая.
\end{thm}

\begin{thm}{Упражнение}\label{ex:saddle-convex}
Пусть $\Sigma$ --- открытая строго седловая поверхность, а $f\:\mathbb{R}^3\z\to\mathbb{R}$ --- гладкая выпуклая функция.
Покажите, что сужение $f$ на $\Sigma$ не имеет точки локального максимума.
\end{thm}

Направление на гладкой поверхности с нулевой нормальной кривизной называется \index{асимптотическое направление и линия}\emph{асимптотическим}.
Гладкая кривая, идущая всё время в асимптотическом направлении, называется
{}\emph{асимптотической линией}.\label{page:asymptotic line}

Напомним, что множество $R$ на плоскости называется \index{звёздное множество}\emph{звёздным}, если существует точка $p\in R$ такая, что для любой точки $x\in R$ отрезок $[p,x]$ принадлежит~$R$.

\begin{thm}{Продвинутое упражнение}\label{ex:panov}
Пусть $\gamma$ --- замкнутая гладкая асимптотическая линия
на графике $z\z=f(x,y)$ гладкой функции~$f$. 
Предположим, что график строго седловой в окрестности~$\gamma$.
Покажите, что область на $(x,y)$-плоскости, ограниченная проекцией $\bar \gamma$ кривой $\gamma$, не может оказаться звёздной. 
\end{thm}

\begin{thm}{Продвинутое упражнение}\label{ex:crosss}
Пусть $p$ --- точка гладкой поверхности $\Sigma$ с $K(p)<0$.
Покажите, что существует окрестность $\Omega$ точки $p$ в $\Sigma$,
такая, что пересечение $\Omega$ с касательной плоскостью $\T_p$ представляет собой объединение двух гладких кривых, \index{трансверсальность}\emph{трансверсально пересекающихся} в точке~$p$;
то есть их векторы скорости линейно независимы в $p$.
\end{thm}

\section{Критерий горбушек}

Обратите внимание, что \textit{замкнутая поверхность не может быть седловой}.
Действительно, пусть $\Sigma$ --- замкнутая поверхность.
Рассмотрим наименьшую сферу, которая окружает~$\Sigma$.
Эта сфера обязана подпирать $\Sigma$ в некоторой точке.
Согласно \ref{prop:surf-support}, в этой точке у $\Sigma$ главные кривизны одного знака.
Следующее более общее утверждение использует похожую идею.

\begin{thm}{Лемма}\label{lem:convex-saddle}
Пусть $\Sigma$ --- компактная седловая поверхность, и её край лежит в замкнутой выпуклой области~$R$.
Тогда вся поверхность $\Sigma$ лежит в~$R$.
\end{thm}

{

\begin{wrapfigure}{r}{46 mm}
\vskip-8mm
\centering
\includegraphics{mppics/pic-73}
\vskip-4mm
\end{wrapfigure}

\parit{Замечание.}
Для строго седловых поверхностей лемму можно вывести из \ref{ex:saddle-convex}.

\parbf{Доказательство.}
Допустим, что существует точка $p\in \Sigma$,  не лежащая в~$R$.
Тогда найдётся плоскость $\Pi$, отделяющая $p$ от $R$ (\ref{lem:separation}).
Обозначим через $\Sigma'$ ту часть $\Sigma$, которая лежит с $p$ на одной стороне от~$\Pi$.

}

Так как $\Sigma$ компактна, её можно окружить сферой;
пусть $\sigma$ --- окружность пересечения этой сферы с $\Pi$.
Рассмотрим наименьший сферический купол $\Sigma_0$ с краем $\sigma$, который окружает~$\Sigma'$.

$\Sigma_0$ подпирает $\Sigma$ в некоторой точке~$q$.
Не умаляя общности, можно предположить, что $\Sigma_0$ и $\Sigma$ сонаправлены в $q$, и у $\Sigma_0$ положительные главные кривизны.
В этом случае, $\Sigma_0$ подпирает $\Sigma$ снаружи,
и из \ref{cor:surf-support}, $K(q)_\Sigma\z\ge K(q)_{\Sigma_0}>0$ --- противоречие.
\qeds

\begin{thm}{Упражнение}\label{ex:proper-saddle}
{\sloppy 
Постройте собственную седловую поверхность, которая не лежит в выпуклой оболочке своего края.
(Согласно \ref{lem:convex-saddle}, она некомпактна.)

}
\end{thm}

\begin{thm}{Упражнение}\label{ex:length-of-bry}
Предположим, что точка $p$ лежит на компактной гладкой седловой поверхности $\Delta$, край которой лежит на единичной сфере с центром в~$p$.
Докажите, что если $\Delta$ является топологическим диском, то $\length(\partial\Delta)\ge 2\cdot\pi$.

Покажите, что это перестаёт быть верным без предположения, что $\Delta$ топологический диск.
\end{thm}

\parit{Замечание.}
На самом деле $\area \Delta\ge \pi$;
то есть единичный круг имеет минимальную возможную площадь среди поверхностей из упражнения.
Это выводится из так называемой \index{формула коплощади}\emph{формулы коплощади}.

\begin{thm}{Упражнение}\label{ex:circular-cone-saddle}
Покажите, что открытая седловая поверхность
не может лежать внутри бесконечного кругового конуса. 
\end{thm}

Топологический диск $\Delta$ на поверхности $\Sigma$ называется \index{горбушка}\emph{горбушкой} поверхности $\Sigma$, если его край $\partial\Delta$ лежит в некоторой плоскости $\Pi$, а $\Delta \setminus \partial \Delta$ лежит строго с одной стороны от $\Pi$.

\begin{thm}{Предложение}\label{prop:hat}
Гладкая поверхность $\Sigma$ является седловой тогда и только тогда, когда на ней нет горбушек.
\end{thm}

Седловая поверхность может содержать замкнутую плоскую кривую.
Например, гиперболоид $x^2+y^2-z^2=1$ пересекается по единичной окружности с плоскостью $(x,y)$.
Однако, согласно предложению (и лемме), плоская кривая на седловой поверхности не может ограничивать топологический диск (как и любое компактное множество).

\parbf{Доказательство.}
Так как плоскость выпукла, необходимость следует из \ref{lem:convex-saddle};
остаётся доказать достаточность.

Допустим, что $\Sigma$ не седловая; то есть у неё есть точка $p$ со строго положительной гауссовой кривизной;
иначе говоря, её главные кривизны $k_1(p)$ и $k_2(p)$ одного знака.
Давайте считать, что $0 < k_1(p) \z\le k_2(p)$, другой случай аналогичен.

Пусть $z=f(x,y)$ --- локальное представление $\Sigma$ в касательно-нормальных координатах при~$p$.
Рассмотрим множество $F_\epsilon$ в плоскости $(x,y)$, определяемое неравенством $f(x,y)\le \epsilon$.
Согласно \ref{lem:gauss+=>convexity}, функция $f$ выпукла в малой окрестности $(0,0)$.
Следовательно, $F_\epsilon$ выпукло при достаточно малых $\epsilon>0$.
В частности, $F_\epsilon$ --- топологический диск.

Отображение $(x,y)\mapsto (x,y,f(x,y))$ есть гомеоморфизм из $F_\epsilon$
в
\[\Delta_\epsilon=\set{(x,y,f(x,y))\in \mathbb{R}^3}{f(x,y)\le \epsilon};\]
следовательно,
$\Delta_\epsilon$ --- топологический диск при любом достаточно малом $\epsilon>0$.
Граничная линия $\Delta_\epsilon$ лежит в плоскости $z=\epsilon$, и весь диск лежит ниже неё;
то есть $\Delta_\epsilon$ является горбушкой~$\Sigma$.
\qeds

Следующее упражнение показывает, что на самом деле $\Delta_\epsilon$ --- гладкий диск.
Из него выводится более сильный вариант критерия горбушек (\ref{prop:hat});
а именно, в определении горбушек можно предполагать, что диск является гладким.

\begin{thm}{Упражнение}\label{ex:disc-hat}
Пусть $f\:\mathbb{R}^2\to\mathbb{R}$ --- гладкая строго выпуклая функция с минимумом в нуле.
Покажите, что для любого $\epsilon>0$, множество $F_\epsilon$ в графике $z=f(x,y)$, определяемое неравенством $f(x,y)\le \epsilon$, есть гладкий диск;
то есть существует диффеоморфизм из
$F_\epsilon$ на единичный диск $\Delta\z=\set{(x,y)\in\mathbb{R}^2}{x^2+y^2\le 1}$.
\end{thm}

\begin{thm}{Упражнение}\label{ex:saddle-linear}
Пусть $L\:\mathbb{R}^3\to\mathbb{R}^3$ --- \index{аффинное преобразование}\emph{аффинное преобразование}; то есть биекция $\mathbb{R}^3\to\mathbb{R}^3$, переводящая любую плоскость в плоскость.
Покажите, что для любой седловой поверхности $\Sigma$ образ $L(\Sigma)$ также седловая поверхность.
\end{thm}


\section{Седловые графики}

Следующая теорема доказана Сергеем Бернштейном \cite{bernstein}.
\index{теорема Бернштейна}

\begin{thm}{Теорема}\label{thm:bernshtein}
Если гладкая функция $f\:\mathbb{R}^2\to\mathbb{R}$ имеет строго седловой график $z=f(x,y)$, то она неограничена;
то есть нету константы $C$, при которой неравенство 
$|f(x,y)|\le C$ выполнялось бы при любых $(x,y)\in\mathbb{R}^2$.
\end{thm}

По теореме, седловой график не может лежать между двумя горизонтальными плоскостями.
Вместе с \ref{ex:saddle-linear} это влечёт, что седловые графики не могут находиться между двумя параллельными плоскостями.
Следующее упражнение показывает, что теорема не применима к седловым поверхностям, не являющимся графиками.

\begin{thm}{Упражнение}\label{ex:between-parallels}
Постройте открытую строго седловую поверхность, лежащую между двумя параллельными плоскостями.
\end{thm}

Бывают функции со строго седловым графиком, ограниченные лишь с одной стороны.
Поскольку $\exp(x-y^2)>0$, следующее упражнение даёт такой пример.
То есть в доказательстве теоремы нам придётся использовать и верхнюю и нижнюю оценку на $f$.

\begin{thm}{Упражнение}\label{ex:one-side-bernshtein}
Покажите, что график
$z=\exp(x-y^2)$
является строго седловым.
\end{thm}

Следующее упражнение даёт усиление теоремы Бернштейна.

\begin{thm}{Продвинутое упражнение}\label{ex:saddle-graph}
Пусть $\Sigma$ --- открытая гладкая строго седловая поверхность в $\mathbb{R}^3$.
Предположим, что существует такое компактное подмножество $K\subset \Sigma$, что дополнение $\Sigma\setminus K$ является графиком $z=f(x,y)$ гладкой функции, определённой в открытой области плоскости $(x,y)$.
Покажите, что вся поверхность $\Sigma$ является графиком.
\end{thm}

Следующая лемма похожа на \ref{ex:length-of-bry};
она потребуется в доказательстве теоремы \ref{thm:bernshtein}.

\begin{thm}{Лемма}\label{lem:region}
Не существует собственной строго седловой гладкой поверхности с краем в плоскости $\Pi$ и расположенной на ограниченном расстоянии от некоторой прямой в $\Pi$.
\end{thm}

\parbf{Доказательство.}
Воспользовавшись \ref{ex:saddle-linear}, лемму можно переформулировать следующим образом:
\textit{не существует собственной строго седловой гладкой поверхности 
с краем в $(x,y)$-плоскости,
содержащейся в бесконечном параллелепипеде следующего вида:}
\[R=\set{(x,y,z)\in\mathbb{R}^3}{0\le z\le r,\ 0\le y\le r}.\]

Предположим обратное, пусть $\Sigma$ --- такая поверхность.
Рассмотрим проекцию $\hat \Sigma$ на $(x,z)$-плоскость.
Она лежит в верхней полуплоскости, но ниже линии $z=r$.

Рассмотрим открытую верхнюю полуплоскость
\[H=\set{(x,z)\in \mathbb{R}^2}{z> 0}.\] 
Пусть $\Theta$ --- компонента связности дополнения $H\setminus \hat \Sigma$, содержащая точки выше линии $z=r$.

Заметим, что $\Theta$ выпукла.
Если это не так, то найдётся отрезок $[p,q]$ с $p,q\in \Theta$, отсекающий от $\hat\Sigma$ компактный кусок.
\begin{figure}[!ht]
\vskip-1mm
\centering
\includegraphics{mppics/pic-74}
\vskip-1mm
\end{figure}
Пусть $\Pi$ --- плоскость через $[p,q]$, перпендикулярная $(x,z)$-плоскости, она отсекает от $\Sigma$ компактную область $\Delta$.
Применив \ref{lem:reg-section}, можно считать, что $\Delta$ --- компактная поверхность с краем в $\Pi$, и остальная часть $\Delta$ лежит строго с одной стороны от $\Pi$.
Поскольку плоскость $\Pi$ выпукла, приходим к противоречию с \ref{lem:convex-saddle}.

Итак, $\Theta$ является открытым выпуклым подмножеством $H$, которое содержит все точки выше $z=r$.
По выпуклости, вместе с любой точкой $w$ множество $\Theta$ содержит все точки на лучах, которые начинаются в $w$ и \textit{направлены вверх}; то есть в направлениях с положительной $z$-координатой. 
Другими словами, вместе с любой точкой $w$
множество $\Theta$ содержит все точки с большей $z$-координатой.
\begin{figure}[!ht]
\vskip-1mm
\centering
\includegraphics{mppics/pic-75}
\vskip-1mm
\end{figure}
Так как~$\Theta$ открыто, оно описывается неравенством $z>r_0$.
Отсюда вытекает, что плоскость $z=r_0$ подпирает $\Sigma$ в некоторой точке (на самом деле во многих точках).
Применив \ref{prop:surf-support}, приходим к противоречию.
\qeds


\parbf{Доказательство \ref{thm:bernshtein}.}
Пусть $\Sigma$ --- график $z=f(x,y)$.
Предположим обратное, а именно, что $\Sigma$ лежит между двумя плоскостями $z=\pm C$. 

Функция $f$ не может быть постоянной.
Следовательно, касательная плоскость $\T_p$ в некоторой точке $p\in\Sigma$ не горизонтальна.

Обозначим через $\Sigma^+$ часть $\Sigma$, лежащую выше $\T_p$.
В ней есть как минимум две связные компоненты, которые приближаются к $p$ с обеих сторон 
вдоль главного направления с положительной главной кривизной.
Иначе нашлась бы кривая в $\Sigma^+$, подходящая к $p$ с обеих сторон.
Она бы отсекла от $\Sigma$ топологический диск, скажем $\Delta$, с краем не ниже $\T_p$ и некоторыми точками строго ниже $\T_p$,
а это противоречит \ref{lem:convex-saddle}. 

Итак, у $\Sigma^+$ как минимум две связные компоненты.

\begin{figure}[!ht]
\vskip-1mm
\centering
\includegraphics{mppics/pic-76}
\caption*{Поверхность $\Sigma$, вид сверху.}
\vskip0mm
\end{figure}

Пусть $\T_p$ задаётся уравнением $z=h(x,y)=a\cdot x+b\cdot y+c$.
Тогда $\Sigma^{+} = \set{(x,y,f(x,y))\in \Sigma}{h(x,y) \leq f(x,y)}$.
Следовательно, $\Sigma^{+}$ содержит связное множество
\[R_-=\set{(x,y,f(x,y))\in\Sigma}{h(x,y)< -C}\] 
и не пересекается с  
\[R_+=\set{(x,y,f(x,y))\in\Sigma}{h(x,y)> C} . \]
Значит одна из связных компонент, скажем, $\Sigma^+_0$, лежит в полосе
\[R_0=\set{(x,y,f(x,y))\in\Sigma}{|h(x,y)|\le  C}.\]
Это множество находится на ограниченном расстоянии от прямой пересечения $\T_p$ с $(x,y)$-плоскостью.

Подвинем $\T_p$ чуть вверх и отсечём от $\Sigma^+_0$ часть выше полученной плоскости, назовём её $\bar\Sigma^+_0$.
Применив \ref{lem:reg-section},
можно считать, что $\bar\Sigma^+_0$ --- гладкая поверхность с краем;
по построению, её край лежит в плоскости.
Расстояние от полученной поверхности $\bar\Sigma^+_0$ до прямой всё ещё ограничено,
что невозможно по \ref{lem:region}.
\qeds

\section{Замечания}

Для нестрого седловых поверхностей, теорема Бернштейна и лемма в её доказательстве не верны;
контрпримеры можно найти среди цилиндрических поверхностей, см. раздел \ref{sec:shape}.
На самом деле других примеров нет;
доказательство основано на той же идее, но сложней.

Согласно критерию горбушек \ref{prop:hat}, седловые поверхности можно определять как гладкие поверхности без горбушек.
Это определение работает для произвольных поверхностей, не обязательно гладких.
Многие из приведённых результатов, например, характеристика седловых графиков, остаются верными для таких обобщённых седловых поверхностей.
Однако этот класс поверхностей пока что недостаточно изучен; см. \cite[глава 4]{alexander-kapovitch-petrunin2019} и ссылки там.

\arxiv{\cleardoublepage
\phantomsection
\AddToShipoutPictureBG*{\includegraphics[width=\paperwidth]{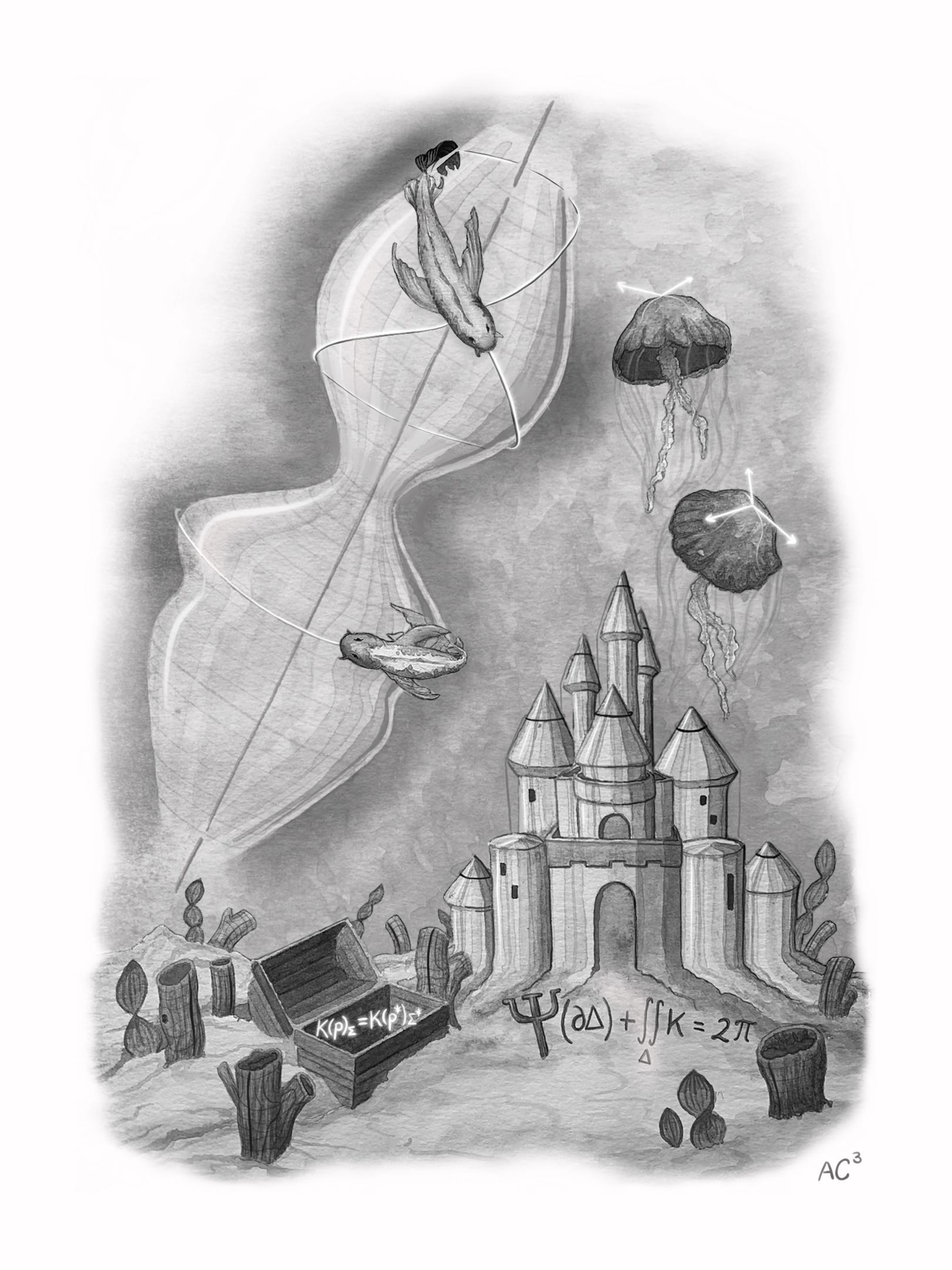}}
\cleardoublepage
\thispagestyle{empty}
\stepcounter{part}
\begin{center}
{\Huge\textbf{Часть \Roman{part}\quad Геодезические}}\
\end{center}
\addcontentsline{toc}{part}{\texorpdfstring{Часть \Roman{part}\quad Геодезические}{Часть \Roman{part} Геодезические}}
\clearpage
}
{\backgroundsetup{
scale=.95,
opacity=1,
angle=0,
vshift=10mm,
contents={%
  \includegraphics[width=\paperwidth
  ]{pics/Geodesics}
  }%
}

\cleardoublepage
\phantomsection
\stepcounter{part}
\addcontentsline{toc}{part}{\texorpdfstring{Часть \Roman{part}\quad Геодезические}{Часть \Roman{part} Геодезические}}
\thispagestyle{empty}
\begin{center}
{\Huge\textbf{Часть \Roman{part}\quad Геодезические}}\
\end{center}
\BgThispage
}

\chapter{Кратчайшие}
\label{chap:shortest}

\section{Внутренняя геометрия}

Геометрию поверхностей разделяют на внутреннюю и внешнюю.
\index{внутренняя геометрия}\emph{Внутренним} называется всё, что можно проверить или померить не вылазя из поверхности;
например, можно мерить длины кривых на поверхности и углы между ними.
Если же существенно используется объемлющее пространство, то мы говорим о чём-то \index{внешняя геометрия}\emph{внешнем};
нормальная кривизна тому пример.

Средняя кривизна, как и гауссова, определяется через главные кривизны, которые являются внешними.
Позднее (\ref{thm:remarkable}) будет показано, что гауссова кривизна является внутренней, то есть её таки можно найти при помощи измерений на самой поверхности.
Средняя же кривизна не является внутренней; например, внутренняя геометрия плоскости не отличается от внутренней геометрии графика $z\z=(x+y)^2$.
Однако, средняя кривизна первой поверхности равна нулю во всех точках, в то время как средняя кривизна второй не равна нулю, например, в начале координат.

Следующее упражнение поможет настроиться на нужный лад.
Может показаться, что это нудная задача по анализу, но на самом деле это забавная задача по геометрии.

\begin{wrapfigure}[5]{r}{33 mm}
\vskip-8mm
\centering
\includegraphics{mppics/pic-77}
\vskip-0mm
\end{wrapfigure}

\begin{thm}{Упражнение}\label{ex:lasso}
Ковбой стоит у подножия ледяной горы, в виде конуса с круглым основанием и углом наклона~$\theta$.
Он забрасывает лассо на вершину конуса, затягивает его и пытается подняться.

Найти максимальный угол $\theta$, при котором ковбой не может забраться на гору.
\end{thm}

\section{Определение}

Пусть $p$ и $q$ --- точки на поверхности~$\Sigma$.
Напомним, что $\dist{p}{q}\Sigma$ обозначает внутреннее расстояние от $p$ до~$q$;
то есть это наименьшая верхняя грань длин путей на $\Sigma$ из $p$ в~$q$.

Кривая $\gamma$ из $p$ в $q$ на $\Sigma$, минимизирующий длину, называется \index{кратчайшая}\emph{кратчайшей}.
Образ такой кратчайшей будет обозначаться $[p,q]$ или $[p,q]_\Sigma$.\index{10aac@$[p,q]$, $[p,q]_\Sigma$ (кратчайшая)}
Если мы пишем $[p,q]_\Sigma$, то предполагается, что кратчайшая существует, и мы выбрали одну из них.
Вообще говоря, таких кратчайших может не быть, также бывает, что их больше одной;
это показано в следующих двух примерах.

\begin{wrapfigure}[9]{r}{28 mm}
\vskip-6mm
\centering
\includegraphics{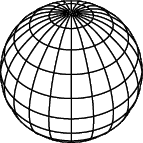}
\bigskip
\includegraphics{mppics/pic-79}
\end{wrapfigure}

\parbf{Неединственность.}
Между полюсами сферы каждый меридиан --- кратчайшая.
Это следует из \ref{obs:S2-length}.

\parbf{Несуществование.}
Пусть $\Sigma$ --- плоскость $z=0$ с удалённым началом координат.
Рассмотрим на ней две точки $p=(1,0,0)$ и $q\z=(-1,0,0)$.
Покажем, что \textit{на~$\Sigma$ нету кратчайших из $p$ в $q$.}

Заметим, что $\dist{p}{q}\Sigma=2$. 
Действительно, пусть $s_\epsilon=(0,\epsilon,0)$ при $\epsilon\z>0$.
Ломаная $ps_\epsilon q$ идёт по $\Sigma$;
её длина $2\cdot\sqrt{1+\epsilon^2}$ стремится к $2$ при $\epsilon\to0$.
Отсюда $\dist{p}{q}\Sigma\le2$.
А поскольку $\dist{p}{q}\Sigma\z\ge \dist{p}{q}{\mathbb{R}^3}=2$, получаем, что $\dist{p}{q}\Sigma=2$.

Следовательно, кратчайшая из $p$ в $q$, если она существует, должна иметь длину~$2$.
По неравенству треугольника, любая кривая длины $2$ из $p$ в $q$ обязана идти вдоль отрезка $[p,q]$;
и в частности, содержать начало координат.
Ну а поскольку в $\Sigma$ нету начала координат, то на ней нет и кратчайшей из $p$ в $q$.

\begin{thm}{Предложение}\label{prop:shortest-paths-exist}
Любые две точки на собственной гладкой поверхности (возможно, с краем) соединимы кратчайшей.
\end{thm}

\parbf{Доказательство.}
Пусть $p$ и $q$ --- точки на собственной гладкой поверхности~$\Sigma$, и $\ell=\dist{p}{q}\Sigma$.

Поскольку $\Sigma$ гладкая, точки $p$ и $q$ соединимы в $\Sigma$ кусочно-гладким путём.
Любой такой путь спрямляем, а значит, $\ell\z<\infty$.

По определению индуцированной внутренней метрики (\ref{sec:Length metric})
на $\Sigma$ существует последовательность путей $\gamma_n$ из $p$ в $q$, таких что
\[\length\gamma_n\to \ell\quad\text{при}\quad n\to \infty.\]

Можно считать, что $\length\gamma_n<\ell+1$ для любого $n$, и каждый $\gamma_n$ параметризован пропорционально длине его дуги.
В частности, каждый путь $\gamma_n\:[0,1]\to\Sigma$ липшицев с константой $(\ell+1)$; 
то есть
\[|\gamma_n(t_0)-\gamma_n(t_1)|\le (\ell+1)\cdot|t_0-t_1|\]
при любых $t_0,t_1\in[0,1]$ и $n$.

При любом $n$, образ $\gamma_n$ лежит в замкнутом шаре $\bar B[p,\ell+1]$.
Следовательно, координатные функции всех $\gamma_n$ равномерно непрерывны и ограничены.
По теореме Арцела --- Асколи (\ref{lem:equicontinuous}),
существует сходящаяся подпоследовательность $\gamma_n$, и её предел, скажем $\gamma_\infty\:[0,1]\to\mathbb{R}^3$, непрерывен.
Кривая $\gamma_\infty$ идёт из $p$ в $q$,
и в частности,
\[\length\gamma_\infty\ge \ell.\]
Поскольку $\Sigma$ --- замкнутое множество, $\gamma_\infty$ лежит в~$\Sigma$.

По полунепрерывности длины (\ref{thm:length-semicont}), $\length\gamma_\infty\le \ell$.
Следовательно, $\length\gamma_\infty\z= \ell$, и $\gamma_\infty$ --- кратчайшая из $p$ в~$q$.
\qeds

\section{Короткая проекция}

Следующая лемма относится к выпуклой геометрии,
но будет играть столь важную роль, что мы её докажем.

\begin{thm}{Лемма}\label{lem:nearest-point-projection}
Пусть $R$ --- замкнутое выпуклое множество в $\mathbb{R}^3$.
Для любой точки $p\in\mathbb{R}^3$ существует единственная точка $\bar p\z\in R$, которая минимизирует расстояние до $R$;
то есть $|p-\bar p|\le |p-x|$ для любой точки $x\in R$.

Более того, отображение $p\mapsto \bar p$ короткое;
то есть
\[|p-q|\ge|\bar p-\bar q| \eqlbl{eq:short-cpp}\]
для любой пары точек $p,q\in \mathbb{R}^3$.
\end{thm}

Отображение $p\mapsto \bar p$ будет называться \label{проекция на ближайшую точку}\index{короткая проекция}\emph{короткой проекцией};
оно отображает евклидово пространство в~$R$.
Отметим, что $\bar p=p$ для любой точки $p\in R$.
В частности, $\bar{\bar p}\equiv\bar p$.

\parbf{Доказательство.}
Выберем $p\in \mathbb{R}^3$.
Пусть 
\[\ell=\inf\set{|p-x|}{x\in R},\]
и $x_n\in R$ такая последовательность, что $|p-x_n|\to \ell$ при $n\to\infty$.

Не умаляя общности, можно считать, что все точки $x_n$ лежат в шаре радиуса $\ell+1$ с центром в~$p$.
Значит, у последовательности $x_n$ найдётся \index{частичный предел}\emph{частичный предел}, скажем $\bar p$;
то есть $\bar p$ --- предел некоторой подпоследовательности $x_n$.
Множество $R$ замкнуто, и, значит, $\bar p\in R$.
По построению 
\begin{align*}
|p-\bar p|&=\lim_{n\to\infty}|p-x_n|=\ell,
\end{align*}
и существование доказано.

{

\begin{wrapfigure}{l}{22 mm}
\vskip-0mm
\centering
\includegraphics{mppics/pic-40}
\vskip-0mm
\end{wrapfigure}

Допустим, что существуют две различные точки $\bar p, \bar p'\in R$, минимизирующие расстояние до~$p$.
Так как $R$ выпукло, их середина $m=\tfrac12\cdot (\bar p+\bar p')$ лежит в~$R$.
Заметим, что $|p-\bar p|\z=|p-\bar p'|=\ell$;
то есть $[p\bar p\bar p']$ --- равнобедренный треугольник.
Следовательно, $[p\bar p m]$ --- прямоугольный треугольник.
Поскольку катет короче гипотенузы, $|p-m|<\ell$ --- противоречие. 

Остаётся проверить \ref{eq:short-cpp}.
Можно считать, что $\bar p\ne\bar q$; иначе и доказывать нечего.

}

Если $\measuredangle \hinge{\bar p}{p}{\bar q}< \tfrac\pi2$, то $\dist{p}{x}{}\z<\dist{p}{\bar p}{}$ для точки $x\in [\bar p,\bar q]$, близкой к $\bar p$,
а это невозможно поскольку $[\bar p,\bar q]\subset R$.

{

\begin{wrapfigure}{r}{37 mm}
\vskip-8mm
\centering
\includegraphics{mppics/pic-41}
\vskip-0mm
\end{wrapfigure}

Значит $\measuredangle \hinge{\bar p}{p}{\bar q}\ge \tfrac\pi2$ или же $p=\bar p$.
В обоих случаях ортогональная проекция $p$ на прямую $\bar p\bar q$ лежит до $\bar p$ или совпадает с $\bar p$.
Аналогично доказывается, что проекция $q$ на прямую $\bar p\bar q$ лежит после $\bar q$ или совпадает с~$\bar q$.
Значит, проекция $[p,q]$ содержит $[\bar p,\bar q]$,
и \ref{eq:short-cpp} следует.
\qeds

}

\section{Кратчайшие на выпуклых поверхностях}

\begin{thm}{Теорема}\label{thm:shorts+convex}
Предположим, что поверхность $\Sigma$ ограничивает замкнутую выпуклую область $R$, и $p,q\in \Sigma$.
Обозначим через $W$ замкнутую внешнюю область относительно $\Sigma$;
то есть $W$ включает $\Sigma$ и дополнение к~$R$.
Тогда 
\[\length\gamma\ge \dist{p}{q}\Sigma\]
для любого пути $\gamma$ в $W$ от $p$ до~$q$.
Более того, если $\gamma$ не лежит на $\Sigma$, то неравенство строгое.
\end{thm}

\parbf{Доказательство.}
Пусть $\bar\gamma$ --- короткая проекция $\gamma$ на $R$.
Кривая $\bar\gamma$ идёт по $\Sigma$ из $p$ в $q$; следовательно, 
\[\length\bar\gamma\ge \dist{p}{q}\Sigma.\]

Давайте покажем, что 
\[\length\gamma\ge\length\bar\gamma.
\eqlbl{bar-gamma=<gamma}\]
Рассмотрим ломаную $\bar p_0\dots \bar p_n$, вписанную в $\bar\gamma$.
Пусть $p_0\dots p_n$ --- соответствующая ломаная, вписанная в $\gamma$;
то есть $p_i=\gamma(t_i)$, если $\bar p_i\z=\bar\gamma(t_i)$.
Согласно \ref{lem:nearest-point-projection}, $|p_i-p_{i-1}|\z\ge|\bar p_i-\bar p_{i-1}|$ для любого~$i$.
Следовательно,
\[\length p_0\dots p_n\ge \length \bar p_0\dots \bar p_n.\]
Применив определение длины, получаем \ref{bar-gamma=<gamma}, и первое утверждение следует;
осталось второе.

\begin{wrapfigure}{o}{37 mm}
\vskip-0mm
\centering
\includegraphics{mppics/pic-82}
\vskip-0mm
\end{wrapfigure}

Допустим, что существует точка $w\z=\gamma(t_1)\z\notin\Sigma$;
заметим, что $w\notin R$.
Тогда существует плоскость $\Pi$, отделяющая $w$ от~$\Sigma$ (\ref{lem:separation}).

Кривая $\gamma$ должна пересечь $\Pi$ в паре точек, до и после $t_1$.
Пусть $a=\gamma(t_0)$ и $b=\gamma(t_2)$ --- эти точки.
Дуга $\gamma$ от $a$ до $b$ длиннее, чем $|a-b|$,
ведь её длина не меньше $|a-w|+|w-b|$, а $|a-w|\z+|w-b|>|a-b|$, так как $w\notin[a,b]$.

Заменим в $\gamma$ дугу от $a$ до $b$ на отрезок $[a,b]$.
Обозначим полученную кривую через $\gamma_1$.
Тогда
\[\length\gamma>\length \gamma_1.\]

Из вышесказанного,
\[\length \gamma_1\ge \dist{p}{q}\Sigma,\]
ибо $\gamma_1$ лежит в~$W$.
И второе утверждение следует.
\qeds

\begin{thm}{Упражнение}\label{ex:length-dist-conv}
Пусть $\Sigma$ --- открытая или замкнутая гладкая поверхность с положительной гауссовой кривизной, и $\Norm$ --- её поле нормалей.
Покажите, что 
\[\dist{p}{q}\Sigma\le 2\cdot \frac{|p-q|}{|\Norm(p)+\Norm(q)|}\]
при любых $p,q\in \Sigma$.
\end{thm}

\begin{wrapfigure}{r}{27 mm}
\vskip-14mm
\centering
\includegraphics{mppics/pic-240}
\end{wrapfigure}

\begin{thm}{Упражнение}\label{ex:hat-convex}
Пусть плоскость $\Pi$ отсекает горбушку $\Delta$ от замкнутой гладкой поверхности~$\Sigma$.
Предположим, что $\Sigma$ ограничивает выпуклую область $R$, и отражение внутренней части $\Delta$ относительно $\Pi$ лежит внутри $R$.
Покажите, что $\Delta$ является \index{выпуклое множество}\emph{выпуклой} относительно внутренней метрики $\Sigma$;
то есть 
если оба конца кратчайшей в $\Sigma$ 
лежат в $\Delta$,
то и вся она лежит в~$\Delta$.
\end{thm}

Определим \index{внутренний диаметр}\emph{внутренний диаметр} поверхности как точную верхнюю грань длин её кратчайших.

\begin{thm}{Упражнение}\label{ex:intrinsic-diameter}
Пусть $\Sigma$ --- замкнутая гладкая поверхность с положительной гауссовой кривизной, лежащая в единичном шаре.

\begin{subthm}{}
Покажите, что внутренний диаметр $\Sigma$ не превышает~$\pi$.
\end{subthm}

\begin{subthm}{}
Покажите, что площадь $\Sigma$ не превосходит $4\cdot \pi$.
\end{subthm}

\end{thm}

\chapter{Геодезические}
\label{chap:geodesics}

\section{Определение}

Гладкая кривая $\gamma$ на гладкой поверхности называется \index{геодезическая}\emph{геодезической}, если её ускорение $\gamma''(t)$ перпендикулярно касательной плоскости $\T_{\gamma(t)}$ при любом~$t$.

\begin{thm}{Упражнение}\label{ex:helix-geodesic}
Покажите, что на цилиндрической поверхности $x^2+y^2=1$,
винтовая линия $\gamma(t)\z\df(\cos t,\sin t, t)$ является геодезической.
\end{thm}

Геодезические описывают траектории частиц, скользящих по $\Sigma$ без трения и сторонних сил.
Ведь если трения нет, то сила, которая удерживает частицу на $\Sigma$, должна быть перпендикулярна~$\Sigma$, и, по второму закону Ньютона, ускорение $\gamma''$ перпендикулярно $\T_{\gamma(t)}$.

Следующую лемму физик вывел бы из закона сохранения энергии.

\begin{thm}{Лемма}\label{lem:constant-speed}
Пусть $\gamma$ --- геодезическая на гладкой поверхности~$\Sigma$. 
Тогда её скорость $|\gamma'|$ постоянна.
Более того, для любого $\lambda\in\mathbb{R}$ кривая 
$\gamma_{\lambda}(t)\df \gamma (\lambda\cdot t)$ также геодезическая. 

\end{thm}

Иными словами, скорость геодезической постоянна, и умножение её параметра на константу оставляет её геодезической.

\parbf{Доказательство.} 
Отметим, что $\gamma''(t)\z\perp\gamma'(t)$,
ибо $\gamma'(t)$ касается поверхности в точке $\gamma(t)$, а $\gamma''(t)$ перпендикулярна к касательной плоскости.
Иначе говоря, $\langle\gamma'',\gamma'\rangle=0$ для любого~$t$.
Отсюда  $\langle\gamma',\gamma'\rangle'\z=2\cdot \langle\gamma'',\gamma'\rangle=0$.
То есть, величина $|\gamma'|^2=\langle\gamma',\gamma'\rangle$ постоянна.

Поскольку $\gamma_{\lambda}''(t) =\lambda^2\cdot \gamma''(\lambda t)$,
получаем вторую часть леммы.
\qeds

{\sloppy

Следующее упражнение описывает \index{соотношение Клеро}\emph{соотношение Клеро};
физик вывел бы его из закона сохранения углового момента.

}

\begin{thm}{Упражнение}\label{ex:clairaut}
Пусть $\gamma$ --- геодезическая на гладкой поверхности вращения,
$r(t)$ --- расстояние от $\gamma(t)$ до оси вращения
и $\theta(t)$ --- угол между $\gamma'(t)$ и параллелью поверхности, проходящей через $\gamma(t)$. 

Покажите, что величина $r(t)\cdot \cos\theta(t)$ не меняется. 
\end{thm}

Напомним, что {}\emph{асимптотическая линия} --- это кривая на поверхности с нулевой нормальной кривизной.

\begin{thm}{Упражнение}\label{ex:asymptotic-geodesic}
Предположим, что кривая $\gamma$ одновременно и геодезическая и асимптотическая линия гладкой поверхности.
Покажите, что $\gamma$ --- прямолинейный отрезок.
\end{thm}

\begin{thm}{Упражнение}\label{ex:reflection-geodesic}
Пусть гладкая поверхность $\Sigma$ пересекается со своей плоскостью симметрии по гладкой кривой~$\gamma$.
Покажите, что~$\gamma$, параметризованная длиной, является геодезической на~$\Sigma$.
\end{thm}

\section{Существование и единственность}

{\sloppy

Следующее предложение говорит, что движение без трения и внешних сил гладко описывается через начальные данные.
В формулировке используется гладкость отображения $w$, которое выдаёт точку в $\mathbb{R}^3$ по точке $p$ поверхности $\Sigma$, касательному в ней вектору $\vec{v}$ и вещественному параметру $t$.
Если выбрать карту $s$ на $\Sigma$, то в ней точка $p$ задаётся парой координат $(u,v)$, а вектор $\vec{v}$ можно представить как сумму $a\cdot s_u+b\cdot s_v$.
Поэтому в локальных координатах $w$ задаётся отображением $(u,v,a,b,t)\mapsto w(p,\vec v,t)$ из подмножества $\mathbb{R}^5$ в $\mathbb{R}^3$, а к нему уже применимо обычное определение гладкости.
Само отображение $(p,\vec v, t)\z\mapsto w(p,\vec v, t)$ считается \index{гладкое отображение}\emph{гладким}, если оно описывается гладким отображением $(u,v,a,b,t)\mapsto w(p,\vec v,t)$ в любых локальных координатах на $\Sigma$.

}

\begin{thm}{Предложение}\label{prop:geod-existence} 
Пусть $\Sigma$ --- гладкая поверхность без края.
Для вектора ${\vec v}$ касательного к $\Sigma$ в точке $p$ существует единственная геодезическая $\gamma\:\mathbb{I}\to \Sigma$, выходящая из $p$ со скоростью~${\vec v}$ (то есть $\gamma(0)=p$ и $\gamma'(0)={\vec v}$) и определённая на максимальном открытом интервале $\mathbb{I}\ni 0$.
Более того,
\begin{subthm}{prop:geod-existence:smooth}
Отображение $w\:(p,{\vec v},t)\mapsto \gamma(t)$ гладкое, и его область определения открыта%
\footnote{То есть если $w$ определено для тройки $p\in \Sigma$, ${\vec v}\in \T_p$ и $t\in \mathbb{R}$,
то оно определено для всякой тройки  $q\in \Sigma$, $\vec u\in \T_q$ и $s\in \mathbb{R}$, где $q$, $\vec u$ и $s$ достаточно близки к $p$, ${\vec v}$ и $t$ соответственно.}%
.
\end{subthm}

\begin{subthm}{prop:geod-existence:whole}
Если $\Sigma$ собственная поверхность, то $\mathbb{I}=\mathbb{R}$; то есть $\gamma$ определена на всей вещественной прямой.
\end{subthm}

\end{thm}

Поверхность, в которой все геодезические бесконечно продолжимы в обе стороны, называется \index{геодезически полная}\emph{геодезически полной}.
Таким образом, согласно \ref{SHORT.prop:geod-existence:whole}, \textit{любая собственная поверхность без края геодезически полна}.
Это часть \index{теорема Хопфа --- Ринова}\emph{теоремы Хопфа --- Ринова} \cite{hopf-rinow}.

В доказательстве мы перепишем определение геодезической через дифференциальное уравнение, а потом применим \ref{thm:ODE-nth-order} и \ref{thm:ODE}.

\begin{thm}{Лемма}\label{lem:geodesic=2nd-order}
Пусть $f$ --- гладкая функция, определённая на открытой области в $\mathbb{R}^2$.
Гладкая кривая $t\mapsto \gamma(t)=(x(t),y(t),z(t))$ является геодезической на графике $z=f(x,y)$ тогда и только тогда, когда $z(t)=f(x(t),y(t))$ для любого $t$, а функции $t\mapsto x(t)$ и $t\mapsto y(t)$ удовлетворяют системе дифференциальных уравнений вида
\[
\begin{cases}
x''=g(x,y,x',y'),
\\
y''=h(x,y,x',y'),
\end{cases}
\]
где $g$ и $h$ --- гладкие функции четырёх переменных, однозначно определяемые функцией~$f$.
\end{thm}

\parbf{Доказательство.}
Первое уравнение $z(t)=f(x(t),y(t))$ означает, что $\gamma(t)$ лежит на графике $z=f(x,y)$.

Далее будем опускать аргументы функций; 
то есть пользоваться сокращениями $x=x(t)$, $f\z=f(x,y)=f(x(t),y(t))$ и так далее.

Сначала выразим $z''$ через $f$, $x$ и $y$.
\[
\begin{aligned}
z''&=f(x,y)''=
\\
&=\left(f_x\cdot x'+ f_y\cdot y'\right)'=
\\
&=
f_{xx}\cdot (x')^2
+
f_x\cdot x''
+2\cdot f_{xy}\cdot x'\cdot y'
+
f_{yy}\cdot (y')^2
+
f_y\cdot y''.
\end{aligned}
\eqlbl{eq:def-geod}
\]

Условие
\[\gamma''(t)\perp\T_{\gamma(t)}\] 
означает, что 
вектор $\gamma''$ перпендикулярен двум базисным векторам в $\T_{\gamma(t)}$, то есть
\[
\begin{cases}
\langle \gamma'',s_x\rangle=0,
\\
\langle \gamma'',s_y\rangle=0,
\end{cases}
\]
где $s(x,y)\df (x,y,f(x,y))$, $x=x(t)$ и $y=y(t)$.
Обратите внимание, что 
$s_x=(1,0,f_x)$ 
и 
$s_y=(0,1,f_y)$.
Поскольку $\gamma''\z=(x'',y'',z'')$, наша система переписывается как
\[
\begin{cases}
x''+ f_x\cdot z''=0,
\\
y''+ f_y\cdot z''=0.
\end{cases}
\]
Остаётся выразить $z''$ из \ref{eq:def-geod}, привести подобные и упростить.
\qeds

\parbf{Доказательство \ref{prop:geod-existence}.}
Пусть $\Sigma$ представляется графиком $z=f(x,y)$ в касательно-нормальных координатах при~$p$.
По лемме~\ref{lem:geodesic=2nd-order}, условие $\gamma''(t)\perp\T_{\gamma(t)}$ записывается системой дифференциальных уравнений второго порядка.
Из \ref{thm:ODE-nth-order} и \ref{thm:ODE} получаем существование и единственность геодезической $\gamma$ на интервале $(-\epsilon,\epsilon)$ для некоторого $\epsilon>0$.

Продолжим геодезическую $\gamma$ до максимального открытого интервала $\mathbb{I}$.
Пусть $\gamma_1$ --- другая геодезическая с теми же начальными данными, определённая на максимальном открытом интервале $\mathbb{I}_1$.
Допустим, что $\gamma_1$ расходится с $\gamma$ при некотором $t_0>0$;
то есть $\gamma_1$ и $\gamma$ совпадают на интервале $[0,t_0)$, но различаются на любом интервале $[0,t_0+\delta)$ при $\delta>0$.
По непрерывности, $\gamma_1(t_0)=\gamma(t_0)$ и $\gamma_1'(t_0)=\gamma'(t_0)$.
Снова применив \ref{thm:ODE-nth-order} и \ref{thm:ODE}, получаем, что $\gamma_1$ совпадает с $\gamma$ в небольшой окрестности $t_0$ --- противоречие.

Аналогичное рассуждение показывает, что $\gamma_1$ не может разойтись с $\gamma$ при $t_0<0$.
Следовательно, $\gamma_1=\gamma$, и в частности, $\mathbb{I}_1=\mathbb{I}$.

Если вся поверхность $\Sigma$ является графиком гладкой функции, то часть \ref{SHORT.prop:geod-existence:smooth} следует из \ref{thm:ODE-nth-order}, \ref{thm:ODE} и леммы.
В этом случае отображение
\[\vec{w}(p,\vec v,t)\df\tfrac{\partial}{\partial t}w(p,\vec v,t)\] также гладкое.
Заметим, что $\vec{w}(p,\vec v,t_0) \in \T_{\gamma(t_0)}$ есть вектор скорости геодезической $\gamma\:t\mapsto w(p,\vec v,t)$ в момент времени $t_0$.

В общем случае пусть $w(p,{\vec v},b)$ определено при $b\ge0$; то есть геодезическая $\gamma\:t\mapsto w(p,{\vec v},t)$ определена в интервале $[0,b]$.
Тогда найдётся такое разбиение $0=t_0<t_1<\dots<t_n=b$ интервала $[0,b]$, что каждая из геодезических $\gamma|_{[t_{i-1},t_i]}$ накрывается картой заданной некоторыми касательно-нормальными координатами.
Пусть $p_i=\gamma(t_i)$ и $\vec v_i=\gamma'(t_i)$ так, что $p_0=p$ и $\vec v_0=\vec v$.
Поскольку $\gamma|_{[t_{i-1},t_i]}$ лежит в графике, по предыдущему рассуждению получаем, что для любого $i$ отображения
$w$ и $\vec w$ определены и гладки в окрестности троек $(p_{i-1},\vec v_{i-1}, t_i-t_{i-1})$.
Заметим, что $p_i=w(p_{i-1},\vec v_{i-1},t_i-t_{i-1})$ и $\vec v_i=\vec w(p_{i-1},\vec v_{i-1},t_i-t_{i-1})$ для любого~$i$.
Поскольку композиция гладких отображений гладкая, отображение $w$ гладко определено в окрестности тройки $(p,\vec v, b)$.

Случай $b\le 0$ аналогичен, и мы получили общий случай в \ref{SHORT.prop:geod-existence:smooth}.

Допустим, что \ref{SHORT.prop:geod-existence:whole} не выполняется;
то есть максимальный интервал $\mathbb{I}$ строго содержится в $\mathbb{R}$.
Не умаляя общности, можно считать, что $b=\sup\mathbb{I}<\infty$.
(Если нет, то обратим параметризацию кривой~$\gamma$.)

{\sloppy

Согласно \ref{lem:constant-speed}, $|\gamma'|$ постоянно; в частности, функция $t\mapsto \gamma(t)$ равномерно непрерывна.
Следовательно, предельная точка
$q\z=\lim_{t\to b}\gamma(t)$
определена, и $q\in \Sigma$, ибо $\Sigma$ собственная.

}

Применив рассуждение выше в касательно-нормальных координатах при $q$, заключаем, что $\gamma$ продолжается геодезической за~$q$.
Следовательно, интервал $\mathbb{I}$ не был максимальным --- противоречие.
\qeds

\begin{thm}{Упражнение}\label{ex:round-torus}
Пусть $\Sigma$ --- гладкий тор вращения; то есть поверхность вращения с гладкой замкнутой образующей кривой.
Докажите, что любая замкнутая геодезическая на $\Sigma$ не стягиваема.

(Другими словами, если $s\:\mathbb{R}^2\to \Sigma$ является естественной би-периодической параметризацией $\Sigma$, то
не существует замкнутой такой кривой $\gamma$ в $\mathbb{R}^2$, что $s\circ\gamma$ является геодезической.)
\end{thm}

\section{Экспоненциальное отображение}\label{sec:exp}

Пусть $p$ --- точка гладкой поверхности $\Sigma$.
Для касательного вектора ${\vec v}\in \T_p$, рассмотрим геодезическую $\gamma_{\vec v}$ на поверхности, выходящую из $p$ с начальной скоростью~$\vec v$;
то есть $\gamma(0)=p$ и $\gamma'(0)={\vec v}$.

Определим \index{экспоненциальное отображение}\emph{экспоненциальное отображение}%
\footnote{Объяснение причины этого термина увело бы нас слишком далеко в сторону.}
в точке $p$ как
\[\exp_p\:\vec v\mapsto \gamma_{\vec v}(1).\]
Согласно \ref{prop:geod-existence}, это отображение гладкое и определено в окрестности нуля касательной плоскости $\T_p$;
более того, если поверхность $\Sigma$ собственная, 
то $\exp_p$ определено на всей плоскости $\T_p$.

Экспоненциальное отображение $\exp_p$ переводит одну поверхность в другую;
первая поверхность --- это касательная плоскость (или её открытое подмножество), а вторая $\Sigma$.
Плоскость $\T_p$ можно отождествить 
со своей касательной плоскостью $\T_0\T_p$, так что дифференциал $d_0(\exp_p)\:\vec v\mapsto D_{\vec v}\exp_p$ отображает $\T_p$ в себя.
Кроме того, по лемме~\ref{lem:constant-speed}, этот дифференциал является тождественным отображением; то есть $(d_0\exp_p)(\vec v)=
\vec v$ для любого $\vec v\in \T_p$.
Получаем следующее.

\begin{thm}{Наблюдение}\label{obs:d(exp)=1}
Пусть $p$ --- точка гладкой поверхности $\Sigma$.

{\sloppy

\begin{subthm}{}
Экспоненциальное отображение $\exp_p$ гладкое, и его область определения $\Dom(\exp_p)$ содержит окрестность нуля в $\T_p$.
Более того, если $\Sigma$ собственная, то $\Dom(\exp_p)=\T_p$.
\end{subthm}

}

\begin{subthm}{}
Дифференциал $d_0(\exp_p)\:\T_p\to \T_p$ является тождественным отображением.
\end{subthm}

\end{thm}

На самом деле легко проверить, что область $\Dom(\exp_p)$ \index{звёздное множество}\emph{звёздная} в $\T_p$;
то есть $\lambda\cdot\vec v\in \Dom(\exp_p)$, если $\vec v\in \Dom(\exp_p)$ и $0\le \lambda\le 1$.

\pagebreak[4]

\section{Радиус инъективности}

\index{радиус инъективности}\emph{Радиусом инъективности} $\inj(p)$ поверхности $\Sigma$ в точке $p$ называется максимальный радиус $r_p\ge 0$ такой, что экспоненциальное отображение $\exp_p$ определено на открытом шаре $B_p\z\df B(0,r_p)_{\T_p}$,
и сужение $\exp_p|_{B_p}$ является гладкой регулярной параметризацией окрестности $p$ в~$\Sigma$.

\begin{thm}{Предложение}\label{prop:exp}
Радиус инъективности положителен в любой точке гладкой поверхности $\Sigma$ (без края).
Более того, он {}\emph{локально отделён от нуля};
то есть для любого $p\in \Sigma$ существует такое $\epsilon>0$, что если $\dist{p}{q}\Sigma<\epsilon$ для некоторого $q\in \Sigma$, то $\inj(q)\ge\epsilon$.
\end{thm}

Предложение будет доказано применением \ref{obs:d(exp)=1} и теоремы об обратной функции (\ref{thm:inverse}).
На самом деле верно, что \textit{функция $\inj\:\Sigma\z\to (0,\infty]$ непрерывна} \cite[5.4]{gromoll-klingenberg-meyer}.

\parbf{Доказательство.}
Пусть $z=f(x,y)$ --- локальное представление $\Sigma$ в касательно-нормальных координатах при~$p$.
В частности, $\T_p$ --- горизонтальная плоскость.

Пусть $h$ --- композиция $\exp_p$ с проекцией $(x,y,z)\z\mapsto (x,y)$.
Согласно \ref{obs:d(exp)=1}, дифференциал $d_0h$ является тождественным;
иными словами, у $h$ единичная матрица Якоби в нуле.
Применив теорему об обратной функции (\ref{thm:inverse}), получим первую часть предложения.

Доказательство второй части аналогично, но более техническое.

Обозначим через $h_q$ композицию $\exp_q$ с ортогональной проекцией $(x,y,z)\mapsto (x,y)$.
Рассмотрим карту $s\:(u,v)\z\mapsto (u,v,f(u,v))$.
Положим 
\[m\:(u,v,a,b)\mapsto h_q(\vec v),\]
где $q=s(u,v)$ и $\vec v=a\cdot s_u+b\cdot s_v$.
По \ref{prop:geod-existence}, $m$ --- гладкое отображение, определённое в окрестности нуля.
Перейдя к меньшей окрестности, можно считать, что первые и вторые частные производные $m$ ограничены.
Из сказанного выше, у отображения $(a,b)\z\mapsto m(0,0,a,b)$ единичная матрица Якоби в нуле.
Значит, при малых $u$ и $v$ матрица Якоби отображения $(a,b)\z\mapsto m(u,v,a,b)$ в нуле близка к единичной.
В частности, к отображению $(a,b)\z\mapsto m(u,v,a,b)$ при малых фиксированных $u$ и $v$
применима вторая часть теоремы об обратной функции (\ref{thm:inverse}), что и завершает доказательство.
\qeds

Идея доказательства следующего предложения приводится в \ref{ex:inj-rad}.

\begin{thm}{Предложение}\label{prop:inj-rad}
Пусть $p$ --- точка гладкой поверхности $\Sigma$ (без края).
Если $\exp_p$ инъективно в $B_p=B(0,r)_{\T_p}$, то сужение $\exp_p|_{B_p}$ есть диффеоморфизм из $B_p$ на свой образ в~$\Sigma$.
\end{thm}

То есть, $\inj(p)$ есть точная верхняя грань таких $r>0$, что сужение $\exp_p|_{B(0,r)_{\T_p}}$ инъективно (отсюда термин \textit{радиус инъективности}).

\section{Кратчайшие и геодезические}

\begin{thm}{Предложение}\label{prop:gamma''}
Любая кратчайшая $\gamma$ на гладкой поверхности $\Sigma$, параметризованная пропорционально длине дуги, является геодезической на~$\Sigma$.
В частности, $\gamma$ --- гладкая кривая.

Обратное верно локально; а именно, для любой точки на $\Sigma$ найдётся окрестность $U$, такая что любая геодезическая, целиком лежащая в $U$, является кратчайшей.
\end{thm}

В частности, достаточно короткий отрезок любой геодезической является кратчайшей.
Если геодезическая является кратчайшей, то её называют \index{минимизирующая геодезическая}\emph{минимизирующей}.
Сейчас мы увидим, что не все геодезические минимизирующие.

\begin{thm}{Упражнение}\label{ex:helix=geodesic}
Пусть $\Sigma$ --- цилиндрическая поверхность, заданная уравнением $x^2\z+y^2=1$.
Покажите, что винтовая линия $\gamma\:[0,2\cdot\pi]\to \Sigma$, определяемая как $\gamma(t)\z\df(\cos t, \sin t, t)$,
является геодезической, но не кратчайшей на~$\Sigma$.
\end{thm}

Полное доказательство предложения будет дано в разделе~\ref{sec:proof-of-gamma''},
но следующее интуитивное объяснение может показаться достаточным.
В предположении гладкости $\gamma$, оно легко переделывается в строгое доказательство.

\parbf{Физическое объяснение.}
Будем думать, что вдоль кратчайшей $\gamma$ натянута резинка, удерживаемая на поверхности реакцией опоры $\vec n$.
Допустим, что  трение отсутствует, а значит, удельная сила $\vec n=\vec n(t)$ пропорциональна нормальному вектору к поверхности в точке~$\gamma(t)$.

Пусть $\tau$ --- натяжение резинки;
оно должно быть одинаковым во всех точках, иначе бы резинка скользила взад-вперёд.

Можно считать, что $\gamma$ параметризована длиной;
тогда равнодействующая сил натяжения дуги $\gamma_{[t_0,t_1]}$ равна $\tau\cdot(\gamma'(t_1)-\gamma'(t_0))$.
Следовательно, её удельная сила при $t_0$ равна
\begin{align*}
\vec f(t_0)&=\lim_{t_1\to t_0}\tau\cdot\frac{\gamma'(t_1)-\gamma'(t_0)}{t_1-t_0}=
\\
&=\tau\cdot\gamma''(t_0).
\end{align*}
По второму закону Ньютона,  
$\vec f+\vec n=0$,
а значит, $\gamma''(t)\perp\T_{\gamma(t)}\Sigma$.
\qeds

\begin{thm}{Следствие}
Пусть $\Sigma$ --- гладкая поверхность, $p\in\Sigma$ и $r\z\le \inj(p)$.
Тогда экспоненциальное отображение $\exp_p$ определяет диффеоморфизм $B(0,r)_{\T_p}\to B(p,r)_\Sigma$.
\end{thm}

\parbf{Доказательство.}
По \ref{prop:inj-rad}, сужение $\exp_p$ на $B_p={B(0,r)_{\T_p}}$ является диффеоморфизмом на его образ $\exp_p(B_p)\subset \Sigma$.

Очевидно, что $B(p,r)_\Sigma\supset\exp_p(B_p)$.
По \ref{prop:gamma''}, $B(p,r)_\Sigma\subset\exp_p(B_p)$, отсюда результат.
\qeds

{\sloppy

Согласно следствию, сужение $\exp_p|_{B(0,r)_{\T_p}}$ допускает обратное отображение, называемое \index{логарифм}\emph{логарифмом};
оно будет обозначается как
\[\log_p\:B(p,r)_\Sigma\to B(0,r)_{\T_p}.\]

}

По предложению выше, любая кратчайшая, параметризованная длиной гладкая.
Это поможет решить следующую пару упражнений.

\begin{thm}{Упражнение}\label{ex:two-min-geod}
{\sloppy
Покажите, что если у двух кратчайших есть две различные общие точки $p$ и $q$, то либо это их общие концы, либо же у кратчайших есть общая дуга от $p$ до~$q$.

}

Постройте пример, геодезических без общих дуг, которые пересекаются в произвольно большом числе точек.
\end{thm}

\begin{thm}{Упражнение}\label{ex:min-geod+plane}
{\sloppy
Предположим, что гладкая поверхность $\Sigma$ зеркально-симметрична относительно плоскости $\Pi$.
Покажите, что кратчайшая на $\Sigma$ не может проходить сквозь $\Pi$ больше раза.
Другими словами, если идти вдоль кратчайшей, то стороны $\Pi$ сменятся не больше одного раза.

}

\end{thm}

{

\begin{thm}{Продвинутое упражнение}\label{ex:milka}\\
Пусть $\gamma\:[0,\ell]\z\to \Sigma$ --- минимизирующая геодезическая с единичной скоростью на 
гладкой замкнутой строго выпуклой поверхности $\Sigma$.

Положим $p\z=\gamma(0)$, $q=\gamma(\ell)$ и 
\[p^s=\gamma(s)-s\cdot\gamma'(s).\]

\begin{wrapfigure}{r}{40 mm}
\vskip-4mm
\centering
\includegraphics{mppics/pic-250}
\vskip-0mm
\end{wrapfigure}

Покажите, что $q$ не видна из $p^s$ при любом $s\in (0,\ell)$;
то есть отрезок $[p^s,q]$ пересекает $\Sigma$ в точке, отличной от~$q$.

Покажите, что утверждение перестаёт быть верным без предположения, что $\gamma$ минимизирующая.
\end{thm}

}

\begin{wrapfigure}[3]{r}{40 mm}
\end{wrapfigure}

\begin{thm}{Упражнение}\label{ex:round-sphere}
Пусть $\Sigma$ --- гладкая замкнутая поверхность.
Предположим, что для любых $p,q\z\in \Sigma$ расстояние $\dist{p}{q}\Sigma$ зависит только от расстояния $\dist{p}{q}{\mathbb{R}^3}$.
Покажите, что $\Sigma$ есть сфера.
\end{thm}

\begin{thm}{Сильно продвинутое упражнение}\label{ex:rad=2}
Пусть \(\Theta\) — сфера радиуса $2$ с центром в $0\in\mathbb{R}^3$,
и пусть \( \Sigma \) — гладкая замкнутая поверхность, содержащаяся в открытом шаре, ограниченном \(\Theta\).
Предположим, что все нормальные кривизны \( \Sigma \) не превышают~$1$ по абсолютной величине.

\begin{subthm}{ex:rad=2:a}
Докажите, что существует диффеоморфизм \(\rho\:\Theta \to \Sigma \) (назовем его \emph{радиальной проекцией}),
который отображает точку \( p \in \Theta \) в единственную точку пересечения \( \Sigma \cap [0,p]_{\mathbb{R}^3} \).
В частности, \( \Sigma \) ограничивает звёздную область.
\end{subthm}

\begin{subthm}{ex:rad=2:b}
Докажите, что \(\rho\: \Theta \to \Sigma \) не увеличивает длины кривых.
\end{subthm}

\begin{subthm}{ex:rad=2:c}
Пусть \( x\in \Sigma \) — точка на максимальном расстоянии от начала координат.
Докажите, что шар с диаметром \( [0, x] \) лежит в области, ограниченной \( \Sigma \).
Выведите отсюда, что эта область содержит единичный шар.
\end{subthm}

\end{thm}

\section{Лемма Либермана}

Вариант следующей леммы использовался Иосифом Либерманом \cite{liberman}.

\begin{thm}{Лемма}
\label{lem:liberman}
\index{лемма Либермана}
Пусть $f$ --- гладкая локально выпуклая функция, определённая на открытом подмножестве плоскости,
и $t\mapsto \gamma(t)\z=(x(t),y(t),z(t))$ --- геодезическая параметризованная длиной на графике $z=f(x,y)$.
Тогда $t\mapsto z(t)$ --- выпуклая функция; то есть $z''(t)\ge 0$ для любого~$t$.
\end{thm}

\parbf{Доказательство.}
Выберем ориентацию графика так, чтобы нормаль $\Norm$ всегда указывала вверх;
то есть в каждой точке $z$-координата $\Norm$ была положительна.
Будем использовать сокращение $\Norm(t)$ для $\Norm(\gamma(t))$.

Поскольку $\gamma$ --- геодезическая, $\gamma''(t)\perp\T_{\gamma(t)}$;
иначе говоря, ускорение $\gamma''(t)$ пропорционально $\Norm(t)$ при любом~$t$.
Более того,
\[\gamma''=k\cdot\Norm,\]
где $k=k(t)$ --- нормальная кривизна графика в точке $\gamma(t)$ и направлении $\gamma'(t)$.

Следовательно,
\[z''=k\cdot\cos\theta,
\eqlbl{eq:z''}\]
где $\theta=\theta(t)$ --- угол между $\Norm(t)$ и осью $z$.

Поскольку $\Norm$ смотрит вверх, $\theta(t)<\tfrac\pi2$, и, значит, $\cos\theta>0$.

Так как $f$ выпукла, касательная плоскость подпирает график снизу в любой точке;
в частности, $k(t)\ge 0$ для любого~$t$.
Значит правая часть в \ref{eq:z''} неотрицательна.
Лемма доказана.
\qeds

\begin{thm}{Упражнение}\label{ex:closed-liberman}
Пусть $\Sigma$ --- это график локально выпуклой функции, определённой на открытом подмножестве плоскости.
Покажите, что $\Sigma$ не содержит замкнутых геодезических.
\end{thm}

\begin{thm}{Упражнение}\label{ex:rho''}
Пусть $\gamma$ --- геодезическая с единичной скоростью на гладкой выпуклой поверхности $\Sigma$, и точка $p$ лежит внутри выпуклого множества, ограниченного~$\Sigma$.
Рассмотрим функцию $\rho(t)=|p-\gamma(t)|^2$.
Покажите, что $\rho''(t)\le 2$ для любого~$t$.
\end{thm}

\section{Полная кривизна геодезической}

Напомним, что $\tc\gamma$ обозначает полную кривизну кривой~$\gamma$, см.~\ref{sec:Total curvature}.

\begin{thm}{Упражнение}\label{ex:tc-spherical-image}
Пусть $\gamma$ --- геодезическая на гладкой поверхности $\Sigma$ 
с полем нормалей $\Norm$.
Покажите, что $\length(\Norm\circ\gamma)\ge \tc\gamma$.
\end{thm}

\begin{thm}{Теорема}\label{thm:usov}
Пусть $\Sigma$ --- график выпуклой $\ell$-липшицевой функции $f$, определённой на открытом множестве в плоскости $(x,y)$.
Тогда полная кривизна любой геодезической на $\Sigma$ не превышает $2\cdot \ell$.
\end{thm}

Эта теорема доказана Владимиром Усовым \cite{usov}.

\parbf{Доказательство.}
Пусть $t\mapsto\gamma(t)=(x(t),y(t),z(t))$ --- геодезическая с единичной скоростью на~$\Sigma$.
По лемме Либермана (\ref{lem:liberman}), функция $t\mapsto z(t)$ выпукла.

Поскольку угловой коэффициент $f$ не превышает $\ell$, 
$|z'(t)|\le \frac{\ell}{\sqrt{1+\ell^2}}$
при любом $t$.
Можно считать, что $\gamma$ определена на интервале $[a,b]$.
Тогда
\[
\begin{aligned}
\int_a^b z''(t) dt&=z'(b)-z'(a)\le 
 2\cdot \frac{\ell}{\sqrt{1+\ell^2}}.
\end{aligned}
\eqlbl{eq:intz''}
\]

Также отметим, что $z''$ является проекцией $\gamma''$ на ось $z$.
При этом угловой коэффициент касательной плоскости $\T_{\gamma (t)} \Sigma$ не превышает $\ell$, для любого~$t$, и
\[|\gamma'' (t)| \le z''(t)\cdot\sqrt{1+ \ell ^2},\]
ибо $\gamma ''$ перпендикулярна этой плоскости.

Из \ref{eq:intz''}, получаем, что
\begin{align*}
\tc\gamma&=\int_a^b|\gamma'' (t)|\cdot dt\le 
\sqrt{1+ \ell ^2}\cdot \int_a^b z''(t)\cdot dt\le 
2\cdot \ell.
\end{align*}
\qedsf

По следующему упражнению, оценка в теореме оптимальна.

\begin{thm}{Упражнение}\label{ex:usov-exact}
Пусть $\Sigma$ --- график $z=\ell\cdot\sqrt{x^2+y^2}$ с удалённым началом координат.
Покажите, что любая бесконечная в обе стороны геодезическая $\gamma$ на $\Sigma$ имеет полную кривизну ровную $2\cdot \ell$.
\end{thm}

\begin{thm}{Упражнение}\label{ex:ruf-bound-mountain}
Предположим, что $f$ является гладкой выпуклой $\tfrac32$-липшицевой функцией, определённой на плоскости $(x,y)$.
Покажите, что любая геодезическая $\gamma$ на графике $z\z=f(x,y)$ является простой; то есть не имеет самопересечений.

{\sloppy

Постройте выпуклую $2$-липшицеву функцию, определённую на плоскости,
с самопересекающейся геодезической $\gamma$ на её графике.

}

\end{thm}

\begin{thm}{Теорема}\label{thm:tc-of-mingeod}
Пусть гладкая поверхность $\Sigma$ ограничивает выпуклое множество $K$, и при этом $B(0,\epsilon)\subset K\subset B(0,1)$.
Тогда полную кривизну любой кратчайшей на $\Sigma$ можно оценить, зная~$\epsilon$.
\end{thm}

{

\begin{wrapfigure}{r}{48 mm}
\vskip4mm
\centering
\includegraphics{mppics/pic-83}
\vskip-0mm
\end{wrapfigure}

\begin{thm}{Доказательство и упражнение}\label{ex:bound-tc}
Пусть $\Sigma$ как в теореме, $\gamma$ --- кратчайшая на~$\Sigma$ с единичной скоростью, $\Norm(t)$ --- вектор нормали к $\Sigma$ в точке $\gamma(t)$, направленный наружу, $\theta(t)$ --- угол между $\Norm(t)$ и направлением к $\gamma(t)$ из начала координат,
и 
$k(t)$ --- кривизна $\gamma$ при~$t$.
Рассмотрим функцию $\rho(t)\z=|\gamma(t)|^2$.

\begin{subthm}{ex:bound-tc:a}
Покажите, что $\cos(\theta(t))\ge \epsilon$ для любого~$t$.
\end{subthm}

\begin{subthm}{ex:bound-tc:b}
Покажите, что $|\rho'(t)|\le 2$ для любого~$t$.
\end{subthm}

\begin{subthm}{ex:bound-tc:c}
Покажите, что
\[\rho''(t)=2-2\cdot k(t)\cdot \cos \theta(t)\cdot |\gamma(t)|\]
для любого~$t$.
\end{subthm}

\begin{subthm}{ex:bound-tc:d}
Воспользуйтесь короткой проекцией из единичной сферы на $\Sigma$, дабы показать, что
\[\length \gamma\le \pi.\]
\end{subthm}

\begin{subthm}{ex:bound-tc:e}
Выведите отсюда, что $\tc\gamma\le 100/\epsilon^2$.
\end{subthm}

\end{thm}

}

\parit{Замечание.}
Полученная оценка идёт к бесконечности при $\epsilon\to 0$,
но известна и оценка, не зависящая от $\epsilon$;
это результат Нины Лебедевой и первого автора \cite{lebedeva-petrunin}.
Алексей Погорелов выдвинул гипотезу, что существует оценка на длину сферического образа кратчайшей \cite{pogorelov}.
Согласно \ref{ex:tc-spherical-image}, эта гипотеза сильнее,
однако к ней и всевозможным её вариантам  нашлись контрпримеры \cite{zalgaller,milka,usov,pach}.

\chapter{Параллельный перенос}
\label{chap:parallel-transport}

\section{Параллельные поля}

Пусть $\gamma\:[a,b]\z\to \Sigma$ --- гладкая кривая на гладкой поверхности $\Sigma$.
Гладкая векторнозначная функция $t\mapsto {\vec v}(t) \in \mathbb{R}^3$ называется \index{касательное!поле}\emph{касательным полем} вдоль $\gamma$, если
вектор ${\vec v}(t)$ лежит в касательной плоскости $\T_{\gamma(t)}\Sigma$ для любого~$t$.

Касательное поле ${\vec v}$ на кривой $\gamma$ называется \index{параллельные поле и перенос}\emph{параллельным}, если ${\vec v}'(t)\perp\T_{\gamma(t)}$ для любого~$t$.

В общем случае семейство касательных плоскостей $\T_{\gamma(t)}\Sigma$ не является параллельным.
Поэтому никак нельзя ожидать наличия \textit{по-настоящему} параллельного поля с ${\vec v}'(t)\equiv 0$ вдоль $\gamma$.
Условие ${\vec v}'(t)\z\perp\T_{\gamma(t)}$ значит, что поле старается вести себя как можно \textit{параллельней} --- ему приходится крутиться вместе с касательной плоскостью, но оно не крутится внутри оной.

По определению геодезической, поле скоростей ${\vec v}(t)\z=\gamma'(t)$ любой геодезической $\gamma$  параллельно вдоль~$\gamma$.

\begin{thm}{Упражнение}\label{ex:parallel}
Пусть $\gamma\:[a,b]\to \Sigma$ --- гладкая кривая на гладкой поверхности $\Sigma$, а ${\vec v}(t)$ и $\vec w(t)$ --- параллельные векторные поля вдоль~$\gamma$.

\begin{subthm}{ex:parallel:a}
Покажите, что величина $|{\vec v}(t)|$ постоянна.
\end{subthm}

\begin{subthm}{ex:parallel:b}
Покажите, что угол $\theta(t)=\measuredangle({\vec v}(t),\vec w(t))$ постоянен.
\end{subthm}

\end{thm}

\section{Параллельный перенос}

\begin{thm}{Предложение и определение}\label{prop:parallel}
Пусть $p\z=\gamma(a)$ и $q\z=\gamma(b)$ --- концы кусочно-гладкой кривой $\gamma\:[a,b]\z\to \Sigma$ на гладкой поверхности~$\Sigma$.

Для данного касательного вектора ${\vec w}\in\T_p$ существует единственное параллельное такое поле ${\vec w}(t)$ вдоль $\gamma$, что ${\vec w}(a)={\vec w}$.

Вектор ${\vec w}(b)\in\T_q$ называется \index{параллельные поле и перенос}\emph{параллельным переносом} вектора ${\vec w}(a)$ вдоль~$\gamma$ на~$\Sigma$.
\end{thm}

Параллельный перенос вдоль $\gamma$ будет обозначаться как $\iota_\gamma$;
то есть можно писать $\vec w(b)=\iota_\gamma({\vec w}(a))$ или $\vec w(b)=\iota_\gamma({\vec w}(a))_\Sigma$, если нужно подчеркнуть, что $\gamma$ лежит на поверхности~$\Sigma$.
Из упражнения~\ref{ex:parallel} следует, что параллельный перенос $\iota_\gamma\:\T_p\z\to\T_q$ является изометрией.
В общем случае, параллельный перенос $\iota_\gamma\:\T_p\z\to\T_q$ зависит от выбора $\gamma$; то есть, для другой кривой $\gamma_1$, идущей из $p$ в $q$, параллельные переносы $\iota_{\gamma_1}$ и $\iota_{\gamma}$ могут различаться.

{\sloppy

\parbf{Идея доказательства.}
Предположим, что кривая $\gamma$ гладкая и она покрыта картой $(u,v)\mapsto s(u,v)$ поверхности $\Sigma$;
тогда $\gamma(t)\z=s(u(t),v(t))$ для гладких функций $t\mapsto u(t)$ и $t\mapsto v(t)$.
Определим 
\[
\vec u(t)=s_u(u(t),v(t)),
\quad
\vec v(t)=s_v(u(t),v(t)),
\quad
\text{и}
\quad
\Norm(t)=\Norm(\gamma(t)).
\]

}

Условия $\vec w(t)\in \T_{\gamma(t)}$ и $\vec w'(t)\perp \T_{\gamma(t)}$ записываются системой уравнений
\[
\begin{cases}
\ \langle\vec w(t), \Norm(t)\rangle=0,
\\
\ \langle\vec w'(t), \vec u(t)\rangle=0,
\\
\ \langle\vec w'(t), \vec v(t)\rangle=0.
\end{cases}
\]
Переписав эту систему через компоненты $\vec u$, $\vec v$, $\vec w$ и $\Norm$, получим систему обыкновенных дифференциальных уравнений на компоненты $\vec w$, и, применив \ref{thm:ODE}, получим решение $\vec w(t)$.
По \ref{ex:parallel}, $|\vec w(t)|$ не меняется;
следовательно, \ref{thm:ODE} влечёт, что решение $\vec w$ определено на всём интервале $[a,b]$.

Если $\gamma$ только кусочно-гладкая и/или не покрыта одной картой, то можно её представить как произведение гладких дуг $\gamma_1,\dots,\gamma_n$ таких, что каждая $\gamma_i$ покрыта одной картой.
Последовательно применив предложение к каждой $\gamma_i$, получим его для $\gamma$.
\qeds

Предположим, что $\gamma_1$ и $\gamma_2$ --- две гладкие кривые на двух гладких поверхностях $\Sigma_1$ и $\Sigma_2$ со сферическими отображениями $\Norm_i\:\Sigma_i\to\mathbb{S}^2$.
Если $\Norm_1\circ\gamma_1(t)= \Norm_2\circ\gamma_2(t)$ для любого $t$, то будем говорить, что кривые $\gamma_1$ и $\gamma_2$ имеют \emph{идентичные нормали} на $\Sigma_1$ и $\Sigma_2$ соответственно.

В этом случае касательная плоскость $\T_{\gamma_1(t)}\Sigma_1$ параллельна касательной плоскости $\T_{\gamma_2(t)}\Sigma_2$ для любого~$t$, и можно отождествить $\T_{\gamma_1(t)}\Sigma_1$ с $\T_{\gamma_2(t)}\Sigma_2$.
В частности, любое поле $\vec v(t)$ касательное вдоль $\gamma_1$ является
касательным вдоль $\gamma_2$.
Более того, $\vec v'(t)\perp \T_{\gamma_1(t)}\Sigma_1$ эквивалентно $\vec v'(t)\perp \T_{\gamma_2(t)}\Sigma_2$; то есть, если $\vec v(t)$ параллельно вдоль $\gamma_1$,
то оно параллельно и вдоль $\gamma_2$.

Рассуждение выше приводит к следующему наблюдению, которое скоро себя проявит.

\begin{thm}{Наблюдение}\label{obs:parallel=}
Пусть $\gamma_1$ и $\gamma_2$ --- гладкие кривые на гладких поверхностях $\Sigma_1$ и $\Sigma_2$.
Предположим, что $\gamma_1$ и $\gamma_2$ имеют идентичные нормали (как кривые на $\Sigma_1$ и $\Sigma_2$ соответственно).
Тогда параллельные переносы $\iota_{\gamma_1}$ и $\iota_{\gamma_2}$ идентичны. 
\end{thm}

\begin{thm}{Упражнение}\label{ex:parallel-transport-support}
Пусть $\Sigma_1$ и $\Sigma_2$ --- две поверхности с общей кривой~$\gamma$.
Предположим, что $\Sigma_2$ содержится в области, ограниченной $\Sigma_1$.
Покажите, что параллельный перенос вдоль $\gamma$ на $\Sigma_1$ 
совпадает с параллельным переносом вдоль $\gamma$ на $\Sigma_2$. 
\end{thm}

\section{Велосипедное колесо и проекции}

Следующие толкования параллельного переноса помогут обзавестись правильной интуицией,
но не помогут писать доказательства.
Первое, основанное на физическом эксперименте с велосипедным колесом, предложено Марком Леви \cite{levi}.
Второе, близко с ним связанное, даётся через ортогональные проекции касательных плоскостей друг на друга.

\begin{figure}[ht!]
\vskip-0mm
\centering
\includegraphics[scale=.3]{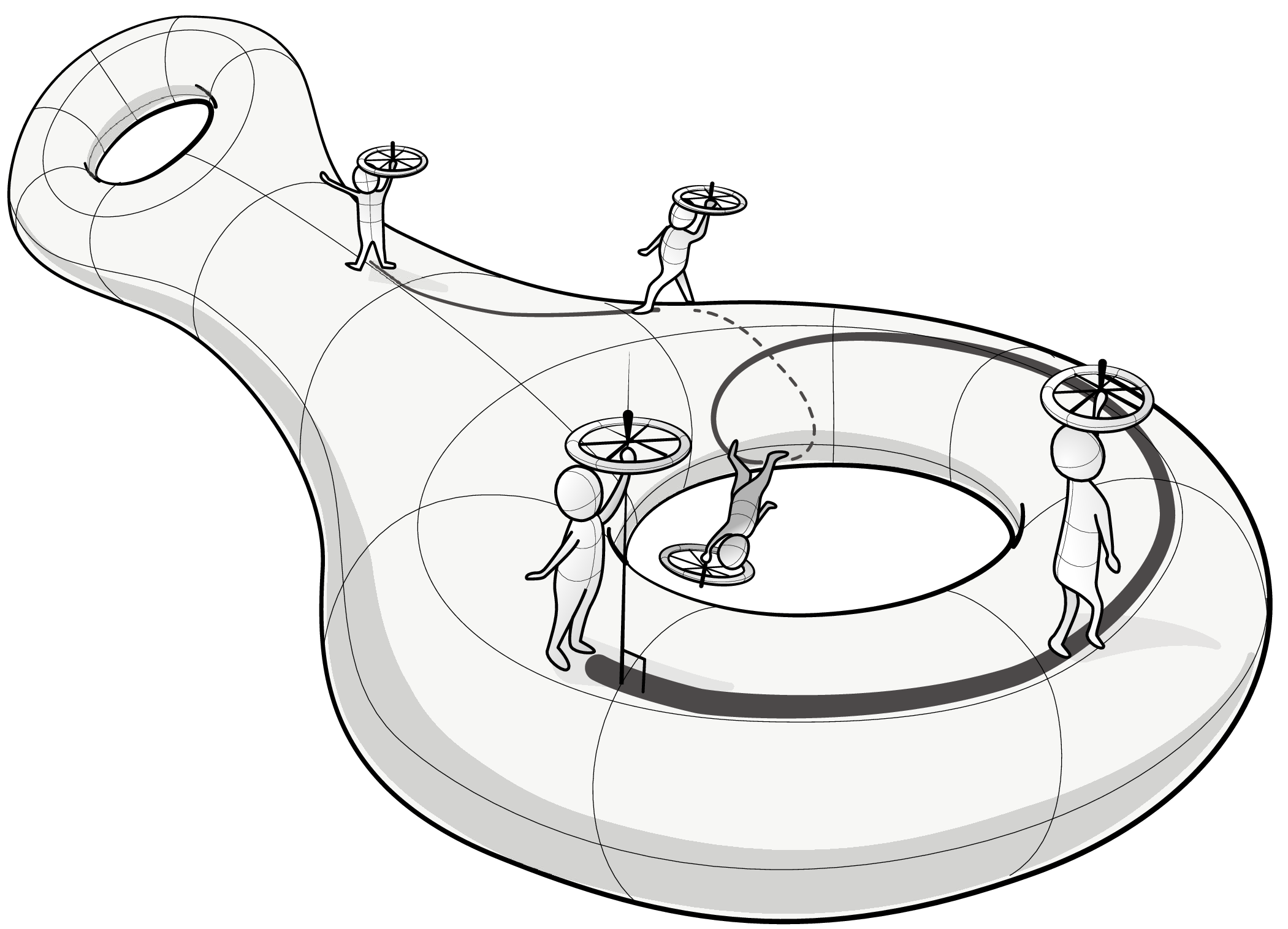}
\end{figure}

Пусть $\gamma\:[a,b]\to\Sigma$ --- гладкая кривая на гладкой поверхности $\Sigma$.
Представьте, что мы идём по поверхности вдоль $\gamma$, неся идеально сбалансированное велосипедное колесо.
При этом ось колеса всё время направлена перпендикулярно $\Sigma$ и мы не трогаем само колесо, держась только за его ось.
Если колесо не вращается в начальной точке $p\z=\gamma(a)$, то оно не будет вращаться после остановки в точке~$q=\gamma(b)$, ведь, прилагая силу к оси, нельзя придать колесу крутящий момент.
Тогда отображение, переводящее начальное положение колеса в конечное, и есть параллельный перенос~$\iota_\gamma$.

Физическая интуиция должна подсказывать, что
\textit{если не менять направление оси,
то не будут меняться и направления спиц};
по сути это переформулировка наблюдения~\ref{obs:parallel=}.

Для второго толкования, выберем разбиение $a=t_0<\dots\z<t_n=b$ отрезка $[a,b]$
и рассмотрим последовательность ортогональных проекций $\phi_i\:\T_{\gamma(t_{i-1})}\to\T_{\gamma(t_i)}$.
Для достаточно мелкого разбиения композиция 
\[\phi_n\circ\dots\circ\phi_1\:\T_p\z\to\T_q\]
становится близкой к $\iota_\gamma$.

Каждая проекция $\phi_i$ не увеличивает длину векторов и то же верно для их композиции.
Однако, для достаточно мелкого разбиения композиция является почти изометрией.
А значит, предел $\iota_\gamma$ --- изометрия $\T_p\z\to\T_q$.

\begin{thm}{Упражнение}\label{ex:holonomy=not0}
Постройте такую петлю $\gamma$ на $\mathbb{S}^2$, что параллельный перенос $\iota_\gamma$ не является тождественным отображением.
\end{thm}

\section{Полная геодезическая кривизна}

Напомним, что геодезическая кривизна определена в~\ref{sec:Darboux}.
Она обозначается через $k_g$ и измеряет, насколько кривая $\gamma$ на ориентированной поверхности отклоняется от геодезической;
если $\gamma$ поворачивает налево, то кривизна положительна, если направо, то отрицательна.
В частности, геодезические линии имеют нулевую геодезическую кривизну.

Пусть $\gamma\:\mathbb{I}\to \Sigma$ --- гладкая кривая с единичной скоростью на ориентированной гладкой поверхности $\Sigma$.
Полная геодезическая кривизна $\gamma$ определяется через интеграл 
\[\tgc\gamma=\tgc\gamma_\Sigma
\df
\int_{\mathbb{I}} k_g(t)\cdot dt.\]

Если $\Sigma$ --- это плоскость, то геодезическая кривизна $\gamma$ совпадает
с её ориентированной кривизной (см.~\ref{sec:Total signed curvature}).
В частности, её полная геодезическая кривизна равна её полной ориентированной кривизне.
По этой причине мы воспользовались тем же обозначением $\tgc\gamma$; если надо указать поверхность, то пишем $\tgc\gamma_\Sigma$.\index{10psi@$\tgc\gamma_\Sigma$ (полная геодезическая кривизна)}

Если $\gamma$ --- кусочно-гладкая кривая на $\Sigma$, то её \index{полная!геодезическая кривизна}\emph{полная геодезическая кривизна} определяется как сумма всех полных геодезических кривизн её дуг плюс сумма ориентированных внешних углов $\gamma$ на стыках.
То есть, если $\gamma$ --- произведение гладких кривых $\gamma_1,\dots,\gamma_n$, то
\[\tgc\gamma
\df
\tgc{\gamma_1}+\dots+\tgc{\gamma_n}+\theta_1+\dots+\theta_{n-1},\]
где $\theta_i$ --- ориентированный внешний угол на стыке $\gamma_i$ и $\gamma_{i+1}$;
он положителен, если $\gamma$ поворачивает влево, отрицателен, если вправо, и не определён, если $\gamma$ делает полный разворот.
Если $\gamma$ замкнута, то 
\[\tgc\gamma
\df
\tgc{\gamma_1}+\dots+\tgc{\gamma_n}+\theta_1+\dots+\theta_n,\]
где $\theta_n$ --- это ориентированный внешний угол на стыке $\gamma_n$ и $\gamma_1$.

Если каждая дуга $\gamma_i$ в произведении является минимизирующей геодезической, то $\gamma$ называется \index{ломаная геодезическая}\emph{ломаной геодезической}.
В этом случае $\tgc{\gamma_i}=0$ для каждого~$i$.
Следовательно, полная геодезическая кривизна $\gamma$ равна сумме её ориентированных внешних углов.

\begin{thm}{Предложение}\label{prop:pt+tgc}
Пусть $\gamma$ --- замкнутая ломаная геодезическая на гладкой ориентированной поверхности $\Sigma$, которая начинается и заканчивается в точке~$p$.
Тогда параллельный перенос $\iota_\gamma\:\T_p\to\T_p$ представляет собой вращение по часовой стрелке плоскости $\T_p$ на угол~$\tgc\gamma$.

Утверждение остается верным для гладких замкнутых кривых и для кусочно-гладких кривых.
\end{thm}

\parbf{Доказательство.}
Пусть $\gamma$ --- циклическое произведение геодезических $\gamma_1,\dots,\gamma_n$.
Выберем касательный вектор ${\vec v}$ в точке $p$ и продолжим его до параллельного векторного поля вдоль~$\gamma$.
Поскольку поле $\tan_i(t)\z=\gamma_i'(t)$ параллельно вдоль $\gamma_i$, ориентированный угол $\phi_i$ от $\tan_i$ до ${\vec v}$ постоянен на каждой дуге~$\gamma_i$.

\begin{wrapfigure}{o}{22 mm}
\vskip-0mm
\centering
\includegraphics{mppics/pic-48}
\vskip-0mm
\end{wrapfigure}

Пусть $\theta_i$ --- внешний угол на стыке $\gamma_{i}$ и $\gamma_{i+1}$.
Тогда 
\[\phi_{i+1}=\phi_i-\theta_i \pmod{2\cdot\pi}.\]
Таким образом, после обхода всех вершин мы получаем, что 
\[\phi_{n+1}-\phi_1=-\theta_1-\dots-\theta_n=-\tgc\gamma \pmod{2\cdot\pi},\]
и первое утверждение следует.

Для гладкой кривой с единичной скоростью $\gamma\:[a,b]\to\Sigma$ доказательство аналогично.
Обозначим через $\phi(t)$ ориентированный угол от $\tan (t)=\gamma'(t)$ до ${\vec v} (t)$ и покажем, что
\[\phi'(t)+k_g(t)\equiv0\eqlbl{eq:phi'+kg}\]

Напомним, что $\mu=\mu(t)$ --- это поворот $\tan \z=\tan(t)$ против часовой стрелки на угол $\tfrac\pi2$ в $\T_{\gamma(t)}$.
Обозначим через $\vec w=\vec w(t)$ поворот $\vec v=\vec v(t)$ против часовой стрелки на угол $\tfrac\pi2$ в $\T_{\gamma(t)}$.
Тогда
\begin{align*}
\tan&=\cos\phi\cdot \vec v-\sin\phi\cdot \vec w,
\\
\mu&=\sin\phi\cdot \vec v+\cos\phi\cdot \vec w.
\end{align*}

Векторные поля $\vec v$ и $\vec w$ параллельны вдоль $\gamma$; то есть $\vec v'(t),\vec w'(t)\z\perp\T_{\gamma(t)}$.
Следовательно, $\langle\vec v',\mu\rangle\z=\langle\vec w',\mu\rangle\z=0$.
Из этого следует, что
\begin{align*}
k_g&=\langle\tan',\mu\rangle=
\\
&=-(\sin^2\phi+\cos^2\phi)\cdot \phi'
\\
&=-\phi',
\end{align*}
что и доказывает \ref{eq:phi'+kg}.

Применив \ref{eq:phi'+kg}, получаем, что 
\begin{align*}
\phi(b)-\phi(a)&=\int_a^b \phi'(t)\cdot dt=
\\
&=-\int_a^b k_g\cdot dt=
\\
&=-\tgc\gamma
\end{align*}

Если $\gamma$ --- кусочно-гладкая, то результат доказывается прямым сочетанием доказательств двух рассмотренных случаев.
\qeds

\chapter{Формула Гаусса --- Бонне}
\label{chap:gauss-bonnet}
\section{Формулировка}

Следующая теорема доказана Карлом Фридрихом Гауссом \cite{gauss}
для геодезических треугольников;
Пьер Бонне и Жак Бине независимо 
обобщили её на произвольные кривые.
\index{формула Гаусса --- Бонне}

\begin{thm}{Теорема}\label{thm:gb}
Пусть $\Delta$ --- топологический диск на гладкой ориентированной поверхности $\Sigma$, 
а его граница $\partial\Delta$ --- кусочно-гладкая кривая так ориентированная, что $\Delta$ лежит слева от неё.
Тогда 
\[\tgc{\partial\Delta}+\iint_\Delta K=2\cdot \pi,\eqlbl{eq:g-b}\]
где $K$ --- гауссова кривизна.
\end{thm}

В этой главе даётся доказательство частного случая этой формулы, его можно обобщить, но полное доказательство будет строиться на другой идее, см.~\ref{sec:gauss--bonnet:formal}.
Сначала, давайте попрактикуемся в приложениях формулы.

\begin{thm}{Упражнение}\label{ex:1=geodesic-curvature}
Пусть $\gamma$ --- простая замкнутая кривая постоянной геодезической кривизны $1$ на гладкой замкнутой поверхности $\Sigma$ с положительной гауссовой кривизной.
Покажите, что $\length\gamma<2\cdot\pi$;
то есть, что $\gamma$ короче единичной окружности.  
\end{thm}

\begin{thm}{Упражнение}\label{ex:GB-hat}
Предположим, что диск $\Delta$ лежит на графике $z=f(x,y)$ гладкой функции,
а его край $\partial\Delta$ лежит в плоскости $(x,y)$,
и $\Delta$ подходит к плоскости под фиксированным углом $\alpha$.
Покажите, что $\iint_\Delta K=2\cdot\pi\cdot(1-\cos\alpha)$.

\end{thm}

\begin{thm}{Упражнение}\label{ex:geodesic-half}
Пусть $\gamma$ --- простая замкнутая геодезическая на гладкой замкнутой поверхности $\Sigma$ с положительной гауссовой кривизной и сферическим отображением $\Norm\:\Sigma\to\mathbb{S}^2$.
Покажите, что кривая $\alpha=\Norm\circ\gamma$ делит сферу на две области равной площади.

Выведите отсюда, что $\length \alpha\ge 2\cdot\pi.$
\end{thm}

\begin{thm}{Упражнение}\label{ex:closed-geodesic}
Пусть $\gamma$ --- замкнутая геодезическая (возможно, с самопересечениями) на гладкой замкнутой поверхности $\Sigma$ с положительной гауссовой кривизной.
Предположим, что $R$ --- одна из областей, которые $\gamma$ вырезает из~$\Sigma$.
Покажите, что $\iint_R K\le 2\cdot\pi$.

Выведите отсюда, что любые две замкнутые геодезические на $\Sigma$ имеют общую точку.
\end{thm}

\begin{thm}{Упражнение}\label{ex:self-intersections}
Пусть $\Sigma$ --- замкнутая гладкая поверхность с положительной гауссовой кривизной. 
Покажите, что замкнутая геодезическая на $\Sigma$ 
ни в какой из карт не может выглядеть как одна из кривых на следующих рисунках.

\begin{figure}[h]
\begin{minipage}{.32\textwidth}
\centering
\includegraphics{mppics/pic-46}
\end{minipage}
\hfill
\begin{minipage}{.32\textwidth}
\centering
\includegraphics{mppics/pic-47}
\end{minipage}
\hfill
\begin{minipage}{.32\textwidth}
\centering
\includegraphics{mppics/pic-471}
\end{minipage}

\medskip

\begin{minipage}{.32\textwidth}
\centering
\caption*{\textit{(лёгкая)}}
\end{minipage}
\hfill
\begin{minipage}{.32\textwidth}
\centering
\caption*{\textit{(сложная)}}
\end{minipage}
\hfill
\begin{minipage}{.32\textwidth}
\centering
\caption*{\textit{(безнадёжная)}}
\end{minipage}
\vskip-5mm
\end{figure}

\end{thm}

Следующее упражнение уточняет оценку в \ref{ex:ruf-bound-mountain}.

\begin{thm}{Упражнение}\label{ex:sqrt(3)}
Предположим, что $f\:\mathbb{R}^2\to\mathbb{R}$ --- гладкая выпуклая $\sqrt{3}$-липшицева функция.
Покажите, что любая геодезическая на графике $z=f(x,y)$ не имеет самопересечений.
\end{thm}

Поверхность $\Sigma$ называется \index{односвязная поверхность}\emph{односвязной}, если любая простая замкнутая кривая на $\Sigma$ ограничивает диск.
Эквивалентно, любую замкнутую кривую на $\Sigma$ можно продеформировать в \index{тривиальная кривая}\emph{тривиальную кривую}; то есть кривую, стоящую всё время в одной точке.

Плоскости и сферы дают примеры односвязных поверхностей;
торы и цилиндры неодносвязны.

\begin{thm}{Упражнение}\label{ex:unique-geod}
Предположим, что $\Sigma$ --- открытая односвязная поверхность с неположительной гауссовой кривизной.

\begin{subthm}{ex:unique-geod:unique}
Покажите, что любые две точки на $\Sigma$ соединены единственной геодезической.
Более того, геодезические на $\Sigma$ не имеют самопересечений.
\end{subthm}

\begin{subthm}{ex:unique-geod:diffeomorphism}
Выведите, что для любой точки $p\in \Sigma$,
экспоненциальное отображение $\exp_p$ есть диффеоморфизм из касательной плоскости $\T_p$ в~$\Sigma$.
В частности, поверхность $\Sigma$ диффеоморфна плоскости.
\end{subthm}

\end{thm}

\section{Аддитивность}

Пусть $\Delta$ --- топологический диск на гладкой ориентированной поверхности $\Sigma$, ограниченный простой кусочно-гладкой кривой $\partial \Delta$.
Как и раньше, мы считаем, что $\partial \Delta$ ориентирована и $\Delta$ лежит от неё слева.
Введём обозначение \index{10gb@$\GB$ (формула Гаусса --- Бонне)}
\[\GB(\Delta)
\df
\tgc{\partial\Delta}+\iint_\Delta K-2\cdot \pi,
\eqlbl{eq:GB}\]
где $K$ --- гауссова кривизна.
Обозначение $\GB$ призвано напоминать, что формулу Гаусса --- Бонне можно записать как
\[\GB(\Delta)=0.\]

{

\begin{wrapfigure}{r}{40 mm}
\vskip-14.5mm
\centering
\includegraphics{mppics/pic-1750}
\vskip-0mm
\end{wrapfigure}

\begin{thm}{Лемма}\label{lem:GB-sum}
Пусть диск $\Delta$ разрезан на два диска $\Delta_1$ и $\Delta_2$ кривой $\delta$.
Тогда
\[
\GB(\Delta)=\GB(\Delta_1)+\GB(\Delta_2).
\]
\end{thm}

}

\parbf{Доказательство.}
Разделим $\partial \Delta$ на две кривые $\gamma_1$ и $\gamma_2$, которые имеют общие концы с $\delta$, так что
$\Delta_i$ ограничен дугами $\gamma_i$ и~$\delta$ при $i=1,2$.

Обозначим через $\phi_1$, $\phi_2$, $\psi_1$ и $\psi_2$ углы между $\delta$ и $\gamma_i$ как на рисунке.
Полагая, что дуги $\gamma_1$, $\gamma_2$ и $\delta$ ориентированы как на рисунке, получаем
\begin{align*}
\tgc{\partial \Delta}&= \tgc{\gamma_1}-\tgc{\gamma_2}+(\pi-\phi_1-\phi_2)+(\pi-\psi_1-\psi_2),
\\
\tgc{\partial \Delta_1}&= \tgc{\gamma_1}-\tgc{\delta}+(\pi-\phi_1)+(\pi-\psi_1),
\\
\tgc{\partial \Delta_2}&= \tgc{\delta}-\tgc{\gamma_2}+(\pi-\phi_2)+(\pi-\psi_2),
\\
\iint_\Delta K&=\iint_{\Delta_1} K+\iint_{\Delta_2} K.
\end{align*}
Остаётся подставить результаты в формулы для $\GB(\Delta)$, $\GB(\Delta_1)$ и $\GB(\Delta_2)$.
\qeds

\enlargethispage{2\baselineskip}

\section{Сферический случай}

Если наша поверхность $\Sigma$ является плоскостью, то её гауссова кривизна равна нулю, и формула Гаусса --- Бонне \ref{eq:g-b} переписывается как 
\[\tgc{\partial\Delta}=2\cdot \pi,\]
а значит, она следует из \ref{prop:total-signed-curvature}.
Другими словами, $\GB(\Delta)=0$ для любого плоского диска~$\Delta$ с кусочно-гладкой границей.

Если $\Sigma$ --- это единичная сфера, то $K\equiv1$;
в этом случае \ref{thm:gb} эквивалентна следующему.

\begin{thm}{Предложение}\label{prop:area-of-spher-polygon}
{\sloppy
Пусть $P$ --- сферический многоугольник, ограниченный простой замкнутой ломаной геодезической $\partial P$.
Предположим, что $\partial P$ ориентирована так, что $P$ находится слева от $\partial P$.
Тогда 
\[\GB(P)=\tgc{\partial P}+\area P-2\cdot \pi=0.\]

}

Более того, та же формула справедлива для любой сферической области, ограниченной кусочно-гладкой простой замкнутой кривой.
\end{thm}

Это предложение будет использовано в следующем разделе.

\parbf{Набросок доказательства.}
Пусть $\Delta$ --- сферический треугольник с углами 
$\alpha$, $\beta$ и~$\gamma$.
Согласно \ref{lem:area-spher-triangle},
\[\area\Delta=\alpha+\beta+\gamma-\pi.\]

Поскольку граница $\partial\Delta$ ориентирована так, что $\Delta$ лежит от неё слева,
её ориентированные внешние углы равны $\pi-\alpha$, $\pi-\beta$ и $\pi-\gamma$.
Следовательно,
\[\tgc{\partial\Delta}=3\cdot\pi-\alpha-\beta-\gamma.\]
Отсюда $\tgc{\partial\Delta}+\area \Delta=2\cdot\pi$, что эквивалентно, $\GB(\Delta)=0$.

Далее, любой сферический многоугольник $P$ можно разбить на треугольники, за несколько шагов, разрезая многоугольник на два по ломаной геодезической на каждом шагу.
По аддитивности (\ref{lem:GB-sum}), получаем, что 
\[\GB(P)=0\]
для любого сферического многоугольника~$P$.

Второе утверждение доказывается через приближения.
Нужно показать, что полная геодезическая кривизна кусочно-гладкой простой кривой
приближается полной геодезической кривизной вписанных в неё ломаных.
Рассуждение схоже с решением \ref{ex:total-curvature=}; мы опускаем подробности.
\qeds

\begin{thm}{Упражнение}\label{ex:half-sphere-total-curvature}
Пусть $\gamma$ --- простая кусочно-гладкая петля на единичной сфере $\mathbb{S}^2$.
Предположим, что $\gamma$ делит $\mathbb{S}^2$ на две области равной площади.
Обозначим через $p$ базовую точку~$\gamma$.
Покажите, что параллельный перенос $\iota_\gamma\:\T_p\mathbb{S}^2\to\T_p\mathbb{S}^2$ является тождественным отображением.
\end{thm}

\section{Наглядное полудоказательство}\label{sec:gb-intuitive-proof}

Следующий частный случай формулы Гаусса --- Бонне является ключевым.
Общий случай можно доказать аналогично, используя ориентированную площадь с учётом кратности, но мы пойдём другим путём, см.~\ref{sec:gauss--bonnet:formal}. 

\parbf{Доказательство \ref{thm:gb} для открытых и замкнутых поверхностей с положительной гауссовой кривизной.}
Пусть $\Norm\:\Sigma\to\mathbb{S}^2$ --- сферическое отображение.
Из \ref{cor:intK},
\[\GB(\Delta)=\tgc{\partial\Delta}+\area(\Norm(\Delta))-2\cdot \pi.
\eqlbl{eq:gb-area}\]

Выберем петлю $\alpha$, которая проходит вдоль $\partial\Delta$ так, чтобы $\Delta$ лежала слева от неё; пусть $p\in \partial\Delta$ --- её базовая точка.
Рассмотрим параллельный перенос $\iota_\alpha\:\T_p\to\T_p$ вдоль $\alpha$.
Согласно \ref{prop:pt+tgc}, $\iota_\alpha$ --- поворот по часовой стрелке на угол $\tgc{\alpha}_\Sigma$.

Пусть $\beta=\Norm\circ\alpha$.
Согласно \ref{obs:parallel=}, $\iota_\alpha=\iota_\beta$, где $\beta$ рассматривается как кривая на единичной сфере.

Далее, $\iota_\beta$ --- поворот по часовой стрелке на угол $\tgc{\beta}_{\mathbb{S}^2}$.
Согласно \ref{prop:area-of-spher-polygon},
\[\GB(\Norm(\Delta))=\tgc{\beta}_{\mathbb{S}^2}+\area(\Norm(\Delta))-2\cdot \pi=0.\]
Значит, 
$\iota_\beta$ --- поворот против часовой стрелки на угол $\area(\Norm(\Delta))$.

То есть, поворот по часовой стрелке на угол $\tgc{\alpha}_\Sigma$ идентичен повороту против часовой стрелки на угол $\area(\Norm(\Delta))$.
Повороты идентичны, если их углы равны по модулю $2\cdot\pi$.
Следовательно, 
\[
\begin{aligned}
\GB(\Delta)&=\tgc{\partial\Delta}_\Sigma+\area(\Norm(\Delta))-2\cdot \pi=
2\cdot n \cdot \pi
\end{aligned}
\eqlbl{eq:sum=2pin}\]
для некоторого целого числа~$n$.

Остаётся показать, что $n=0$.
Согласно \ref{prop:total-signed-curvature}, это верно для топологического диска на плоскости. 
В общем случае, диск $\Delta$ можно рассматривать как результат непрерывной деформации плоского диска. 
Целое число $n$ не может измениться в процессе деформации, ибо левая часть в \ref{eq:sum=2pin} меняется непрерывно,
и раз $n=0$ в начале, то $n=0$ и в конце деформации.
\qeds

\enlargethispage{2\baselineskip}

\section{Простая геодезическая}

Следующая теорема даёт интересное приложение формулы Гаусса --- Бонне; она доказана Стефаном Кон-Фоссеном \cite[Satz 9]{convossen}.

\begin{thm}{Теорема}\label{thm:cohn-vossen}
Любая открытая гладкая поверхность с положительной гауссовой кривизной имеет простую двусторонне бесконечную геодезическую.
\end{thm}

\parbf{Доказательство.}
Пусть $\Sigma$ --- открытая поверхность с положительной гауссовой кривизной.
Выберем двусторонне бесконечную геодезическую $\gamma$ на~$\Sigma$.

Если у кривой $\gamma$ есть самопересечения, то она содержит простую петлю;
то есть для некоторого интервала $[a,b]$,
сужение $\ell=\gamma|_{[a,b]}$ представляет собой простую петлю.

Согласно \ref{ex:convex-proper-plane}, $\Sigma$ параметризуется открытой выпуклой областью $\Omega$ на плоскости.
По теореме Жордана (\ref{thm:jordan}), $\ell$ ограничивает топологический диск на $\Sigma$; обозначим его через~$\Delta$.
Если $\phi$ --- внутренний угол в базовой точке петли, то по формуле Гаусса --- Бонне
\[\iint_\Delta K=\pi+\phi.\] 

Напомним, что 
\[\iint_\Sigma K\le 2\cdot\pi,
\eqlbl{intK=<2pi+}\]
см. \ref{ex:intK:2pi}.
Следовательно, $0<\phi<\pi$; то есть $\gamma$ не имеет вогнутых простых петель.

Допустим, что у $\gamma$ есть две простые петли, $\ell_1$ и $\ell_2$;
они ограничивают диски скажем $\Delta_1$ и $\Delta_2$.
Тогда диски $\Delta_1$ и $\Delta_2$ должны пересекаться;
иначе кривизна $\Sigma$ превысила бы $2\cdot\pi$, что невозможно по \ref{intK=<2pi+}.

Значит, после выхода из $\Delta_1$ геодезической $\gamma$ придётся войти в него снова, прежде чем образовать новую петлю.
\begin{figure}[h!]
\vskip-0mm
\centering
\includegraphics{mppics/pic-1550}
\end{figure}
Рассмотрим момент, когда $\gamma$ снова входит в $\Delta_1$;
на рисунке показаны два возможных сценария.
На левом рисунке мы получаем два неперекрывающихся диска, что, как мы знаем, невозможно.
Правый рисунок также невозможен --- в этом случае мы получаем вогнутую простую петлю.

\enlargethispage{2\baselineskip}

Следовательно, $\gamma$ содержит только одну простую петлю.
Эта петля вырезает из $\Sigma$ диск и обходит его справа или слева.
Так, все самопересекающиеся геодезические 
делятся на два типа: {}\emph{правые} и {}\emph{левые}.

Если геодезическая $t\mapsto \gamma(t)$ правая, то обращение параметризации $t\mapsto \gamma(-t)$ делает её левой.
Выпустим геодезическую в каждом направлении из точки $p=\gamma(0)$.
Это даёт однопараметрическое семейство геодезических $\gamma_s$ для $s\in[0,\pi]$, соединяющее геодезическую $t\mapsto \gamma(t)$ с $t\mapsto \gamma(-t)$; то есть $\gamma_0(t)\z=\gamma(t)$, а $\gamma_\pi(t)=\gamma(-t)$.

Геодезическая $\gamma_s$ является правой (или левой) при $s$ из открытого множества в $[0,\pi]$.
То есть, если $\gamma_s$ правая, то правые и все $\gamma_t$ при $t$, близких к~$s$.%
\footnote{Неформально говоря, это означает, что самопересечение не может исчезнуть вдруг. Попробуйте в этом убедиться.}

Поскольку интервал $[0,\pi]$ связен, он не разбивается на два открытых множества.
Значит, при каком-то $s$ геодезическая $\gamma_s$ не является ни правой, ни левой;
то есть $\gamma_s$ не имеет самопересечений.
\qeds

{

\begin{wrapfigure}{r}{17 mm}
\vskip-0mm
\centering
\includegraphics{mppics/pic-1575}
\end{wrapfigure}

\begin{thm}{Упражнение}\label{ex:cohn-vossen}
Пусть $\Sigma$ --- открытая гладкая поверхность с положительной гауссовой кривизной,
и $\alpha\:[0,1]\z\to \Sigma$ --- такая гладкая петля, что $\alpha'(0)\z=-\alpha'(1)$.
Покажите, что найдётся простая двухсторонне бесконечная геодезическая, которая касается $\alpha$.
\end{thm}

}

{

\begin{wrapfigure}[3]{r}{43 mm}
\vskip-6mm
\centering
\includegraphics{mppics/pic-1577}
\end{wrapfigure}

\begin{thm}{Упражнение}\label{ex:3-curves}
Какие из кривых на рисунке могут задавать геодезические на картах открытых гладких поверхностей положительной гауссовой кривизны и почему?
\end{thm}

}

\enlargethispage{2\baselineskip}

\section{Области общего вида}
\index{формула Гаусса --- Бонне}

Следующее обобщение формулы Гаусса --- Бонне получено Вальтером фон Диком \cite{dyck}.

\begin{thm}{Теорема}\label{thm:GB-generalized}
Пусть $\Lambda$ --- компактная область на гладкой поверхности.
Предположим, что $\Lambda$ ограничена конечным (возможно, пустым) набором замкнутых кусочно-гладких кривых $\gamma_1,\dots,\gamma_n$, и каждая $\gamma_i$ ориентирована так, что $\Lambda$ лежит слева от неё.
Тогда
\[\iint_\Lambda K=2\cdot \pi\cdot \chi-\tgc{\gamma_1}-\dots-\tgc{\gamma_n}\eqlbl{eq:g-b++}\]
для целого числа $\chi=\chi(\Lambda)$.

Более того, если граф с $v$ вершинами и $e$ рёбрами разбивает $\Lambda$ на $f$ дисков и содержит все $\gamma_i$, то $\chi=v-e+f$.
\end{thm}

\begin{wrapfigure}{r}{29 mm}
\vskip-8mm
\centering
\includegraphics{mppics/pic-1580}
\end{wrapfigure}

Число $\chi=\chi(\Lambda)$ называется \index{эйлерова характеристика}\emph{эйлеровой характеристикой} области $\Lambda$. 
Она не зависит от выбора разбиения, ведь остальные члены в формуле \ref{eq:g-b++} от него не зависят.
Формула \ref{eq:g-b++} выводится из стандартной формулы Гаусса --- Бонне (\ref{thm:gb}).
Геометрия та же, что в \ref{lem:GB-sum}, но комбинаторика сложней.

Прежде чем перейти к доказательству, попробуйте вывести формулу для разбиения кольца $A$ на рисунке следуя рассуждениям в \ref{lem:GB-sum}.
У этого графа $4$ вершины и $6$ рёбер (одно из них --- петля), и он разбивает $A$ на два диска $\Delta_1$ и $\Delta_2$.
Таким образом, $\chi(A)=4-6+2=0$, и
\[\iint_A K=-\tgc{\gamma_1}-\tgc{\gamma_2}.\]

\parbf{Доказательство.}
Пусть граф с $v$ вершинами и $e$ рёбрами разбивает область $\Lambda$ на $f$ дисков $\Delta_1,\dots,\Delta_f$.
Применим формулу Гаусса --- Бонне к каждому диску и сложим результаты:
\[
\begin{aligned}
\iint_\Lambda K&=\iint_{\Delta_1} K+\dots+\iint_{\Delta_f} K=
2\cdot f\cdot \pi-\tgc{\partial\Delta_1}-\dots-\tgc{\partial\Delta_f}.
\end{aligned}
\]
Остаётся показать, что  
\[\tgc{\gamma_1}+\dots+\tgc{\gamma_n}-\tgc{\partial\Delta_1}-\dots-\tgc{\partial\Delta_f}
=
2\cdot\pi\cdot(v-e).
\eqlbl{eq:GB-sum}\]
Для этого мы вычислим левую часть, суммируя отдельно вклад каждого ребра и каждой вершины.

Пусть $\sigma$ --- ребро графа.
Если $\sigma$ не является частью какой-либо $\gamma_i$,
то оно встречается дважды на краях дисков, скажем, на $\partial \Delta_i$ и $\partial \Delta_j$.
В этом случае можно считать, что $\Delta_i$ лежит слева от $\sigma$, а $\Delta_j$ --- справа, поэтому 
$\sigma$ добавляет $\tgc\sigma$ в $\tgc{\partial\Delta_i}$ и $-\tgc\sigma$ в $\tgc{\partial\Delta_j}$; следовательно, $\sigma$ ничего не добавляет в левую часть формулы \ref{eq:GB-sum}.
Может случиться, что $i=j$ (как для одного из рёбер на рисунке выше);
в этом случае у нас один и тот же диск с обеих сторон от~$\sigma$, но всё равно оно ничего не добавляет в \ref{eq:GB-sum}.
Если же ребро $\sigma$ идёт по краю, то оно встречается один раз на краю какого-то диска, скажем, на $\partial \Delta_j$.
Можно считать, что и $\Lambda$, и $\Delta_j$ лежат слева от $\sigma$,
поэтому оно вносит $\tgc\sigma$ как в $\tgc{\gamma_i}$, так и в $\tgc{\partial\Delta_j}$, а 
эти вклады компенсируют друг друга в левой части формулы \ref{eq:GB-sum}.

Итак, в левой части формулы \ref{eq:GB-sum} вклады геодезических кривизн рёбер компенсируют друг друга.

\enlargethispage{2\baselineskip}

Теперь займёмся внешними углами при вершинах.
Выберем вершину $p$, пусть $d$ --- её \index{степень вершины}\emph{степень},
то есть число рёбер исходящих из~$p$.

\begin{wrapfigure}{r}{23 mm}
\vskip-3mm
\centering
\includegraphics{mppics/pic-1585}
\end{wrapfigure}

Предположим, что $p$ лежит во внутренней части~$\Lambda$.
Пусть $\delta_1,\dots,\delta_d$ --- внутренние углы при $p$ дисков, прилегающих к $p$,
и $\phi_{i}=\pi-\delta_{i}$ --- соответствующие внешние углы.
Тогда вклад $p$ в сумму равен 
$-\phi_1-\dots-\phi_d$.
Поскольку $\delta_1+\dots+\delta_d=2\cdot\pi$, вклад $p$ в левую часть формулы \ref{eq:GB-sum} равен
\[-\phi_1-\dots-\phi_d = (\delta_1+\dots+\delta_d) - d\cdot \pi=(2-d)\cdot \pi.\]

\begin{wrapfigure}{r}{23 mm}
\vskip-0mm
\centering
\includegraphics{mppics/pic-1590}
\end{wrapfigure}

Если же $p$ лежит на краю $\Lambda$, то она является вершиной $(d-1)$-го внутреннего угла
$\delta_1,\dots,\delta_{d-1}$,
и $\phi_{i}\z=\pi-\delta_{i}$ --- соответствующие внешние углы.
Заметим, что
\[\delta_1+\dots+\delta_{d-1}\z=\pi-\theta,\]
где $\theta\in(-\pi,\pi)$ --- внешний угол области $\Lambda$ при $p$.
А значит, вклад $p$ в левую часть формулы \ref{eq:GB-sum} опять равен
\[\theta-\sum\phi_{i}=(2-d)\cdot \pi.\]

Итак, если $p_1,\dots,p_v$ --- вершины графа, а $d_1,\dots,d_v$ --- их степени,
то общий вклад от внешних углов в левую часть формулы \ref{eq:GB-sum} равен
\[2\cdot v\cdot \pi-(d_1+\dots+d_v)\cdot\pi.
\eqlbl{eq:GB-sum-d}\]
Так как рёбра не внесли ничего, левая часть формулы \ref{eq:GB-sum} равна \ref{eq:GB-sum-d}.

Остаётся проверить, что $d_1+\dots+d_v=2\cdot e$.
И действительно, $d_1+\z\dots+d_v$ --- это число концов всех рёбер в графе, а их $2\cdot e$,
ибо у каждого ребра по два конца.
Отсюда получаем \ref{eq:GB-sum} и \ref{eq:g-b++}.
\qeds

\begin{thm}{Упражнение}\label{ex:g-b-chi}
{\sloppy
Найдите интеграл гауссовой кривизны по каждой из следующих поверхностей:

}

\setlength{\columnseprule}{0.4pt}
\begin{multicols}{2}

\begin{subthm}{ex:g-b-chi:torus}
Тор.
\end{subthm}

\begin{Figure}
\vskip-0mm
\centering
\includegraphics{mppics/pic-1595}
\end{Figure}

\begin{subthm}{ex:g-b-chi:moebius}
Лента Мёбиуса с геодезическим краем.
\end{subthm}

\begin{Figure}
\vskip-0mm
\centering
\includegraphics{mppics/pic-1605}
\end{Figure}

\begin{subthm}{ex:g-b-chi:pair-of-pants}
Пара штанов с геодезическими компонентами края.
\end{subthm}
\begin{Figure}
\vskip-0mm
\centering
\includegraphics{mppics/pic-1600}
\end{Figure}

\begin{subthm}{ex:g-b-chi:two-handles}
Сфера с двумя ручками.
\end{subthm}

\begin{Figure}
\vskip-0mm
\centering
\includegraphics{mppics/pic-1610}
\end{Figure}

\end{multicols}

\begin{subthm}{ex:g-b-chi:cylinder}
Цилиндр, окрестности краёв которого лежат в плоскостях.
\begin{Figure}
\vskip-0mm
\centering
\includegraphics{mppics/pic-1620}
\end{Figure}
\end{subthm}

\end{thm}

\chapter{Полугеодезические координаты}
\label{chap:semigeodesic}

Это вычислительная глава, в ней выводятся несколько утверждений, обсуждавшихся выше, включая
эквивалентное определение радиуса инъективности (\ref{prop:inj-rad}),
то, что кратчайшие являются геодезическими (\ref{prop:gamma''})
и формулу Гаусса --- Бонне (\ref{thm:gb}).
Кроме того, мы обсудим внутренние изометрии между поверхностями и докажем 
замечательную теорему Гаусса о том, что \textit{гауссова кривизна является внутренним инвариантом}.

\section{Полярные координаты}\label{sec:Polar coordinates}

{\sloppy

Свойство экспоненциального отображения, описанное в \ref{prop:exp}, можно использовать для задания \index{полярные координаты}\emph{полярных координат} на гладкой поверхности.

}

А именно, пусть $p$ --- точка гладкой поверхности $\Sigma$ и
$(r,\theta) \z\in \mathbb{R}_{\ge0} \times \mathbb{S}^1$ --- полярные координаты на касательной плоскости $\T_p$.
Если $\vec v\in \T_p$ имеет координаты $(r,\theta)$,
то будем говорить, что $s(r,\theta)\z=\exp_p\vec v$ --- это точка на $\Sigma$ с полярными координатами $(r,\theta)$.

Из точки $p$ в данную точку $x$ может идти много геодезических, а может их вовсе не быть.
Поэтому одна точка может иметь несколько представлений в полярных координатах, а может и не иметь их вовсе.
Однако из~\ref{prop:exp} получаем следующее.

\begin{thm}{Наблюдение}\label{obs:polar}
Пусть $s\:(r,\theta)\mapsto s(r,\theta)$ задаёт полярные координаты на гладкой поверхности~$\Sigma$ с началом в точке $p$.
Тогда существует такое $r_0>0$, что $s$ является регулярным для любой пары $(r,\theta)$ при $0<r<r_0$.

Более того, если $0\le r_1,r_2<r_0$, то $s(r_1,\theta_1) \z= s(r_2,\theta_2)$ тогда и только тогда, когда
$r_1=r_2=0$ или $r_1=r_2$ и $\theta_1\z=\theta_2+2\cdot n\cdot\pi$ для некоторого целого числа~$n$.
\end{thm}

Далее мы всегда будем предполагать, что полярные координаты на поверхности определены только при $r<r_0$,
а значит, они ведут себя обычным образом.

Следующее утверждение сыграет основную роль в строгом доказательстве того, что кратчайшие являются геодезическими, см.~\ref{sec:proof-of-gamma''}.

\begin{thm}{Лемма Гаусса}\label{lem:palar-perp}
Пусть $(r,\theta)\mapsto s(r,\theta)$ --- полярные координаты на гладкой поверхности с началом в точке $p$.
Тогда
$s_\theta\perp s_r$
при любых $r$ и~$\theta$.
\end{thm}

\parbf{Доказательство.}
Выберем $\theta \in \mathbb{S}^1$.
По определению экспоненциального отображения, кривая $\gamma(t)\z=s(t,\theta)$ --- геодезическая с единичной скоростью, исходящая из $p$.
\begin{enumerate}[(i)]
\item Поскольку у $\gamma$ единичная скорость, $|s_r|=|\gamma'|=1$, и в частности,
 \[
 \tfrac{\partial}{\partial \theta}
 \langle s_r,s_r\rangle=0.\]
\item Поскольку $\gamma$ --- геодезическая, $s_{rr}(r,\theta)=\gamma''(r)\perp\T_{\gamma(r)}$,
и следовательно, 
\[
\langle s_\theta, s_{rr}\rangle=0.\]
\end{enumerate}
Отсюда вытекает, что
\[
\begin{aligned}
\tfrac{\partial}{\partial r}
\langle s_\theta, s_r\rangle
&=
\langle s_{\theta r},s_r\rangle
+
\cancel{\langle s_\theta,s_{rr}\rangle}=
\\
&=
\tfrac12
\cdot 
\tfrac{\partial}{\partial \theta}
\langle s_r, s_r\rangle=
\\
&=0.
\end{aligned}
\eqlbl{eq:<s',s'>'=0}
\]

Далее, $s_\theta(0,\theta)=0$, ведь $s(0,\theta)=p$ для любого $\theta$.
В частности,
$\langle s_\theta, s_r\rangle=0$
при $r=0$.
Согласно \ref{eq:<s',s'>'=0}, значение 
$\langle  s_\theta, s_r\rangle$ не зависит от~$r$ при фиксированном~$\theta$.
А значит,
\[\langle s_\theta, s_r\rangle=0\]
при любых $r$ и $\theta$.
\qeds

\section{Снова кратчайшие и геодезические}
\label{sec:proof-of-gamma''}

В этом разделе мы воспользуемся полярными координатами и леммой Гаусса (\ref{lem:palar-perp}) для доказательства предложения~\ref{prop:gamma''}.

\parbf{Доказательство \ref{prop:gamma''}.}
Пусть $\gamma\:[0,\ell]\to\Sigma$ --- кратчайшая, параметризованная длиной, и $p=\gamma(0)$.
Предположим, что длина $\ell=\length\gamma$ достаточно мала, так что $\gamma$ описывается в полярных координатах с началом в $p$;
скажем $\gamma(t)=s(r(t),\theta(t))$ для функций $t\mapsto \theta(t)$ и $t\mapsto r(t)$, где $r(0)=0$.

По правилу дифференцирования сложной функции,
\[\gamma'= s_\theta\cdot \theta'+ s_r\cdot r'
\eqlbl{eq:chain(gamma)}\]
если левая часть определена и $t>0$.
По лемме Гаусса \ref{lem:palar-perp}, $s_\theta\perp s_r$, и по определению полярных координат, $|s_r|=1$.
Следовательно, из \ref{eq:chain(gamma)} следует
\[|\gamma'(t)|\ge r'(t).\eqlbl{eq:|gamma'|=r'}\]
для любого $t>0$, где $\gamma'(t)$ определена.

Поскольку $\gamma$ параметризована длиной, 
\[|\gamma(t_2)-\gamma(t_1)|\le |t_2-t_1|.\]
В частности, $\gamma$ липшицева, и, по теореме Радемахера (\ref{thm:rademacher}), производная $\gamma'$ определена почти везде.
Из \ref{adex:integral-length:a} получаем
\begin{align*}
\length\gamma&=\int_0^\ell|\gamma'(t)|\cdot dt\ge
\\
&\ge\int_0^\ell r'(t)\cdot dt=
\\
&=r(\ell).
\end{align*}

По построению полярных координат, найдётся геодезическая длины $r(\ell)$ из $p=\gamma(0)$ в $q=\gamma(\ell)$.
Поскольку $\gamma$ --- кратчайшая, $r(\ell)=\ell$, и, более того, $r(t)=t$ для любого~$t$.
Это равенство выполняется, тогда и только тогда, когда в \ref{eq:|gamma'|=r'} достигается равенство при почти всех~$t$.
Последнее означает, что $\gamma$ --- геодезическая.

Остаётся проверить верность почти обратного.

Зафиксируем точку $p\in\Sigma$.
Пусть $\epsilon>0$, как в \ref{prop:exp}.
Предположим, что короткая геодезическая $\gamma$ (короче чем $\epsilon$) от $p$ до $q$ не минимизирует длину между своими концами.
Тогда существует кратчайшая от $p$ до $q$ отличная от $\gamma$.
Если пропараметризовать её длиной, то получим другую геодезическую.
То есть, существуют две геодезические от $p$ до $q$ длины меньше $\epsilon$.
Иными словами, найдутся такие два вектора ${\vec v},\vec w\in\T_p$, что $|{\vec v}|<\epsilon$, $|\vec w|<\epsilon$ и 
$q=\exp_p\vec v\z=\exp_p\vec w$.
Однако, по \ref{prop:exp}, экспоненциальное отображение $\T_p \to \Sigma$  инъективно в $\epsilon$-окрестности нуля --- противоречие.
\qeds

\section{Гауссова кривизна}\label{sec:jacobi-formula}

Пусть $s$ --- гладкое отображение из (возможно, неограниченного) координатного прямоугольника на плоскости $(u,v)$ в гладкую поверхность~$\Sigma$.
Отображение $s$ называется \index{полугеодезические}\emph{полугеодезическим}, если для любого фиксированного $v$ отображение $u\mapsto s(u,v)$ является геодезической, параметризованной длиной, и $s_u\perp s_v$ при любых $(u,v)$.

По лемме Гаусса (\ref{lem:palar-perp}), полярные координаты на $\Sigma$ задаются полугедезическим отображением.

Пусть $\Norm=\Norm(u,v)$ --- нормаль к $\Sigma$ в точке $s(u,v)$.
Векторы $\Norm$, $\vec u\z=s_u$ и $\vec v\z=\Norm\times \vec u$ образуют ортонормированный базис для каждой пары $(u,v)$.
Напомним, что $s_v\perp \vec u$ и $s_v\perp \Norm$, ведь вектор $s_v(u,v)$ касателен к $\Sigma$ в $s(u,v)$. 
Следовательно, $s_v=b\cdot\vec v$ для некоторой гладкой функции $(u,v)\z\mapsto b(u,v)$.%
\footnote{Для фиксированного значения $v_0$ векторное поле $s_v=b\cdot\vec v$ описывает разницу между $\gamma_0$ и \textit{инфинитезимально близкой} геодезической $\gamma_1\:u\z\mapsto s(u,v_1)$.
Поля с таким свойством называются \index{поле Якоби}\emph{полями Якоби} вдоль $\gamma_0$.}

\begin{thm}{Предложение}\label{prop:jaccobi}
Пусть $(u,v)\mapsto s(u,v)$ --- полугеодезическое отображение на гладкую поверхность $\Sigma$, для которого $\Norm$, $\vec u$, $\vec v$, и $b$ определены выше.
Тогда 
\[b\cdot K+b_{uu}=0,\]
где $K=K(u,v)$ --- гауссова кривизна $\Sigma$ в точке $s(u,v)$.

Более того, 
\[
\langle\vec u_u,\vec u\rangle=
\langle\vec u_u,\vec v\rangle=
\langle\vec u_v,\vec u\rangle=0,
\quad\text{и}\quad
\langle\vec u_v,\vec v\rangle=b_u.
\]

\end{thm}

Доказательство проводится длинным, но несложным вычислением.

\parbf{Доказательство.}
Пусть $\ell=\ell(u,v)$, $m=m(u,v)$ и $n=n(u,v)$ --- компоненты матрицы, описывающей оператор формы в базисе $\vec u, \vec v$;
то есть
\[
\begin{aligned}
\Shape\vec u&=\ell\cdot \vec u+ m\cdot \vec v,
&
\Shape\vec v&=m\cdot \vec u+ n\cdot \vec v.
\end{aligned}
\eqlbl{eq:Shape(u,v)}
\]
Напомним, что (см.~\ref{sec:More curvatures})
\[K=\ell\cdot n-m^2.\]

Сначала давайте выведем предложение из следующие тождеств
\[
\begin{aligned}
\vec u_u&=\ell\cdot \Norm,
&
\vec u_v&=\phantom{-}b_u\cdot \vec v+b\cdot m\cdot\Norm,
\\
\vec v_u&=m\cdot \Norm,
&
\vec v_v&=-b_u\cdot \vec u+b\cdot n\cdot\Norm.
\end{aligned}
\eqlbl{eq:uu-vv}
\]
Действительно, 
\begin{align*}
b\cdot K&=b\cdot (\ell\cdot n-m^2)=
\\
&=\langle\vec u_u,\vec v_v\rangle-\langle\vec u_v,\vec v_u\rangle=
\tag{по \ref{eq:uu-vv}}
\\
&= 
\left(
\tfrac{\partial}{\partial v}
\langle\vec u_u,\vec v\rangle
-
\cancel{\langle\vec u_{uv},\vec v\rangle}
\right)-
\left(
\tfrac{\partial}{\partial u}
\langle\vec u_v,\vec v\rangle
-
\cancel{\langle\vec u_{uv},\vec v\rangle}
\right)=
\\
&=0-b_{uu}.
\tag{по \ref{eq:uu-vv}}
\end{align*}
Остальные тождества в предложении прямо следуют из \ref{eq:uu-vv}.
Остаётся доказать четыре тождества в \ref{eq:uu-vv}.

\parit{Вывод $\vec u_u=\ell\cdot \Norm$.}
Поскольку базис $\Norm$, $\vec u$ и $\vec v$ ортонормирован, это векторное тождество можно переписать в виде  трёх скалярных
\[
\begin{aligned}
\langle\vec u_u,\vec u\rangle&=0,
&
\langle\vec u_u,\vec v\rangle&=0,
&
\langle\vec u_u,\Norm\rangle&=\ell.
\end{aligned}
\]
Поскольку $u\mapsto s(u,v)$ --- геодезическая, $\vec u_u=s_{uu}(u,v)\perp\T_{s(u,v)}$.
Отсюда следуют первые два тождества.
Далее,
\begin{align*}
\langle\vec u_u,\Norm\rangle
&=\langle s_{uu},\Norm\rangle=
\tag{по \ref{thm:shape-chart}}
\\
&=    \langle \Shape s_u,s_u\rangle=
\tag{поскольку $\vec u=s_u$}
\\
&=\langle \Shape \vec u,\vec u\rangle=
\tag{по \ref{eq:Shape(u,v)}}
\\
&=\ell.
\end{align*}

\parit{Вывод $\vec u_v= b_u\cdot \vec v+b\cdot m\cdot\Norm$.}
Перепишем его в виде скалярных тождеств:
\[
\begin{aligned}
\langle\vec u_v,\vec u\rangle&=0,
&
\langle\vec u_v,\vec v\rangle&=b_u,
&
\langle\vec u_v,\Norm\rangle&=b\cdot m.
\end{aligned}
\eqlbl{eq:uu-vv:2}
\]

Поскольку $\langle\vec u,\vec u\rangle=1$, имеем
$0=\tfrac{\partial}{\partial v}\langle\vec u,\vec u\rangle=2\cdot\langle\vec u_v,\vec u\rangle$.
Отсюда вытекает первое тождество в \ref{eq:uu-vv:2}.
Далее,
\begin{align*}
\langle\vec u_v,\vec v\rangle&=\langle s_{vu},\vec v\rangle=
\tag{$s_v=b\cdot \vec v$}
\\
&=\langle \tfrac{\partial}{\partial u} (b\cdot\vec v),\vec v\rangle =
\\
&=b_u\cdot \langle \vec v,\vec v\rangle+b\cdot \langle \vec v_u,\vec v\rangle=
\tag*{($\langle \vec v,\vec v\rangle=1$ и $\qquad\qquad\qquad$} 
\\
&=b_u. 
\tag*{$0=\tfrac{\partial}{\partial u}\langle \vec v,\vec v\rangle=2\cdot\langle \vec v_u,\vec v\rangle$)}
\end{align*}
Получили второе тождество в \ref{eq:uu-vv:2}.
Наконец,
\begin{align*}
\langle\vec u_v,\Norm\rangle
&=\langle s_{uv},\Norm\rangle=
\tag{по \ref{thm:shape-chart}}
\\
&=\langle \Shape s_u,s_v\rangle=
\tag{$\vec u=s_u$ и $s_v=b\cdot \vec v$}
\\
&=\langle \Shape \vec u,b\cdot \vec v\rangle=
\tag{по \ref{eq:Shape(u,v)}}
\\
&=b\cdot m.
\end{align*}

\parit{Вывод $\vec v_u=m\cdot \Norm$ и $\vec v_v=-b_u\cdot \vec u+b\cdot n\cdot\Norm$.}
Напомним, что $\vec v=\Norm\times \vec u$.
Следовательно,
\[
\vec v_u=\Norm_u\times \vec u+\Norm\times \vec u_u,
\qquad\qquad
\vec v_v=\Norm_v\times \vec u+\Norm\times \vec u_v.
\eqlbl{eq:uu-vv:3+4}
\]
Выражения для $\vec u_u$ и $\vec u_v$ в \ref{eq:uu-vv} уже доказаны.
Далее,
\begin{align*}
-\Norm_u&=\Shape s_u=
&
-\Norm_v&=\Shape s_v=
\\
&=\Shape \vec u=
&
&=b\cdot\Shape \vec v=
\\
&=\ell\cdot\vec u+m\cdot\vec v,
&
&=b\cdot(m\cdot \vec u+ n\cdot \vec v),
\end{align*}
Остаётся подставить в \ref{eq:uu-vv:3+4} выражения для $\vec u_u$, $\vec u_v$, $\Norm_u$ и~$\Norm_v$. 
\qeds

Карта $(u,v)\mapsto s(u,v)$ и её локальные координаты называются \index{полугеодезические}\emph{полугеодезическими}, если отображение $(u,v)\z\mapsto s(u,v)$ является полугеодезическим.
Заметим, что функция $b=b(u,v)$ в полугеодезических координатах имеет постоянный знак.
Поэтому, обратив знак $\Norm$, можно (и должно) считать, что $b>0$;
иными словами, $b=|s_v|$.

\begin{thm}{Упражнение}\label{ex:semigeodesc-chart}
Покажите, что любую точку $p$ на гладкой поверхности $\Sigma$ можно покрыть полугеодезической картой.
\end{thm}

\begin{thm}{Упражнение}\label{ex:inj-rad}
{\sloppy
Пусть $p$ --- точка гладкой поверхности~$\Sigma$.
Предположим, что отображение $\exp_p$ инъективно в шаре $B\z=B(0,r_0)_{\T_p}$.
Пусть полугеодезическое отображение $(r,\theta)\mapsto s(r,\theta)$ задаёт полярные координаты с началом в $p$, а функция $(r,\theta)\mapsto b(r,\theta)$ та же, что выше.

}

Докажите следующие утверждения:

\begin{subthm}{ex:inj-rad:sign}
$b(r,\theta)$ не меняет знак при $0\z\le r\z<r_0$.
\end{subthm}

\begin{subthm}{ex:inj-rad:0}
$b(r,\theta)\ne0$, если $0< r<r_0$.
\end{subthm}

\begin{subthm}{ex:inj-rad:prop:inj-rad}
Примените \ref{SHORT.ex:inj-rad:sign} и \ref{SHORT.ex:inj-rad:0}, чтобы доказать \ref{prop:inj-rad}.
\end{subthm}
 
\end{thm}

Карта $(u,v)\mapsto s(u,v)$ называется \index{ортогональная карта}\emph{ортогональной}, если $s_u\perp s_v$ для любого $(u,v)$.
Например, любая полугеодезическая карта ортогональна.

Решение следующего упражнения похоже на \ref{prop:jaccobi}.

\begin{thm}{Упражнение}\label{lem:K(orthogonal)}
Пусть $(u,v)\mapsto s(u,v)$ --- ортогональная карта гладкой поверхности~$\Sigma$.
Обозначим через $K=K(u,v)$ гауссову кривизну $\Sigma$ в точке $s(u,v)$.
Определим 
\begin{align*}
a=a(u,v)&\df|s_u|,&
b=b(u,v)&\df|s_v|,\\
\vec u=\vec u(u,v)&\df\tfrac{s_u}a,&
\vec v=\vec v(u,v)&\df\tfrac{s_v}b.
\end{align*}
Пусть $\Norm=\Norm(u,v)$ --- единичный нормальный вектор в точке $s(u,v)$.

\begin{subthm}{lem:K(orthogonal):uu-vv}
Покажите, что 
\begin{align*}
\vec u_u
&=
-\tfrac1{b}\cdot a_v
\cdot
\vec v 
+
a\cdot \ell\cdot \Norm
,
&
\vec v_u
&=
\tfrac1{b}\cdot a_v
\cdot \vec u
+
a\cdot m\cdot \Norm
\\
\vec u_v
&=
\tfrac1{a}\cdot b_u\cdot\vec v
+
b\cdot m\cdot \Norm
,
&
\vec v_v
&=
-\tfrac1{a}\cdot b_u\cdot\vec u
+
b\cdot n\cdot \Norm,
\end{align*}
где $\ell=\ell(u,v)$, $m=m(u,v)$ и $n=n(u,v)$ --- компоненты матрицы, описывающей оператор формы в базисе $\vec u, \vec v$.
\end{subthm}

\begin{subthm}{lem:K(orthogonal):K}
Покажите, что
\[K=-\frac1{a\cdot b}\cdot
\left(
\frac{\partial}{\partial u}
\left(\frac{b_u}a \right)
+
\frac{\partial}{\partial v}
\left(\frac{a_v}b\right)
\right).\]
\end{subthm}
\end{thm}

\begin{thm}{Упражнение}\label{ex:conformal}
Предположим, что $(u,v)\mapsto s(u,v)$ --- это \index{конформная карта}\emph{конформная карта};
то есть существует функция $(u,v)\mapsto b(u,v)$ такая, что $b=|s_u|=|s_v|$ и $s_u\perp s_v$ для любого $(u,v)$.
(Функция $b$ называется {}\emph{конформным множителем} карты.)

Используя \ref{lem:K(orthogonal)}, покажите, что  
\[K=-\frac{\triangle (\ln b)}{b^2},\]
где $\triangle$ --- \index{лапласиан}\emph{лапласиан}; то есть $\triangle=\tfrac{\partial^2}{\partial u^2}+\tfrac{\partial^2}{\partial v^2}$, и 
 $K=K(u,v)$ --- гауссова кривизна $\Sigma$ в точке $s(u,v)$.
\end{thm}

Полезно знать, что \textit{любую точку гладкой поверхности можно покрыть конформной картой},
соответствующие локальные координаты называются \index{изотермические координаты}\emph{изотермическими}.

\section{Вращение векторного поля}

Пусть $\gamma\:[0,1]\to \Sigma$ --- простая петля на гладкой ориентированной поверхности~$\Sigma$,
и $\vec u$ --- поле единичных касательных векторов на $\Sigma$, определённое в окрестности~$\gamma$.
Обозначим через $\vec v$ поворот $\vec u$ против часовой стрелки на угол $\tfrac{\pi}2$ в касательной плоскости в каждой точке; поле $\vec v$ можно также определить как $\vec v\df\Norm\times\vec u$, где $\Norm$ --- поле нормалей к $\Sigma$.
Тогда \index{вращение}\emph{вращение} $\vec u$ вдоль $\gamma$ определяется как интеграл
\[\rot_\gamma\vec u
\df
\int_0^1\langle\vec u'(t),\vec v(t)\rangle\cdot dt.\]

\begin{thm}{Лемма}\label{lem:rotation-parallel}
Пусть $\gamma\:[0,1]\to \Sigma$ --- простая петля с базовой точкой $p$ на гладкой ориентированной поверхности~$\Sigma$, а $\vec u$ --- поле касательных единичных векторов к $\Sigma$, определённое в окрестности~$\gamma$.
Тогда параллельный перенос $\iota_\gamma\:\T_p\to\T_p$ есть поворот по часовой стрелке на угол $\rot_\gamma\vec u$.

В частности, вращения различных векторных полей вдоль $\gamma$ могут различаться только на числа кратные $2\cdot\pi$.
\end{thm}

\parbf{Доказательство.}
Как и ранее, пусть $\vec v=\Norm\times\vec u$. 
Будем обозначать через $\vec u(t)$ и $\vec v(t)$ векторы полей $\vec u$ и $\vec v$ в точке $\gamma(t)$.

Пусть $t\mapsto \vec x(t)\in \T_{\gamma(t)}$ --- это параллельное векторное поле вдоль $\gamma$ с  $\vec x(0)\z=\vec u(0)$, и $\vec y\df\Norm\times\vec x$.

Заметим, что существует непрерывная функция $t\mapsto \phi(t)$ такая, что 
$\vec u(t)$ --- это поворот $\vec x(t)$ против часовой стрелки на угол $\phi(t)$.
Так как $\vec x(0)=\vec u(0)$, можно считать, что $\phi(0)=0$.
Тогда
\begin{align*}
\vec u&=\cos\phi\cdot \vec x+\sin\phi\cdot \vec y,
\\
\vec v&=-\sin\phi\cdot \vec x+\cos\phi\cdot \vec y.
\end{align*}
Отсюда
\begin{align*}
\langle\vec u',\vec v\rangle
=\phi'\cdot\biggl(&(\sin \phi)^2\cdot \langle\vec x,\vec x\rangle+(\cos \phi)^2\cdot \langle\vec y,\vec y\rangle
\biggr)=
\phi'.
\end{align*}

Следовательно,
\begin{align*}
\rot_\gamma\vec u&=\int_0^1\langle\vec u'(t),\vec v(t)\rangle\cdot dt=
\\
&=\int_0^1\phi'(t)\cdot dt=
\\
&=\phi(1).
\end{align*}

Заметим, что 
\begin{itemize}
\item $\iota_\gamma(\vec x(0))=\vec x(1)$,

\item  $\vec x (0) = \vec u (0) = \vec u (1),$ 

\item $\vec u(1)$ --- поворот $\vec x(1)$ против часовой стрелки на угол $\phi(1)\z=\rot_\gamma\vec u$,

\end{itemize}
Отсюда вытекает, что $\vec x(1)$ --- это \textit{поворот по часовой стрелке} $\vec x(0)$ на угол $\rot_\gamma\vec u$, и результат следует.

Последнее утверждение следует из \ref{prop:pt+tgc}.
\qeds

Следующая лемма пригодится при доказательстве формулы Гаусса --- Бонне.

\begin{thm}{Лемма}\label{lem:rotation-semigeoesic}
Пусть $(u,v)\mapsto s(u,v)$ --- полугеодезическая карта на гладкой поверхности~$\Sigma$.
Предположим, что простая замкнутая кривая $\gamma$ ограничивает диск $\Delta$, который полностью покрывается картой~$s$.
Тогда 
\[\rot_\gamma\vec u+\iint_\Delta K=0,\]
где $\vec u=s_u$, а $K$ обозначает гауссову кривизну поверхности~$\Sigma$.
\end{thm}

Вычисления ниже используют так называемую \index{формула Грина}\emph{формулу Грина}, формулирующейся следующим образом.

\textit{Пусть $D$ --- компактная область на $(u,v)$-плоскости, ограниченная кусочно-гладким простым замкнутым путём $\alpha\:t\mapsto (u(t),v(t))$.
Предположим, что $\alpha$ ориентирована так, что $D$ находится слева от неё.
Тогда}
\[\iint_D (Q_u- P_v)\cdot du\cdot dv=\int_\alpha (P\cdot du+Q\cdot dv)\df \int_0^1 (P\cdot u'+Q\cdot v')\cdot dt\]
\textit{для любых двух гладких функций $P$ и $Q$, определённых на $D$.}

Формулы Грина и Гаусса --- Бонне похожи, обе связывают интеграл по области с интегралом по его краю, и неудивительно, что одна поможет доказать другую.

\parbf{Доказательство.}
Пусть $\vec u$, $\vec v$ и $b$ определены, как в~\ref{sec:jacobi-formula}.
Запишем $\gamma$ в $(u,v)$-координатах: $\gamma(t)=s(u(t),v(t))$.
Тогда
\begin{align*}
\rot_\gamma \vec u&=\int_0^1\langle\vec u',\vec v\rangle\cdot dt=
\tag{поскольку $\vec u'=\vec u_u\cdot u'+\vec u_v\cdot v'$}
\\
&=\int_0^1[\langle\vec u_u,\vec v\rangle\cdot u'+\langle\vec u_v,\vec v\rangle\cdot v']\cdot dt=\tag{по \ref{prop:jaccobi}}
\\
&=\int_0^1b_u\cdot v'\cdot dt=\int_{s^{-1}\circ\gamma}b_u\cdot dv=
\tag{по формуле Грина}
\\
&=\iint_{s^{-1}(\Delta)}b_{uu}\cdot du\cdot dv=
\tag{$s_u\perp s_v \Longrightarrow \jac s=|s_u|\cdot|s_v|=b$}
\\
&=\iint_\Delta\frac{b_{uu}}{b}=
\tag{по \ref{prop:jaccobi}}
\\
&=-\iint_{\Delta}K.
\end{align*}
\qedsf

\section{Вывод формулы Гаусса --- Бонне}\label{sec:gauss--bonnet:formal}

{\sloppy

Напомним, что формула Гаусса --- Бонне записывается как $\GB(\Delta)\z=0$, где 
\[\GB(\Delta)
\df
\tgc{\partial\Delta}+\iint_\Delta K-2\cdot \pi,\]
$\Delta$ --- топологический диск на гладкой ориентированной поверхности,
ограниченный кусочно-гладкой кривой $\partial \Delta$ ориентированной так, что $\Delta$ лежит от неё слева.

}

\parbf{Вывод формулы Гаусса --- Бонне (\ref{thm:gb}).}
Предположим, что $\Delta$ покрывается полугеодезической картой.
Из \ref{prop:pt+tgc}, \ref{lem:rotation-parallel} и \ref{lem:rotation-semigeoesic},
\[\GB(\Delta)
=
2\cdot n\cdot \pi,
\eqlbl{eq:gb(n)}\]
где $n=n(\Delta)$ --- целое число.

Согласно \ref{ex:semigeodesc-chart}, полугеодезической картой можно покрыть любую точку поверхности.
Таким образом, применив аддитивность $\GB$ (\ref{lem:GB-sum}) конечное число раз, получим \ref{eq:gb(n)} для любого диска $\Delta$ на~$\Sigma$.
Точнее, мы можем разрезать $\Delta$ гладкой кривой, которая проходит от края до края,
и повторить такое разбиение рекурсивно для полученных дисков;
см. рисунок.
\begin{figure}[!ht]
\vskip-0mm
\centering
\includegraphics{mppics/pic-1700}
\vskip-0mm
\end{figure}
После нескольких таких шагов каждый малый диск покроется полугеодезической картой.
В частности, \ref{eq:gb(n)} выполняется для каждого малого диска.
Затем, применив \ref{lem:GB-sum} несколько раз, получим \ref{eq:gb(n)} для исходного диска.

Остаётся показать, что $n=0$.
Предположим, что $\Delta$ лежит в локальной реализации нашей поверхности графиком $z = f(x,y)$.
Рассмотрим однопараметрическое семейство графов $z = t \cdot f(x,y)$;
обозначим через $\Delta_t$ соответствующий диск на $\Sigma_t$, так что $\Delta_1 = \Delta$, а $\Delta_0$ --- его проекция на плоскость $(x,y)$.
Функция $h \: t \mapsto \GB (\Delta_t )$ непрерывна.
Согласно \ref{eq:gb(n)}, $h(t)$ --- целое кратное $2 \cdot \pi$ для любого $t$.
Значит $h$ постоянна.
Следовательно,
\[\GB(\Delta)=\GB(\Delta_0)=0;\]
последнее равенство следует из \ref{prop:total-signed-curvature}.

Итак, мы доказали, что 
\[\GB (\Delta)=0
\eqlbl{eq:GB=0}\]
если $\Delta$ лежит на графике $z = f(x,y)$ в некоторой системе координат $(x,y,z)$.
Любая точка поверхности~$\Sigma$ имеет окрестность, которая покрывается таким графиком, и, применив аддитивность $\GB$ (\ref{lem:GB-sum}), как выше, получаем, что \ref{eq:GB=0} выполняется для любого диска~$\Delta$ на~$\Sigma$.
\qeds

\section{Сравнение Рауха}

Следующее предложение является частным случаем так называемой \index{теорема сравнения Рауха}\emph{теоремы сравнения Рауха}.

\begin{thm}{Предложение}\label{prop:rauch}
Пусть $p$ --- точка гладкой поверхности $\Sigma$, и $r\le \inj(p)$.
Далее, пусть $\tilde\gamma$ --- кривая в $r$-окрестности нуля касательной плоскости $\T_p$,
и $\gamma$ --- кривая на~$\Sigma$, определяемая как
\[\gamma=\exp_p\circ\tilde\gamma,
\quad
\text{что эквивалентно}
\quad
\log_p\circ\gamma=\tilde\gamma.\]

\begin{subthm}{prop:rauch:K=<0}
Если поверхность $\Sigma$ имеет неположительную гауссову кривизну, то логарифмическое отображение $\log_p$ не увеличивает длину в $r$-окрестности точки $p$ на~$\Sigma$;
то есть
\[\length \gamma\ge \length \tilde\gamma\]
для любой кривой $\gamma$ в открытом шаре $B(p,r)_{\Sigma}$.
\end{subthm}

{\sloppy

\begin{subthm}{prop:rauch:K>=0}
Если поверхность $\Sigma$ имеет неотрицательную гауссову кривизну, то экспоненциальное отображение $\exp_p$ не увеличивает длину в $r$-окрестности нуля в $\T_p$;
то есть
\[\length \gamma\le \length \tilde\gamma\]
для любой кривой $\tilde\gamma$ в открытом шаре $B(0,r)_{\T_p}$.
\end{subthm}

}

\end{thm}

\parbf{Доказательство.}
Пусть $(r(t),\theta(t))$ --- полярные координаты $\tilde\gamma(t)$.
Тогда $(r(t),\theta(t))$ --- полярные координаты $\gamma(t)$ с началом в $p$ на~$\Sigma$;
то есть $\gamma(t)\z=s(r(t),\theta(t))$, где $s$ как в \ref{sec:Polar coordinates}.

Докажем, что $b(0,\theta)=0$ и $b_r(0,\theta)=1$ для любого $\theta$.
Можно предположить, что $\theta=0$.
Выберем стандартный базис $\vec v,\vec w$ в $\T_p$, с $\vec v$ направленным по кривой $t\mapsto s(t,0)$.
Заметим, что%
\footnote{Напомним, что $(D_{\vec w}\exp_p)(r\cdot \vec v)\df h'(0)$, где $h(t)=\exp_p(r\cdot \vec v+t\cdot \vec w)$; см.~\ref{sec:dirder}.}
\[b(r,0)=r\cdot |(D_{\vec w}\exp_p)(r\cdot \vec v)|.\]
В частности, $b(0,0)=0$.
Согласно \ref{obs:d(exp)=1}, $|(D_{\vec w}\exp_p)(0)|=1$,
и, значит, $|(D_{\vec w}\exp_p)(r\cdot \vec v)|\z\ne 0$ при малых $r$.
Отсюда следует, что функция $r\mapsto|(D_{\vec w}\exp_p)(r\cdot \vec v)|$ дифференцируема при $r=0$.
Взяв частную производную от выражения для $b(r,0)$, получаем $b_r(0,0)=1$.

Пусть $b(r,\theta)\df|s_\theta|$.
Согласно \ref{prop:jaccobi},
\[b_{rr}=-K\cdot b.\]
Если $K\ge 0$, то  для фиксированного $\theta$ функция $r\mapsto b(r,\theta)$ вогнута,
а если $K\le 0$, то $r\mapsto b(r,\theta)$ выпукла.
Поскольку $b(0,\theta)=0$ и $b_r(0,\theta)=1$,
\[
\begin{aligned}
b(r,\theta)\ge r\quad\text{если}\quad K&\le 0;
\\
b(r,\theta)\le r\quad\text{если}\quad K&\ge 0.
\end{aligned}
\eqlbl{eq:b-K}
\]

Можно считать, что $\tilde\gamma\:[a,b]\to \T_p$ параметризована длиной;
в частности, она липшицева.
Заметим, что
\begin{align*}
\length\tilde\gamma&=\int_a^b\sqrt{r'(t)^2+r(t)^2\cdot\theta'(t)^2}\cdot dt.
\shortintertext{Применив \ref{lem:palar-perp}, получим}
\length\gamma&=\int_a^b\sqrt{r'(t)^2+b(r(t),\theta(t))^2\cdot\theta'(t)^2}\cdot dt.
\end{align*}
А теперь утверждения в \ref{SHORT.prop:rauch:K=<0} и \ref{SHORT.prop:rauch:K>=0} следуют из \ref{eq:b-K}.
\qeds

\section{Внутренние изометрии}

Пусть $\Sigma$ и $\Sigma^{*}$ --- две гладкие поверхности.
Будем говорить, что отображение $f\:\Sigma\to \Sigma^{*}$ \index{сохраняющее длину}\emph{сохраняет длины}, если для любой кривой $\gamma$ на $\Sigma$ кривая $\gamma^{*}=f\circ\gamma$ на $\Sigma^{*}$ имеет ту же длину. 
Если $f$ ещё и гладкая биекция, то она будет называться \index{внутренняя изометрия}\emph{внутренней изометрией}.

Пример сохраняющего длину отображения можно получить, сворачивая плоскость в цилиндр при помощи отображения $s\:\mathbb{R}^2\to\mathbb{R}^3$, заданного как $s(x,y)=(\cos x,\sin x,y)$.

\begin{thm}{Упражнение}\label{ex:isom(geod)}
Покажите, что любая внутренняя изометрия между гладкими поверхностями переводит геодезические в геодезические.
\end{thm}

\begin{thm}{Упражнение}\label{ex:K=0}
Пусть гауссова кривизна гладкой поверхности $\Sigma$ равна нулю.
Покажите, что $\Sigma$ является \index{локально плоская поверхность}\emph{локально плоской} поверхностью;
то есть некоторая окрестность любой точки на $\Sigma$ допускает внутреннюю изометрию на открытое множество евклидовой плоскости.  
\end{thm}

\begin{thm}{Упражнение}\label{ex:K=1}
Пусть гауссова кривизна гладкой поверхности $\Sigma$ равна $1$ в каждой точке.
Покажите, что малая окрестность любой точки на $\Sigma$ допускает внутреннюю изометрию на открытое подмножество единичной сферы.
\end{thm}

\begin{thm}{Упражнение}\label{ex:deformation}
Для данного $a>0$, покажите, что существует гладкая кривая с единичной скоростью 
$\gamma(t)=(x(t),y(t))$, такая что $y(t) = a\cdot \cos t$ и $y>0$.
Найдите её интервал определения.

Пусть $\Sigma_a$ --- поверхность вращения кривой $\gamma$ вокруг оси $x$.
\begin{figure}[h!]
\vskip-0mm
\centering
\begin{lpic}[t(-0mm),b(6mm),r(0mm),l(0mm)]{asy/deformation(1.2)}
\lbl[t]{8,-.5;$a=2$}
\lbl[t]{24,3;$a=\sqrt{2}$}
\lbl[t]{41,4;$a=1$}
\lbl[t]{57,7;$a=\tfrac1{\sqrt{2}}$}
\lbl[t]{73,8;$a=\tfrac12$}
\end{lpic}
\vskip-6mm
\end{figure}
Покажите, что поверхность $\Sigma_a$ имеет единичную гауссову кривизну в каждой точке.

Воспользуйтесь этим построением и \ref{ex:K=1}, чтобы построить гладкую сохраняющую длину деформацию малого диска $\Delta$ на $\mathbb{S}^2$;
то есть однопараметрическое семейство $\Delta_t$ поверхностей с краем, для которого $\Delta_0=\Delta$ и $\Delta_t$ не конгруэнтна $\Delta_0$ для $t\ne0$.%
\footnote{На самом деле, любой диск на $\mathbb{S}^2$ допускает гладкую сохраняющую длину деформацию.
Однако если диск больше полусферы, то доказательство требует дополнительных усилий;
оно выводится из двух результатов Александра Александрова: теоремы о склеивании и теоремы о существовании выпуклой поверхности с абстрактно заданной метрикой \cite[с. 44]{pogorelov}.
}
\end{thm}

Следующее упражнение иллюстрирует заключительный шаг в доказательстве того, что \textit{любая открытая поверхность с нулевой гауссовой кривизной является цилиндрической поверхностью}.
См. обсуждение после \ref{ex:flat-plane}.

\begin{thm}{Продвинутое упражнение}\label{ex:line-cylinder} 
Пусть $(u,v)\mapsto f(u,v)$ задаёт внутреннюю изометрию из плоскости со стандартными $(u,v)$-координатами на поверхность $\Sigma$ в $\mathbb{R}^3$.
Предположим, что $f$ изометрически отображает $v$-ось на $z$-ось.
Покажите, что $\Sigma$ является цилиндрической поверхностью;
точнее, $\Sigma$ есть объединение семейства прямых, параллельных $z$-оси.
\end{thm}

\section{Замечательная теорема}

\begin{thm}{Теорема}\label{thm:remarkable}
Предположим, что $f\:\Sigma\to \Sigma^{*}$ является внутренней изометрией между двумя гладкими поверхностями; $p\in \Sigma$ и $p^{*}\z=f(p)\in \Sigma^{*}$.
Тогда 
\[K(p)_{\Sigma}=K(p^{*})_{\Sigma^{*}};\]
то есть в точке $p$ гауссова кривизна у $\Sigma$ та же, что в точке $p^{*}$ у~$\Sigma^{*}$.
\end{thm}

Напомним, что гауссова кривизна определяется как произведение главных кривизн, которые могут быть различными в точках $p$ и $p^*$; однако, согласно теореме, их произведения одинаковы.
Другими словами, гауссова кривизна является \textit{внутренним инвариантом}.
Эта теорема была доказана Карлом Фридрихом Гауссом \cite{gauss}, и справедливо названа {}\emph{замечательной} ({}\emph{Theorema Egregium}).

На самом деле, кривизну $K(p)$ можно получить из \textit{внутренних} измерений.
Например, она появляется в следующей формуле для длины $c(r)$ геодезической окружности с центром в точке $p$ на поверхности:
\[c(r)=2\cdot\pi\cdot r-\tfrac\pi3\cdot K(p)\cdot r^3+o(r^3).\]

Из теоремы следует, например, что не существует гладкого сохраняющего длину отображения, которое отправляет открытую область на единичной сфере в плоскость.%
\footnote{Гладкость существенна --- существует множество негладких отображений из сферы в плоскость сохраняющих длины; см. \cite{petrunin-yashinski} и ссылки в там.}
Это следует из того, что гауссова кривизна плоскости равна нулю, а гауссова кривизна единичной сферы равна $1$.
В частности, любая географическая карта обязана иметь искажения.

\parbf{Доказательство.}
Выберем карту $(u,v)\mapsto s(u,v)$ на $\Sigma$ и пусть
$s^{*}\z=f\circ s$.
Заметим, что $s^{*}$ --- карта на $\Sigma^{*}$, и 
\begin{align*}
\langle s_u,s_u\rangle
&=
\langle s_u^{*}, s_u^{*}\rangle,
&
\langle s_u, s_v\rangle
&=
\langle s_u^{*}, s_v^{*}\rangle,
&
\langle s_v, s_v\rangle
&=
\langle s_v^{*}, s_v^{*}\rangle
\end{align*}
в любой точке $(u,v)$.
Действительно, поскольку $f$ сохраняет длины координатных линий $\gamma\:t\mapsto s(t,v)$ и  $\gamma\:t\z\mapsto s(u,t)$, мы получаем первое и третье равенства.
Теперь, поскольку $f$ сохраняет длины кривых $\gamma\:t\z\mapsto s(t,c-t)$ для любой константы~$c$, первое и третье равенства влекут второе.

Из \ref{prop:gamma''}, если $s$ --- полугеодезическая карта, то и $s^{*}$ таковая.
Остаётся применить \ref{prop:jaccobi} и \ref{ex:semigeodesc-chart}.
\qeds

\chapter{Теоремы сравнения}
\label{chap:comparison}

Теоремы сравнения --- мощный инструмент, дающий возможность применять евклидову интуицию в дифференциальной геометрии.

\section{Треугольники и створки}

Напомним, что $[x,y]$ и $[x,y]_\Sigma$ обозначают кратчайшую между точками $x$ и $y$ на поверхности $\Sigma$, а
$\dist{x}{y}\Sigma$ --- \index{внутренняя метрика}\emph{внутреннее расстояние} от $x$ до $y$ на~$\Sigma$.

\index{геодезический треугольник}\emph{Геодезический треугольник} на поверхности $\Sigma$ определяется как тройка точек $x,y,z\z\in \Sigma$ с выбранными кратчайшими $[x,y]_\Sigma$, $[y,z]_\Sigma$ и $[z,x]_\Sigma$.
Точки $x,y,z$ называются {}\emph{вершинами} треугольника,
а кратчайшие $[x,y]_\Sigma$, $[y,z]_\Sigma$ и $[z,x]_\Sigma$ --- его {}\emph{сторонами};
сам треугольник обозначается как $[xyz]$ или $[xyz]_\Sigma$;
как всегда, последнее обозначение используется, если хочется подчеркнуть, что треугольник лежит на поверхности~$\Sigma$.\index{10aad@$[xyz]$, $[xyz]_\Sigma$ (геодезический треугольник)}

Треугольник $[\tilde x\tilde y\tilde z]$ на плоскости $\mathbb{R}^2$ называется {}\emph{модельным треугольником} треугольника $[xyz]$,
если их соответствующие стороны равны;
то есть
\[\dist{\tilde x}{\tilde y}{\mathbb{R}^2}=\dist{x}{y}\Sigma,
\quad
\dist{\tilde y}{\tilde z}{\mathbb{R}^2}=\dist{y}{z}\Sigma,
\quad
\dist{\tilde z}{\tilde x}{\mathbb{R}^2}=\dist{z}{x}\Sigma.
\]
В этом случае пишут $[\tilde x\tilde y\tilde z]=\tilde\triangle xyz$.
\index{10aae@$\tilde\triangle$ (модельный треугольник)}

Пара кратчайших $[x,y]$ и $[x,z]$, исходящих из одной точки $x$, называется \index{створка}\emph{створкой} и обозначается как $\hinge xyz$.\index{10aac@$\hinge yxz$ (створка)}
Угол между этими геодезическими в точке $x$ обозначается как $\measuredangle\hinge xyz$.
Величина угла $\measuredangle\hinge {\tilde x}{\tilde y}{\tilde z}$ модельного треугольника $[\tilde x\tilde y\tilde z]=\tilde\triangle xyz$ обозначается через $\modangle xyz$,
\index{10aab@$\modangle yxz$ (модельный угол)}
и называется \index{модельный треугольник и угол}\emph{модельным углом} треугольника $[xyz]$ при~$x$.

Согласно признаку равенства треугольников по трём сторонам,
модельный треугольник $[\tilde x\tilde y\tilde z]$ определён с точностью до конгруэнтности.
Поэтому модельный угол $\tilde\theta=\modangle xyz$ однозначно определён.
По теореме косинусов, 
\[\cos \tilde\theta=\frac{a^2+b^2-c^2}{2\cdot a \cdot b},\]
где $a=\dist{x}{y}{\Sigma}$, $b=\dist{x}{z}{\Sigma}$, и $c=\dist{y}{z}{\Sigma}$.

\begin{thm}{Упражнение}\label{ex:wide-hinges}
Пусть $[x_ny_nz_n]$ --- последовательность треугольников на гладкой поверхности~$\Sigma$ 
со сторонами $a_n=\dist{x_n}{y_n}{\Sigma}$,
$b_n=\dist{x_n}{z_n}{\Sigma}$,
$c_n=\dist{y_n}{z_n}{\Sigma}$, и пусть $\tilde\theta_n=\modangle {x_n}{y_n}{z_n}$.
Предположим, что последовательности $a_n$ и $b_n$ отделена от нуля;
то есть $a_n>\epsilon$ и $b_n>\epsilon$ для фиксированного $\epsilon>0$ и любого~$n$.
Покажите, что
\[(a_n+b_n-c_n)\to 0\qquad\iff\qquad \tilde\theta_n\to \pi\]
при $n\to\infty$.
\end{thm}

\section{Формулировки}

Часть \ref{SHORT.thm:comp:cat} следующей теоремы называется {}\emph{теоремой Картана--Адамара};
она доказана 
Гансом фон Мангольдтом \cite{mangoldt} и обобщена 
Эли Картаном \cite{cartan} и
Жаком Адамаром \cite{hadamard}.
Часть \ref{SHORT.thm:comp:toponogov} называется {}\emph{теоремой сравнения Топоногова} и иногда {}\emph{теоремой сравнения Александрова};
она доказана Паоло Пиццетти \cite{pizzetti}, переоткрыта Александром Александровым \cite{aleksandrov}, и 
обобщена Виктором Топоноговым \cite{toponogov1957}.

Напомним, что поверхность называется \index{односвязная поверхность}\emph{односвязной}, если любая простая замкнутая кривая на ней ограничивает диск.

\begin{thm}{Теоремы сравнения}
\label{thm:comp}
\index{теорема сравнения}
\ 

\begin{subthm}{thm:comp:cat}
Пусть $\Sigma$ --- открытая односвязная гладкая поверхность с неположительной гауссовой кривизной.
Тогда 
\[\measuredangle\hinge {x}{y}{z}\le\modangle xyz\]
для любого геодезического треугольника $[xyz]$ на $\Sigma$.
\end{subthm}

\begin{subthm}{thm:comp:toponogov}
Пусть $\Sigma$ --- замкнутая (или открытая) гладкая поверхность с неотрицательной гауссовой кривизной.
Тогда 
 \[\measuredangle\hinge {x}{y}{z}\ge\modangle xyz\]
для любого геодезического треугольника $[xyz]$ на $\Sigma$.
\end{subthm}

\end{thm}

Доказательства частей \ref{SHORT.thm:comp:cat} и \ref{SHORT.thm:comp:toponogov} даны в разделах~\ref{sec:nonpos-comp} и~\ref{sec:nonneg-comp} соответственно.

\begin{thm}{Упражнение}\label{ex:thm:comp:cat:nsc}
Покажите, что неравенство в \ref{thm:comp:cat} не выполняется для гиперболоида $\set{(x,y,z)\in\mathbb{R}^3}{x^2+y^2-z^2=1}$.
В частности, \ref{thm:comp:cat} неверно без предположения об односвязности.
\end{thm}

Давайте сравним формулу Гаусса --- Бонне с теоремами сравнения.
Рассмотрим диск $\Delta$ ограниченный геодезическим треугольником $[xyz]$ с внутренними углами $\alpha$, $\beta$ и~$\gamma$.
По формуле Гаусса --- Бонне, 
\[\alpha+\beta+\gamma-\pi=\iint_\Delta K;\]
в частности, у обеих сторон равенства тот же знак:
\begin{itemize}
\item если $K_\Sigma\ge 0$, то $\alpha+\beta+\gamma\ge\pi$, и
\item если $K_\Sigma\le 0$, то $\alpha+\beta+\gamma\le\pi$.
\end{itemize}

\begin{wrapfigure}{r}{35mm}
\centering
\vskip-10mm
\includegraphics{mppics/pic-2307}
\end{wrapfigure}

Теперь пусть 
$\hat\alpha\z=\measuredangle\hinge {x}{y}{z}$,
$\hat\beta\z=\measuredangle\hinge {y}{z}{x}$,
и $\hat\gamma\z=\measuredangle\hinge {z}{x}{y}$.
Заметим, что $\hat\alpha,\hat\beta,\hat\gamma\in[0,\pi]$.
Так как сумма углов любого плоского треугольника равна $\pi$, то из теорем сравнения вытекает, что
\begin{itemize}
\item если $K_\Sigma\ge 0$, то $\hat\alpha+\hat\beta+\hat\gamma\ge\pi$, и
\item если $K_\Sigma\le 0$, то $\hat\alpha+\hat\beta+\hat\gamma\le\pi$.
\end{itemize}

Получается, что формула Гаусса --- Бонне и теоремы сравнения тесно связаны,
однако эта связь не столь прямолинейна.

Например, допустим, что $K\ge 0$.
Тогда формула Гаусса --- Бонне не запрещает внутренним углам $\alpha$, $\beta$ и $\gamma$ одновременно быть близкими к $2\cdot\pi$.
Однако $\hat \alpha=\alpha$, если $\alpha\le \pi$, а если нет, то $\hat \alpha=\pi-\alpha$;
то есть
\begin{align*}
\hat \alpha&=\min\{\,\alpha,2\cdot\pi-\alpha\,\},
&
\hat\beta &=\min\{\,\beta,2\cdot\pi-\beta\,\},
&
\hat\gamma&=\min\{\,\gamma,2\cdot\pi-\gamma\,\}.
\end{align*}
Значит, если $\alpha$, $\beta$ и $\gamma$ близки к $2\cdot\pi$, то $\hat\alpha$, $\hat\beta$ и $\hat\gamma$ будут близки к $0$,
ну а последнее невозможно по теореме сравнения.

\begin{thm}{Упражнение}\label{ex:diam-angle}
Пусть точки $p$ и $q$ максимизируют расстояние на гладкой замкнутой выпуклой поверхности $\Sigma$;
то есть $\dist{p}{q}\Sigma\z\ge\dist{x}{y}\Sigma$ при любых $x,y\in \Sigma$.
Покажите, что $\measuredangle\hinge xpq\ge \tfrac\pi3$
для любой точки $x\z\in \Sigma\setminus\{p,q\}$.
\end{thm}

\begin{thm}{Упражнение}\label{ex:sum=<2pi}
Пусть $\Sigma$ --- замкнутая (или открытая) поверхность неотрицательной гауссовой кривизной.
Покажите, что 
\[\modangle pxy+\modangle pyz+\modangle pzx\le2\cdot \pi\]
для любых четырёх различных точек $p,x,y,z$ на~$\Sigma$.
\end{thm}

\section{Локальные теоремы}\label{sec:loc-comp}

Следующая теорема --- первый шаг в доказательстве теорем сравнения (\ref{thm:comp});
она будет выведена из теоремы сравнения Рауха (\ref{prop:rauch}).

\begin{thm}{Теорема}\label{thm:loc-comp}
Теоремы сравнения (\ref{thm:comp}) выполняются в малой окрестности любой точки.

Более того, если $\Sigma$ --- гладкая поверхность без края,
то для любой точки $p\in \Sigma$ существует такое $r>0$, что если $\dist{p}{x}\Sigma<r$, то $\inj(x)_\Sigma>r$, и справедливы следующие утверждения:

\begin{subthm}{thm:loc-comp:cba}
Если $\Sigma$ имеет неположительную гауссову кривизну, то
\[\measuredangle\hinge {x}{y}{z}\le\modangle xyz\]
для любого геодезического треугольника $[xyz]$ в $B(p,\tfrac r4)_\Sigma$.
\end{subthm}

\begin{subthm}{thm:loc-comp:cbb}
Если $\Sigma$ имеет неотрицательную гауссову кривизну, то 
\[\measuredangle\hinge {x}{y}{z}\ge\modangle xyz\]
для любого геодезического треугольника $[xyz]$ в $B(p,\tfrac r4)_\Sigma$.
\end{subthm}

\end{thm}

\parbf{Доказательство.}
Существование $r>0$ вытекает из \ref{prop:exp}.
Пусть $[xyz]$ --- геодезический треугольник в $B(p,\tfrac{r}4)$.

Поскольку $r<\inj(x)_\Sigma$, найдутся такие векторы $\vec v,\vec w\in\T_x$, что 
\begin{align*}
y&=\exp_x\vec v,
& 
z&=\exp_x\vec w,
\\
|\vec v|_{\T_x}&=\dist{x}{y}\Sigma,
&
|\vec w|_{\T_x}&=\dist{x}{z}\Sigma,
&
\measuredangle\hinge 0{\vec v}{\vec w}_{\T_x}&=\measuredangle\hinge xyz_\Sigma;
\end{align*}
в частности, $|\vec v|, |\vec w|< \tfrac r2$.

\parit{\ref{SHORT.thm:loc-comp:cba}.}
Пусть $\gamma$ --- минимизирующая геодезическая из $y$ в $z$.
Поскольку $\dist{x}{y}{\Sigma},\dist{x}{z}{\Sigma}\z<\tfrac r2$, по неравенству треугольника, $\gamma$ лежит в $r$-окрестности~$x$.
В частности, кривая
$\tilde \gamma\df\log_x\circ\gamma$ определена и лежит в $r$-окрестности нуля в $\T_x$.
Заметим, что $\tilde\gamma$ соединяет $\vec v$ и $\vec w$ в~$\T_x$.

По теореме сравнения Рауха (\ref{prop:rauch:K=<0}),
\[\length \tilde \gamma\le \length\gamma.\]
Так как $\dist{\vec v}{\vec w}{\T_x}\le\length\tilde \gamma$ и $\length\gamma=\dist{y}{z}{\Sigma}$, получаем 
\[\dist{\vec v}{\vec w}{\T_x}\le \dist{y}{z}\Sigma.\]
По монотонности угла (\ref{lem:angle-monotonicity}), 
\[\measuredangle\hinge 0{\vec v}{\vec w}_{\T_x}\le \modangle xzy,\]
и результат следует.

\parit{\ref{SHORT.thm:loc-comp:cbb}.}
Рассмотрим отрезок прямой $\tilde \gamma$, соединяющий $\vec v$ и $\vec w$ в касательной плоскости $\T_x$, и пусть $\gamma\df\exp_x\circ\tilde \gamma$.
По теореме Рауха (\ref{prop:rauch:K>=0}), 
\[\length \tilde \gamma\ge\length\gamma.\]
Поскольку $\dist{\vec v}{\vec w}{\T_x}=\length\tilde \gamma$ и $\length\gamma\ge\dist{y}{z}\Sigma$, 
\[\dist{\vec v}{\vec w}{\T_x}\ge \dist{y}{z}\Sigma.\]
По монотонности угла (\ref{lem:angle-monotonicity}), 
\[\measuredangle\hinge 0{\vec v}{\vec w}_{\T_x}\ge\modangle xzy,\]
отсюда результат.
\qeds

\section{Неположительная кривизна}\label{sec:nonpos-comp}

\parbf{Доказательство \ref{thm:comp:cat}.}
Так как $\Sigma$  односвязна, из \ref{ex:unique-geod} получаем, что 
\[\inj(p)_\Sigma=\infty\]
для любой точки $p\in\Sigma$.
Следовательно, \ref{thm:loc-comp:cba} влечёт \ref{thm:comp:cat}.
\qeds

\section{Неотрицательная кривизна}\label{sec:nonneg-comp}

Сейчас мы докажем \ref{thm:comp:toponogov} для компактных поверхностей.
Общий случай требует лишь незначительных изменений; они приводятся в упражнении \ref{ex:open-comparison} в конце раздела.
Доказательство взято из \cite{alexander-kapovitch-petrunin2027}, похожее доказательство независимо нашли Урс Лэнг и Виктор Шрёдер \cite{lang-schroeder}.

\parbf{Доказательство \ref{thm:comp:toponogov} в компактном случае.}\label{proof(thm:comp:toponogov)}
Предположим, что $\Sigma$ компактна.
По \ref{thm:loc-comp} найдётся такое $\epsilon>0$, что неравенство 
\[\measuredangle\hinge {x}{p}{q}\ge\modangle xpq\]
выполняется, если
$\dist{x}{p}\Sigma+\dist{x}{q}\Sigma<\epsilon$.
Следующая лемма утверждает, что в этом случае то же верно и если $\dist{x}{p}\Sigma+\dist{x}{q}\Sigma<\tfrac32\cdot\epsilon$.
Применив эту лемму несколько раз, получаем, сравнение для любой створки, что и доказывает \mbox{\ref{thm:comp:toponogov}}.
\qeds

\begin{thm}{Основная лемма}\label{key-lem:globalization}
Пусть $\Sigma$ --- гладкая открытая или замкнутая поверхность.
Предположим, что на $\Sigma$ выполняется сравнение
\[\measuredangle\hinge x y z
\ge\modangle x y z\eqlbl{eq:key-lem:globalization}\]
если 
$\dist{x}{y}\Sigma+\dist{x}{z}\Sigma
<
\frac{2}{3}\cdot\ell$.
Тогда сравнение \ref{eq:key-lem:globalization} выполняется и в случае, если $\dist{x}{y}\Sigma+\dist{x}{z}\Sigma<\ell$.
\end{thm}

{

\begin{wrapfigure}{r}{35mm}
\centering
\vskip-0mm
\includegraphics{mppics/pic-2308}
\end{wrapfigure}

Для данной створки $\hinge x p q$,
рассмотрим треугольник на плоскости с углом $\measuredangle\hinge x p q$ и двумя прилежащими сторонами $\dist{x}{p}\Sigma$ и $\dist{x}{q}\Sigma$.
Обозначим через $\side \hinge x p q$ третью сторону этого треугольника; она будет называться \index{модельная сторона}\emph{модельной стороной} створки.

Следующее вычислительное упражнение будет использовано в доказательстве леммы.

}

\begin{thm}{Упражнение}\label{ex:s-r}
Предположим, что створки $\hinge xpq$ и $\hinge xpy$ имеют общую сторону $[x,p]$, и $[x,y]\subset [x,q]$.
Покажите, что
\[\frac{\dist{x}{p}{}+\dist{x}{q}{}-\side\hinge xpq}{\dist{x}{q}{}}
\le
\frac{\dist{x}{p}{}+\dist{x}{y}{}-\side\hinge xpy}{\dist{x}{y}{}}.\]
\end{thm}

\parbf{Доказательство.} 
По монотонности угла (\ref{lem:angle-monotonicity})
\[\measuredangle\hinge x p q\ge \modangle x p q\quad\iff\quad\side \hinge x p q
\ge\dist{p}{q}\Sigma.\]
Следовательно, достаточно доказать, что
\[\side \hinge x p q
\ge\dist{p}{q}\Sigma
\eqlbl{eq:thm:=def-loc}\]
при условии, что $\dist{x}{p}\Sigma+\dist{x}{q}\Sigma<\ell$.

Опишем построение новой створки $\hinge{x'}p q$ по данной створке $\hinge x p q$, для которой 
\[\tfrac{2}{3}\cdot\ell \le\dist{p}{x}\Sigma+\dist{x}{q}\Sigma< \ell.\]

\begin{wrapfigure}{r}{32mm}
\vskip-6mm
\centering
\includegraphics{mppics/pic-2310}
\end{wrapfigure}

Будем считать, что $\dist{x}{q}\Sigma\ge\dist{x}{p}\Sigma$; иначе поменяем $p$ и $q$ ролями.
Выберем $x'\in [x, q]$ так, что 
\[\dist{p}{x}\Sigma+3\cdot\dist{x}{x'}\Sigma
=\tfrac{2}{3}\cdot\ell
\eqlbl{3|xx'|}\]
Далее выберем геодезическую $[x', p]$ и рассмотрим створку $\hinge{x'}p q$, образованную $[x',p]$ и $[x',q]\z\subset [x,q]$.

По неравенству треугольника,
\[
\dist{p}{x}\Sigma+\dist{x}{q}\Sigma\ge\dist{p}{x'}\Sigma+\dist{x'}{q}\Sigma.
\eqlbl{eq:thm:=def-loc-fourstar}\]
Покажем, что
\[\side \hinge x p q
\ge\side \hinge{x'}p q
\eqlbl{eq:thm:=def-loc-fivestar}\]

Из \ref{3|xx'|},
\[
\begin{aligned}
\dist{p}{x}{\Sigma}\z+\dist{x}{x'}{\Sigma}&<\tfrac{2}{3}\cdot\ell,
\\
\dist{p}{x'}{\Sigma}\z+\dist{x'}{x}{\Sigma}&<\tfrac{2}{3}\cdot\ell.
\end{aligned}
\]
Значит, по предположению 
\[\begin{aligned}
\measuredangle\hinge x p{x'}
\ge\modangle x p{x'}
\quad\text{и}\quad
\measuredangle\hinge {x'}p x
\ge\modangle {x'}p x.
  \end{aligned}
\eqlbl{eq:thm:=def-loc-threestar}
\]

Рассмотрим модельный треугольник $[\tilde x\tilde x'\tilde p]\z=\modtrig xx'p$. 
Выберем такую точку $\tilde q$ на продолжении отрезка $[\tilde x,\tilde x']$ за $x'$, что $\dist{\tilde x}{\tilde q}\Sigma=\dist{x}{q}\Sigma$, и, значит, $\dist{\tilde x'}{\tilde q}\Sigma\z=\dist{x'}{q}\Sigma$.

Из \ref{eq:thm:=def-loc-threestar}
\[\measuredangle\hinge x p q
=\measuredangle\hinge x p{x'}\ge\modangle x p{x'}.\]
Следовательно,
\[
\side\hinge x q p
\ge
\dist{\tilde p}{\tilde q}{\mathbb{R}^2}.
\]
Далее, так как $\measuredangle\hinge{x'}p x+\measuredangle\hinge{x'}p q= \pi$,
неравенства в \ref{eq:thm:=def-loc-threestar} означают, что
\[
\pi
-\modangle{x'}p x
\ge
\pi-\measuredangle\hinge{x'}p x
\ge
\measuredangle\hinge{x'}p q.
\]
Следовательно,
\[\dist{\tilde p}{\tilde q}{\mathbb{R}^2}\ge\side \hinge{x'}q p \]
и из этого следует \ref{eq:thm:=def-loc-fivestar}.

Пусть $x_0=x$; применяя построение выше, получаем последовательность створок $\hinge{x_n}p q$, где $x_{n+1}=x_n'$.
Согласно \ref{eq:thm:=def-loc-fivestar} и \ref{eq:thm:=def-loc-fourstar}, обе последовательности
\[s_n=\side \hinge{x_n}pq\quad\text{и}\quad r_n=\dist{p}{x_n}\Sigma+\dist{x_n}{q}\Sigma\]
невозрастающие.

Последовательность может остановиться на $x_n$ только если $r_n\z< \tfrac{2}{3}\cdot\ell$.
Тогда, по предположению леммы, 
\[s_n=\side \hinge{x_n}p q\ge\dist{p}{q}\Sigma.\]
А так как последовательность $s_n$ не возрастает, получаем, что
\[\side \hinge{x}p q=s_0\ge s_n\ge\dist{p}{q}\Sigma;\]
откуда следует \ref{eq:thm:=def-loc}.

\begin{figure}[!ht]
\centering
\includegraphics{mppics/pic-2315}
\end{figure}

Осталось доказать \ref{eq:thm:=def-loc} в случае если последовательность $x_n$ бесконечна.
Из \ref{3|xx'|}
\[
\dist{x_n}{x_{n-1}}\Sigma
\ge 
\tfrac1{100}\cdot \ell.
\eqlbl{eq:|x-x|><l}
\]
Согласно \ref{eq:thm:=def-loc-fourstar}, $\dist{x_n}{p}{},\dist{x_n}{q}{}<
\ell$ для любого~$n$.
В случае, если $x_{n+1}\z\in [x_n,q]$, применим \ref{ex:s-r} к створкам $\hinge{x_n}pq$ и $\hinge{x_n}p{x_{n+1}}$.
Согласно \ref{eq:thm:=def-loc-threestar}, $\dist{p}{x_{n+1}}{}\le \side \hinge{x_n}{x_{n+1}}{p}$.
Следовательно,
\[r_n-s_n\le 100\cdot (r_n-r_{n+1})\eqlbl{eq:r-s<100(r-r)}\]
В случае, если $x_{n+1}\in [x_n,p]$, неравенство \ref{eq:r-s<100(r-r)} следует, если применить \ref{ex:s-r} к створкам $\hinge{x_n}pq$ и $\hinge{x_n}{x_{n+1}}q$.

Последовательности $r_n$ и $s_n$ являются невозрастающими и неотрицательными;
поэтому они обязаны сходиться.
В частности, $(r_n\z-r_{n+1})\to0$ при $n\to \infty$.
Следовательно, из \ref{eq:r-s<100(r-r)} следует, что
\[\lim_{n\to\infty}s_n=\lim_{n\to\infty}r_n.\]
По неравенству треугольника, $r_n\ge \dist{p}{q}\Sigma$ для любого~$n$.
Поскольку последовательность $s_n$ невозрастающая, получаем
\[\side \hinge{x}p q=s_0\ge \lim_{n\to\infty}s_n=\lim_{n\to\infty}r_n\ge \dist{p}{q}\Sigma,\]
что завершает доказательство \ref{eq:thm:=def-loc}.
\qeds

\begin{thm}{Упражнение}\label{ex:open-comparison}
Пусть $\Sigma$ --- открытая поверхность с неотрицательной гауссовой кривизной.
Для данного $p\in\Sigma$ обозначим $R_p$ 
(\emph{радиус сравнения} в точке $p$) 
максимальное значение (возможно, $\infty$), при котором сравнение 
\[\measuredangle\hinge x p y
\ge\modangle x p y\]
выполняется для любой створки $\hinge x p y$, если $\dist{p}{x}\Sigma+\dist{x}{y}\Sigma<R_p$.

\begin{subthm}{ex:open-comparison:positive}
Покажите, что для любого компактного подмножества $K\z\subset \Sigma$ существует такое $\epsilon>0$, что $R_p>\epsilon$ для любой точки $p\in K$.
\end{subthm}

\begin{subthm}{ex:open-comparison:almost-min}
Используя часть \ref{SHORT.ex:open-comparison:positive}, покажите, что 
существует точка $p\z\in\Sigma$, такая что 
\[R_q>(1-\tfrac1{100})\cdot R_p,\]
для любой точки $q\in B(p,100\cdot R_p)_\Sigma$.
\end{subthm}

\begin{subthm}{ex:open-comparison:proof}
Объясните как, используя \ref{SHORT.ex:open-comparison:almost-min}, можно распространить доказательство \ref{thm:comp:toponogov} (страница \pageref{proof(thm:comp:toponogov)}) на открытые поверхности. 
(То есть доказать, что $R_p=\infty$ для любой точки $p\in\Sigma$.) 
\end{subthm}

\end{thm}

\section{Лемма Александрова}
\index{лемма Александрова}

Следующая лемма (\ref{lem:alex-reformulation}) позволит получить несколько полезных переформулировок теорем сравнения.

\begin{thm}{Лемма}
\label{lem:alex}
Пусть $pxyz$ и $p'x'y'z'$ --- два четырехугольника на евклидовой плоскости с равными соответствующими сторонами.
Предположим, что стороны $[x',y']$ и $[y',z']$ являются продолжениями друг друга; то есть точка $y'$ лежит на отрезке $[x',z']$.
Тогда следующие выражения имеют одинаковые знаки:
\begin{enumerate}[(i)]
 \item $|p-y|-|p'-y'|$;
 \item $\measuredangle\hinge xpy-\measuredangle\hinge {x'}{p'}{y'}$;
 \item $\pi-\measuredangle\hinge ypx-\measuredangle\hinge ypz$.
\end{enumerate}
\end{thm}

\parbf{Доказательство.} 
Рассмотрим точку $\bar z$ на продолжении отрезка 
$[x,y]$ за точкой $y$ так, что $\dist{y}{\bar z}{}=\dist{y}{z}{}$ (и, следовательно, $\dist{x}{\bar z}{}=\dist{x'}{z'}{}$).

\begin{figure}[!ht]
\vskip-0mm
\centering
\includegraphics{mppics/pic-50}
\vskip-0mm
\end{figure}

По монотонности угла (\ref{lem:angle-monotonicity}), следующие выражения одного знака:
\begin{enumerate}[(i)]
\item $|p-y|-|p'-y'|$;
\item $\measuredangle\hinge{x}{y}{p}-\measuredangle\hinge{x'}{y'}{p'}=\measuredangle\hinge{x}{\bar z}{p}-\measuredangle\hinge{x'}{z'}{p'}$;
\item $|p-\bar z|-|p'-z'| = | p - \bar z | - | p-z | $;
\item $\measuredangle\hinge{y}{\bar z}{p}-\measuredangle\hinge{y}{z}{p}$.
\end{enumerate}
Поскольку
\[\measuredangle\hinge{y'}{z'}{p'}+\measuredangle\hinge{y'}{x'}{p'}=\pi
\quad\text{и}\quad
\measuredangle\hinge{y}{\bar z}{p}+\measuredangle\hinge{y}{x}{p}=\pi,\]
утверждение следует.
\qeds

\section{Переформулировки}

Для любого треугольника $[xyz]$ на поверхности $\Sigma$ и его модельного треугольника $[\tilde x \tilde y \tilde z]$, существует естественное отображение $p\mapsto \tilde p$, которое изометрически отображает геодезические $[x,y]$, $[y,z]$, $[z,x]$ на отрезки $[\tilde x,\tilde y]$, $[\tilde y, \tilde z]$, $[\tilde z, \tilde x]$ соответственно.
Треугольник $[xyz]$ называется \index{толстый/тонкий треугольник}\emph{толстым} (или {}\emph{тонким}), если неравенство
\[\dist{p}{q}{\Sigma}\ge |\tilde p- \tilde q|_{\mathbb{R}^2}\qquad \text{(или, соответственно,}\quad \dist{p}{q}{\Sigma}\le |\tilde p- \tilde q|_{\mathbb{R}^2})\]
выполняется для любых двух точек $p$ и $q$ на сторонах треугольника $[xyz]$.

\begin{thm}{Предложение}\label{prop:comp-reformulations}
Пусть $\Sigma$ --- открытая или замкнутая гладкая поверхность.
Тогда следующие три условия эквивалентны:

\begin{subthm}{mang>angk}
Для любого геодезического треугольника $[xyz]$ на $\Sigma$ выполняется неравенство
 \[\measuredangle\hinge{x}{y}{z}\ge\modangle xyz.\]
\end{subthm}

\begin{subthm}{angk>angk}
Для любого геодезического треугольника $[pxz]$ на $\Sigma$ и точки $y$ на стороне $[x,z]$ выполняется неравенство
 \[\modangle xpy \ge \modangle xpz.\]
 
\end{subthm}

\begin{subthm}{fat}
Любой геодезический треугольник на $\Sigma$ является толстым.
\end{subthm}

\medskip

Аналогично, следующие три условия эквивалентны:

\begin{subthmA}{mang<angk}
Для любого геодезического треугольника $[xyz]$ на $\Sigma$ выполняется неравенство
 \[\measuredangle\hinge{x}{y}{z}\le\modangle xyz.\]
\end{subthmA}

\begin{subthmA}{angk<angk}
Для любого геодезического треугольника $[pxz]$ на $\Sigma$ и точки $y$ на стороне $[x,z]$ выполняется неравенство
 \[\modangle xpy \le \modangle xpz.\]
\end{subthmA}

\begin{subthmA}{thin}
Любой геодезический треугольник на $\Sigma$ является тонким.
\end{subthmA}

\end{thm}

Давайте перепишем лемму Александрова (\ref{lem:alex}) на языке модельных треугольников и углов.

\begin{thm}{Переформулировка леммы Александрова}\label{lem:alex-reformulation}
Пусть $[pxz]$ --- треугольник на поверхности $\Sigma$, 
и $y\in[x,z]$. 
Рассмотрим его модельный треугольник $[\tilde p\tilde x\tilde z]\z=\tilde\triangle pxz$, и пусть $\tilde y$ будет соответствующей точкой на стороне $[\tilde x,\tilde z]$.

\begin{wrapfigure}{r}{25mm}
\vskip-2mm
\centering
\includegraphics{mppics/pic-2305}
\end{wrapfigure}

Тогда следующие выражения имеют одинаковые знаки:
\begin{enumerate}[(i)]
 \item $\dist{p}{y}\Sigma-\dist{\tilde p}{\tilde y}{\mathbb{R}^2}$;
 \item $\modangle xpy-\modangle {x}{p}{z}$;
 \item $\pi-\modangle ypx-\modangle ypz$;
\end{enumerate}
\end{thm}

\parbf{Доказательство \ref{prop:comp-reformulations}.}
Мы докажем последовательность импликаций \ref{SHORT.mang>angk}$\Rightarrow$\ref{SHORT.angk>angk}$\Rightarrow$\ref{SHORT.fat}$\Rightarrow$\ref{SHORT.mang>angk}.
Импликации \ref{SHORT.mang<angk}$\Rightarrow$\ref{SHORT.angk<angk}$\Rightarrow$\ref{SHORT.thin}$\Rightarrow$\ref{SHORT.mang<angk} доказываются также, заменяя знаки в неравенствах.

\parit{\ref{SHORT.mang>angk}$\Rightarrow$\ref{SHORT.angk>angk}.}
Заметим, что $\measuredangle\hinge ypx+\measuredangle\hinge ypz=\pi$.
По \ref{SHORT.mang>angk}, 
\[\modangle ypx+\modangle ypz\le \pi.\]
Остаётся применить лемму Александрова (\ref{lem:alex-reformulation}).

\parit{\ref{SHORT.angk>angk}$\Rightarrow$\ref{SHORT.fat}.}
Применив \ref{SHORT.angk>angk} дважды, сначала для $y\in [x,z]$, а затем для $w\in [p,x]$, получим
\[\modangle xwy \ge \modangle xpy \ge \modangle xpz.\]
Следовательно,
\[\dist{w}{y}\Sigma\ge \dist{\tilde w}{\tilde y}{\mathbb{R}^2},\]
где $\tilde w$ и $\tilde y$ --- точки, соответствующие $w$ и $y$ на сторонах модельного треугольника. 

\parit{\ref{SHORT.fat}$\Rightarrow$\ref{SHORT.mang>angk}.}
Поскольку треугольник толстый, имеем
\[\modangle xwy \ge \modangle xpz\]
для любого $w\in [x,p]\setminus \{x\}$ и $y\in [x,z]\setminus \{x\}$.
Заметим, что $\modangle xwy\z\to \measuredangle\hinge xpz$ при $w,y\to x$,
и импликация следует.
\qeds

\begin{thm}{Упражнение}\label{ex:geod-convexity}
Пусть $\gamma$ --- геодезическая с единичной скоростью на гладкой открытой поверхности
$\Sigma$, и $p\in\Sigma$.

Рассмотрим функцию
\[h(t)=\dist{p}{\gamma(t)}\Sigma^2-t^2.\]

\begin{subthm}{}
Покажите, что если $\Sigma$ односвязна и её гауссова кривизна неположительна, то функция $h$ является выпуклой.
\end{subthm}

\begin{subthm}{}
Покажите, что если гауссова кривизна $\Sigma$ неотрицательна, то функция $h$ является вогнутой.
\end{subthm}

\end{thm}

\begin{thm}{Упражнение}\label{ex:midpoints}
Пусть $\bar x$ и $\bar y$ --- середины кратчайших $[p,x]$ и $[p,y]$ на гладкой открытой поверхности~$\Sigma$.

\begin{subthm}{}
Покажите, что если $\Sigma$ является односвязной и её гауссова кривизна неположительна, то 
\[2\cdot \dist{\bar x}{\bar y}\Sigma\le \dist{x}{y}\Sigma.\]
\end{subthm}

\begin{subthm}{}
Покажите, что если гауссова кривизна $\Sigma$ неотрицательна, то 
 \[2\cdot \dist{\bar x}{\bar y}\Sigma\ge \dist{x}{y}\Sigma.\]
\end{subthm}

\end{thm}

\begin{thm}{Упражнение}\label{ex:convex-dist}
Пусть $\gamma_1$ и $\gamma_2$ --- две геодезические на односвязной открытой гладкой поверхности $\Sigma$ с неположительной гауссовой кривизной.
Покажите, что функция $h(t)\df\dist{\gamma_1(t)}{\gamma_2(t)}\Sigma$
является выпуклой.
\end{thm}

\begin{thm}{Упражнение}\label{ex:disc+}
Пусть $\Sigma$ --- открытая или замкнутая гладкая поверхность с неотрицательной гауссовой кривизной.
Докажите, что площадь любого $R$-шара во внутренней метрике $\Sigma$ не превышает $\pi\cdot R^2$.
\end{thm}

\begin{thm}{Упражнение}\label{ex:disc-}
Пусть $\Delta$ --- $R$-шар во внутренней метрике открытой, односвязной, гладкой поверхности $\Sigma$ неположительной гауссовой кривизной.

\begin{subthm}{ex:disc-:kg}
Докажите, что граница $\Delta$ является гладкой кривой с геодезической кривизной не менее $\tfrac1R$.
\end{subthm}

\begin{subthm}{ex:disc-:area}
Докажите, что площадь $\Delta$ не меньше $\pi\cdot R^2$.
\end{subthm}

\end{thm}

Следующее упражнение обобщает задачу о луне в луже (\ref{thm:moon-orginal}).

\begin{thm}{Продвинутое упражнение}\label{ex:moon-}
Пусть $\Delta$ --- топологический диск на гладкой поверхности $\Sigma$ неположительной гауссовой кривизны.
Предположим, что $\Delta$ ограничен гладкой кривой $\gamma$ с геодезической кривизной, не превышающей $1$ по абсолютному значению.
Докажите, что $\Delta$ содержит единичный круг во внутренней метрике $\Sigma$.

Выведите отсюда, что площадь $\Delta$ не меньше $\pi$.
\end{thm}

{
\appendix
\bookmarksetupnext{startatroot,level=0}
\chapter[Приложение]{}
\chaptermark{}

\vskip5mm

{\footnotesize

Это приложение лучше использовать как справочник.
В основном здесь даются формулировки нужных результатов и ссылки на их доказательства.

\section{Метрические пространства}\label{sec:metric-spcaes}

В этом разделе мы вводим обозначения, которые будем использовать.
Все эти темы подробно обсуждаются в вводной части книги Дмитрия Бураго, Юрия Бураго и Сергея Иванова \cite{burago-burago-ivanov}.
Предполагается, что читатель знаком с метрикой на евклидовом пространстве.

\spell{\begin{multicols}{2}}{}

\subsection*{Определения}

\emph{Метрика} --- это функция, которая возвращает вещественное значение $\Dist(x,y)$ для любой пары элементов $x,y$ данного множества $\spc{X}$ и удовлетворяет следующим свойствам для любой тройки $x,y,z\in\spc{X}$: \label{page:def:metric}
\begin{enumerate}[(a)]
\item\label{def:metric-space:a}
Положительность: 
$$\Dist(x,y)\ge 0.$$
\item\label{def:metric-space:b}
$x=y$ тогда и только тогда, когда 
$$\Dist(x,y)=0.$$
\item\label{def:metric-space:c}
Симметричность: $$\Dist(x, y) = \Dist(y, x).$$
\item\label{def:metric-space:d}
Неравенство треугольника: 
$$\Dist(x, z) \le \Dist(x, y) + \Dist(y, z).$$
\end{enumerate}

{\sloppy

Множество с выбранной метрикой называется \index{метрическое пространство}\emph{метрическим пространством}, а элементы множества называются его \index{точка}\emph{точками}.

}

\subsection*{Обозначение метрики}

Как правило, мы рассматриваем только одну метрику на множестве, поэтому допустимо обозначать метрическое пространство и его подлежащее множество одной и той же буквой, скажем, $\spc{X}$.
В этом случае мы также используем сокращённые обозначения $\dist{x}{y}{}$ или $\dist{x}{y}{\spc{X}}$ для \emph{расстояния} $\Dist(x,y)$ от $x$ до $y$ в $\spc{X}$.\index{10aaa@$\dist{x}{y}{}$, $\dist{x}{y}{\spc{X}}$ (расстояние)}
Например, неравенство треугольника можно записать как
$$\dist{x}{z}{\spc{X}}\le \dist{x}{y}{\spc{X}}+\dist{y}{z}{\spc{X}}.$$

В евклидовом пространстве (как и на плоскости и на вещественной прямой) можно думать, что расстояние $\dist{x}{y}{}$ от $x$ до $y$ есть норма вектора $x-y$.
И хоть в метрическом пространстве разница точек $x-y$ не имеет смысла, мы  обозначаем расстояние как $\dist{x}{y}{}$.

\subsection*{Ещё примеры}

Обычно, если мы говорим \emph{плоскость} или \emph{пространство}, имеется в виду \emph{евклидова} плоскость или пространство.
Однако плоскость (как и пространство) допускает множество других метрик; например, так называемую \emph{метрику городских кварталов} из следующего упражнения.

\begin{thm}{Упражнение}\label{ex:ell-infty}
Рассмотрим функцию
$$\Dist(p,q)=|x_1-x_2|+|y_1-y_2|,$$
где $p=(x_1,y_1)$ и $q=(x_2,y_2)$ --- точки на координатной плоскости~$\mathbb{R}^2$.
Докажите, что это метрика на $\mathbb{R}^2$.
\end{thm}

Ещё один пример: \emph{дискретное пространство} --- произвольное непустое множество $\spc{X}$ с метрикой, определяемой как $\dist{x}{y}{\spc{X}}=0$, если $x=y$, и $\dist{x}{y}{\spc{X}}=1$ в противном случае.

\subsection*{Подпространства}

{\sloppy

Любое подмножество метрического пространства само окажется метрическим пространством, если сузить исходную метрику на это подмножество;
полученное метрическое пространство называется \emph{подпространством}.
В частности, все подмножества евклидова пространства являются метрическими пространствами.

}

\subsection*{Шары}

Для данной точки $p$ метрического пространства $\spc{X}$ и вещественного числа $R\ge 0$ множество точек $x$, находящихся на расстоянии меньше чем $R$ (не больше чем $R$) от $p$, называется \index{открытый!шар}\emph{открытым} (соответственно, \index{замкнутый!шар}\emph{замкнутым}) \emph{шаром} радиуса $R$ с центром в~$p$.
Открытый шар обозначается как $B(p,R)$ или $B(p,R)_{\spc{X}}$;
второе обозначение используется, если хочется подчеркнуть, что наш шар живёт в $\spc{X}$.
Строго говоря,
\begin{align*}
B(p,R)&=B(p,R)_{\spc{X}}=
\\
&=\set{x\in \spc{X}}{\dist{x}{p}{\spc{X}}< R}.
\end{align*}
\index{10b@$B(p,R)_{\spc{X}}$, $\bar B[p,R]_{\spc{X}}$ (шар)}
Замкнутый шар обозначается через $\bar B[p,R]$ или $\bar B[p,R]_{\spc{X}}$, и
\begin{align*}
\bar B[p,R]&=\bar B[p,R]_{\spc{X}}=
\\
&=\set{x\in \spc{X}}{\dist{x}{p}{\spc{X}}\le R}.
\end{align*}

\begin{thm}{Упражнение}\label{ex:B2inB1}

\begin{subthm}{ex:B2inB1:a}
Пусть $p$ и $q$ --- точки метрического пространства $\spc{X}$.
Докажите, что если $\bar B[p,2]\z\subset \bar B[q,1]$, то $\bar B[p,2]\z=\bar B[q,1]$.
\end{subthm}

\begin{subthm}{ex:B2inB1:b}
Постройте метрическое пространство $\spc{X}$ и два шара $B(p,\tfrac32)$ и $B(q,1)$ в нём, со строгим включением
$B(p,\tfrac32)\z\subsetneq B(q,1)$.
\end{subthm}

\end{thm}

\subsection*{Изометрии}

Пусть $\spc{X}$ и $\spc{Y}$ --- метрические пространства.
Отображение $f\:\spc{X} \z\to \spc{Y}$ \index{сохраняющее расстояние}\emph{сохраняет расстояние}, если 
$$\dist{f(x)}{f(y)}{\spc{Y}}
 = \dist{x}{y}{\spc{X}}$$
при любых $x,y\in {\spc{X}}$.

\index{изометрия}\emph{Изометрия} определяется как биективное отображение, сохраняющее расстояние.
Два метрических пространства называются \index{изометричные пространства}\emph{изометричными}, если найдётся изометрия между ними.

\begin{thm}{Упражнение}\label{ex:dist-preserv=>injective}
Покажите, что если отображение $f\:\spc{X}\to\spc{Y}$ сохраняет расстояния, то оно \index{инъективное отображение}\emph{инъективно};
то есть $f(x)\ne f(y)$ для любой пары различных точек $x, y\in \spc{X}$.
\end{thm}

\subsection*{Непрерывность}

\begin{thm}{Определение}
Пусть ${\spc{X}}$ --- метрическое пространство.
Последовательность точек $x_1, x_2, \ldots$ в ${\spc{X}}$ \emph{сходится}, 
если существует точка $x_\infty\in {\spc{X}}$, такая что $\dist{x_\infty}{x_n}{}\to 0$ при $n\to\infty$.  
То есть для любого $\epsilon > 0$ существует натуральное число $N$, такое что 
\[
\dist{x_\infty}{x_n}{\spc{X}}
<
\epsilon
\]
при всех $n \ge N$.

В этом случае говорят, что последовательность $x_n$ \emph{сходится} к $x_\infty$, 
или что $x_\infty$ --- \emph{предел} последовательности $x_n$,
и пишут $x_n\to x_\infty$ при $n\to\infty$ или $x_\infty=\lim_{n\to\infty} x_n$.
\end{thm}

\begin{thm}{Определение}\label{def:continuous}
Пусть $\spc{X}$ и $\spc{Y}$ --- метрические пространства.
Отображение $f\:\spc{X}\to \spc{Y}$ называется \index{непрерывный}\emph{непрерывным}, если из сходимости $x_n\z\to x_\infty$ в ${\spc{X}}$
следует сходимость $f(x_n) \z\to f(x_\infty)$ в $\spc{Y}$.

Эквивалентно, $f\:\spc{X}\to \spc{Y}$ непрерывно, если для любого $x\in {\spc{X}}$ и любого $\epsilon>0$
существует $\delta>0$ такое, что 
$$\dist{x}{y}{\spc{X}}<\delta\quad\Longrightarrow\quad \dist{f(x)}{f(y)}{\spc{Y}}<\epsilon.$$

\end{thm}

\begin{thm}{Упражнение}\label{ex:shrt=>continuous}
Пусть $f\:\spc{X}\z\to \spc{Y}$ --- \emph{короткое} отображение между метрическими пространствами, то есть
\[\dist{f(x)}{f(y)}{\spc{Y}}\le \dist{x}{y}{\spc{X}}\]
при любых $x,y\in \spc{X}$.
Покажите, что $f$ непрерывно.
\end{thm}

\spell{\end{multicols}}{}

\section{Топология}\label{sec:topology}

Материал этого раздела обсуждается в любом вводном тексте по топологии;
например, в классическом учебнике Чеса Коснёвского \cite{kosniowski}.

\spell{\begin{multicols}{2}}{}

\subsection*{Замкнутые и открытые множества}

\begin{thm}{Определение}
Подмножество $C$ метрического пространства $\spc{X}$ называется \index{замкнутое множество}\emph{замкнутым}, если всякая сходящаяся в $\spc{X}$ последовательность точек из $C$ имеет предел в $C$.

Множество $\Omega \subset \spc{X}$ называется \index{открытое!множество}\emph{открытым}, если для любой точки $z\in \Omega$ 
существует $\epsilon>0$ такое, что $B(z,\epsilon)\subset\Omega$.
\end{thm}

\begin{thm}{Упражнение}\label{ex:close-open}
Пусть $Q$ --- подмножество метрического пространства $\spc{X}$.
Покажите, что $Q$ замкнуто тогда и только тогда, когда его дополнение $\Omega=\spc{X}\setminus Q$ открыто.
\end{thm}

Открытое множество $\Omega$, содержащее заданную точку $p$, называется \index{окрестность}\emph{окрестностью}~$p$.
Замкнутое подмножество $C$, содержащее точку $p$ вместе с её какой-то её окрестностью, называется \emph{замкнутой окрестностью}~$p$.

Точка $p$ лежит на \index{граница}\emph{границе} множества $Q$ (это обозначается как $p\in\partial Q$), если любая окрестность точки $p$ содержит точки как из $Q$, так и из его дополнения.

\subsection*{Компактные множества}

{\sloppy
Подмножество $K$ метрического пространства называется \index{компактное множество}\emph{компактным}, если в любой последовательности точек из $K$ найдётся сходящаяся в $K$ подпоследовательность.

}

Из определения следует пара свойств.

\begin{itemize}
\item Замкнутое подмножество компактного множества компактно.
\item Для непрерывного отображения, образ компактного множества компактен.
\end{itemize}

{\sloppy

\begin{thm}{Лемма Гейне --- Бореля}\label{thm:Heine--Borel}
Подмножество евклидова пространства компактно тогда и только тогда, когда оно замкнуто и ограничено.
\end{thm}

}

\subsection*{Гомеоморфизмы\\ и вложения}

Биекция $f\:\spc{X}\to\spc{Y}$ между метрическими пространствами называется \index{гомеоморфизм}\emph{гомеоморфизмом}, если $f$ и $f^{-1}$ непрерывны.
Гомеоморфизм на свой образ называется \index{вложение}\emph{вложением}.

Если существует гомеоморфизм $f\:\spc{X}\to \spc{Y}$,
то мы говорим, что ${\spc{X}}$ {}\emph{гомеоморфно} $\spc{Y}$,
или что $\spc{X}$ и $\spc{Y}$ {}\emph{гомеоморфны}.

{\sloppy

Если метрическое пространство гомеоморфно известному пространству, например, плоскости, сфере, диску, окружности и так далее,
то его можно назвать \index{топологическая поверхность}\emph{топологической} плоскостью, сферой, диском, окружностью и так далее.

}

Следующая теорема характеризует гомеоморфизмы между компактными пространствами:

\begin{thm}{Теорема}\label{thm:Hausdorff-compact}
Непрерывная биекция $f$ между компактными метрическими пространствами имеет непрерывную обратную.
В частности, справедливо следующее.

\begin{subthm}{}
Любая непрерывная биекция между компактными метрическими пространствами
является гомеоморфизмом.
\end{subthm}

{\sloppy

\begin{subthm}{}
Любое непрерывное инъективное отображение компактного метрического пространства в другое метрическое пространство
является вложением.
\end{subthm}

}

\end{thm}

\subsection*{Связные множества}

Напомним, что непрерывное отображение $\alpha$ из отрезка $[0,1]$ в евклидово пространство называется \index{путь}\emph{путём}.
Если $p=\alpha (0)$ и $q = \alpha (1)$, то мы говорим, что $\alpha$ соединяет $p$ с~$q$.

Непустое множество $X$ евклидова пространства называется \index{линейно связное множество}\emph{линейно связным}, если любые две точки $x,y\z\in X$ можно соединить путём, лежащим в~$X$.

Непустое множество $X$ евклидова пространства называется \index{связное множество}\emph{связным}, если его нельзя покрыть двумя непересекающимися открытыми множествами $V$ и $W$ такими, что оба пересечения $X\cap V$ и $X\cap W$ непусты.

Обратите внимание, что связные и линейно связные пространства по определению не пусты. 

\begin{thm}{Утверждение}
Любое линейно связное множество связно.

Более того, любое открытое связное множество евклидова пространства (или плоскости) линейно связно.
\end{thm}

Для данной точки $x\in X$ максимальное связное подмножество $X$, содержащее $x$, называется {}\emph{связной компонентой} $x$ в~$X$.

\subsection*{Теорема Жордана}
\index{теорема Жордана}

Первая часть следующей теоремы доказана Камилем Жорданом, вторая --- Артуром Шёнфлисом.

\begin{thm}{Теорема}\label{thm:jordan}
Дополнение любой простой замкнутой кривой $\gamma$ в $\mathbb{R}^2$ имеет ровно две связные компоненты. 
Более того, существует гомеоморфизм $h\:\mathbb{R}^2\to \mathbb{R}^2$, отображающий единичную окружность в~$\gamma$.
В частности, $\gamma$ --- граница топологического диска.
\end{thm}

Эта теорема известна простотой формулировки и сложностью доказательства.
Короткое доказательство первого утверждения, основанное на несколько продвинутой технике, дано Патриком Дойлом \cite{doyle}, но есть и вполне элементарные доказательства, одно из таких найденно Алексеем Филипповым~\cite{filippov}.

В основном мы будем пользоваться этой теоремой для гладких кривых.
Этот случай проще, любопытное его доказательство найдено Григорием Чамберсом и Евгением Лиокумовичем \cite{chambers-liokumovich}.

\spell{\end{multicols}}{}

\section{Элементарная геометрия}

\spell{\begin{multicols}{2}}{}

\subsection*{Внутренние углы}

Многоугольник определяется как компактное множество, ограниченное замкнутой ломаной. 
Напомним, что внутренний угол многоугольника $P$ при вершине $v$
определяется как угловая мера пересечения $P$ и маленькой окружности с центром в~$v$.

\begin{thm}{Теорема}\label{thm:sum=(n-2)pi}
Сумма всех внутренних углов $n$-угольника равна $(n\z-2)\cdot\pi$.
\end{thm}

Честное доказательство этой теоремы приводится например в \cite{meisters}.
Оно использует индукцию по $n$ и основано на следующем не вполне тривиальном факте.

\begin{thm}{Факт}
Пусть $P$ --- $n$-угольник с $n\ge 4$.
Тогда одна из диагоналей $P$ полностью лежит в~$P$.
\end{thm}

\subsection*{Монотонность углов}

{}\emph{Мера} угла со сторонами $[p,x]$ и $[p,y]$ будет обозначаться $\measuredangle\hinge pxy$\index{10aab@$\measuredangle\hinge yxz$ (мера угла)};
она принимает значение в интервале $[0,\pi]$.

Следующее следствие теоремы косинусов простое и полезное.
Оно говорит, что угол треугольника монотонно зависит от противолежащей стороны, если две другие стороны не меняют длины.

\begin{thm}{Лемма}\label{lem:angle-monotonicity}
Предположим, что для точек $x$, $y$, $z$, $x^{*}$, $y^{*}$ и $z^{*}$ выполнено $\dist{x}{y}{}\z=|x^{*}-y^{*}|>0$ и $|y-z|\z=|y^{*}-z^{*}|>0$.
Тогда 
\[\measuredangle\hinge yxz
\ge
\measuredangle\hinge {y^{*}}{x^{*}}{z^{*}}
\ \Longleftrightarrow\
|x-z|\ge |x^{*}-z^{*}|.\]
\end{thm}

\subsection*{Сферическое неравенство треугольника}

Следующая теорема утверждает, что неравенство треугольника выполняется для углов между лучами с общим началом.
В частности, из этого следует, что сфера с угловой метрикой является метрическим пространством.

\begin{thm}{Теорема}\label{thm:spherical-triangle-inq}
Для любых трёх отрезков $[o,a]$, $[o,b]$ и $[o,c]$ евклидова пространства выполняется неравенство
\[\measuredangle\hinge oab
+
\measuredangle\hinge obc
\ge
\measuredangle\hinge oac.\]

\end{thm}

Эту теорему часто используют вовсе без упоминания, но её доказательство не такое уж простое.
Его можно найти в классическом учебнике евклидовой геометрии Андрея Киселёва \cite[§ 47]{kiselyov}.

\subsection*{Площадь сферического треугольника}

\begin{thm}{Лемма}\label{lem:area-spher-triangle}
Пусть $\Delta$ --- сферический треугольник,
то есть пересечение трёх замкнутых полусфер на единичной сфере $\mathbb{S}^2$.
Тогда 
\[\area\Delta=\alpha+\beta+\gamma-\pi,\eqlbl{eq:area(Delta)}\]
где $\alpha$, $\beta$ и $\gamma$ --- углы треугольника $\Delta$.
\end{thm}

Величина $\alpha+\beta+\gamma-\pi$ называется \index{избыток треугольника}\emph{избытком} треугольника $\Delta$;
то есть можно сказать, что площадь сферического треугольника равна его избытку.

Эта лемма сыграет роль в интуитивном объяснении формулы Гаусса --- Бонне.
Поэтому мы приведём её доказательство.

\begin{wrapfigure}{r}{16 mm}
\vskip-2mm
\centering
\includegraphics{mppics/pic-43}
\vskip2mm
\end{wrapfigure}

\parbf{Доказательство.}
Напомним, что 
\[\area\mathbb{S}^2=4\cdot\pi.\eqlbl{eq:area(S2)}\]

Заметим, что площадь сферического сектора $S_\alpha$ между двумя меридианами, пересекающимися под углом $\alpha$, пропорциональна $\alpha$.
Поскольку $S_\pi$ --- это полусфера, из \ref{eq:area(S2)} получаем $\area S_\pi\z=2\cdot\pi$.
Следовательно, коэффициент равен $2$; то есть
\[\area S_\alpha=2\cdot \alpha
\eqlbl{eq:area(Sa)}\]
при любом $\alpha$.

Продлив стороны $\Delta$, получим $6$ секторов: два $S_\alpha$, два $S_\beta$ и два $S_\gamma$.
Они покрывают почти всю сферу один раз,
но треугольник $\Delta$ и его центрально-симметричная копия $\Delta^{*}$ покрываются трижды.
Отсюда следует, что
\begin{align*}
2\cdot \area S_\alpha &+2\cdot \area S_\beta+2\cdot \area S_\gamma=
\\
&=\area\mathbb{S}^2+4\cdot\area\Delta.
\end{align*}
Остаётся применить \ref{eq:area(S2)} и \ref{eq:area(Sa)}.
\qeds

\spell{\end{multicols}}{}

\section{Выпуклая геометрия}

\spell{\begin{multicols}{2}}{}

Множество $X$ евклидова пространства называется \index{выпуклое множество}\emph{выпуклым}, если для любых двух точек $x,y\in X$ любая точка $z$ между ними лежит в~$X$.
Далее, $X$ называется {}\emph{строго выпуклым}, если для любых двух точек $x,y\in X$ любая точка $z$ между $x$ и $y$ лежит во внутренности~$X$.

Из определения видно, что пересечение произвольного семейства выпуклых множеств выпукло. 
Пересечение всех выпуклых множеств, содержащих данное множество $X$, называется его \index{выпуклая!оболочка}\emph{выпуклой оболочкой};
это минимальное выпуклое множество, содержащее~$X$.

Эти определения и следующие утверждения найдутся на первых страницах любого вводного текста по выпуклой геометрии;
например, в книге Курта Лайхтвайса \cite{leichtweiss} 

\subsection*{Разделяющие и опорные плоскости}

Следующая лемма --- частный случай так называемой \index{теорема о разделяющей гиперплоскости}\emph{теоремы о разделяющей гиперплоскости}.

\begin{thm}{Лемма}\label{lem:separation}
Пусть $K\subset \mathbb{R}^3$ --- замкнутое выпуклое множество.
Тогда для любой точки $p\notin K$ существует плоскость $\Pi$, которая разделяет $K$ и $p$;
то есть $K$ и $p$ лежат в противоположных открытых полупространствах, разделённых плоскостью $\Pi$.

Более того, для любой граничной точки $p\in\partial K$ существует \index{опорная!плоскость}\emph{опорная плоскость}%
\footnote{Можно также говорить, что $\Pi$ подпирает $K$ в точке $p$.} $\Pi$ к $K$ в точке $p$;
то есть $\Pi\ni p$ и $K$ лежит в замкнутом полупространстве, ограниченном плоскостью $\Pi$.
\end{thm}

\spell{\end{multicols}}{}

\begin{figure*}[h!]
\begin{minipage}{.48\textwidth}
\centering
\includegraphics{mppics/pic-3540}
\end{minipage}\hfill
\begin{minipage}{.48\textwidth}
\centering
\includegraphics{mppics/pic-3542}
\end{minipage}
\end{figure*}

\section{Линейная алгебра}

\spell{\begin{multicols}{2}}{}

Следующая теорема найдётся в любом учебнике по линейной алгебре;
подойдёт, например, книга Сергея Трейла \cite{treil}.

{\sloppy

\begin{thm}{Спектральная теорема}\label{thm:spectral}
Любая симметричная матрица диагонализируется ортогональной матрицей.
\end{thm}

}

Нам потребуется только матрицы $2{\times}2$.
В этом случае теорему можно переформулировать следующим образом:
Рассмотрим функцию
\begin{align*}
f(x,y)&=
\begin{pmatrix}
x&y
\end{pmatrix}
\cdot
\begin{pmatrix}
\ell&m
\\
m&n
\end{pmatrix}
\cdot
\begin{pmatrix}
x\\y
\end{pmatrix}=
\\
&=\ell\cdot x^2+2\cdot m\cdot x\cdot y+n\cdot y^2,
\end{align*}
определенную на координатной плоскости.
Тогда после подходящего поворота координат
выражение для $f$ в новых $(x,y)$-координатах будет иметь вид
\begin{align*}
\bar f(x,y)&=
\begin{pmatrix}
x&y
\end{pmatrix}
\cdot
\begin{pmatrix}
k_1&0
\\
0&k_2
\end{pmatrix}
\cdot
\begin{pmatrix}
x\\y
\end{pmatrix}=
\\
&=k_1\cdot x^2+k_2\cdot y^2.
\end{align*}

\spell{\end{multicols}}{}

\section{Анализ}\label{sec:analysis}

Следующие теоремы обсуждаются в любом вводном курсе по вещественному анализу.
Например, в классическом учебнике Уолтера Рудина \cite{rudin}.

\spell{\begin{multicols}{2}}{}

\subsection*{Измеримые функции}

Функция называется \index{измеримая функция}\emph{измеримой}, если прообраз любого борелевского множества борелевский.
Практически все функции, которые естественным образом появляются в геометрии, являются измеримыми.

Следующая теорема позволяет распространить многие утверждения с непрерывных функций на измеримые.

\begin{thm}{Теорема Лузина}\label{thm:lusin}
Пусть $\phi\:[a,b]\to \mathbb{R}$ --- измеримая функция.
Для любого $\epsilon>0$ найдётся непрерывная функция $\psi_\epsilon\:[a,b]\z\to \mathbb{R}$, совпадающая с $\phi$ вне множества меры, не превышающей $\epsilon$.
Более того, если $\phi$ ограничена сверху и/или снизу, то можно считать, что то же верно и для~$\psi_\epsilon$.  
\end{thm}

\subsection*{Условие Липшица}

Напомним, что функция $f$ между метрическими пространствами называется \index{липшицева функция}\emph{липшицевой}, если найдётся константа $L$ такая, что 
\[\dist{f(x)}{f(y)}{}\le L\cdot\dist{x}{y}{}\]
для всех значений $x$ и $y$ в области определения~$f$.

Следующая теорема позволяет распространить многие результаты с гладких функций на липшицевы.
Напомним, что {}\emph{почти все} означает \textit{все, за исключением множества нулевой меры Лебега}.

\begin{thm}{Теорема Радемахера}\label{thm:rademacher}
Пусть $f\:[a,b]\to\mathbb{R}$ --- липшицева функция.
Тогда её производная $f'$ --- ограниченная измеримая функция, определённая почти везде на $[a,b]$.
Более того, к ней применима формула Ньютона --- Лейбница; то есть верно равенство
\[f(b)-f(a)=\int_a^b f'(x)\cdot dx,\]
если интеграл понимать в смысле Лебега.
\end{thm}

\subsection*{Равномерная непрерывность и сходимость}

Пусть $f\:{\spc{X}}\to \spc{Y}$ --- отображение между метрическими пространствами.
Если для любого $\epsilon>0$ существует $\delta\z>0$, такое что
\[\dist{x_1}{x_2}{\spc{X}}<\delta\ \Longrightarrow\ \dist{f(x_1)}{f(x_2)}{\spc{Y}}<\epsilon,\]
то функция $f$ называется \index{равномерная непрерывность}\emph{равномерно непрерывной}.

Ясно, что любая равномерно непрерывная функция непрерывна;
обратное неверно.
Например, функция $f(x)=x^2$ непрерывна, но не является равномерно непрерывной.

{\sloppy

\begin{thm}{Теорема}
Любая непрерывная функция на компактном метрическом пространстве равномерно непрерывна.
\end{thm}

}

Если условие выше выполняется для любой функции $f_n$ в последовательности, и $\delta$ зависит только от $\epsilon$,
то такая последовательность $f_n$ называется \index{равностепенно непрерывная последовательность}\emph{равностепенно непрерывной}.
Более точно, 
последовательность функций $f_n:{\spc{X}}\to \spc{Y}$ называется равностепенно непрерывной, если 
для любого $\epsilon>0$ существует $\delta>0$, такое что 
\[\dist{x_1}{x_2}{\spc{X}}<\delta\ \Rightarrow\ \dist{f_n(x_1)}{f_n(x_2)}{\spc{Y}}<\epsilon\]
для любого~$n$.

Мы говорим, что последовательность функций $f_i\: {\spc{X}} \to \spc{Y}$ \index{равномерная сходимость}\emph{равномерно сходится} к функции $f_{\infty}\: {\spc{X}} \to \spc{Y}$, если для любого 
$\epsilon>0$ существует натуральное число $N$ такое, что $\dist{f_{\infty}(x)}{f_n (x)}{}<\epsilon$ для всех $n \ge N$.

{\sloppy

\begin{thm}{Теорема Арцела --- Асколи}\label{lem:equicontinuous}
Пусть $\spc{X}$ и $\spc{Y}$ --- компактные метрические пространства. 
Тогда в любой равностепенно непрерывной последовательности функций $f_n\:\spc{X}\z\to \spc{Y}$ найдётся подпоследовательность, равномерно сходящаяся к некоторой непрерывной функции $f_\infty\:\spc{X}\z\to \spc{Y}$. 
\end{thm}

}

\subsection*{Срезки и сглаживания}

Здесь мы обсудим как строить гладкие функции, которые ведут себя как определённые модельные функций.
Эти построения используются при сглаживании модельных геометрических объектов.

Для начала рассмотрим следующие две функции
\begin{align*}
h(t)&\df
\begin{cases}
0&\text{если}\ t\le 0,
\\
t&\text{если}\ t> 0;
\end{cases}
\\
f(t)&\df
\begin{cases}
0&\text{если}\ t\le 0,
\\
\frac{t}{e^{1\!/\!t}}&\text{если}\ t> 0.
\end{cases}
\end{align*}
Они ведут себя похоже ---
обе равны нулю при $t\le 0$ и возрастают до бесконечности при положительных $t$.
Функция $h$ не является гладкой --- у неё неопределена производная в нуле.
При этом функция $f$ гладкая.
Действительно, существование всех производных $f^{(n)}(x)$ при $x\ne 0$ очевидно, в то время как прямые вычисления дают, что $f^{(n)}(0)=0$ для всех~$n$.

Другой пример: \index{колокольчик}\emph{колокольчик} --- гладкая функция, которая положительна в $\epsilon$-окрестности нуля и обращающаяся в ноль за её пределами.
Такие функции можно получить из функции $f$, построенной выше, например как
\[b_\epsilon(t)\df c\cdot f(\epsilon^2-t^2);\]
обычно константу $c$ подбирают так, чтобы $\int b_\epsilon=1$.

Ещё один полезный пример: \index{сигмоид}\emph{сигмоид} --- неубывающая функция, которая равна нулю при $t\le -\epsilon$ и единице при $t\ge \epsilon$.
Она имитирует ступенчатую функцию; её можно задать как \label{page:sigma-function}
\[\sigma_\epsilon(t)
\df 
\int_{-\infty}^t b_\epsilon(x)\cdot dx.\]

\spell{\end{multicols}}{}

\begin{figure*}[h!]
\begin{minipage}{.48\textwidth}
\centering
\includegraphics{mppics/pic-320}
\end{minipage}\hfill
\begin{minipage}{.48\textwidth}
\centering
\includegraphics{mppics/pic-321}
\end{minipage}
\end{figure*}

\begin{figure*}[h!]
\begin{minipage}{.48\textwidth}
\centering
\includegraphics{mppics/pic-325}
\end{minipage}\hfill
\begin{minipage}{.48\textwidth}
\centering
\includegraphics{mppics/pic-326}
\end{minipage}
\end{figure*}

\section{Векторный анализ}\label{sec:Multivariable calculus}

Следующие теоремы обсуждаются в любом курсе по векторному анализу;
например, в уже упомянутом учебнике Уолтера Рудина \cite{rudin}.

\spell{\begin{multicols}{2}}{}

\subsection*{Регулярные значения}

Пусть $\Omega\subset\mathbb{R}^m$ --- открытое множество.
Любое отображение $\bm{f}\:\Omega\z\to\mathbb{R}^n$ можно рассматривать как набор координатных функций
\[f_1,\dots,f_n\:\Omega\to \mathbb{R};\]
$\bm{f}$ называется \index{гладкое отображение}\emph{гладким}, если каждая функция $f_i$ гладкая;
то есть всевозможные частные производные функции $f_i$ определены в области $\Omega$.

\emph{Матрица Якоби} отображения $\bm{f}$ в точке $\bm{x}\in\mathbb{R}^m$ определяется как \index{10j@$\Jac$ (матрица Якоби)}
\[\Jac_{\bm{x}}\bm{f}=
\begin{pmatrix}
\dfrac{\partial f_1}{\partial x_1} & \cdots & \dfrac{\partial f_1}{\partial x_m}\\
\vdots & \ddots & \vdots\\
\dfrac{\partial f_n}{\partial x_1} & \cdots & \dfrac{\partial f_n}{\partial x_m} \end{pmatrix};\]
предполагается, что правая часть вычислена в точке $\bm{x}=(x_1,\dots,x_m)$.

Если матрица Якоби определяет сюръективное линейное отображение $\mathbb{R}^m\to\mathbb{R}^n$ (то есть если $\rank(\Jac_{\bm{x}}\bm{f})\z=n$), то мы говорим, что
$\bm{x}$ --- это \index{регулярная точка}\emph{регулярная точка}~$\bm{f}$.

Если $\bm{x}$ является регулярной точкой всякий раз, когда $\bm{f}(\bm{x})=\bm{y}$,
то $\bm{y}$ называется \index{регулярное значение}\emph{регулярным значением}~$\bm{f}$.
Следующая лемма утверждает, что \textit{большинство} точек в $\mathbb{R}^n$ являются регулярными значениями $\bm{f}$.

\begin{thm}{Лемма Сарда}\label{lem:sard}
Для гладкого отображения $\bm{f}\colon \Omega \z\to \mathbb{R}^n$, определённого на открытом множестве $\Omega\subset \mathbb{R}^m$, почти все значения в $\mathbb{R}^n$ регулярны.
\end{thm}

Слова \index{почти все}\emph{почти все} означают все, за исключением множества нулевой меры Лебега.
В частности, если выбрать случайное значение, равномерно распределённое в произвольно малом шаре $B\z\subset \mathbb{R}^n$, то оно окажется регулярным $\bm{f}$ с вероятностью~$1$.

Обратите внимание, что если $m\z<n$, то любая точка $\bm{y}=\bm{f}(\bm{x})$ \textit{не} является регулярным значением $\bm{f}$.
Таким образом, регулярные значения $\bm{f}$ не являются значениями $\bm{f}$, то есть принадлежат дополнению образа $\Im \bm{f}$.
В этом случае теорема утверждает, что почти все точки в $\mathbb{R}^n$ \textit{не} принадлежат $\Im \bm{f}$.

\subsection*{Теорема об обратной функции}

{\sloppy

\index{теорема об обратной функции}\emph{Теорема об обратной функции} даёт условие, при котором отображение $\bm{f}$ обратимо в окрестности заданной точки $\bm{x}$.
Это условие формулируется в терминах $\Jac_{\bm{x}}\bm{f}$ --- матрицы Якоби отображения $\bm{f}$ в точке~$\bm{x}$.

}

\index{теорема о неявной функции}\emph{Теорема о неявной функции} --- близкая ей родственница и следствие.
Эта теорема понадобится при переходе от параметрического описания кривых и поверхностей к их неявному описанию и наоборот.

Обе теоремы сводят существование отображения, удовлетворяющего определённому уравнению, к задаче линейной алгебры.
Нам они потребуются только при в размерностях $\le 3$.

\begin{thm}{Теорема}\label{thm:inverse}
{\sloppy
Пусть $\bm{f}=(f_1,\z\dots,f_n)\:\Omega\to\mathbb{R}^n$ --- гладкое отображение, определённое на открытом множестве $\Omega\subset \mathbb{R}^n$.
Предположим, что
$\Jac_{\bm{x}}\bm{f}$
обратима в точке $\bm{x}\in \Omega$.
Тогда существует гладкое отображение $\bm{h}\:\Phi\to\mathbb{R}^n$, определённое в некоторой открытой окрестности $\Phi$ точки ${\bm{y}}\z=\bm{f}(\bm{x})$, которое является {}\emph{локально обратным} к $\bm{f}$ в точке $\bm{x}$;
то есть существует окрестность $\Psi\ni \bm{x}$ такая, что
$\bm{f}$ определяет гомеоморфизм $\Psi\leftrightarrow \Phi$, и
$\bm{h} \circ \bm{f}$ является тождественным отображением на~$\Psi$.

}

Более того, если $|\det[\Jac_{\bm{x}}\bm{f}]|\z>\epsilon\z>0$, область $\Omega$ содержит $\epsilon$-окрестность точки $\bm{x}$,
а первые и вторые частные производные $\tfrac{\partial f_i}{\partial x_j}$, $\tfrac{\partial^2 f_i}{\partial x_j\partial x_k}$ ограничены константой $C$ для всех $i$, $j$ и $k$, то можно предположить, что $\Phi$ является $\delta$-окрестностью точки $\bm{y}$, где $\delta>0$ зависит только от $\epsilon$ и~$C$. 
\end{thm}

\begin{thm}{Теорема о неявной функции}\label{thm:imlicit}
Пусть $\bm{f}=(f_1,\dots,f_n)\:\Omega\z\to\mathbb{R}^n$ --- гладкое отображение, определённое на открытом подмножестве $\Omega\subset\mathbb{R}^{n+m}$, где
$m,n\z\ge 1$.
Рассмотрим $\mathbb{R}^{n+m}$ как декартово произведение $\mathbb{R}^n\times \mathbb{R}^m$ с координатами
$(\bm{x},\bm{y})\z=(x_1,\dots,x_n,y_1,\dots,y_m)$.
Далее, рассмотрим следующую матрицу
\[
M=\begin{pmatrix}
\dfrac{\partial f_1}{\partial x_1} & \cdots & \dfrac{\partial f_1}{\partial x_n}\\
\vdots & \ddots & \vdots\\
\dfrac{\partial f_n}{\partial x_1} & \cdots & \dfrac{\partial f_n}{\partial x_n} \end{pmatrix}\]
состоящую из первых $n$ столбцов матрицы $\Jac_{(\bm{x},\bm{y})}\bm{f}$.
Предположим, что $M$ обратима в точке $(\bm{x},\bm{y})\in \Omega$, и $\bm{f}(\bm{x},\bm{y})=\bm{0}$.
Тогда существует гладкая функция $\bm{h}\:\Phi\z\to\mathbb{R}^n$, определённая на окрестности $\Phi\ni \bm{y}$ в $\mathbb{R}^m$
и окрестность $\Psi\ni (\bm{x},\bm{y})$ в $\mathbb{R}^n\times \mathbb{R}^m$, такие что
для любых $(x_1,\dots,x_n,y_1,\dots, y_m)\z\in \Psi$ равенство
\[\bm{f}(x_1,\dots,x_n,y_1,\dots, y_m)=\bm{0}\]
выполняется тогда и только тогда, когда
\[(x_1,\dots, x_n)=\bm{h}(y_1,\dots, y_m).\]

\end{thm}

\subsection*{Кратные интегралы}

Пусть $\bm{f}\:\mathbb{R}^n\to\mathbb{R}^n$ --- гладкое отображение (возможно, частично определённое).
Определим
\[\jac_{\bm{x}}\bm{f}\df|\det[\Jac_{\bm{x}}\bm{f}]|;
\index{10j@$\jac$ (определитель Якоби)}\]
то есть $\jac_{\bm{x}}\bm{f}$ --- модуль определителя матрицы Якоби $\bm{f}$ в точке~$\bm{x}$;
будем называть его \index{якобиан}\emph{якобианом}.

Следующая теорема предоставляет формулу замены переменных под кратным интегралом.

\index{борелевские подмножества}\emph{Борелевские подмножества} определяются как класс подмножеств, который порождается открытыми множествами с помощью последовательного применения следующих операций: счётного объединения, счётного пересечения и дополнения.
Так как дополнение замкнутого множества является открытым, и наоборот, эти множества могут также быть получены из всех замкнутых множеств.
Этот класс множеств включает практически все множества, которые естественным образом появляются в геометрии, но не включает патологические примеры, создающие проблемы с интегрированием.

\begin{thm}{Теорема}\label{thm:mult-substitution} 
Пусть $h\:K\to\mathbb{R}$ --- непрерывная функция на борелевском подмножестве $K\subset \mathbb{R}^n$.
Предположим, что $\bm{f}\:\Omega\to \mathbb{R}^n$ --- инъективное гладкое отображение, определённое на открытом множестве $\Omega\supset K$.
Тогда
\[\idotsint_{\bm{x}\in K} \!h(\bm{x})\cdot \jac_{\bm{x}}\bm{f}
=
\idotsint_{\bm{y}\in \bm{f}(K)} \!h\circ \bm{f}^{-1}(\bm{y}).\]

\end{thm}

\subsection*{Выпуклые функции}

Следующие утверждения потребуются только при $n\le 3$.

Функция $f\:\mathbb{R}^n\to \mathbb{R}$ называется \index{выпуклая!функция}\emph{(строго) выпуклой} (соответственно \index{вогнутая!функция}\emph{вогнутой}), если её эпиграф $z\z\ge f(\bm{x})$ (соответственно подграф $z\z\le f(\bm{x})$) является (строго) выпуклым множеством в $\mathbb{R}^n\times \mathbb{R}$.

Пусть $f\:\mathbb{R}^n\to \mathbb{R}$ --- гладкая функция (возможно, частично определённая).
Выберем вектор $\vec w\in \mathbb{R}^n$.
Для заданной точки $p\in\mathbb{R}^n$ рассмотрим функцию $\phi(t)=f(p+t\cdot \vec w)$.
Тогда производная $(D_{\vec w}f)(p)$ функции $f$ в точке $p$ \index{производная по направлению}\emph{по направлению} вектора $\vec w$ определяется как
\[(D_{\vec w}f)(p)=\phi'(0).\]

\begin{thm}{Теорема}\label{thm:Jensen}
Гладкая функция $f\:K\to \mathbb{R}$, определённая на выпуклом подмножестве $K\subset\mathbb{R}^n$, является выпуклой тогда и только тогда, когда выполняется одно из следующих эквивалентных условий:

\begin{subthm}{}
Вторая производная $f$ в любой точке по любому направлению неотрицательна, то есть
\[(D_{\vec w}^2f)(p)\ge 0\]
при любых $p\in K$ и $\vec w\in\mathbb{R}^n$.
\end{subthm}

\begin{subthm}{}
Так называемое \index{неравенство Йенсена}\emph{неравенство Йенсена}
\begin{align*}
f (&(1-t)\cdot x_0 + t\cdot x_1 ) \le
\\
&\le (1-t)\cdot f(x_0)+ t\cdot f(x_1)
\end{align*}
выполняется при любых $x_0,x_1\z\in K$ и $t\z\in[0,1]$.

\end{subthm}

\begin{subthm}{}
Для любых $x_0,x_1\in K$ выполняется неравенство
\[f \left (\frac{x_0 + x_1}2 \right ) \le \frac{f(x_0) + f(x_1)}2.\]
\end{subthm}

\end{thm}

\spell{\end{multicols}}{}

\section{Дифференциальные уравнения}

Следующий материал обсуждается в начале любого курса по обыкновенным дифференциальным уравнениям;
например в классическом учебнике Владимира Арнольда \cite{arnold}.

\spell{\begin{multicols}{2}}{}

\subsection*{Уравнения первого порядка}

Следующая теорема гарантирует существование и единственность решения задачи Коши для системы обыкновенных дифференциальных уравнений первого порядка
\[
\begin{cases}
x_1'&=f_1(x_1,\dots,x_n,t),
\\
&\,\,\vdots
\\
x_n'&=f_n(x_1,\dots,x_n,t),
\end{cases}
\]
где каждая функция $t\mapsto x_i=x_i(t)$ определена на вещественном интервале $\mathbb{J}$ и принимает вещественные значения, а каждая функция $f_i$ гладкая и определена на открытом подмножестве $\Omega\subset \mathbb{R}^n\times \mathbb{R}$.

Набор функций $(f_1,\dots,f_n)$ можно объединить в одну векторозначную функцию $\bm{f}\:\Omega\to \mathbb{R}^n$; аналогично, набор функций $(x_1,\dots,x_n)$ можно объединить в векторозначную функцию $\bm{x}\:\mathbb{J}\to\mathbb{R}^n$.
Таким образом, систему переписывается как одно векторное уравнение
\[\bm{x}'=\bm{f}(\bm{x}, t).\]

\begin{thm}{Теорема}\label{thm:ODE}
Пусть $\bm{f}\:\Omega\to \mathbb{R}^n$ --- гладкая функция, определённая на открытом подмножестве $\Omega\z\subset \mathbb{R}^n\z\times \mathbb{R}$.
Тогда для любых начальных данных $\bm{x}(t_0)=\bm{u}$ таких, что $(\bm{u},t_0)\in\Omega$, дифференциальное уравнение
\[\bm{x}'=\bm{f}(\bm{x},t)\]
имеет единственное решение $t\mapsto \bm{x}(t)$, определённое на максимальном интервале $\mathbb{J}$, содержащем $t_0$.
Более того,
\begin{enumerate}[(a)]
\item если $\mathbb{J}\ne \mathbb{R}$, то есть один из концов $\mathbb{J}$, скажем $b$, конечен, тогда $(\bm{x}(t),t)$ не имеет предельной точки в $\Omega$ при $t\to b$;%
\footnote{Другими словами, если $\bm{x}(t_n)$ сходится для последовательности $t_n\to b$, то её предел не лежит в $\Omega$.}

\item функция $w\:(\bm{u},t_0,t)\mapsto \bm{x}(t)$ имеет открытую область определения в $\Omega\times \mathbb{R}$, содержащую все точки вида $(\bm{u},t_0,t_0)$ при $(\bm{u},t_0)\in\Omega$,
и $w$ является гладкой функцией в этой области.
\end{enumerate}

\end{thm}

\subsection*{Высшие порядки}

Рассмотрим обыкновенное дифференциальное уравнение порядка~$k$
\[\bm{x}^{(k)}=\bm{f}(\bm{x},\bm{x}',\dots,\bm{x}^{(k-1)},t),\eqlbl{eq:nth-order}\]
где $t\mapsto\bm{x}=\bm{x}(t)$ --- функция, заданная на вещественном интервале, и принимающая значения в $\mathbb{R}^n$.

Это уравнение можно переписать как $k$ уравнений первого порядка, добавив $(k-1)$-у новую векторную переменную
$\bm{y}_1\z=\bm{x}'$,
$\bm{y}_2\z=\bm{x}'',\z\dots,\bm{y}_{k-1}\z=\bm{x}^{(k-1)}$:
\[
\begin{cases}
\bm{x}'(t)&=\bm{y}_1(t),
\\
\bm{y}_1'(t)&=\bm{y}_2(t),
\\
&\,\,\vdots
\\
\bm{y}_{k-2}'(t)&=\bm{y}_{k-1}(t),
\\
\bm{y}_{k-1}'(t)&=\bm{f}(\bm{x},\bm{y}_{1},\dots,\bm{y}_{k-1},t).
\end{cases}
\eqlbl{eq:nth-order-new}
\]

Итак, мы вывели следующую теорему.

\begin{thm}{Теорема}\label{thm:ODE-nth-order}
Уравнение $k$-го порядка \ref{eq:nth-order} эквивалентно системе \ref{eq:nth-order-new} из $k$ уравнений первого порядка.
\end{thm}

Этот трюк сводит вопросы существования и единственности для обыкновенных дифференциальных уравнений высших порядков к соответственным вопросам для уравнений первого порядка.
То есть теорема \ref{thm:ODE} обобщается на уравнения высших порядков, надо только предположить, что $\Omega\z\subset \mathbb{R}^{n\cdot k}\times \mathbb{R}$, а начальные данные состоят из $\bm{x}(t_0)$, $\bm{x}'(t_0),\dots,\bm{x}^{(k-1)}(t_0)$.

\spell{\end{multicols}}{}
}

}

\backmatter

\arxiv{\newgeometry{top=0.83in, bottom=0.83in,inner=0.52in, outer=0.52in}}{\newgeometry{top=0.955in, bottom=0.955in,inner=0.52in, outer=0.645in}}

{
\scriptsize

\stepcounter{chapter}
\setcounter{eqtn}{0}

\spell{\begin{multicols}{2}}{}

\chapter[Подсказки]{\vspace*{-15mm} Подсказки\vspace*{-10mm}}


\stepcounter{chapter}
\setcounter{eqtn}{0}

\parbf{\ref{ex:9}}; \ref{SHORT.ex:9:compact}.
Воспользуйтесь тем, что непрерывное инъективное отображение, определённое на компактном множестве, является вложением (\ref{thm:Hausdorff-compact}).

\parit{\ref{SHORT.ex:9:9}.}
Образ $\gamma$ может выглядеть как цифра $8$ или $9$.

\parbf{\ref{ex:mono}.}
Сдавите окрестности концов в их концы.

\parbf{\ref{aex:simple-curve}.}
Пусть $\alpha$ --- путь, соединяющий $p$ с~$q$.
Можно предположить, что $\alpha(t)\ne p,q$ для $t\ne0,1$.

Назовём открытое множество $\Omega$ в $(0,1)$ {}\emph{подходящим},
если $\alpha(a)=\alpha(b)$ для любой компоненты связности $(a,b)$ в $\Omega$.
Покажите, что объединение вложенных подходящих множеств, является подходящим.
Следовательно, найдётся максимальное подходящее множество~$\hat \Omega$.

Определим $\beta(t)=\alpha(a)$ для любого $t$ в компоненте связности $(a,b)\subset\hat \Omega$, и $\beta (t) \z= \alpha (t) $ для $t\notin\hat{\Omega}$.
Заметьте, что для любого $x\in [0,1]$ множество $\beta^{-1}\{\beta(x)\}$ связно.

Остаётся получить простой путь, репараметризовав $\beta$.
Для этого нужно построить такую монотонную функцию $\tau\:[0,1]\z\to[0,1]$, что 
$\tau(t_1)\z=\tau(t_2)$ тогда и только тогда, когда существует компонента связности $(a,b)\subset\hat \Omega$ такая, что $t_1,t_2\z\in [a,b]$.

Функция $\tau$ похожа на так называемую {}\emph{чёртову лестницу};
разберитесь с её построением и воспользуйтесь им, чтобы построить $\tau$.

\parbf{\ref{ex:L-shape}.}
Обозначим объединение двух полуосей через~$L$.

Заметьте, что $f(t)\to\infty$ при $t\to \infty$.
Так как $f(0)=0$, по теореме о промежуточных значениях, $f(t)$ принимает все неотрицательные значения при $t\ge 0$.
Воспользуйтесь этим, чтобы показать, что $L$ есть образ~$\alpha$.

Далее, покажите, что функция $f$ строго возрастает при $t> 0$.
Выведите отсюда, что отображение $t\mapsto \alpha(t)$ инъективно.

Итак, $\alpha$ --- гладкая параметризация~$L$.

Теперь начнём с произвольной гладкой параметризации~$L$, скажем $\beta\:t\z\mapsto (x(t),y(t))$.
Можно считать, что $x(0)\z=y(0)=0$.
Заметим, что $x(t)\ge 0$ для любого $t$, следовательно, $x'(0)=0$.
Точно также получаем, что $y'(0)\z=0$.
То есть, $\beta'(0)=0$.
Значит, $L$ не допускает гладкой \textit{регулярной} параметризации.

\parbf{\ref{ex:cycloid}.}
Примените определения.
В \ref{SHORT.ex:cycloid:regular} нужно проверить, что $\gamma'_\ell\ne 0$.
В \ref{SHORT.ex:cycloid:simple} нужно проверить, что $\gamma_\ell(t_0)\z=\gamma_\ell(t_1)$ только если $t_0=t_1$.

\parbf{\ref{ex:nonregular}.}
Заметьте, что параметризация $t\mapsto (t,t^3)$ гладкая и регулярная.
Измените её так, чтобы скорость в одной из точек была равна нулю.

\begin{wrapfigure}{r}{20 mm}
\vskip-2mm
\centering
\includegraphics{mppics/pic-270}
\vskip-2mm
\end{wrapfigure}

\parbf{\ref{ex:y^2=x^3}.}
Это \index{полукубическая парабола}\emph{полукубическая парабола}; она показана на рисунке.
Попробуйте рассуждать аналогично \ref{ex:L-shape}.

\parbf{\ref{ex:viviani}.}
При $\ell=0$ система описывает пару точек $(0,0,\pm1)$, так что можно считать, что $\ell\z\ne 0$.

Первое уравнение описывает единичную сферу с центром в начале координат, а второе --- цилиндр над окружностью в плоскости $(x,y)$ с центром в точке $(-\tfrac\ell2,0)$ и радиусом~$|\tfrac\ell2|$.

\begin{Figure}
\centering
\vskip-0mm
\begin{lpic}[t(2mm),b(0mm),r(0mm),l(0mm)]{asy/viviani(1)}
\lbl[r]{-.5,18;$x$}
\lbl[l]{41,22;$y$}
\lbl[r]{18,54;$z$}
\end{lpic}
\end{Figure}

Найдите градиенты $\nabla f$ и $\nabla h$ функций
\begin{align*}
 f(x,y,z)&=x^2+y^2+z^2-1,
 \\
 h(x,y,z)&=x^2+\ell\cdot x+y^2.
\end{align*}
Покажите, что при $\ell\ne 0$ градиенты линейно зависимы только на оси $x$.
Выведите отсюда, что для $\ell\ne\pm 1$ каждая компонента связности множества решений является гладкой кривой.

Покажите, что 
\begin{itemize}
\item если $|\ell|<1$, то множество имеет две связные компоненты с $z>0$ и $z<0$.
\item если же $|\ell|\ge1$, то множество связно.
\end{itemize}

Линейная независимость градиентов предоставляет лишь достаточное условие.
Поэтому случай $\ell=\pm1$ придётся проверить руками.
В этом случае окрестность точки $(\pm1,0,0)$ не допускает гладкой регулярной параметризации --- попробуйте это доказать.
Случай $\ell=1$ показан на рисунке.

\parit{Замечание.}
При $\ell=\pm1$ получается так называемая \index{кривая Вивиани}\emph{кривая Вивиани}.
Она допускает следующую гладкую регулярную параметризацию с самопересечением в точке $(\pm1,0,0)$
\[t\mapsto(\pm(\cos t)^2,\cos t\cdot\sin t,\sin t).\]

\parbf{\ref{ex:open-curve}.}
Допустим, что $|\gamma(t)|\z\to\infty$ при $t\to\pm\infty$.
Выберем компактное множество $K\subset \mathbb{R}^3$.
Покажите, что прообраз $\gamma^{-1}(K)$ ограничен и замкнут в $\mathbb{R}$.
Далее примените лемму Гейне --- Бореля (\ref{thm:Heine--Borel}).
Выведите, что $\gamma$ собственная.

Теперь предположим, что $\gamma(t_n)$ сходится для какой-то последовательности $t_n\to \pm \infty$; пусть $p$ --- её предел и $K$ --- замкнутый шар с центром в~$p$.
Убедитесь, что прообраз $\gamma^{-1}(K)$ некомпактен.
Выведите отсюда, что $\gamma$ несобственная.

\parbf{\ref{ex:proper-closed}.}
Докажите, что множество $C\subset \mathbb{R}^3$ является замкнутым тогда и только тогда, когда пересечение $K\cap C$ компактно для любого компактного $K\subset \mathbb{R}^3$,
и воспользуйтесь этим.

\parbf{\ref{ex:proper-curve}.}
Не умаляя общности, можно считать, что начало координат не лежит на кривой.

Покажите, что инверсия плоскости $(x,y)\z\mapsto (\tfrac{x}{x^2+y^2},\tfrac{y}{x^2+y^2})$ отображает нашу кривую в замкнутую кривую с удалённым началом координат.
Примените теорему Жордана для полученной кривой и снова воспользуйтесь инверсией.


\stepcounter{chapter}
\setcounter{eqtn}{0}

\parbf{\ref{ex:integral-length-0}.}
Покажите, что если взять точную верхнюю грань в \ref{def:length} для всех последовательностей
$a=t_0\le t_1\le\z\dots\le t_k=b$, то результат будет тем же самым.

Предположим, что $\gamma_2$ --- репараметризация $\gamma_1$ с помощью $\tau\:[a_1,b_1]\to [a_2,b_2]$; не умаляя общности, можно считать, что $\tau$ неубывающая.
Пусть $\theta_i=\tau(t_i)$.
Заметьте, что $a_2=\theta_0\z\le\theta_1\le \z\dots\le\theta_k=b_2$ тогда и только тогда, когда 
$a_1=t_0\le t_1\le\z\dots\le t_k=b_1$.
Сделайте последний шаг.

\parbf{\ref{ex:length-chain}.}
Покажите, что для любой вписанной ломаной $\beta$ и любого $\epsilon>0$ выполняется
\[\length\beta_n>\length\beta-\epsilon\]
при всех достаточно больших $n$.
Выведите отсюда, что
\[\liminf_{n\to\infty}\length\beta_n\ge \length \gamma.\]
Используя определение длины, покажите, что 
\[\limsup_{n\to\infty}\length\beta_n\le \length \gamma.\]
Убедитесь, что из полученных неравенств следует всё, что надо.

\parbf{\ref{ex:length-image}.}
Для данного разбиения $0=t_0<\z\dots <t_n\z=1$ отрезка $[0,1]$, положим $\tau_0=0$ и 
\[\tau_i=\max\set{\tau \in[0,1]}{\beta(\tau)=\gamma(t_i)}\]
при $i>0$.
Покажите, что $(\tau_i)$ является разбиением $[0,1]$;
то есть $0=\tau_0<\tau_1<\z\dots<\tau_n=1$.

По построению 
\begin{align*}
&|\gamma(t_0)-\gamma(t_1)|+|\gamma(t_1)-\gamma(t_2)|+\dots
\\
&\qquad\qquad\dots+|\gamma(t_{n-1})-\gamma(t_n)|=
\\
&=
|\beta(\tau_0)-\beta(\tau_1)|+|\beta(\tau_1)-\beta(\tau_2)|+\dots
\\
&\qquad\qquad\dots+|\beta(\tau_{n-1})-\beta(\tau_n)|.
\end{align*}
Поскольку разбиение $(t_i)$ произвольно, получаем 
\[\length \beta\ge \length \gamma.\]

\parit{Замечания.}
Полезно сравнить это упражнение с \ref{obs:S2-length}.

Неравенство может оказаться строгим.
Такое происходит, если $\beta$ ходит вдоль $\gamma$ туда-сюда.
В этом случае разбиение $(\tau_i)$ выше нельзя выбрать произвольно.

В предположении $\beta([0,1])\z\supset\gamma([0,1])$, задача решается применением следующего неравенства \cite[2.6.1+2.6.2]{burago-burago-ivanov}:
\[h\le \length\gamma,\] где $h$ --- одномерная мера Хаусдорфа образа $\gamma([0,1])$.
Более того, равенство достигается тогда и только тогда, когда $\gamma$ простая.

\parbf{\ref{ex:integral-length}}; \ref{SHORT.ex:integral-length>}.
Примените формулу Ньютона --- Лейбница к каждому отрезку разбиения.

\parit{\ref{SHORT.ex:integral-length<}.}
Рассмотрите разбиение, для которого вектор скорости $\alpha'(t)$ почти постоянен на каждом из отрезков.

\parbf{\ref{adex:integral-length}.}
Воспользуйтесь теоремами Радемахера и Лузина (\ref{thm:rademacher} и \ref{thm:lusin}).

\parbf{\ref{ex:nonrectifiable-curve}}; \ref{SHORT.ex:nonrectifiable-curve:a}.
Посмотрите на рисунок и догадайтесь как запараметризовать дугу снежинки отрезком $[0,1]$.
Расширьте параметризацию на всю снежинку.
Убедитесь, что она действительно описывает вложение окружности в плоскость.

\begin{Figure}
\vskip-0mm
\centering
\includegraphics{mppics/pic-226}
\vskip0mm
\end{Figure}

\parit{\ref{SHORT.ex:nonrectifiable-curve:b}.}
Пусть $\gamma\:[0,1]\to\mathbb{R}^2$ --- спрямляемая кривая, и $\gamma_k$ --- её гомотетия с коэффициентом $k>0$;
то есть $\gamma_k(t)=k\cdot\gamma(t)$ для любого~$t$.
Покажите, что 
\[\length\gamma_k=k\cdot\length \gamma.\]

Теперь допустим, что дуга $\gamma$ снежинки Коха, показанная на рисунке, спрямляема,
и $\ell$ --- её длина.
Заметьте, что $\gamma$ можно разделить на $4$ дуги, каждая из которых конгруэнтна гомотетии $\gamma$ с коэффициентом~$\tfrac13$.
Следовательно, $\ell=\tfrac43\cdot\ell$,
но ведь $\ell>0$ --- противоречие.

\parbf{\ref{ex:cont-length}.}
Применив \ref{thm:length-semicont}, покажите, что $s$ полунепрерывна снизу;
то есть если $t_n\to t_\infty$ при $n\to\infty$, тогда 
\[\liminf_{n\to\infty} s(t_n)\ge s(t_\infty).\]
Убедитесь, что
\[s(t)=s(b)-\length(\gamma|_{[t,b]}).\]
Применив это равенство с \ref{thm:length-semicont}, покажите, что $s$ полунепрерывна сверху.
Выведите отсюда, что функция $s$ непрерывна.
Наконец, покажите, что $s$ неубывающая, 
и воспользуйтесь этим.

\parbf{\ref{ex:arc-length-helix}.} 
Можно предположить, что $a\ne 0$ или $b\ne 0$,
иначе задача тривиальна.

Покажите, что $|\gamma'(t)|\equiv \sqrt{a^2+b^2}$;
в частности, скорость постоянна.
Значит, $s\z=t/\sqrt{a^2+b^2}$ --- параметр длины.
Остаётся подставить $s\cdot \sqrt{a^2+b^2}$ вместо~$t$.

\parbf{\ref{ex:convex-hull}.}
Пусть $p_1\dots p_n$ --- замкнутая ломаная, вписанная в $\beta$.
По \ref{cor:convex=>rectifiable}, можно считать, что её длина сколь угодно близка к длине $\beta$;
то есть
\[\length (p_1\dots p_n)>\length\beta-\epsilon\]
для любого наперёд заданного $\epsilon>0$.

Убедитесь, что при всём при этом можно предположить, что каждая точка $p_i$ лежит на $\alpha$.

Поскольку $\alpha$ простая, точки $p_1,\dots,p_n$ появляются на $\alpha$ в одном и том же циклическом порядке;
то есть ломаная $p_1\dots p_n$ также вписана в $\alpha$.
В частности,
\[\length\alpha\ge \length (p_1\dots p_n).\]
Следовательно, 
\[\length\alpha>\length\beta-\epsilon\]
для любого $\epsilon>0$.
Отсюда
\[\length\alpha\ge\length\beta.\]

\begin{wrapfigure}{r}{25 mm}
\vskip-0mm
\centering
\includegraphics{mppics/pic-275}
\vskip0mm
\end{wrapfigure}

Если у $\alpha$ есть самопересечения, то точки $p_1,\dots, p_n$ могут располагаться на $\alpha$ в другом циклическом порядке, скажем, $p_{i_1},\dots,p_{i_n}$.
Примените неравенство треугольника, чтобы показать, что
\[\length(p_{i_1}\dots p_{i_n})\ge \length (p_1\dots p_n)\]
и воспользуйтесь этим, чтобы обобщить доказательство.

\parbf{\ref{ex:convex-croftons}.} 
Обозначим через $\ell_{\vec u}$ отрезок, 
полученный ортогональной проекцией $\gamma$ на прямую в направлении ${\vec u}$.
Поскольку $\gamma_{\vec u}$ проходит вдоль $\ell_{\vec u}$ туда и обратно, получаем
\[\length\gamma_{\vec u}\ge 2\cdot\length\ell_{\vec u}.\]
По формуле Крофтона, 
\[\length\gamma\ge \pi\cdot \overline{\length\ell_{\vec u}}.\]

В случае равенства кривая $\gamma_{\vec u}$ проходит точно туда и обратно вдоль $\ell_{\vec u}$ без дополнительных зигзагов;
это должно происходить для почти всех (а следовательно, и для всех) направлений~${\vec u}$.

Пусть $K$ --- замкнутое множество, ограниченное~$\gamma$.
Последнее утверждение означает, что если прямая пересекает $K$, то пересечение отрезок или точка.
Отсюда следует, что $K$ выпукло.

\parbf{\ref{adex:more-croftons}.}
Доказательство то же, что у обычной формулы Крофтона.
Чтобы найти коэффициенты, достаточно проверить её на единичном интервале,
и это делается интегрированием
\begin{align*}
\frac1{k_1}&=\frac{1}{\area \mathbb{S}^2}\cdot\iint_{\mathbb{S}^2} |x|;
\\
\frac1{k_2}&=\frac{1}{\area \mathbb{S}^2}\cdot\iint_{\mathbb{S}^2} \sqrt{1-x^2}.
\end{align*}
Ответы: $k_1=2$ и $k_2=\tfrac4\pi$.

\parbf{\ref{ex:induced-is-length}.}
Пусть $d$ — исходная метрика на $\spc{X}$.
Пусть $\gamma\:[a,b]\to (\spc{X},\ell)$ непрерывно.
Докажите, что $\gamma\:[a,b]\to (\spc{X},d)$ также непрерывно.

Выберем кривую $\gamma$ в $(\spc{X},\ell)$.
Пусть $r$ и $s$ --- её $d$-длина и $\ell$-длина соответственно.
Докажите, что $s\ge r\ge s-\epsilon$ для любого $\epsilon>0$.
Выведите отсюда, что $d$-длина и $\ell$-длина любой кривая в $(\spc{X},\ell)$ равны.
Сделайте ещё шаг.

\parit{Замечание.} Обратное к первому утверждению в общем случае неверно.
Рассмотрите, например, следующую метрику на верхней полуплоскости:
\[d(x,y)=\min\{|x|+|y|,\bigl||x|-|y|\bigr|+\sqrt{\measuredangle(x,y)}\}.\]

\parbf{\ref{ex:intrinsic-convex}.}
Необходимость очевидна.
Чтобы доказать достаточность, допустим, что $A$ не выпукло;
то есть найдутся точки $x,y\in A$ и точка $z\notin A$, которая лежит между $x$ и~$y$.

Так как $A$ замкнуто, его дополнение открыто.
Иными словами, существует шар $B(z,\epsilon)$, который не пересекает $A$ для некоторого $\epsilon>0$.

Покажите, что найдётся такое $\delta>0$, что любая кривая из $x$ в $y$ длины не более $\dist{x}{y}{\mathbb{R}^3}+\delta$ задевает $B(z,\epsilon)$.
Отсюда $\dist{x}{y}A\z\ge \dist{x}{y}{\mathbb{R}^3}+\delta$; 
в частности, $\dist{x}{y}A\ne \dist{x}{y}{\mathbb{R}^3}$.

\begin{wrapfigure}{r}{23 mm}
\vskip-0mm
\centering
\includegraphics{mppics/pic-280}
\vskip0mm
\end{wrapfigure}

\parbf{\ref{ex:antipodal}.}
Сферическая кривая, показанная на рисунке, не имеет антиподальных пар точек.
Однако на одной из её сторон лежат три точки $x,y,z$ с одного экватора, а на другой --- их антиподы $-x,-y,-z$.
(Мы предполагаем, что точки $x,-y,z,-x,y,-z$ лежат в этом же порядке на экваторе.)

Покажите, что такая кривая не может лежать ни в одной полусфере.

\parbf{\ref{ex:bisection-of-S2}.}
Допустим, что $\gamma$ лежит в открытой полусфере (\ref{lem:hemisphere}).
В частности, она не делит $\mathbb{S}^2$ на две области равной площади --- противоречие.

\parbf{\ref{ex:flaw}.}
Уже первое предложение неверно --- \textit{недостаточно} показать, что диаметр не превышает~$2$.
Например, если у равностороннего треугольника радиус описанной окружности чуть больше $1$,
то его диаметр (который определяется как максимальное расстояние между его точками) будет чуть больше $\sqrt3$, так что он меньше $2$.

С другой стороны, доказательство леммы о полусфере (\ref{lem:hemisphere}) можно приспособить для получения правильного решения.
А именно: (1) выберите две точки $p$ и $q$ на $\gamma$, которые делят её на две дуги одинаковой длины;
(2) пусть $z$ --- средина между $p$ и $q$;
(3) покажите, что $\gamma$ лежит в единичном круге с центром в~$z$.

\parbf{\ref{adex:crofton}}; \ref{SHORT.adex:crofton:crofton}.
Рассуждение схоже с доказательством обычной формулы Крофтона (\ref{sec:crofton}).

\parit{\ref{SHORT.adex:crofton:hemisphere}.}
Предположим, что $\length \gamma<2\cdot\pi$.
Согласно \ref{SHORT.adex:crofton:crofton},
\[\overline{\length \gamma^*_{\vec u}}<2\cdot\pi.\]
Следовательно, 
\[\length \gamma^*_{\vec u}<2\cdot\pi\]
для некоторого ${\vec u}$. 

Заметьте, что $\gamma^*_{\vec u}$ идёт по некоторой полуокружности~$h$.
Следовательно, $\gamma$ лежит в полусфере, построенной на $h$ как на диаметре.


\stepcounter{chapter}
\setcounter{eqtn}{0}

\parbf{\ref{ex:zero-curvature-curve}.}
Пусть кривая $\gamma$ имеет единичную скорость и нулевую кривизну. 
Тогда $\gamma''\equiv 0$, и, значит, вектор скорости $\vec v=\gamma'$ не меняется.
Отсюда, $\gamma(t)\z=p+(t-t_0)\cdot \vec v$, где $p=\gamma(t_0)$.

\parbf{\ref{ex:scaled-curvature}.} 
Заметьте, что $\alpha(t)\df\gamma_{\lambda}(t/\lambda)$ --- это параметризация кривой $ \gamma_{\lambda}$ с единичной скоростью,
и примените дважды правило дифференцирования сложной функции.

\parbf{\ref{ex:curvature-of-spherical-curve}.}
Продифференцируйте тождество $\langle\gamma(s),\gamma(s)\rangle=1$ несколько раз.

\parbf{\ref{ex:curvature-formulas}.} 
Пусть $\tan=\tfrac{\gamma'}{|\gamma'|}$.
Выведите  тождества
\begin{align*}
\gamma''-(\gamma'')^\perp&=\tan\cdot\langle\gamma'',\tan\rangle,
&
|\gamma'|&=\sqrt{\langle \gamma',\gamma'\rangle},
\end{align*}
и воспользуйтесь ими.

\parbf{\ref{ex:curvature-graph}.} 
Примените \ref{ex:curvature-formulas:a} к параметризации $t\z\mapsto (t,f(t))$.

\parbf{\ref{ex:approximation-const-curvature}.}
Не умаляя общности, можно предположить, что у $\gamma$ единичная скорость.

Заметьте, что касательная индикатриса $\tan(s)\z=\gamma'(s)$ --- сферическая кривая и $|\tan'|\le 1$.
Воспользуйтесь этим для построения последовательности сферических кривых с единичной скоростью $\tan_n\:\mathbb{I}\to\mathbb{S}^2$ таких, что $\tan_n(s)\to \tan(s)$ при $n\to\infty$ для любого~$s$.
Покажите, что следующая последовательность кривых даёт решение
\[\gamma_n(s)=\gamma(a)+\int_a^s\tan_n(t)\cdot dt.\]

\parbf{\ref{ex:no-parallel-tangents}.}
Воспользуйтесь построением в \ref{ex:antipodal}, дабы получить замкнутую гладкую простую сферическую кривую $\tan\:\mathbb{S}^1\z\to\mathbb{S}^2$, которая не содержит пар антиподальных точек и имеет нулевое среднее значение.
Затем постройте кривую с касательной индикатрисой $\tan$.

\parit{Источник:}
Эту кривую построил Беньям\'{и}но С\'{е}гре \cite{segre}.

\parbf{\ref{ex:helix-curvature}.}
Покажите, что $\gamma_{a,b}''\perp \gamma'_{a,b}$, и примените \ref{ex:curvature-formulas:a}.

\parbf{\ref{ex:length>=2pi}.}
Примените теорему Фенхеля.

\parbf{\ref{ex:gamma/|gamma|}.}
Можно предположить, что у $\gamma$ единичная скорость.
Тогда $\langle\tan,\tan\rangle=1$ и  $\tan'\perp \tan$.
Пусть $\theta(s)=\measuredangle(\gamma(s),\gamma'(s))$ и, значит, $\langle \tan,\sigma\rangle=\cos\theta$.
Поэтому
\begin{align*}
\kur\cdot \sin\theta
&=|\tan'|\cdot \sin\theta\ge
-\langle \tan',\sigma\rangle=
\\
&=
\langle \tan,\sigma'\rangle-\langle \tan,\sigma\rangle'=
\\
&=
(|\sigma'|+\theta')\cdot \sin\theta.
\end{align*}
Следовательно, $\kur\ge |\sigma'|+\theta'$,
если $\theta\ne0,\pi$.
Остаётся проинтегрировать это неравенство и показать, что множество, для которого $\theta\z\ne0$ или $\pi$, не создаёт проблем.

\parit{Другое решение} строится на основе \ref{prop:inscribed-total-curvature}.

\parbf{\ref{ex:DNA}}; \ref{SHORT.ex:DNA:c''c>=k}.
Поскольку $|\gamma(s)|\le 1$, 
\begin{align*}
\langle\gamma''(s),\gamma(s)\rangle&\ge -|\gamma''(s)|\cdot|\gamma(s)|\ge-\kur(s)
\end{align*}
при всех~$s$.

\parit{\ref{SHORT.ex:DNA:int>=length-tc}.}
Поскольку $\gamma$ параметризована длиной, $\ell\z=\length\gamma$ и $\langle\gamma',\gamma'\rangle\equiv1$.
Следовательно,
\begin{align*}
\int_0^\ell\langle\gamma(s),\gamma'(s)\rangle'\cdot ds&=\\
=\int_0^\ell\langle\gamma'(s),\gamma'(s)\rangle\cdot ds&+\int_0^\ell\langle\gamma(s),\gamma''(s)\rangle\cdot ds\ge
\\
&\ge \length\gamma-\tc\gamma.
\end{align*}

\parit{\ref{SHORT.ex:DNA:end}.}
Согласно формуле Ньютона --- Лейбница,
\begin{align*}
\int_0^\ell\langle\gamma(s),\gamma'(s)\rangle'\cdot ds
&=\langle\gamma(s),\gamma'(s)\rangle\bigg|_0^\ell.
\end{align*}
Поскольку $\gamma(0)=\gamma(\ell)$ и $\gamma'(0)=\gamma'(\ell)$, правая часть равенства обращается в ноль.

\columnbreak

Можно считать, что кривая в \ref{thm:DNA} описывается петлёй $\gamma\:[0,\ell]\to\mathbb{R}^3$, параметризованной длиной и ещё, что центр шара находится в начале координат, то есть $|\gamma|\le 1$.
Поскольку $\gamma$ гладкая и замкнутая, 
$\gamma'(0)=\gamma'(\ell)$ и $\gamma(0)=\gamma(\ell)$.
Следовательно, \ref{SHORT.ex:DNA:int>=length-tc} и \ref{SHORT.ex:DNA:end} влекут \ref{thm:DNA}.

\parbf{\ref{ex:tangent-support}.}
Покажите, что никакая прямая, отличная от $\ell$, не может подпирать $F$ в точке~$p$. 
Примените \ref{lem:separation}, чтобы показать, что какая-то прямая подпирает $F$ в точке $p$.
Сделайте вывод.

{

\begin{wrapfigure}{r}{22 mm}
\vskip-6mm
\centering
\includegraphics{mppics/pic-255}
\vskip0mm
\end{wrapfigure}

\parbf{\ref{ex:anti-bow}.}
Начнём с кривой $\gamma_1$, изображённой на рисунке.
Чтобы получить $\gamma_2$, слегка распрямите пунктирную дугу (то есть уменьшите её кривизну).

}

\parbf{\ref{ex:bow'}}; \ref{SHORT.ex:bow'+}.
Не умаляя общности, можно предположить, что $p=\gamma_1(a)\z=\gamma_2(a)$ --- начало координат, 
вектор $\vec u=\tan_1(a)\z=\tan_2(a)$ направлен вдоль оси $x$,
а точки $q_1=\gamma_1(b)$ и $q_2=\gamma_2(b)$ лежат в верхней полуплоскости $(x,y)$.

Допустим, что $\alpha_1<\alpha_2$.
Найдите такую точку $x$, что
$\dist{x}{q_1}{}<\dist{x}{q_2}{}$, и
$\gamma_1$ вкупе с отрезками $[p,x]$ и $[x,q_1]$ ограничивает выпуклую область в $(x,y)$-плоскости.
(Здесь понадобится, что $\beta_1\le\tfrac\pi2$.)

Проверьте, что доказательство леммы о луке (\ref{lem:bow}) работает для произведений $[x,p]$ с $\gamma_1$ и $\gamma_2$, и придите к противоречию.

\parit{\ref{SHORT.ex:bow'-}.}
Посмотрите на рисунок и подумайте.

\begin{Figure}
\vskip-0mm
\centering
\includegraphics{mppics/pic-257}
\vskip0mm
\end{Figure}

\parit{Замечание.}
Это упражнение полезно сравнить с леммой о хорде (\ref{lem:chord}).

\parbf{\ref{ex:length-dist}}; \ref{SHORT.ex:length-dist:>}.
Выберите значение $s_0\in[a,b]$, которое делит полную кривизну $\gamma$ пополам.
Заметьте, что $\measuredangle(\gamma'(s_0),\gamma'(s))\le \theta$ для любого~$s$.
Используйте это неравенство так же, как в доказательстве леммы о луке.

\begin{wrapfigure}[4]{r}{22 mm}
\vskip-3mm
\centering
\includegraphics{mppics/pic-290}
\vskip-0mm
\end{wrapfigure}

\parit{\ref{SHORT.ex:length-dist:self-intersection:>pi}.}
Воспользуйтесь \ref{SHORT.ex:length-dist:>} и посмотрите на рисунок.

\parit{\ref{SHORT.ex:length-dist:=}.}
Начните с ломаной из двух звеньев равной длины и с внешним углом $2\cdot\theta$ и сгладьте вершину, используя срезки и сглаживания (\ref{sec:analysis}).

\parbf{\ref{ex:schwartz}.}
Пусть $\ell=\length\gamma$.
Допустим, что $\ell_1\z<\ell<\ell_2$, и пусть $\gamma_1$ --- дуга единичной окружности длины $\ell$.

Покажите, что расстояние между концами $\gamma_1$ больше, чем $|p-q|$, и примените лемму о луке (\ref{lem:bow}).

\parit{Источник:} Это утверждение приписывается Герману Шварцу \cite{shur}.

\parbf{\ref{ex:loop}.}
Допустим, что $\length\gamma<2\cdot\pi$. 
Примените лемму о луке (\ref{lem:bow}) к $\gamma$ и дуге единичной окружности той же длины.

\parbf{\ref{ex:bow-upper}.}
Выберите гладкую сферическую кривую $\alpha$, которая всё время идёт близко к одной точке;
например, маленькую сферическую окружность с центром в этой точке.
Рассмотрите кривую $\tan$, которая идёт по $\alpha$ со скоростью $\kappa(s)$ при любом $s\in [0,\ell]$.
Покажите, что кривая с касательной индикатрисой $\tan$ решает задачу.

\parbf{\ref{ex:gromov-twist}.}
Можно предположить, что $\gamma$ имеет единичную скорость.
Допустим, что неравенство не выполняется при $t=0$.
Можно предположить, что $\alpha(0)\le\tfrac\pi2$;
в противном случае развернём параметризацию.

В плоскости, натянутой на вектора $\gamma(0)$ и $\gamma'(0)$, выберите дугу окружности (или отрезок прямой) $\sigma$ с единичной скоростью от $0$ до $\gamma(0)$, которая приходит в $\gamma(0)$ в направлении, противоположном $\gamma'(0)$.
Рассмотрите полуокружность $\tilde\gamma$ с единичной скоростью и кривизной $2$, которая начинается в $\gamma(0)$ в направлении $\gamma'(0)$ так, что произведение $\sigma*\tilde\gamma$ является дугой выпуклой плоской кривой; см. рисунок.

\begin{Figure}
\vskip-1mm
\centering
\includegraphics{mppics/pic-282}
\vskip-1mm
\end{Figure}

Покажите, что если $|\gamma(0)|> \sin(\alpha(0))$, то $\tilde\gamma$ выходит за пределы единичного шара; то есть $|\tilde\gamma(t_0)|>1$ для некоторого $t_0$.

Произведения $\sigma*\gamma$ и $\sigma*\tilde\gamma$ не являются гладкими на стыке, но они дифференцируемы в этой точке.
Проверьте, что доказательство леммы о луке работает и в этом более общем случае.

Применив этот вариант леммы о луке к $\sigma*\gamma$ и $\sigma*\tilde\gamma$, получите, что $|\gamma(t_0)|\ge |\tilde\gamma(t_0)|$, и придите к противоречию.

\parit{Источник:} Это утверждение использовалось первым автором \cite{petrunin2023}.


\stepcounter{chapter}
\setcounter{eqtn}{0}

\parbf{\ref{ex:chord-lemma-optimal}.} 
Пусть $\alpha=\measuredangle(\vec u,\vec w)$ 
и $\beta=\measuredangle(\vec w,\vec v)$.
Попробуйте догадаться до примера по рисунку.

\begin{Figure}
\vskip-1mm
\centering
\includegraphics{mppics/pic-285}
\vskip-1mm
\end{Figure}

Показанная кривая разделена на три дуги: I, II и III. 
Дуга I поворачивает от $\vec u$ к $\vec w$;
её полная кривизна равна $\alpha$.
Аналогично, дуга III поворачивает от $\vec w$ к $\vec v$ и имеет полную кривизну $\beta$. 
Дуга II проходит очень близко и почти параллельно хорде $pq$, и её полную кривизну можно сделать произвольно малой.

\parbf{\ref{ex:monotonic-tc}.}
Используйте, что внешний угол треугольника равен сумме внутренних углов в двух других вершинах.
Во второй части примените индукцию по числу вершин.

\parbf{\ref{ex:sef-intersection}.}
Допустим, что $x$ --- это точка самопересечения.
Покажите, что можно выбрать такие две точки $y$ и $z$ на $\gamma$ между самопересечениями, чтобы треугольник $xyz$ оказался невырожденным.
В частности, 
$\measuredangle\hinge yxz
\z+
\measuredangle\hinge zyx
<\pi$, или, что эквивалентно, $\tc{xyzx}\z>\pi$ для \textit{незамкнутой} вписанной ломаной $xyzx$.
Остаётся применить~\ref{prop:inscribed-total-curvature}.

\parbf{\ref{ex:quadrisecant}.}
Рассмотрим замкнутую ломаную $acbd$.
Убедитесь, что $\tc{acbd}\z=4\cdot\pi$, и примените \ref{prop:inscribed-total-curvature}.

\parit{Замечания.}
Прямые, пересекающие кривую, как в этом упражнении, называются \index{четырёхкратная секущая}\emph{альтернированными четырёхкратными секущими}.
Известно, что такая секущая найдётся у любого {}\emph{нетривиального узла} \cite{denne};
согласно упражнению, это доказывает так называемую {}\emph{теорему Фари --- Милнора} --- \textit{полная кривизна любого узла больше, чем~$4\cdot \pi$}; см. \cite{petrunin-stadler} и тамошние ссылки.

\parbf{\ref{ex:total-curvature=}.}
Согласно \ref{prop:inscribed-total-curvature}, $\tc\gamma\ge \tc\beta$;
остаётся показать, что
$\tc\gamma\le\sup\{\tc\beta\}$.
Другими словами, 
для любого $\epsilon>0$ и ломаной $\sigma=\vec u_0\dots \vec u_k$, вписанной в касательную индикатрису $\tan$ кривой $\gamma$, 
нам нужно построить такую ломаную $\beta$, вписанную в $\gamma$, что
\[\length\sigma<\tc\beta+\epsilon.
\eqlbl{eq:tc=<tc}\]

Предположим, что $\vec u_i=\tan(s_i)$.
Выберите вписанную ломаную $\beta=p_0\dots p_{2\cdot k+1}$ так, чтобы точки $p_{2\cdot i}$ и $p_{2\cdot i+1}$ лежали достаточно близко к $\gamma(s_i)$; этим можно добиться того, что направление вектора $p_{2\cdot i+1}-p_{2\cdot i}$ достаточно близко к вектору $\vec u_i$ для каждого~$i$.
Покажите, что \ref{eq:tc=<tc} выполняется для построенной ломаной~$\beta$.

\parbf{\ref{ex:tc-rectifiable}.}
Покажите, что для любой ломаной $\beta$, лежащей в шаре радиуса $R$, выполняется неравенство
\[\tc{\beta}+2\cdot \pi\ge\frac{\length\beta}R.\]
Убедитесь, что $\gamma$ лежит в каком-то шаре, и примените это неравенство.

Пример для второй части можно найти среди логарифмических спиралей.


\stepcounter{chapter}
\setcounter{eqtn}{0}

\parbf{\ref{ex:helix-torsion}.} 
Параметр длины $s$ уже найден в \ref{ex:arc-length-helix}.
Остаётся вычислить базис Френе, кривизну и кручение.

{

\begin{wrapfigure}{r}{25 mm}
\vskip-3mm
\centering
\begin{lpic}[t(-0mm),b(0mm),r(0mm),l(0mm)]{asy/helix(1)}
\lbl[br]{8,24;$\norm$}
\lbl[b]{2,26;$\bi$}
\lbl[wl]{15,25;$\tan$}
\end{lpic}
\vskip-0mm
\end{wrapfigure}

Должно получиться
\begin{align*}
\tan(t)&=\tfrac{(-a\cdot\sin t, a\cdot\cos t,b)}{\sqrt{a^2+b^2}},
\\
\norm(t)&=(-\cos t,-\sin t,0),
\\
\bi(t)&=\tfrac{(b\cdot \sin t,-b\cdot \cos t, a)}{\sqrt{a^2+b^2}},
\\
\kur(t) &\equiv \tfrac{a}{a^2+b^2},
\\
\tor(t) &\equiv \tfrac{b}{a^2+b^2}.
\end{align*}

Остаётся показать, что отображение $(a,b) \z\mapsto (\frac{a}{a^2+b^2}, \frac{b}{a^2+b^2})$ биективно отображает полуплоскость $a>0$ на себя.

}

\parbf{\ref{ex:beta-from-tau+nu}.}
По правилу произведения получаем
\begin{align*}
\bi'&=(\tan\times \norm)'=
\tan'\times \norm+\tan\times\norm'.
\end{align*}
Остаётся подставить значения из \ref{eq:frenet-tau} и \ref{eq:frenet-nu} и упростить.

\parbf{\ref{ex:torsion=0}.}
Воспользуйтесь уравнением $\bi' = - \tor\cdot \norm $.

\parbf{\ref{ex:+B}.}
Можно считать, что $\gamma_0$ параметризована длиной.
Покажите, что
$|\gamma_1'|=|\tan-\tor\cdot \norm|\ge 1$ и используйте это.

\parbf{\ref{ex:frenet}.}
$\tfrac{\gamma'\times\gamma''}{|\gamma'\times\gamma''|}$ --- единичный вектор, перпендикулярный плоскости, натянутой на $\gamma'$ и $\gamma''$.
Поэтому, с точностью до знака, он равен $\bi$.
Остаётся найти знак.

\parbf{\ref{ex:moment-curve}.}
Примените \ref{ex:curvature-formulas:b} и \ref{ex:frenet}.

\parbf{\ref{ex:bow-converse}.}
Можно считать, что $t_0=0$ и $\gamma_i(0)\z=0$.
Рассмотрим функции
\begin{align*}
\rho_i(t)&\df|\gamma_i(t)|^2=\langle \gamma_i(t),\gamma_i(t)\rangle.
\end{align*}
Заметьте, что $\rho_1\ge \rho_2$ и $\rho_i(0)=0$.
Убедитесь, что 
\begin{align*}
\rho_i'(0)&=0,
&
\rho_i''(0)&=2,
\\
\rho_i'''(0)&=0,
&
\rho_i''''(0)&=-2\cdot\kur(0)^2_{\gamma_i}.
\end{align*}
Сделайте вывод. (Сравните с \ref{ex:const-dist}.)

\parbf{\ref{ex:torsion-indicatrix}.}
Допустим, что касательная индикатриса не имеет самопересечений.
Покажите, что она лежит в открытой полусфере, и рассуждайте, как в теореме Фенхеля.

\parbf{\ref{ex:lancret}}; \ref{SHORT.ex:lancret:a}.
Заметьте, что 
$\langle \vec w,\tan\rangle'=0$.
Выведите отсюда, что $\langle \vec w, \norm\rangle =0$.
Далее, покажите, что $\langle \vec w, \tan\rangle^2+\langle \vec w, \bi\rangle^2=\langle \vec w, \vec w\rangle$, и воспользуйтесь этим.
При доказательстве последнего тождества примените формулу Френе для~$\norm'$.

\parit{\ref{SHORT.ex:lancret:b}.}
Покажите, что $\vec w'=0$;
это означает, что $\langle \vec w,\tan\rangle\z=\tfrac\tor\kur$.
В частности, вектор скорости образует постоянный угол с $\vec w$; то есть $\gamma$ --- линия откоса.

\parbf{\ref{ex:evolvent-constant-slope}.}
Пусть $\alpha$ --- эвольвента кривой $\gamma$, а $\vec w$ --- фиксированный вектор.
Покажите, что скалярное произведение $\langle \vec w,\alpha\rangle$ постоянно, если $\gamma$ идёт под постоянным углом к $\vec w$, и воспользуйтесь этим.

\parbf{\ref{ex:spherical-frenet}}.
Часть \ref{SHORT.ex:spherical-frenet:tau} следует из того, что $(\tan,\norm,\bi)$ является ортонормированным базисом.
Для \ref{SHORT.ex:spherical-frenet:nu} возьмите первую и вторую производные от тождества $\langle\gamma,\gamma\rangle=1$ и упростите их, пользуясь формулами Френе.
Часть \ref{SHORT.ex:spherical-frenet:beta} вытекает из \ref{SHORT.ex:spherical-frenet:nu} и формул Френе.
Согласно \ref{SHORT.ex:spherical-frenet:beta}, $\int\tfrac\tor\kur=0$, следовательно, выполняется \ref{SHORT.ex:spherical-frenet:beta+}.
Часть \ref{SHORT.ex:spherical-frenet:kur-tor} доказывается алгебраическими манипуляциями.
В \ref{SHORT.ex:spherical-frenet:f},
воспользовавшись формулами Френе, покажите, что $(\gamma+\tfrac1\kur\cdot \norm+\tfrac{\kur'}{\kur^2\cdot\tor}\cdot\bi)'=0$.

\parit{Замечание.}
Из \ref{SHORT.ex:spherical-frenet:beta} получаем, что \textit{$\int\tfrac\tor\kur=0$ для любой замкнутой гладкой сферической кривой}.
Известно, что это свойство характеризует сферы и плоскости; это доказано Беньямино Сегре \cite{segre}.

\parbf{\ref{ex:cur+tor=helix}.}
Воспользуйтесь второй частью \ref{ex:helix-torsion}.

\parbf{\ref{ex:const-dist}.}
По условию,
\begin{align*}
\rho(\ell)&=|\gamma(t+\ell)-\gamma(t)|^2=
\\
&=\langle \gamma(t+\ell)-\gamma(t),\gamma(t+\ell)-\gamma(t)\rangle
\end{align*}
--- гладкая функция, не зависящая от~$t$.
Выразите скорость, кривизну и кручение $\gamma$ через производные $\rho^{(n)}(0)$.
Будьте терпеливы, потребуются две производные для скорости,
четыре для кривизны,
и шесть для кручения.
После этого, примените \ref{ex:cur+tor=helix}.


\stepcounter{chapter}
\setcounter{eqtn}{0}

\parbf{\ref{ex:bike}.}
Можно считать, что $\gamma_0$ параметризована длиной.
Тогда
\begin{align*}
|\gamma_1'|&=|\gamma_0'+\tan'|=|\tan+\kur\cdot\norm|.
\end{align*}
Отсюда 
\[|\gamma_1'(t)|\ge|\gamma_0'(t)|
\quad\text{и}\quad
|\gamma_1'(t)|\ge|\kur(t)_{\gamma_0}|
\]
для любого $t\in[a,b]$.
Проинтегрируйте эти неравенства и примените 
\ref{ex:integral-length}.

\parbf{\ref{ex:trochoids}.}
Заметьте, что 
\[\gamma'_a(t)=(1+a\cdot \cos t, -a\cdot \sin t);\]
то есть $\gamma'_a$ движется по часовой стрелке по окружности с центром в точке $(1,0)$ и радиусом $\vert a \vert$.

\parit{Случай $|a|>1$.}
Обратите внимание, что $\tan_a(t)\z=\gamma'_a/|\gamma'_a|$ движется по часовой стрелке и совершает полный оборот за время $2\cdot\pi$.
Следовательно, $\tgc{\gamma_a}=-2\cdot\pi$ и $\tc{\gamma_a}\z=|\tgc{\gamma_a}|=2\cdot\pi$.

\parit{Случай $|a|<1$.}
Пусть $\theta_a=\arcsin a$.
Покажите, что $\tan_a(t)=\gamma'_a/|\gamma'_a|$ начинается с горизонтального направления $\tan_a(0)=(1,0)$, затем монотонно поворачивается до угла $\theta_a$, затем монотонно до $-\theta_a$, и монотонно возвращается к $\tan_a(2\cdot\pi)\z=(1,0)$.
Отсюда, 
$\tgc{\gamma_a}=0$ и $\tc{\gamma_a}=4\cdot\theta_a$.

\parit{Случай $a=-1$.}
Скорость $\gamma'_{-1}(t)$ обращается в ноль при $t=0$ и $2\cdot\pi$.
Тем не менее кривая допускает гладкую регулярную параметризацию --- найдите её.
Должно получиться $\tc{\gamma_{-1}}\z=-\tgc{\gamma_{-1}}=\pi$.

\parit{Случай $a=1$.}
Скорость $\gamma'_1(t)$ обращается в ноль при $t=\pi$.
При $t=\pi$ у кривой точка возврата;
то есть $\gamma_1$ разворачивается в ней назад.
Значит, у $\gamma_1(t)$ ориентированная полная кривизна не определена.
Кривая состоит из двух гладких дуг с внешним углом $\pi$, и
полная кривизна каждой дуги равна $\tfrac\pi2$, таким образом, 
$\tc{\gamma_{1}}=\tfrac\pi2+\pi+\tfrac\pi2=2\cdot\pi$.

\parbf{\ref{ex:zero-tsc}.}
Ответы на рисунке.
На последнем рисунке предполагается, что у отмеченных точек параллельные касательные. 

\begin{Figure}
\begin{minipage}{.27\textwidth}
\centering
\includegraphics{mppics/pic-260}
\end{minipage}\hfill
\begin{minipage}{.42\textwidth}
\centering
\includegraphics{mppics/pic-261}
\end{minipage}
\hfill
\begin{minipage}{.27\textwidth}
\centering
\includegraphics{mppics/pic-262}
\end{minipage}
\end{Figure}

\parbf{\ref{ex:length'}}; \ref{SHORT.ex:length':reg}.
Покажите, что
\[
\gamma_\ell'(t)=(1-\ell\cdot\skur(t))\cdot \gamma'(t).
\]
Поскольку $\gamma$ регулярна, $\gamma'\ne0$.
Следовательно, если $\gamma_\ell'(t)=0$, то $\ell\cdot \skur(t)=1$.

\parit{\ref{SHORT.ex:length':formula}.}
Можно предположить, что $\gamma$ параметризована длиной, так что $\gamma'(t)=\tan(t)$.
Предположим, что $|\ell|<\frac{1}{\kur(t)}=\frac{1}{|\skur(t)|}$ для любого~$t$.
Тогда 
\[
|\gamma_\ell'(t)|=(1-\ell\cdot\skur(t)).
\]
Следовательно,
\begin{align*}
L(\ell)
&=
\int_a^b(1-\ell\cdot\skur(t))\cdot dt=
\\&=
\int_a^b1\cdot dt-\ell\cdot \int_a^b\skur(t)\cdot dt=
\\
&=
L(0)-\ell\cdot\tgc\gamma.
\end{align*}

\parit{\ref{SHORT.ex:length':antiformula}.}
Рассмотрите единичную окружность $\gamma(t)\z=(\cos t,\sin t)$ при $t\in[0,2\cdot\pi]$ и $\gamma_\ell$ при $\ell=2$.

\parbf{\ref{ex:inverse}.}
Воспользуйтесь определением соприкасающейся окружности через порядок касания и тем, что инверсии переводят окружности в окружности или прямые.

\columnbreak

\parbf{\ref{ex:evolute}.}
Предположим, что $\gamma$ параметризована длиной.
Покажите, что $\omega'=\tfrac{\skur'}{\skur^2}\cdot \norm$. 
Выведите отсюда, что у $\omega$ базис Френе либо $(\norm,-\tan)$, либо $(-\norm,\tan)$, и её кривизна равна $|\tfrac{\skur^3}{\skur'}|$.

\parbf{\ref{ex:3D-spiral}.}
Начните с плоской спирали, показанной на рисунке.
Применяя \ref{thm:fund-curves}, увеличивайте кручение пунктирной дуги, не изменяя кривизну, до появления самопересечения.

\begin{wrapfigure}{r}{17 mm}
\vskip-5mm
\centering
\includegraphics{mppics/pic-296}
\vskip-0mm
\end{wrapfigure}

Можно считать, что пунктирная дуга очень короткая, и касательный вектор $\tan$ почти не меняется на этой дуге. 
Увеличение кручения может вращать нормальный вектор $\norm$ произвольно вокруг $\tan$.
Пересечение появится, если $\norm$ будет повёрнут на некоторый угол, близкий к $\pi$.
(Сравните с \ref{ex:approximation-const-curvature}.)

\parbf{\ref{ex:double-tangent}.} 
Заметьте, что если прямая или окружность касается кривой $\gamma$, то она касается и соприкасающейся окружности в той же точке, и примените лемму о спирали (\ref{lem:spiral}).

\parbf{\ref{ex:spherical-spiral}.}
Убедитесь, что соприкасающиеся окружности к сферической кривой лежат на сфере.
Докажите аналог \ref{lem:spiral} для этих окружностей.
(Сравните с \ref{ex:spherical-frenet:beta+}.)


\stepcounter{chapter}
\setcounter{eqtn}{0}

\parbf{\ref{ex:vertex-support}.}
Примените лемму о спирали (\ref{lem:spiral}).

\parit{Вычислительное решение.} 
Будем считать, что кривизна не обнуляется в точке $p$, оставшийся случай проще.
Можно предположить, что $\gamma$ параметризована длиной, $p=\gamma(0)$, а начало координат --- центр кривизны при $p$;
другими словами, $\kur(0)\cdot\gamma(0)\z+ \norm(0)=0$.

Рассмотрим функцию $f\:t\mapsto \langle\gamma(t),\gamma(t)\rangle$.
Прямые вычисления дают 
\begin{align*}
f'
&=\langle\gamma,\gamma\rangle'
=2\cdot\langle\tan,\gamma\rangle,
\\
f''
&=2\cdot\langle\tan,\gamma\rangle'
=2\cdot\kur\cdot\langle\norm,\gamma\rangle+2,
\\
f'''
&=2\cdot\kur'\cdot\langle\norm,\gamma\rangle
+2\cdot\kur\cdot\langle\norm',\gamma\rangle+\cancel{2\cdot\kur\cdot\langle\norm,\tan\rangle}.
\end{align*}

Заметьте, что $\norm'(0)\perp\gamma(0)$.
Следовательно, $f'(0)=0$, $f''(0)=0$, и $f'''(0)\z=-2\cdot\kur'/\kur$.

Если соприкасающаяся окружность подпирает $\gamma$ в точке $p$,
то функция $f$ имеет локальный максимум или минимум в точке $0$.
Следовательно, $f'''(0)=0$ и $\kur'=0$.

\parbf{\ref{ex:support}.}
Выберем систему координат с началом в точке $p$ и осью $x$ касательной к $\gamma_1$ и $\gamma_2$.
Можно предположить, что $\gamma_1$ и $\gamma_2$ определены на $(-\epsilon,\epsilon)$ для малого $\epsilon>0$,
и они идут почти горизонтально.

Для данного $t\in[0,1]$ рассмотрим кривую $\gamma_t$, касательную и сонаправленную с осью $x$ в точке $\gamma_t(0)=p$ и с ориентированной кривизной $\skur_t(s)\z\df(1\z-t)\cdot\skur_0(s)+t\cdot\skur_1(s)$.
По \ref{thm:fund-curves-2D}, такая существует.

Выберите $s\approx 0$.
Рассмотрите кривую $\alpha_s\:t\z\mapsto \gamma_t(s)$.
Покажите, что $\alpha_s$ движется почти вертикально вверх.
Воспользовавшись тем, что $\gamma_t$ движется почти горизонтально, покажите, что в малой окрестности точки $p$ кривая $\gamma_1$ лежит выше $\gamma_0$,
откуда утверждение следует.

\parbf{\ref{ex:in-circle}.}
Двигайте прямо единичную окружность, пока её центр не достигнет базовой точки петли, затем сжимайте радиус окружности до нуля.
В момент, когда окружность впервые дотронется до петли, она подопрёт её в точке, отличной от базовой.
Примените~\ref{prop:supporting-circline}.

\parbf{\ref{ex:between-parallels-1} $\bm{+}$ \ref{ex:in-triangle} $\bm{+}$ \ref{ex:lens}.}
Убедитесь, что одна из дуг кривизны $1$ в семействе, показанном на рисунке, подопрёт $\gamma$, и примените \ref{prop:supporting-circline}.

\begin{Figure}
\begin{minipage}{.35\textwidth}
\centering
\includegraphics{mppics/pic-265}
\end{minipage}
\hfill
\begin{minipage}{.3\textwidth}
\centering
\includegraphics{mppics/pic-266}
\end{minipage}
\hfill
\begin{minipage}{.25\textwidth}
\centering
\includegraphics{mppics/pic-267}
\end{minipage}
\end{Figure}

Вторую часть в \ref{ex:between-parallels-1} можно свести к \ref{ex:in-circle}, используя показанное семейство и другое семейство дуг, изогнутых в противоположном направлении.

\parit{Замечание.}
Сравните с \ref{ex:moon-rad}.

\parbf{\ref{ex:convex small}.}
Достаточно рассмотреть случай когда $\gamma$ ограничивает выпуклую фигуру, скажем~$F$.
В противном случае, по \ref{prop:convex}, её кривизна меняет знак.
Следовательно, $\gamma$ имеет нулевую кривизну в некоторой точке.

\begin{wrapfigure}{r}{32 mm}
\vskip-5mm
\centering
\includegraphics{mppics/pic-268}
\vskip0mm
\end{wrapfigure}

Выберите внутри $F$ такие две точки $x$ и $y$, что $|x\z-y|\z>2$.
Рассмотрите максимальную линзу, ограниченную двумя дугами с общей хордой $xy$, которая лежит в~$F$.
Примените признак опорной (\ref{prop:supporting-circline}).

\parbf{\ref{ex:convex-lens}.}
Применив лемму о линзе (\ref{lem:lens}), покажите, что $\gamma$ лежит с одной стороны от прямой $pq$.
Выведите отсюда, что объединение дуги $\gamma$ и её хорды $[p,q]$ образует простую замкнутую кривую;
по теореме Жордана она ограничивает фигуру, назовём её~$F$.

Предположите, что $F$ не выпукла, и придите к противоречию, как в \ref{prop:convex}.

\parbf{\ref{ex:diameter-of-simple-curve}.}
Выберите две точки $p,q\in\gamma$ на максимальном расстоянии.
Примените \ref{ex:convex-lens} к дуге $\gamma$ от $p$ до $q$.
После этого используйте лемму о луке (\ref{lem:bow}) или рассуждайте как в \ref{ex:convex small}.

Для второй части упражнения рассмотрите наименьшую окружность $\sigma$, окружающую $\gamma$, и повторите рассуждение выше для одной из дуг $\gamma$ между общими точками с $\sigma$.

\parit{Замечание.}
Если разрешить самопересечения, то диаметр кривой можно сделать сколь угодно большим.
Догадайтесь до примера по рисунку.
На самом деле диаметр можно сделать произвольно близким к длине кривой.
\begin{Figure}
\vskip-0mm
\centering
\includegraphics{mppics/pic-269}
\vskip-0mm
\end{Figure}

\parbf{\ref{ex:moon-rad}.}
Докажите, что $\gamma$ содержит простую петлю $\gamma_1$.
Примените \ref{thm:moon-orginal} к $\gamma_1$.

\parbf{\ref{ex:2-squares}.}
Можно считать, что меньший квадрат лежит слева от $\gamma$.
Воспользовавшись \ref{thm:4-vert-supporting}, покажите, что на $\gamma$ есть две точки с ориентированной кривизной $\ge 1$ и $\le 1$ соответственно.
Применив теорему о промежуточных значениях к ориентированной кривизне, покажите, что существует точка с кривизной $1$.

\parit{Более сложный вопрос:} Каково минимальное число таких точек?

\parbf{\ref{ex:moon-area}.} 
Применив гомотетию, можно предположить, что $a=\pi$; то есть наша петля ограничивает фигуру равновеликую единичному кругу.
Применив \ref{thm:moon}, покажите, что кривизна в какой-то точке не меньше $1$.
Рассуждая как в конце доказательства \ref{thm:4-vert-supporting}, покажите, что кривизна в какой-то точке не больше $1$.

Применив теорему о промежуточных значениях, покажите, что в какой-то точке кривизна ровно $1$.

\begin{wrapfigure}{r}{20 mm}
\vskip-6mm
\centering
\includegraphics{mppics/pic-305}
\vskip-2mm
\end{wrapfigure}

\parbf{\ref{ex:curve-crosses-circle}.}
Повторите доказательство \ref{thm:moon} для каждого циклического произведения дуги $\gamma$ от $p_i$ до $p_{i+3}$ с дугой окружности, на которой лежат точки $p_{i+4},\z\dots,p_{i-1}$.

До примера во второй части можно догадаться по рисунку.




\parbf{\ref{ex:berk}.}
Необходимость очевидна; докажем достаточность.
Выберем такие две точки $p,q\in\gamma$, что $\dist{p}{q}{}\ge 4$.
Обозначим через $\sigma_p$, $\sigma_q$, $r_p$, $r_q$ вписанные окружности и их радиусы при точках $p$ и $q$.

Предположим, что $r_p\ge 1$ и $r_q\ge 1$.
Покажите, что два единичных диска, касающиеся $\gamma$ изнутри в точках $p$ и $q$, решают задачу.

Предположим, что $r_p<1$ и $r_q<1$.
Убедитесь, что $\sigma_p$ и $\sigma_q$ не пересекаются.
Заметьте, что мы можем выбрать две непересекающиеся дуги $\gamma_p$ и $\gamma_q$ кривой $\gamma$, концы которых лежат на $\sigma_p$ и $\sigma_q$ соответственно.
Начните с этих двух дуг и рассуждайте, как в \ref{thm:moon}.
Покажите, что полученные опорные соприкасающиеся окружности не пересекаются; сделайте последний шаг.

В оставшемся случае $r_p\ge 1$ и $r_q<1$ объедините эти два рассуждения.

\parit{Источник:} Это задача Берка Сейлана \cite{ceylan}.
Неизвестно, \textit{ограничивает ли любая простая замкнутая плоская кривая с кривизной не более $1$ и длиной не менее $4\cdot\pi$ два непересекающихся единичных диска.}

\parbf{\ref{ex:4x0-torsion}.}
Убедитесь, что касательная индикатриса $\tan$ кривой $\gamma$ является гладкой простой замкнутой сферической кривой.
Покажите, что каждая точка на кривой $\gamma$ с нулевым кручением соответствует \emph{точке перегиба} у $\tan$;
то есть точке с нулевой геодезической кривизной (см.~\ref{sec:Darboux}).

Попробуйте доказать, что если у индикатрисы $\tan$ меньше $4$ точек перегиба, то она лежит в открытой полусфере,
или прочитайте доказательство теоремы о теннисном мячике \cite[§ 20]{arnold1994}.
После этого рассуждайте, как в \ref{thm:fenchel}.


\stepcounter{chapter}
\setcounter{eqtn}{0}

\parbf{\ref{ex:hyperboloids}.}
Пусть $\Sigma_\ell$ --- множество уровня.

Покажите, что $\nabla_p f=0$ тогда и только тогда, когда $p\z=(0,0,0)$.
Воспользовавшись \ref{prop:implicit-surface}, покажите, что при $\ell\ne 0$ каждая связная компонента в $\Sigma_\ell$ являет собой гладкую поверхность.

Докажите, что $\Sigma_\ell$ связно тогда и только тогда, когда $\ell\ge 0$.
Из этого вытекает, что $\Sigma_\ell$ --- гладкая поверхность при $\ell>0$, а при $\ell<0$ --- нет.

Случай $\ell=0$ нужно разобрать вручную --- хотя достаточное условие в \ref{prop:implicit-surface} не выполнено, это ещё не означает, что $\Sigma_0$ не является гладкой поверхностью.

Покажите, что любая окрестность начала координат в $\Sigma_0$ не описывается графиком ни в какой системе координат.
Значит по определению (см.~\ref{sec:def-smooth-surface}) $\Sigma_0$ не является гладкой поверхностью.

\parbf{\ref{ex:9-surf}.}
Пример похож на \ref{ex:9:9}.

\parbf{\ref{ex:smooth-fun(surf)}.}
При доказательстве достаточности воспользуйтесь тем, что композиция гладких функций гладкая.

Далее, пусть график $z\z=f(x,y)$ задаёт $\Sigma$ в окрестности точки $p$.
Заметьте, что $g$ является гладкой в этой окрестности тогда и только тогда, когда функция $\hat g\:(x,y)\mapsto g(x,y,f(x,y))$ является гладкой в своей области определения.
Определите $h(x,y,z)\z\df \hat g(x,y)$, отсюда достаточность.

Для последней части рассмотрите функцию $g\:(x,y)\mapsto\tfrac1{x^2+y^2}$, определённую на $(x,y)$-плоскости с удалённым началом координат.
Заметьте, что $g$ не продолжается до непрерывной функции на $\mathbb{R}^3$.

\parit{Замечание.}
Если $\Sigma$ --- гладкая собственная поверхность,
то \textit{любая гладкая функция $g\:\Sigma\to\mathbb{R}$ продолжается до гладкой функции на $\mathbb{R}^3$}.
Доказательство использует так называемое \emph{разбиение единицы};
почитайте о нём и докажите данное утверждение.

\parbf{\ref{ex:inversion-chart}.} 
Проверьте, что образ $s$ лежит на единичной сфере с центром в точке $(0,0,1)$;
то есть $
\left(\tfrac{2\cdot u}{1+u^2+v^2}\right)^2
\z+
\left(\tfrac{2\cdot v}{1+u^2+v^2}\right)^2
\z+
\left(\tfrac{2}{1+u^2+v^2}-1\right)^2\z=1$
при любых $u$ и~$v$.

Покажите, что отображение 
$(x,y,z)\z\mapsto \tfrac{2\cdot(x,y)}{x^2+y^2+z^2}$
--- обратное к $s$, и оно является непрерывным вне начала координат.
В частности, $s$ является вложением, которое покрывает всю сферу, кроме начала координат.

Покажите, что $s$ регулярно; то есть $s_u$ и $s_v$ линейно независимы при всех $u$ и $v$.

\parit{Замечание.}
Можно думать, что $s$ задаёт \index{стереографическая проекция}\emph{стереографическую проекцию}
\[(u,v,1)\mapsto (\tfrac{2\cdot u}{1+u^2+v^2},\tfrac{2\cdot v}{1+u^2+v^2},\tfrac{2}{1+u^2+v^2})\]
из плоскости $z=1$ на единичную сферу с центром в точке $(0,0,1)$.
Заметьте, что точка $(u,v,1)$ и её образ лежат на одном и том же луче, исходящем из начала координат.

\begin{Figure}
\vskip-0mm
\centering
\includegraphics{mppics/pic-750}
\vskip0mm
\end{Figure}

\parbf{\ref{ex:revolution}.}
Рассмотрите отображение
\[s\:(t,\theta)\mapsto (x(t), y(t)\cdot\cos\theta,y(t)\cdot\sin\theta).\]
Покажите, что оно регулярно; то есть $s_t$ и $s_\theta$ линейно независимы.
(Здесь пригодится, что $s_t\z\perp s_\theta$).

Убедитесь, что $s$ является вложением;
то есть любая точка $(t_0,\theta_0)$ допускает такую окрестность $U$ на $(t,\theta)$-плоскости, что сужение $s|_U$ имеет непрерывное обратное отображение.
Остаётся применить \ref{cor:reg-parmeterization}.

\parbf{\ref{ex:inv-diffeomorphism}.}
Примените теорему об обратной функции (\ref{thm:inverse}) в картах для обеих поверхностей.

\parbf{\ref{ex:star-shaped-disc}.} 
Решения этих упражнений основаны на следующей общей конструкции, известной как \index{трюк Мозера}\emph{трюк Мозера}.

Предположим, что $\vec u_t$ --- это гладкое векторное поле на плоскости, зависящее от времени.
Рассмотрим обыкновенное дифференциальное уравнение $x'(t)=\vec u_t(x(t))$ и отображение $\iota\:x(0)\mapsto x(1)$, где $t\z\mapsto x(t)$ --- решение уравнения.
Отображение $\iota$ называется \index{поток}\emph{потоком} векторного поля $\vec u_t$ во временном интервале $[0,1]$.
Согласно \ref{thm:ODE}, оно гладко в своей области определения.
Более того, то же самое верно для его обратного $\iota^{-1}$;
действительно, $\iota^{-1}$ является потоком векторного поля $-\vec u_{1-t}$.
То есть, $\iota$ является диффеоморфизмом из своей области определения на свой образ.

Таким образом, для построения диффеоморфизма из одного открытого подмножества плоскости в другое достаточно построить такое гладкое векторное поле, что его поток отображает одно множество на другое;
такое отображение автоматически является диффеоморфизмом.

\parit{\ref{SHORT.ex:plane-n}.}
Рассмотрим два множества $\mathbb{R}^2\z\setminus\{p_1,\dots,p_n\}$ и $\mathbb{R}^2\setminus\{q_1,\dots,q_n\}$.
Выберем гладкие траектории $\gamma_i\:[0,1]\z\to \mathbb{R}^2$ так, чтобы $\gamma_i(0)=p_i$,
$\gamma_i(1)=q_i$, и $\gamma_i(t)\ne \gamma_j(t)$ при $i\ne j$.

Выберите такое гладкое векторное поле $\vec v_t$, что $\vec v_t(\gamma_i(t))=\gamma'_i(t)$ для любого $i$ и~$t$.
Мы можем дополнительно предположить, что $\vec v_t$ обращается в ноль за пределами достаточно большого диска; это можно организовать, умножив векторное поле на подходящую функцию;
например, $\sigma_1(R-|x|)$, где $R$ велико, а $\sigma_1$ --- сигмоид (стр. \pageref{page:sigma-function}).

Остаётся применить трюк Мозера к построенному векторному полю.

\parit{\ref{SHORT.ex:star-shaped-disc:smooth}--\ref{SHORT.ex:star-shaped-disc:star-shaped}.}
В \ref{SHORT.ex:star-shaped-disc:smooth}--\ref{SHORT.ex:star-shaped-disc:nonsmooth}, можно считать, что начало координат принадлежит обоим множествам, а в \ref{SHORT.ex:star-shaped-disc:star-shaped}, что фигуры звёздны относительно начала координат.

В каждом случае покажите, что существует векторное поле $\vec v$, определённое на $\mathbb{R}^2$, которое отображает одну поверхность на другую.
На самом деле, можно выбрать радиальные поля такого типа,
но будьте осторожны с задачами \ref{SHORT.ex:star-shaped-disc:nonsmooth} и \ref{SHORT.ex:star-shaped-disc:star-shaped} --- они сложней чем могут показаться.


\stepcounter{chapter}
\setcounter{eqtn}{0}

\parbf{\ref{ex:tangent-normal}.}
Пусть $\gamma$ --- гладкая кривая на~$\Sigma$.
Заметьте, что $f\circ\gamma(t)\equiv 0$.
Продифференцируйте это тождество и примените определение касательного вектора (\ref{def:tangent-vector}).

\parbf{\ref{ex:vertical-tangent}.}
Предположим, что окрестность точки $p$ на $\Sigma$ задана графиком $z=f(x,y)$.
В этом случае, $s\:(u,v)\z\mapsto (u,v,f(u,v))$ --- гладкая карта, накрывающая~$p$.
Покажите, что плоскость, натянутая на $s_u$ и $s_v$, не является вертикальной;
вместе с \ref{def:tangent-plane} это доказывает достаточность.

Пусть
$s\:(u,v)\z\mapsto(x(u,v),y(u,v),z(u,v))$ --- карта.
Примените теорему об обратной функции к $(u,v)\z\mapsto(x(u,v),y(u,v))$ и выведите необходимость.

\parbf{\ref{ex:tangent-single-point}.}
Выберите $(x,y,z)$-координаты так, чтобы $\Pi$ была $(x,y)$-плоскостью, а $p$ --- началом координат.
Пусть $(u,v)\mapsto s(u,v)$ --- такая карта на $\Sigma$, что $p=s(0,0)$,
и пусть $\vec k$ --- единичный вектор в направлении $z$-оси.

Убедитесь, что можно считать, что $\langle s(u,v),\vec k\rangle>0$ в проколотой окрестности точки $0$.
Выведите отсюда, что $s_u\perp \vec k$ и $s_v\perp \vec k$ в
$0$; отсюда упражнение следует.

\parbf{\ref{ex:lin-ind-chart}.}
Воспользовавшись \ref{thm:ODE}, постройте такую кривую $\alpha$ на $\Sigma$, что $\alpha(0)=p$ и $\alpha'(x)\z=\vec x_{\alpha(x)}$.
Точно также постройте кривую $\beta_x$ на $\Sigma$, для которой $\beta_x(0)=\alpha(x)$ и $\beta_x'(y)=\vec y_{\beta_x(y)}$.
Применив \ref{thm:inverse}, покажите, что отображение $(x,y)\mapsto \beta_x(y)$ описывает карту малой окрестности $W$ точки $p$.

Покажите, что $u\:\beta_x(y)= x$ удовлетворяет требуемым условиям в $W$.
Тем же способом постройте $v$.
Останется снова применить \ref{thm:inverse}.

\parbf{\ref{ex:const-normal}.}
Покажите, что для любой гладкой кривой $\gamma$ на $\Sigma$ функция $t\mapsto \langle\nu_0,\gamma(t)\rangle$ постоянна, и воспользуйтесь этим.

\parbf{\ref{ex:implicit-orientable}.}
По \ref{ex:tangent-normal}, $\Norm=\tfrac{\nabla h}{|\nabla h|}$ --- поле нормалей.

\parbf{\ref{ex:plane-section}.}
Воспользовавшись срезками и сглаживаниями (\ref{sec:analysis}), постройте такую гладкую неотрицательную функцию $f$ на плоскости $(x,y)$, что $f(x,y)=0$ тогда и только тогда, когда $(x,y)\in A$.
График $z=f(x,y)$ описывает требуемую поверхность.


\stepcounter{chapter}
\setcounter{eqtn}{0}

\parbf{\ref{ex:line-of-curvature}.}
Возьмём точку $p$ на~$\gamma$.
Убедитесь, что $\T_p\perp \Pi$, поскольку $\Sigma$ зеркально-симметрична относительно $\Pi$.

Выберите $(x,y)$-координаты на $\T_p$ с осью $x$ идущей по пересечению $\Pi\cap \T_p$.
Предположим, что соприкасающийся параболоид описывается графиком 
$z=\tfrac12\cdot(\ell\cdot x^2+2\cdot m\cdot x\cdot y+n\cdot y^2)$.
Раз $\Sigma$ зеркально-симметрична, то и параболоид зеркально-симметричен;
то есть замена $y$ на $(-y)$ не меняет 
$\ell\cdot x^2+2\cdot m\cdot x\cdot y+n\cdot y^2$.
Иначе говоря, $m=0$, то есть ось $x$ указывает в направлении кривизны.

\parbf{\ref{ex:gauss+orientable}.}
Заметьте, что в каждой точке главные кривизны одного знака.
Поэтому нормаль $\Norm$ в каждой точке можно выбрать так, чтобы обе главные кривизны были положительными.
Покажите, что это определяет глобальное поле на поверхности.

\parbf{\ref{ex:re-scale-surface-curvature}.}
Вычислите ряд Тейлора функции $g(x,y)= \lambda \cdot f( x/ \lambda , y/\lambda)$.

\parbf{\ref{ex:self-adjoint}.}
Примените \ref{thm:shape-chart} к отображению $s$, такому что $s_u(0,0)=\vec u$ и $s_v(0,0)=\vec v$.

\parbf{\ref{ex:normal-curvature=const}}; \ref{SHORT.ex:normal-curvature=const:a}.
Заметьте, что $\Sigma$ имеет единичную матрицу Гессе в каждой точке, и примените определение оператора формы.

\parit{\ref{SHORT.ex:normal-curvature=const:b}.}
Выберите карту $s$ на~$\Sigma$.
Покажите, что
\[\tfrac{\partial }{\partial u}(s+\Norm)
=
\tfrac{\partial }{\partial v}(s+\Norm)
=
0.\]
Сделайте последний шаг.

\parbf{\ref{ex:normal-curvature=0}.}
Задача кажется тривиальной пока не осознаешь, что связные множества бывают весьма причудливы.
Например, в них может не быть нетривиальных кривых.

\medskip

Применив лемму Сарда (\ref{lem:sard}), покажите, что $\Norm$ постоянна на $Z_0$.
Пусть $\Norm_0$ --- значение нормали на $Z_0$.
Снова применив лемму Сарда, покажите, что функция $x\mapsto \langle \Norm_0,x\rangle$ постоянна на $Z_0$.
Сделайте вывод.

\parit{Замечание.}
Этот результат был использован Ричардом Сакстедером \cite[Lemma 6]{sacksteder}.  
Попробуйте решить следующую чуть более сложную задачу:  
\textit{Пусть $\Sigma$ --- гладкая поверхность с ориентацией, заданной полем нормалей $\Norm$, и пусть $Z_1\subset \Sigma$ --- связное множество с единичным оператором формы.  
Покажите, что $Z_1$ лежит на единичной сфере.}

\parbf{\ref{ex:shape-curvature-line}.} 
Можно считать, что $\gamma$ параметризована длиной.
Обозначим через $\Norm_1(s)$ и $\Norm_2(s)$ единичные нормали к $\Sigma_1$ и $\Sigma_2$ в точке $\gamma(s)$.

По предположению, $\langle \Norm_1,\Norm_2\rangle= c$ для некоторой константы $c$.
Мы знаем, что
$\Norm_1'$ пропорциональна $\gamma'$, а показать нужно, что $\Norm_2'$ пропорциональна $\gamma'$.
Если $c=\pm1$, то $\Norm_1\equiv\pm \Norm_2$, и отсюда всё следует.

Убедитесь, что $\langle\Norm_1',\Norm_2\rangle=0$.
Воспользовавшись равенством $c'=0$, покажите, что $\langle\Norm_1,\Norm_2'\rangle=0$.
Далее покажите, что $\langle\Norm_2,\Norm_2'\rangle=0$, и выведите отсюда, что производная $\Norm_2'$ пропорциональна $\gamma'$.

\parit{Источник:}
Этот результат получен Фердинандом Йоахимшталем \cite{joachimsthal}, и обобщён Оссианом Бонне \cite{bonnet}.

\parbf{\ref{ex:equidistant}};
\ref{SHORT.ex:equidistant:smooth}.
Для данного $t$, рассмотрим  отображение $f_t\:p\mapsto p+t\cdot\Norm(p)$; оно отображает $\Sigma$ в $\Sigma_t$.

Применив определение оператора формы, покажите, что $d_pf_t(\vec v)=\vec v-t\cdot\Shape_p(\vec v)$.
Так как $\Sigma$ замкнута, норма оператора $\Shape$ ограничена.
Значит $f_t$ регулярно при $t\approx 0$.

Далее покажите, что $f_t$ инъективно при $t\approx0$;
выведите отсюда, что $f_t$ --- гладкое вложение.

\parit{\ref{SHORT.ex:equidistant:area}.}
По формуле площади (\ref{prop:surface-integral}), 
\[\area\Sigma_t=\int_{p\in \Sigma}\jac_pf_t.\]
Покажите, что $\jac_pf_t=1-t\cdot H+t^2\cdot K$, и воспользуйтесь этим.

\parbf{\ref{ex:flat-plane};} \ref{SHORT.ex:flat-plane:orthonormal}.
Воспользуйтесь тем, что главные направления ортогональны; см. \ref{sec:Principal curvatures}.

\parit{\ref{SHORT.ex:flat-plane:depend}.} 
Применив \ref{cor:Shape(ij)}, покажите, что $\Norm_u=0$, и, воспользовавшись этим, покажите, что $\Norm$ зависит только от~$v$.

Убедитесь, что $\Norm_{uv}=\Norm_{vu}=0$.
Воспользовавшись \ref{cor:Shape(ij)}, покажите, что $\Norm_v\z=-k \cdot s_v$.
Выведите отсюда, что $\vec v_u=0$, и, следовательно, $\vec v$ зависит только от~$v$.

Примените \ref{SHORT.ex:flat-plane:orthonormal}, чтобы показать, что $\vec u$ зависит только от~$v$, и
завершите доказательство.

\parit{\ref{SHORT.ex:flat-plane:depend-u}.}
Первая часть вытекает из \ref{SHORT.ex:flat-plane:depend}.

Согласно \ref{SHORT.ex:flat-plane:orthonormal}, $\langle s_u,s_v\rangle=0$.
Поскольку $s_{uu}$ пропорционален $s_u$, мы также получаем $\langle s_{v},s_{uu}\rangle=0$.
Следовательно,
\begin{align*}
\tfrac{\partial}{\partial v}\langle s_u,s_u\rangle&=2\cdot \langle s_{vu},s_u\rangle=
\\
&=2\cdot \tfrac{\partial}{\partial u}\langle s_v,s_u\rangle-2\cdot \langle s_{v},s_{uu}\rangle=0.
\end{align*}
Сделайте последний шаг.

\parit{\ref{SHORT.ex:flat-plane:linear}.} 
Поскольку $\langle s_u,s_u\rangle=1$,
\[0=\tfrac{\partial}{\partial u}\langle s_u,s_u\rangle=2\cdot\langle s_{uu},s_u\rangle.\]
По \ref{SHORT.ex:flat-plane:depend-u}, $s_{uu}\parallel s_u$; следовательно, $s_{uu}=0$.

Выведите отсюда, что $s_v$ зависит линейно от $u$.
Воспользуйтесь этим и уравнениями
\[\Norm_v=-k\cdot s_v,\quad \Norm_{uv}=0.\]

\parit{Замечания.}
Это доказательство приводится в лекции Сергея Иванова \cite[3-й сем. Лекция 13]{ivanov}.
Другое доказательство получается применением формул Петерсона --- Кодацци.

\stepcounter{chapter}
\setcounter{eqtn}{0}

\parbf{\ref{ex:mean-curvature}.}
Примените \ref{obs:k1-k2} и определение средней кривизны.

\parbf{\ref{ex:average}.}
Вспомните, что $H=k_1+k_2$ и $K=k_1\cdot k_2$, где $k_1$ и $k_2$ --- главные кривизны.
По тождеству Эйлера, нам нужно найти среднее значение
\[(k_1\cdot (\cos\theta)^2+k_2\cdot (\sin\theta)^2)^2.\]
Покажите, что $\tfrac38$, $\tfrac38$ и $\tfrac18$ --- средние значения $(\cos\theta)^4$, $(\sin\theta)^4$ и $(\cos\theta\cdot \sin\theta)^2$ соответственно, и воспользуйтесь этим.

\parbf{\ref{ex:meusnier}.}
Воспользовавшись теоремой Мёнье (\ref{thm:meusnier}), найдите центр и радиус кривизны $\gamma$ через её нормальную кривизну в точке $p$.
Сделайте последний шаг.

\parit{Источник:}
Это утверждение, как и \ref{thm:meusnier}, доказано Жаном Батистом Мёнье \cite{meusnier}.

\parbf{\ref{ex:principal-revolution}.}
Воспользуйтесь \ref{ex:line-of-curvature} и теоремой Мёнье (\ref{thm:meusnier}).

\parbf{\ref{ex:catenoid-is-minimal}.}
Используйте \ref{ex:principal-revolution:a},  \ref{ex:curvature-graph} и \ref{thm:meusnier}.

\parbf{\ref{ex:helicoid-is-minimal}.}
Применив теорему Мёнье (\ref{thm:meusnier}), покажите, что координатные кривые $\alpha_v\:u\z\mapsto s(u,v)$ и $\beta_u\:v\z\mapsto s(u,v)$ асимптотические; то есть их нормальные кривизны обнуляются.

Убедитесь, что эти два семейства ортогональны друг другу.
Примените \ref{ex:mean-curvature}.

\parbf{\ref{ex:rev(sin)}.}
Согласно \ref{ex:principal-revolution:a}, параллели и меридианы являются линиями кривизны.
Одна из главных кривизн --- это кривизна образующей.
Воспользовавшись \ref{ex:curvature-formulas:a}, покажите, что она не превышает $a\cdot \sin x$.
Заметьте, что кривизна параллели равна $\tfrac1{a\cdot \sin x}$;
применив теорему Мёнье, покажите, что главная кривизна в её направлении не может быть больше.
Сделайте последний шаг.

\parbf{\ref{ex:rev(lin)}.}
Примените \ref{ex:principal-revolution:formula}.

\parbf{\ref{ex:moon-revolution}.}
Используйте \ref{ex:line-of-curvature} и \ref{thm:moon-orginal}.

\parbf{\ref{ex:lagunov-genus4}.}
Просверлите ещё пару дырок или объедините два примера вместе.
Во второй части объедините $V_2$ and $V_3$.

\parit{Замечание.}
Можно построить такой пример с торической границей, но это очень продвинутое упражнение.
Тело, ограниченное сферой, обязано содержать шар радиуса чуть большего, чем $r_3=\sqrt{3/2}-1>r_2$; см. \cite{lagunov-1960, lagunov-1961, lagunov-fet-1963, lagunov-fet-1965}.


\stepcounter{chapter}
\setcounter{eqtn}{0}

\parbf{\ref{ex:supp>tan}.}
Можно считать, что $\Sigma_1$ задаётся неявно гладкой функцией $f$ в малой окрестности $U$ точки $p$.
Поскольку $\Sigma_1$ локально подпирает $\Sigma_2$ в точке $p$, можно считать, что $f(x)\ge 0$ для всех $x\in \Sigma_2\cap U$
и при этом $f(p)=0$.

Докажите, что $(f\circ\gamma_2)'(0)=0$ для любой гладкой кривой $\gamma_2\:(-\varepsilon,\varepsilon)\to\Sigma_2$ такой, что $\gamma_2(0)=p$.
Выведите отсюда, что вектор скорости $\gamma_2'(0)$ касателен к $\Sigma_1$ в $p$ и сделайте последний шаг.

\parbf{\ref{ex:surf-support}.}
Выберите кривизны так, чтобы 
\[k_2(p)_{\Sigma_1}\z>k_2(p)_{\Sigma_2}> k_1(p)_{\Sigma_1}> k_1(p)_{\Sigma_2},\]
и предположите, что первое главное направление $\Sigma_1$ совпадает со вторым главным направлением $\Sigma_2$ и наоборот.

\parbf{\ref{ex:positive-gauss-0} $\bm{+}$ \ref{ex:positive-gauss}.}
Рассуждайте как в \ref{ex:between-parallels-1}--\ref{ex:lens} с использованием семейств сфер.

\parbf{\ref{ex:convex-surf}.}
Покажите, что любая касательная плоскость $\T_p$ подпирает $\Sigma$ в точке~$p$,
и примените \ref{prop:surf-support}.

\parbf{\ref{ex:convex-lagunov}.}
Предположим, что максимальный шар в $R$ имеет радиус $r$ и касается границы $R$ в точках $p$ и~$q$.

\begin{wrapfigure}{r}{15 mm}
\vskip-2mm
\centering
\includegraphics{mppics/pic-1050}
\vskip-0mm
\end{wrapfigure}

Пусть $S$ --- проекция области $R$ на плоскость, проходящую через $p$, $q$ и центр шара.
Покажите, что $S$ ограничена гладкой замкнутой выпуклой кривой с кривизной не более $1$.
Рассуждая как в \ref{ex:between-parallels-1}, покажите, что $r\ge 1$.

\parbf{\ref{ex:section-of-convex}.}
Выберем плоскость $\Pi$.
Пусть точка $p$ лежит на пересечении $\Pi\cap\Sigma$.

Покажите, что если $\Pi$ касается поверхности в точке $p$, то $p$ является изолированной точкой пересечения $\Pi\cap\Sigma$.

Отсюда следует, что если $\gamma$ --- компонента связности пересечения $\Pi\cap\Sigma$, которая не является изолированной точкой, то $\Pi$ пересекает $\Sigma$ \index{трансверсальность}\emph{трансверсально} вдоль $\gamma$; то есть $\T_p\Sigma\ne\Pi$ для любой точки $p \in \gamma$. Примените теорему о неявной функции, чтобы показать, что $\gamma$ --- гладкая кривая.

Наконец, кривизна $\gamma$ не меньше нормальной кривизны $\Sigma$ в том же направлении.
Следовательно, у $\gamma$ нет точек с нулевой кривизной.
Поэтому её ориентированная кривизна имеет постоянный знак.

\parbf{\ref{ex:surrounds-disc}.}
Предположим противное.
По \ref{thm:convex-embedded}, наша поверхность выпукла.
Покажите, что поверхность подпирается изнутри сферическим куполом с нашей единичной окружностью в качестве края.
Придите к противоречию с \ref{cor:surf-support}.

\parbf{\ref{ex:small-gauss}.}
Мы можем предположить, что гауссова кривизна поверхности положительна; в противном случае утверждение очевидно.
По \ref{thm:convex-embedded} поверхность ограничивает выпуклую область, содержащую отрезок длины~$\pi$.

\begin{Figure}
\vskip-0mm
\centering
\includegraphics{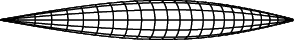}
\vskip-0mm
\end{Figure}

По \ref{ex:rev(sin)} гауссова кривизна поверхности вращения графика $y=a\cdot \sin x$ для $x\in(0,\pi)$ не превосходит $1$.
Попробуйте подпереть поверхность $\Sigma$ изнутри поверхностью вращения указанного типа.
(Сравните с \ref{ex:convex small}.)

\parit{Замечание.}
Упражнение можно вывести из следующего более глубокого результата: \textit{если гауссова кривизна замкнутой поверхности $\Sigma$ не меньше $1$,
то внутренний диаметр $\Sigma$ не больше $\pi$}.
Это означает, что любые две точки на $\Sigma$ можно соединить путём на $\Sigma$, длины не более~$\pi$.
Эту теорему доказали Хайнц Хопф и Вилли Ринов \cite{hopf-rinow} и она названа в честь, обобщившего её, Самнера Майерса \cite{myers}.

\parbf{\ref{ex:convex-proper-sphere}}; \ref{SHORT.ex:convex-proper-sphere:single}.
Воспользуйтесь выпуклостью~$R$.

\parit{\ref{SHORT.ex:convex-proper-sphere:smooth}.}
Заметьте, что отображение $\Sigma\z\to\mathbb{S}^2$ гладкое и регулярное.
Затем, применив теорему об обратной функции, покажите, что его обратное отображение также гладко.

\parbf{\ref{ex:convex-proper-plane}}; \ref{SHORT.ex:convex-proper-plane:a}
(Рассуждение похоже на \ref{prop:convex-monotone:open}.)
Можно предположить, что начало координат лежит на~$\Sigma$.
Рассмотрим последовательность точек $x_n\in \Sigma$, таких что $|x_n|\z\to \infty$ при $n\to \infty$.
Обозначим через $\vec u_n$ единичный вектор в направлении $x_n$; то есть $\vec u_n\z=\tfrac{x_n}{|x_n|}$.

Поскольку единичная сфера компактна, найдётся подпоследовательность $x_n$, такая что $\vec u_n$ сходится к единичному вектору, скажем к $\vec u$.
Покажите, что луч из начала координат в направлении $\vec u$ можно взять за~$\ell$.

\begin{wrapfigure}{r}{29 mm}
\vskip-0mm
\centering
\includegraphics{mppics/pic-1182}
\vskip-3mm
\end{wrapfigure}

\parit{\ref{SHORT.ex:convex-proper-plane:b} $+$ \ref{SHORT.ex:convex-proper-plane:c} $+$ \ref{SHORT.ex:convex-proper-plane:d}.}
Так как $R$ выпукло, то выпукла и его проекция $\Omega$.

Убедитесь, что для любого $q\in \Sigma$ направления $\vec v_n=\tfrac{x_n-q}{|x_n-q|}$ также сходятся к $\vec u$, и воспользуйтесь этим.

Покажите, что у $\Sigma$ нет вертикальных касательных плоскостей.
Выведите отсюда, что проекция из $\Sigma$ на $(x,y)$-плоскость регулярна.
Воспользовавшись теоремой об обратной функции, покажите, что множество $\Omega$ открыто.

\parit{\ref{SHORT.ex:convex-proper-plane:e}.}
Рассуждая от противного, допустим, что для некоторой последовательности $(x_n,y_n)\z\to(x_\infty,y_\infty)\in \partial\Omega$ последовательность $f(x_n,y_n)$ остаётся ограниченной сверху.
Тогда найдётся подпоследовательность, что $f(x_n,y_n)$ сходится к конечному значению, скажем $z_\infty$, либо расходится к $-\infty$.

В первом случае покажите, что точка $(x_\infty, y_\infty,z_\infty)$ не лежит на $\Sigma$, но имеет произвольно близкие точки на~$\Sigma$.
То есть $\Sigma$ несобственная --- противоречие.

Если $f(x_n,y_n)\to -\infty$, используйте выпуклость $f$, чтобы показать, что $f(\tfrac{x_n}2 ,\tfrac{y_n}2)\z\to -\infty$.
Заметьте, что начало координат принадлежит $\Omega$;
используйте это, чтобы показать, что $(\tfrac{x_\infty}2, \tfrac{y_\infty}2)\in\Omega$,
и придите к противоречию.

\parbf{\ref{ex:open+convex=plane}.}
Согласно \ref{ex:convex-proper-plane:d}, $\Sigma$ параметризуется открытой выпуклой областью $\Omega$ плоскости.
Остаётся показать, что $\Omega$-й можно запараметризовать всю плоскость.

Можно предположить, что начало координат плоскости лежит в $\Omega$.
Убедитесь, что в этом случае граница $\Omega$ записывается в полярных координатах как $(\theta,f(\theta))$, где $f\:\mathbb{S}^1\to\mathbb{R}$ --- положительная непрерывная функция.
Тогда гомеоморфизм из $\Omega$ на плоскость можно описать в полярных координатах, изменяя только радиальную координату;
например, как 
$(\theta,r)\z\mapsto (\theta,
\tfrac{r}{1-r/f(\theta)})$.

Во второй части воспользуйтесь \ref{ex:star-shaped-disc:nonsmooth}.

\parbf{\ref{ex:circular-cone}.}
Выберите координаты так, чтобы плоскость $(x,y)$ подпирала $\Sigma$ снизу в начале координат. 

Покажите, что для некоторого $\epsilon>0$ любая прямая, выходящая из начала координат под наклоном не больше $\epsilon$, может пересекать $\Sigma$ только в единичном шаре с центром в начале координат;
можно предположить, что $\epsilon$ мало, скажем $\epsilon<1$.
Рассмотрите конус, образованный лучами из начала координат под наклоном $\epsilon$ и сдвинутый вниз на $1$.
Убедитесь, что вся поверхность лежит внутри этого конуса.

\parbf{\ref{ex:intK}}.
Выберите различные точки $p,q\in\Sigma$.
Применив \ref{thm:convex-embedded}, покажите, что угол $\measuredangle(\Norm(p),p\z-q)$ острый, а $\measuredangle(\Norm(q),p\z-q)$ тупой.
Выведите отсюда, что $\Norm(p)\ne\Norm(q)$;
то есть отображение $\Norm\:\Sigma\to\mathbb{S}^2$ инъективно.

\parit{\ref{SHORT.ex:intK:4pi}.}
Для заданного единичного вектора $\vec u$ рассмотрим точку $p\in \Sigma$, которая максимизирует скалярное произведение $\langle p,\vec u\rangle$.
Покажите, что $\Norm(p)=\vec u$.
Выведите отсюда, что отображение $\Norm\:\Sigma\to\mathbb{S}^2$ сюръективно, и, следовательно, биективно.

Применив \ref{cor:intK}, получите, что наш интеграл равен $4\cdot\pi=\area\mathbb{S}^2$.

\parit{\ref{SHORT.ex:intK:2pi}.}
Выберите $(x,y,z)$-координаты, полученные в \ref{ex:convex-proper-plane:d}.
Убедитесь, что нормаль $\Norm(p)$ образует тупой угол с осью $z$ для любой точки~$p$.
Следовательно, образ $\Norm(\Sigma)$ лежит в южном полушарии.

Применив \ref{cor:intK}, получите, что наш интеграл не превосходит $2\cdot\pi=\tfrac12\cdot\area\mathbb{S}^2$.

\parit{Замечание.}
Это упражнение напоминает следствие~\ref{cor:fenchel=convex}.

\stepcounter{chapter}
\setcounter{eqtn}{0}

\parbf{\ref{ex:convex-revolution}.}
Примените \ref{ex:principal-revolution}.

\parbf{\ref{ex:ruled=>saddle}.}
Докажите, что каждая точка поверхности имеет направление с нулевой нормальной кривизной, и воспользуйтесь этим.

\parbf{\ref{ex:saddle-convex}.}
Допустим, что $p\in \Sigma$ --- точка локального максимума~$f$.
Покажите, что $\Sigma$ подпирается своей касательной плоскостью в~$p$,
и придите к противоречию.

\parbf{\ref{ex:panov}.}
Пусть $\Pi_t$ --- аффинная касательная плоскость к $\Sigma$ в точке $\gamma(t)$, а через $\ell_t$ --- касательная прямая к $\gamma$ в момент времени~$t$.

Заметим, что $\Pi_t$ --- график функции, скажем, $h_t$, определённой на $(x, y)$-плоскости.
Обозначим через $\bar\ell_t$ проекцию $\ell_t$ на $(x, y)$-плоскость.
Покажите, что производная $\tfrac{d}{dt}h_t(w)$ обращается в ноль в точке $w$ тогда и только тогда, когда $w\in \bar\ell_t$, и эта производная меняет знак, если $w$ переходит с одной стороны $\bar\ell_t$ на другую.

Если $\bar\gamma$ звёздна относительно точки $w$, то $w$ не может пересечь $\bar\ell_t$.
Поэтому функция $t\mapsto h_t(w)$ монотонна на $\mathbb{S}^1$.
Убедитесь, что эта функция непостоянна, и придите к противоречию.

\parit{Источник:}
Это результат Галины Ковалёвой \cite{kovaleva}, переоткрытый Дмитрием Пановым~\cite{panov-curves}.

\parbf{\ref{ex:crosss}.}
Это частный случай так называемой \index{лемма Морса}\emph{леммы Морса}.
Мы приведём выжимку из её доказательства.
Более концептуальное доказательство \cite{abraham-marsden-ratiu} строится на трюке Мозера; см. решение упражнения \ref{ex:star-shaped-disc}.

\medskip

Выберите касательно-нормальные координаты при $p$ с координатными осями вдоль главных направлений.
Предположим, что $\Sigma$ локально задаётся графиком $z=f(x,y)$.
Надо показать, что решение $f(x,y)=0$ есть объединение двух гладких кривых, пересекающихся в $p$.

Докажите следующее утверждение:
\textit{Пусть $x\z\mapsto h(x)$ --- гладкая функция, определённая в открытом интервале $\mathbb{I}\ni0$, такая что $h(0)=h'(0)=0$ и $h''(0)>0$.
Тогда, на меньшем интервале $\mathbb{J}\ni0$ существует единственная гладкая функция $a\:\mathbb{J}\to\mathbb{R}$ такая, что $h=a^2$, $a(0)=0$ и $a'(0)> 0$.}

Заметьте, что, сузив область до малого квадрата  $|x|,|y|<\epsilon$, можно считать, что $f_{xx}\z>\epsilon$ и $f_{yy}<-\epsilon$.
Покажите, что если $\epsilon$ мало, то для каждого $x$ существует единственное $y(x)$ такое, что $f_x(x,y(x))=0$; более того, функция $x\z\mapsto y(x)$ является гладкой.

Пусть $h(x)=f(x,y(x))$.
Убедитесь, что $h(0)\z=h'(0)=0$ и $h''>0$.
Применив утверждение выше, получите такую функцию $a$, что $h=a^2$, где $a(0)=0$ и $a'(0)>0$.

Убедитесь, что $g(x,y)=h(x)-f(x,y)\z\ge 0$, $g_y(x,y(x))\z=g(x,y(x))\z=0$, и $g_{yy}>0$.
Применяя утверждение к каждой функции $y\mapsto g(x,y)$ с фиксированным $x$, мы получаем, что $g(x,y)\z=b(x,y)^2$ для некоторой гладкой функции $b$, такой что
$b(x,y(x))=0$ и $b_y(x,y(x))>0$.

Отсюда следует, что 
\begin{align*}
f(x,y)&=a(x)^2-b(x,y)^2=
\\
&=
(a(x)-b(x,y))\cdot (a(x)+b(x,y)).
\end{align*}
То есть $f(x,y)=0$ значит, что $a(x)=\pm b(x,y)$.

Остаётся проверить, что обе функции $g_\pm(x,y)=a(x)\pm b(x,y)$ имеют различные ненулевые градиенты в нуле, и, значит, каждое уравнение $a(x)\pm b(x,y) =0$ определяет гладкую кривую в окрестности $p$;
см.~\ref{sec:implicit-curves}.

\parbf{\ref{ex:proper-saddle}.}
Посмотрите на псевдосферу в \ref{ex:principal-revolution:pseudosphere}.

\parbf{\ref{ex:length-of-bry}.}
Воспользуйтесь \ref{lem:convex-saddle} и леммой о полусфере (\ref{lem:hemisphere}).

Во второй части постройте тонкую трубку, ограниченную двумя замкнутыми сферическими кривыми.

\parbf{\ref{ex:circular-cone-saddle}.}
Пусть $\Sigma$ --- открытая седловая поверхность, лежащая в конусе $K$.
Покажите, что существует плоскость $\Pi$, которая пересекает $\Sigma$ и отсекает от $K$ компактную область.
Выведите отсюда, что $\Pi$ также отсекает от $\Sigma$ компактную область.

Воспользовавшись \ref{lem:reg-section}, сдвиньте немного плоскость $\Pi$ так, чтобы она отсекала от $\Sigma$ компактную поверхность с краем, и примените \ref{lem:convex-saddle}.

\parbf{\ref{ex:disc-hat}.}
Рассмотрим радиальную проекцию $F_\epsilon$ на сферу $\Sigma$ с центром в точке $p=(0,0,\epsilon)$;
то есть точка $q\in F_\epsilon$ переходит в точку $s(q)$ на сфере, лежащую на луче $pq$.

Покажите, что $s$ задаёт диффеоморфизм из $F_\epsilon$ на южную полусферу в~$\Sigma$.
Остаётся заметить, что единичный диск диффеоморфен полусфере.

\parbf{\ref{ex:saddle-linear}.} 
По основной теореме аффинной геометрии, аффинное преобразование является гладким.
Остаётся применить критерий горбушек (\ref{prop:hat}).

\parbf{\ref{ex:between-parallels}.}
Поищите пример среди поверхностей вращения;
воспользуйтесь \ref{ex:principal-revolution}.

\parbf{\ref{ex:one-side-bernshtein}.}
Рассмотрите сечения графика плоскостями, параллельными плоскостям $(x,y)$ и $(x,z)$, и примените теорему Мёнье (\ref{thm:meusnier}).

\parbf{\ref{ex:saddle-graph}.}
Допустим, что ортогональная проекция $\Sigma$ на $(x,y)$-плоскость не инъективна.
Покажите, что существует точка $p\in\Sigma$ с вертикальной касательной плоскостью;
то есть $\T_p\Sigma$ перпендикулярна плоскости $(x,y)$.

Пусть $\Gamma_p$ --- связная компонента точки $p$ в пересечении $\Sigma$ с аффинной касательной плоскостью $\T_p\Sigma$.
Воспользовавшись \ref{ex:crosss}, покажите, что $\Gamma_p$ представляет собой объединение гладких кривых, которые могут пересекаться друг с другом трансверсально.
Более того, две из этих кривых проходят через $p$, и $\Gamma_p$ не ограничивает компактную область на~$\Sigma$.

Покажите, что у множества $\Gamma_p$ должно быть как минимум $4$ направления, в которых оно уходит на бесконечность.
С другой стороны, $\Sigma$ является графиком вне компактного множества $K$.
Значит, $\Gamma_p\setminus K$ --- график вещественной функции одной переменной, который имеет только два направления, в которых он уходит на бесконечность --- противоречие.


\stepcounter{chapter}
\setcounter{eqtn}{0}

\parbf{\ref{ex:lasso}.}
Разрежьте боковую поверхность горы по отрезку от ковбоя до вершины и
разверните её на плоскости (см. рисунок).
Подумайте, как может выглядеть образ натянутого лассо.

\begin{wrapfigure}{r}{24 mm}
\vskip-4mm
\centering
\includegraphics{mppics/pic-1250}
\vskip-4mm
\end{wrapfigure}

\parit{Источник:}
Эту задачу мы узнали от Джоэля Файна, который приписывает её Фредерику Буржуа \cite{fine};
см. также \cite[Problem 12]{khesin-tabachnikov}.

\parbf{\ref{ex:length-dist-conv}.} 
Согласно \ref{thm:convex-embedded}, $\Sigma$ ограничивает строго выпуклую область.
Поэтому можно предположить, что $\Norm(p)\ne\Norm(q)$; в противном случае, $p=q$, и неравенство очевидно.

Ещё можно предположить, что $\Norm(p)\z+\Norm(q)\z\ne 0$; в противном случае, правая часть неопределена.

В оставшемся случае касательные плоскости $\T_p$ и $\T_q$ пересекаются по прямой, скажем, $\ell$.
Положим $\alpha=\tfrac12\cdot\measuredangle(\Norm(p),\Norm(q))$.
Убедитесь, что $2\cdot\cos\alpha\z= |\Norm(p)+\Norm(q)|$.
Пусть точка $x\in \ell$ минимизирует сумму $|p-x|\z+|x-q|$.
Покажите, что $\measuredangle\hinge xpq\ge \pi-2\cdot\alpha$.
Выведите отсюда, что 
\[
|p-x|+|x-q|\le \tfrac{|p-q|}{\cos\alpha}.
\]
Применив \ref{thm:shorts+convex}, покажите, что
$\dist{p}{q}\Sigma\z\le |p-x|\z+|x-q|$.

\parbf{\ref{ex:hat-convex}.}
Допустим, что нашлась кратчайшая $\gamma\z{\not\subset}\Delta$ с концами $p$ и $q$ в $\Delta$.

Можно предположить, что только одна дуга $\gamma$ лежит вне $\Delta$.
Отражение этой дуги относительно $\Pi$ вкупе с оставшейся частью $\gamma$ образует другую кривую $\hat\gamma$ из $p$ в $q$;
часть её идёт по $\Sigma$, 
а часть снаружи $\Sigma$,
но она никогда не заходит внутрь.
Заметьте, что
\[\length\hat\gamma=\length\gamma.\]

Обозначим через $\bar\gamma$ ближайшую проекцию точки $\hat\gamma$ на~$\Sigma$.
Кривая $\bar\gamma$ лежит в $\Sigma$, 
имеет те же концы, что и $\gamma$,
и по \ref{thm:shorts+convex}
\[\length\bar\gamma<\length\gamma.\]
Это означает, что $\gamma$ не была кратчайшей --- 
противоречие.

\parbf{\ref{ex:intrinsic-diameter}.} 
Докажите, что короткая проекция $\mathbb{S}^2\to\Sigma$ сюръективна.
Примените \ref{lem:nearest-point-projection} и \ref{thm:area-axioms:monotonicity}.


\stepcounter{chapter}
\setcounter{eqtn}{0}

\parbf{\ref{ex:helix-geodesic}.}
Докажите, что $\gamma''(t)$ пропорциональна $\nabla_{\gamma(t)} f$, где $f=x^2+y^2$. 
Примените \ref{ex:tangent-normal}.

\parbf{\ref{ex:clairaut}.} 
Мы можем предположить, что начало координат лежит на оси вращения, а вектор $\vec i$ направлен вдоль этой оси.
Воспользовавшись \ref{lem:constant-speed}, покажите, что достаточно доказать, что 
значение $\langle\gamma'\times \gamma,\vec i\rangle$ постоянно.

Поскольку $\gamma''(t)\perp\T_{\gamma(t)}$, три вектора $\vec i$, $\gamma$, и $\gamma''$ лежат в одной плоскости.
В частности, $\langle\gamma''\times \gamma,\vec i\rangle=0$.
Следовательно,
\[
\langle\gamma'\times \gamma,\vec i\rangle'
=
\langle\gamma'\times \gamma',\vec i\rangle+\langle\gamma''\times \gamma,\vec i\rangle =0
.\]

\parbf{\ref{ex:asymptotic-geodesic}.}
По лемме \ref{lem:constant-speed}, можно считать, что $\gamma$ параметризована длиной.
По определению геодезической, $\gamma''(s)\perp\T_{\gamma(s)}$. 
Следовательно, 
\[\gamma''(s)=k_n(s)\cdot\Norm(\gamma(s)),\]
где $k_n(s)$ --- нормальная кривизна кривой $\gamma$ в момент времени $s$.
Поскольку $\gamma$ --- асимптотическая линия, $k_n(s)\equiv 0$;
то есть $\gamma''(s)\equiv 0$.
Следовательно, $\gamma'$ постоянна, и, значит, $\gamma$ идёт вдоль прямой; см. \ref{ex:zero-curvature-curve}.

\parbf{\ref{ex:reflection-geodesic}.}
Обозначим через $\mu$ единичный вектор, перпендикулярный плоскости симметрии $\Pi$.
Раз уж $\gamma$ лежит в $\Pi$, её ускорение $\gamma''$ направлено вдоль $\Pi$, иначе говоря, $\gamma''\perp \mu$.
Из \ref{prop:a'-pertp-a''}, $\gamma''\perp\gamma'$, ведь $\gamma$ параметризована длиной.

Поскольку $\Sigma$ симметрична относительно плоскости $\Pi$,
касательная плоскость $\T_{\gamma(t)}\Sigma$ также симметрична относительно $\Pi$.
Отсюда вытекает, что $\T_{\gamma(t)}\Sigma$ натянута на $\mu$ и $\gamma'(t)$.
Следовательно, $\gamma''\perp \mu$ и $\gamma''\perp\gamma'$ подразумевают $\gamma''\perp\T_{\gamma(t)}\Sigma$;
то есть $\gamma$ является геодезической.

\parbf{\ref{ex:round-torus}.}
Согласно \ref{ex:reflection-geodesic}, любой меридиан $\Sigma$ является замкнутой геодезической.
Рассмотрим произвольную замкнутую геодезическую~$\gamma$.

Если $\gamma$ касается меридиана в какой-то точке, то по единственности в \ref{prop:geod-existence}, кривая $\gamma$ идёт вдоль этого меридиана;
в частности, она не стягиваема.

В оставшемся случае $\gamma$ может пересекать меридианы только трансверсально.
Следовательно, долгота $\gamma$ монотонна, и снова $\gamma$ не стягиваема.

\parbf{\ref{ex:helix=geodesic}.}
Проверьте, что $\Norm(\gamma(t))\z=(\cos t,\sin t,0)$.
Далее, вычислите ускорение $\gamma''(t)$ и убедитесь, что оно пропорционально $\Norm(\gamma(t))$,
а значит, $\gamma''(t)\z\perp\T_{\gamma(t)}$. 

Убедитесь, что отрезок от $\gamma (0)$ до $\gamma (2{\cdot}\pi)$ содержится в~$\Sigma$,
и воспользуйтесь этим.

\parbf{\ref{ex:two-min-geod}.}
Допустим, что две кратчайшие $\alpha$ и $\beta$ имеют пару общих точек $p$ и~$q$.
Обозначим через $\alpha_1$ и $\beta_1$ дуги $\alpha$ и $\beta$ от $p$ до~$q$.
Предположим, что $\alpha_1$ не совпадает с $\beta_1$.

Заметьте, что $\alpha_1$ и $\beta_1$ --- кратчайшие с теми же концами;
в частности, они одинаковой длины.
Замена $\alpha_1$ в $\alpha$ на $\beta_1$ даёт кратчайшую, скажем, $\hat\alpha$, отличную от $\alpha$.
По \ref{prop:gamma''}, $\hat\alpha$ --- геодезическая.

\begin{Figure}
\vskip-0mm
\centering
\includegraphics{mppics/pic-308}
\vskip0mm
\end{Figure}

Предположим, что $\alpha_1$ --- собственная поддуга $\alpha$;
то есть $\alpha_1\ne\alpha$, иначе говоря, $p$ или $q$ не конец кривой $\alpha$.
Тогда у $\alpha$ и $\hat\alpha$ есть общая точка и общий вектор скорости в этой точке.
По \ref{prop:geod-existence}, $\alpha$ совпадает с $\hat\alpha$ --- противоречие.

Следовательно, $p$ и $q$ --- концы~$\alpha$.
Аналогично, $p$ и $q$ --- концы~$\beta$.

Для второй части можно рассмотреть две различные геодезические вида 
\[ \gamma_b(t) = ( \cos t , \sin t , b\cdot t ) , t \in \mathbb{R} \]
на цилиндрической поверхности $x^2 + y^2 =1$.

\parbf{\ref{ex:min-geod+plane}.}
Допустим, что кратчайшая $\alpha$ проходит сквозь $\Pi$ два раза или больше.
В этом случае найдётся дуга $\alpha_1$ в $\alpha$, которая лежит строго на одной стороне от $\Pi$, только её концы лежат на $\Pi$, и эти концы отличаются от концов $\alpha$.

\begin{Figure}
\vskip-1mm
\centering
\includegraphics{mppics/pic-310}
\vskip-1mm
\end{Figure}

Удалим из $\alpha$ дугу $\alpha_1$ и заменим её отражением $\alpha_1$ в $\Pi$.
Полученная кривая, назовём её $\beta$, лежит на поверхности;
она имеет ту же длину, что и $\alpha$, и соединяет ту же пару точек, скажем, $p$ и~$q$.
Следовательно, $\beta$ --- это кратчайшая от $p$ до $q$ на $\Sigma$, отличная от~$\alpha$.

Мы можем предположить, что и $\alpha$, и $\beta$ параметризованы длиной.
По \ref{prop:gamma''}, $\alpha$ и $\beta$ --- геодезические.
Поскольку $\alpha$ и $\beta$ имеют общую дугу, у них есть общая точка и общий вектор скорости в этой точке;
по \ref{prop:geod-existence} $\alpha$ совпадает с $\beta$ --- противоречие.

\parbf{\ref{ex:milka}.}
Пусть $W$ --- замкнутая область снаружи~$\Sigma$.
Убедитесь, что $\dist{p^s}{q}W$ не зависит от~$s$.
Иначе говоря, произведение отрезка $[p^s,\gamma(s)]$ и дуги $\gamma|_{[s,\ell]}$ есть кратчайшая от $p^s$ до $q$ в~$W$.

Так как $\Sigma$ строго выпукла, 
\[ \dist{\gamma (s)}{q}W > \dist{\gamma(s)}{q}{\mathbb{R}^3} \]
для всех $s < \ell$.
Следовательно,
\begin{align*}
\dist{p^s}{q}W&=\dist{p^s}{\gamma(s)}W+\dist{\gamma (s)}{q}W 
> 
\\
&>\dist{p^s}{\gamma(s)}{\mathbb{R}^3} + \dist{ \gamma (s)}{q}{\mathbb{R}^3} 
\ge
\\
&\ge\dist{p^s}{q}{\mathbb{R}^3}. 
\end{align*}
Однако, если отрезок $[p^s,q]$ целиком содержится в $W$, то $\dist{p^s}{q}W = \dist{p^s}{q}{\mathbb{R}^3}$, и это приводит к противоречию.

\parbf{\ref{ex:round-sphere}.}
Рассуждая как в \ref{ex:const-dist}, докажите, что все геодезические на поверхности  имеют постоянную кривизну.
Выведите отсюда, все нормальные кривизны поверхности равны между собой, и рассуждайте как в \ref{ex:normal-curvature=const}.

\parit{Источник:} Задача Али Тагави \cite{taghavi}.

\parbf{\ref{ex:rad=2}}; \ref{SHORT.ex:rad=2:a}.
Пусть\( \gamma \) --- геодезическая на \( \Sigma \).
Убедитесь, что её кривизна не превышает \( 1 \).
Рассуждая как в \ref{ex:gromov-twist}, покажите, что существует \( \epsilon > 0 \), такое что
\[\epsilon
<
|\gamma(t)|
<
2\cdot\sin(\measuredangle(\gamma'(t),\gamma(t)))\]
для любого \( t \).
Выведите отсюда, что найдётся такое поле нормалей \(\Norm\) на \(\Sigma\), что
\[\epsilon<|x|<2\cdot\cos(\measuredangle(\Norm(x),x)) \eqlbl{eq:|gamma(t)|}\]
для любой точки \( x \in \Sigma \).

Воспользовавшись \ref{eq:|gamma(t)|}, покажите, что отображение \( \Sigma \to \Theta \), заданное как \( x\mapsto 2\cdot\tfrac x{|x|} \), локально обратимо.
Сделайте вывод, что оно обратимо.
Обратное к нему отображение и есть радиальная проекция \( \rho\:\Theta\to\Sigma \).

\parit{\ref{SHORT.ex:rad=2:b}.}
Предположим, что \( x=\rho(p) \) и \( \vec v \in \T_p\Theta \).
Покажите, что
\[2\cdot\cos(\measuredangle(\Norm(x),x))\cdot|d_p\rho(\vec v)|
\le
|x|\cdot|\vec v|.\]
Применив \ref{eq:|gamma(t)|}, убедитесь, что \( |d_p\rho(\vec v)| \le|\vec v| \), и заключите, что \( \rho \) не увеличивает длину.

\parit{\ref{SHORT.ex:rad=2:c}.}
Предположим, что \( y\in \Sigma \) минимизирует расстояние до \( \tfrac x2 \).
Пусть \( x=\rho(p) \) и \( y=\rho(q) \).

Рассуждая от противного, предположим, что
\[\dist{\tfrac x2}{y}{}\z<|\tfrac x2|.\]
Убедитесь, что \(\Norm\: \Sigma\to \mathbb{S}^2 \) неувеличивает длины кривых.
Используя всё это и \ref{SHORT.ex:rad=2:b}, выведите следующее:
\begin{align*}
\measuredangle(\Norm(x),\Norm(y))&= \measuredangle(x,y-\tfrac x2)>
\\
&>2\cdot\measuredangle\hinge 0pq=
\\
&=\dist{p}{q}{\Theta}\ge
\\
&\ge\dist{x}{y}{\Sigma}\ge
\\
&\ge\measuredangle(\Norm(x),\Norm(y));
\end{align*}
что даёт противоречие.

Наконец, заметим, что поверхность \( \Sigma \) можно сдвинуть так, чтобы точка \( x \) стала сколь угодно близка к \( \Theta \).
Применив рассуждение выше, получим, что наша область содержит единичный шар, касающийся \(\Sigma\) в точке \(x\).

\parit{Источник:} Задача из статьи Хонгды Цю \cite{qiu2025}.

\parbf{\ref{ex:closed-liberman}.}
Пусть $t\mapsto \gamma(t)=(x(t),y(t),z(t))$ --- замкнутая геодезическая с единичной скоростью на графике $z=f(x,y)$.
По лемме Либермана, функция $t\mapsto z'(t)$ монотонна.
Выведите отсюда, что координата $z$ постоянна на $\gamma$, то есть $\gamma$ лежит в горизонтальной плоскости.
Докажите, что $\gamma$ --- прямая, и придите к противоречию.

\parbf{\ref{ex:rho''}.}
Снабдите $\Sigma$ полем нормалей $\Norm$, направленным внутрь.
Обозначим через $k(t)$ нормальную кривизну $\Sigma$ в точке $\gamma(t)$ в направлении $\gamma'(t)$.
Поскольку $\Sigma$ выпукла, $k(t)\ge 0$ для любого~$t$.

Поскольку $\gamma$ --- геодезическая с единичной скоростью, $\gamma''(t)=k(t)\cdot\Norm(\gamma(t))$ и $\langle\gamma'(t),\gamma'(t)\rangle\z=1$ для любого~$t$.
Можно предположить, что $p$ --- начало координат в~$\mathbb{R}^3$.
Поскольку $p$ находится внутри $\Sigma$, получаем, что $\langle\gamma(t),\Norm(\gamma(t))\rangle\le 0$ и
\[
\langle\gamma''(t),\gamma(t)\rangle=k(t)\cdot \langle\gamma(t),\Norm(\gamma(t))\rangle\le 0
\]
для любого~$t$.
Выведите отсюда, что
\begin{align*}
&\rho''(t)
=\langle\gamma(t),\gamma(t)\rangle''=
\\
&=2\cdot\langle\gamma''(t),\gamma(t)\rangle+2\cdot\langle\gamma'(t),\gamma'(t)\rangle\le 2.
\end{align*}

\parbf{\ref{ex:tc-spherical-image}.}  
Можно считать, что $\gamma$ параметризована длиной.
Определим $\Norm (t) \df \Norm ( \gamma (t))$.
Тогда $\langle \gamma'(t), \Norm(t) \rangle \equiv 0$,
и поскольку $\gamma$ геодезическая, $\gamma''(t) \parallel \Norm (t)$.
Отсюда
\[
\begin{aligned}
|\gamma''(t)|
&=|\langle\gamma''(t),\Norm(t)\rangle|=
\\
&=|\langle\gamma'(t),\Norm(t)\rangle'-\langle\gamma'(t),\Norm'(t)\rangle|=
\\
&=|\langle\gamma'(t),\Norm'(t)\rangle|\le
\\
&\le|\gamma'(t)|\cdot|\Norm'(t)|=|\Norm'(t)|.
\end{aligned}
\]
Проинтегрируйте по $t$ и получите 
\[\tc\gamma\le\length(\Norm\circ\gamma).\]

\parbf{\ref{ex:usov-exact}.} 
Пусть $\gamma(t)=(x(t),y(t),z(t))$. 
Покажите, что
\[
|\gamma'' (t)| = z''(t)\cdot\sqrt{1+ \ell ^2}\eqlbl{eq:gamma''=z''}
\]
для любого~$t$.

Убедитесь, что $z'(t)\to\pm \tfrac\ell{\sqrt{1+ \ell ^2}}$ при $t\z\to\pm\infty$.
Выведите отсюда, что 
\[
\int_{-\infty}^{+\infty}z''(t)\cdot dt
=
\frac{2\cdot\ell}{\sqrt{1+ \ell ^2}}.\eqlbl{eq:int z''}
\]
Из \ref{eq:gamma''=z''} и \ref{eq:int z''},
\begin{align*}
\tc\gamma&=\int_{-\infty}^{+\infty}|\gamma''(t)|\cdot dt=
\\
&=
\sqrt{1+ \ell ^2}\cdot \int_{-\infty}^{+\infty}z''(t)\cdot dt=
\\
&=2\cdot \ell.
\end{align*}

\parbf{\ref{ex:ruf-bound-mountain}.}
Воспользуйтесь \ref{thm:usov} и \ref{ex:sef-intersection}.

Во второй части рассмотрите геодезическую на конусе с угловым коэффициентом $2$, и сгладьте его вершину.

\parit{Замечание.}
Утверждение остаётся верным для $\sqrt{3}$-липшицевых функций, и $\sqrt{3}=\tg\tfrac\pi3$ --- оптимальная константа; см. \ref{ex:sqrt(3)}.
Это тот же угловой коэффициент, что и в задаче про ковбоя и лассо (\ref{ex:lasso}).

\parbf{\ref{ex:bound-tc}}; \ref{SHORT.ex:bound-tc:a}.
Убедитесь, что 
\[|\gamma|\le 1
\qquad
\text{и}
\qquad
\langle\Norm,\epsilon\cdot \Norm-\gamma\rangle\le0,\]
и воспользуйтесь этим.

\parit{\ref{SHORT.ex:bound-tc:b}.}
Покажите, что  
\[|\gamma|\le1,\quad |\gamma'|=1
\quad
\text{и}
\quad \rho'=2\cdot \langle\gamma,\gamma'\rangle,\]
и воспользуйтесь этим.

\parit{\ref{SHORT.ex:bound-tc:c}.}
Покажите, что
\begin{align*}
|\gamma'|&=1,
&
\langle\gamma,\Norm\rangle&=|\gamma|\cdot\cos\theta,
\\
\gamma''&=-k\cdot \Norm,
&
\rho''&=2\cdot \langle\gamma',\gamma'\rangle+2\cdot \langle\gamma,\gamma''\rangle,
\end{align*}
и воспользуйтесь этим.

\parit{\ref{SHORT.ex:bound-tc:d}.}
Пусть $x,y$ --- точки на единичной сфере, проектирующиеся в концы $\gamma$. 
Соедините $x$ с $y$ дугой длины $\le \pi$ на сфере, спроектируйте её на $\Sigma$.
Примените \ref{lem:nearest-point-projection} и воспользуйтесь тем, что $\gamma$ --- кратчайшая.

\parit{\ref{SHORT.ex:bound-tc:e}.}
Пусть $\gamma\:[0,\ell]\to\Sigma$ --- параметризация длиной.
Воспользовавшись \ref{SHORT.ex:bound-tc:d}, покажите, что $\ell\z\le \pi$.
Воспользовавшись \ref{SHORT.ex:bound-tc:a}, \ref{SHORT.ex:bound-tc:b}, \ref{SHORT.ex:bound-tc:c} и тем, что $|\gamma|\ge \epsilon$, выведите, что 
\[2\cdot \ell-2\cdot \epsilon^2\cdot \tc{\gamma}
\ge
\int_0^\ell\rho''(t)\cdot dt.\]
Из \ref{SHORT.ex:bound-tc:b} выведите, что $\int\rho''\ge -4$.
Сделайте последний шаг.

 
\stepcounter{chapter}
\setcounter{eqtn}{0}

\parbf{\ref{ex:parallel}}; \ref{SHORT.ex:parallel:a}.
Покажите, что $\langle\vec v(t),\vec v'(t)\rangle=0$, и воспользуйтесь этим.

\parit{\ref{SHORT.ex:parallel:b}}
Покажите, что $|\vec v(t)|$, $|\vec w(t)|$, и
$\langle\vec v(t),\vec w(t)\rangle$
не меняются; это можно сделать так же, как в \ref{SHORT.ex:parallel:a}.
Затем воспользуйтесь тем, что 
$\langle\vec v(t),\vec w(t)\rangle\z=|\vec v(t)|\cdot|\vec w(t)|\cdot\cos\theta$.

\parbf{\ref{ex:parallel-transport-support}.}
Заметьте, что $\Sigma_1$ подпирает $\Sigma_2$ в каждой точке кривой~$\gamma$.
Выведите отсюда, что у $\gamma$ как кривой на $\Sigma_1$ и на $\Sigma_2$ идентичные нормали.
Примените \ref{obs:parallel=}.

\parbf{\ref{ex:holonomy=not0}.}
Почти любая замкнутая кривая решает задачу.
Например, рассмотрите прямоугольный сферический треугольник, который вырезается октантом $\mathbb{R}^3$ из сферы.
Покажите, что параллельный перенос вокруг него поворачивает касательную плоскость на угол $\tfrac\pi 2$.

\stepcounter{chapter}
\setcounter{eqtn}{0}

\parbf{\ref{ex:1=geodesic-curvature}.}
По \ref{ex:convex-proper-sphere}, $\Sigma$ --- гладкая сфера.
По теореме Жордана (\ref{thm:jordan}), кривая $\gamma$ делит $\Sigma$ на два диска.
Обозначим через $\Delta$ диск, который лежит слева от~$\gamma$.

Убедитесь, что $\tgc\gamma=\length\gamma$, и примените формулу Гаусса --- Бонне (\ref{thm:gb}) к~$\Delta$.

\parbf{\ref{ex:GB-hat}.}
Покажите, что $\skur=\cos\alpha\cdot\skur_0$,
где $\skur$ --- геодезическая кривизна $\partial \Delta$ в $\Delta$,
а $\skur_0$ --- геодезическая кривизна $\partial \Delta$ в $(x,y)$-плоскости.
Примените \ref{prop:total-signed-curvature} и формулу Гаусса --- Бонне (\ref{thm:gb}).

\parit{Замечание.}
Согласно \ref{ex:shape-curvature-line}, $\partial \Delta$ --- это линия кривизны в графике.

\parbf{\ref{ex:geodesic-half}.}
Воспользуйтесь \ref{cor:intK}, \ref{ex:intK:4pi} и формулой Гаусса --- Бонне (\ref{thm:gb}).

В последней части примените \ref{ex:bisection-of-S2}.

\parbf{\ref{ex:closed-geodesic}.}
В первой части примените формулу Гаусса --- Бонне (\ref{thm:gb}).

Во второй части, если бы две замкнутые геодезические $\gamma_1$ и $\gamma_2$ не пересекались, то 
$\gamma_2$ лежала бы в одной из областей, скажем $R_1$, которую $\gamma_1$ вырезает из~$\Sigma$.
Аналогично, $\gamma_1$ лежала бы в одной из областей, скажем $R_2$, которую $\gamma_2$ вырезает из~$\Sigma$.

Заметьте, что $R_1$ и $R_2$ покрывают $\Sigma$ с перекрытием.
Из первой части получаем, что интеграл гауссовой кривизны по $\Sigma$ строго меньше $4\cdot\pi$, а это противоречит \ref{ex:intK:4pi}.

\parbf{\ref{ex:self-intersections}}; \textit{(легкая)}.
Рассмотрите $4$ области, ограниченные петлями.
Применив формулу Гаусса --- Бонне (\ref{thm:gb}), покажите, что интеграл гауссовой кривизны по каждой из этих областей превышает $\pi$.
Придите к противоречию с \ref{ex:intK:4pi}.

\parit{(сложная)}.
Пусть $\alpha$, $\beta$ и $\gamma$ --- углы треугольника.
Применив формулу Гаусса --- Бонне (\ref{thm:gb}), покажите, что три петли окружают области, на которых интеграл гауссовой кривизны равен $\pi+\alpha$, $\pi+\beta$ и $\pi+\gamma$.

Применив формулу Гаусса --- Бонне к треугольной области, убедитесь, что $\alpha+\beta\z+\gamma>\pi$.
Воспользуйтесь \ref{ex:intK:4pi} и придите к противоречию.

\parit{(безнадёжная)}.
Это правда сложно.
Решение на основе \ref{thm:comp:toponogov} приведено в \cite{petrunin2021}.

\parbf{\ref{ex:sqrt(3)}.}
Достаточно показать, что на поверхности нет геодезических петель.
Допустим нашлась петля, оцените интегралы гауссовой кривизны по всей поверхности и по диску, окружённому геодезической петлёй.

\parbf{\ref{ex:unique-geod}}; \ref{SHORT.ex:unique-geod:unique}.
Из \ref{prop:shortest-paths-exist} и \ref{prop:gamma''} вытекает, что любые две точки на $\Sigma$ можно соединить геодезической.
Допустим, что точки $p$ и $q$ можно соединить двумя различными геодезическими $\gamma_1$ и $\gamma_2$.
Перейдя к их дугам, можно считать, что у $\gamma_1$ и $\gamma_2$ нет общих точек кроме концов.
Поскольку поверхность односвязна, $\gamma_1$ и $\gamma_2$ вместе ограничивают диск, скажем, $\Delta\subset\Sigma$.
Остаётся применить формулу Гаусса --- Бонне к $\Delta$ и сделать вывод.

\parit{\ref{SHORT.ex:unique-geod:diffeomorphism}.} 
Воспользуйтесь частью \ref{SHORT.ex:unique-geod:unique} и \ref{prop:inj-rad}.

\parbf{\ref{ex:half-sphere-total-curvature}.}
Примените \ref{prop:pt+tgc} и \ref{prop:area-of-spher-polygon}.

\parbf{\ref{ex:cohn-vossen}.}
Рассуждайте как в конце доказательства \ref{thm:cohn-vossen} для однопараметрического семейства геодезических $\gamma_\tau$, заданных условием $\gamma_\tau(0)=\alpha(\tau)$ и $\gamma'_\tau(0)=\alpha'(\tau)$.

\parbf{\ref{ex:3-curves}.}
Первая кривая соответсвует геодезической на сглаженном конусе.
Убедитесь, что вторая и третья запрещаются рассуждением в доказательстве теоремы.

\parbf{\ref{ex:g-b-chi}.}
Разбейте поверхность на диски, подсчитайте число рёбер, вершин и дисков в полученном разбиении.
Вычислите эйлерову характеристику и примените \ref{thm:GB-generalized}.
Для ленты Мёбиуса и пары штанов, воспользуйтесь тем, что граничные кривые имеют нулевую геодезическую кривизну.
Для цилиндра, воспользуйтесь тем, что граничные кривые плоски и замкнуты и примените \ref{prop:total-signed-curvature}.


\stepcounter{chapter}
\setcounter{eqtn}{0}

\parbf{\ref{ex:semigeodesc-chart}.}
По лемме Гаусса (\ref{lem:palar-perp}), полярные координаты с началом в точке $q$ задают полугеодезические карты в близлежащих точках.
Таким образом, достаточно найти точку $q\ne p$, чтобы полярные координаты на $\Sigma$ с началом в $q$ покрывали $p$.
Согласно \ref{prop:exp}, подойдёт любая точка $q$ достаточно близкая к $p$.

\parbf{\ref{ex:inj-rad}}; \ref{SHORT.ex:inj-rad:sign}.
Обратите внимание, что ориентацию на $\Sigma$ можно выбрать так, что $b_r(0,\theta)\z=1$ для любого $\theta$.
После этого, можно считать, что $b(r,\theta)>0$ для всех малых $r>0$.

Допустим, что $b(r_1,\theta_1)<0$ при некоторой паре $(r_1,\theta_1)$ с $0<r_1<r_0$.
Убедитесь, что если $\theta_2$ достаточно близко к $\theta_1$, то радиальные кривые $r\z\mapsto b(r,\theta_1)$ и $r\z\mapsto b(r,\theta_2)$, определённые на интервале $(0,r_1)$, пересекаются.
Следовательно, $\exp_p$ не является инъективным в $B$ --- противоречие.

\parit{\ref{SHORT.ex:inj-rad:0}.}
Допустим, что $b(r_1,\theta_1)=0$.
Примените \ref{SHORT.ex:inj-rad:sign}, чтобы показать, что $b_r(r_1,\theta_1)=0$.

Применив \ref{prop:jaccobi}, выведите, что $b(r,\theta_1)=0$ для любого $r$.
Это противоречит тому, что $b_r(0,\theta_1)=1$.

\parit{\ref{SHORT.ex:inj-rad:prop:inj-rad}.}
Нужно показать, что $\exp_p$ регулярно в $B$.
Предположим, что вектор $\vec v\in B$ имеет полярные координаты $(r,\theta)$ для некоторого $r>0$.
Покажите, что $\exp_p$ регулярно в точке $\vec v$, если $b(r,\theta)\ne 0$.
Выведите, что $\exp_p$ регулярно в $B\setminus \{0\}$.

По \ref{obs:d(exp)=1}, $\exp_p$ регулярно в нуле.
Значит сужение $\exp_p|_B$ --- регулярное инъективное гладкое отображение.
Воспользовавшись теоремой об обратной функции (\ref{thm:inverse}), покажите, что сужение $\exp_p|_B$ есть диффеоморфизм на свой образ.

\parbf{\ref{lem:K(orthogonal)}}; \ref{SHORT.lem:K(orthogonal):uu-vv}.
Так как базис $\vec u, \vec v,\Norm$ ортонормирован,
первые два векторных тождества эквивалентны следующим шести вещественным:
\[
\begin{aligned}
\langle\vec u_u,\vec u\rangle
&=0,
&
\langle\vec v_u,\vec v\rangle
&=0,
\\
\langle\vec u_u,\vec v\rangle
&=-\tfrac1{b}\cdot a_v,
&
\langle\vec v_u,\vec u\rangle
&=
\tfrac1{b}\cdot a_v,
\\
\langle\vec u_u,\Norm\rangle
&=a\cdot \ell,
&
\langle\vec v_u,\Norm\rangle
&=
a\cdot m.
\end{aligned}
\eqlbl{eq:key-orthogonal/2}
\]

Взяв частные производные от тождеств
$\langle\vec u,\vec u\rangle=1$ и
$\langle\vec v,\vec v\rangle=1$ по $u$,
мы получаем первые два тождества в \ref{eq:key-orthogonal/2}.

Кроме того, заметим, что
\[
\begin{aligned}
\vec v_u
&=
\tfrac{\partial}{\partial v}
(\tfrac1b\cdot s_v)=\tfrac1b\cdot s_{uv}
+
\tfrac{\partial}{\partial u}(\tfrac1b)
\cdot
 s_v.
\end{aligned}
\eqlbl{eq:dv/du}
\]
Так как $s_u\perp s_v$, получаем, что
\begin{align*}
\langle\vec v_u,\vec u\rangle
&=
\tfrac1{a\cdot b}\cdot \langle s_{vu}, s_u\rangle=
\\
&=
\tfrac1{2\cdot a\cdot b}\cdot \tfrac{\partial}{\partial v}\langle s_u, s_u\rangle=
\\
&=
\tfrac1{2\cdot a\cdot b}\cdot \tfrac{\partial a^2}{\partial v}=
\\
&=
\tfrac1{b}\cdot a_v.
\end{align*}
Взяв частную производную $\langle\vec u,\vec v\rangle=0$ по $u$, получаем
\begin{align*}
\langle\vec v_u,\vec u\rangle+
\langle\vec v,\vec u_u\rangle
&=0.
\end{align*}
Получаем ещё пару тождеств из \ref{eq:key-orthogonal/2}.

Так как $\vec u$ и $\vec v$ ортонормированы, из \ref{thm:shape-chart} следует, что
\[
\begin{aligned}
\tfrac1{a^2}
\cdot
\langle s_{uu},\Norm\rangle
&=\ell,
&
\tfrac1{a\cdot b}
\cdot
\langle s_{uv},\Norm\rangle
&=m,
\\
\tfrac1{a\cdot b}
\cdot
\langle s_{vu},\Norm\rangle
&=m,
&
\tfrac1{b^2}
\cdot
\langle s_{vv},\Norm\rangle
&=n.
\end{aligned}
\eqlbl{eq:shape-lmn}
\]

Применив \ref{eq:dv/du}, \ref{eq:shape-lmn} и то, что $s_v\perp\Norm$, получаем
\begin{align*}
\langle\vec u_u,\Norm\rangle
&=
\tfrac1{a}\cdot \langle s_{uu},\Norm\rangle
=a\cdot \ell,
\\
\langle\vec v_u,\Norm\rangle
&=
\tfrac1{a}\cdot \langle s_{uv},\Norm\rangle
=a\cdot m.
\end{align*}
Мы получили последние два тождества в \ref{eq:key-orthogonal/2}.

Итак, мы доказали первые два тождества в \ref{SHORT.lem:K(orthogonal):uu-vv};
остальные доказываются также.

\parit{\ref{SHORT.lem:K(orthogonal):K}.}
Напомним, что гауссова кривизна равна определителю матрицы $
(\begin{smallmatrix}
\ell&m
\\
m&n
\end{smallmatrix}
)
$;
то есть $K\z=\ell\cdot n-m^2$.
Следовательно, 
\begin{align*}
a\cdot b\cdot K
&=
a\cdot b\cdot(\ell\cdot n-m^2)
=
\\
&=
\langle\vec u_u,\vec v_v\rangle 
-
\langle\vec u_v,\vec v_u\rangle= 
\\
&= 
\left(
\tfrac{\partial}{\partial v}
\langle\vec u_u,\vec v\rangle
-
\cancel{\langle\vec u_{uv},\vec v\rangle}
\right)-
\\
&-
\left(
\tfrac{\partial}{\partial u}
\langle\vec u_v,\vec v\rangle
-
\cancel{\langle\vec u_{uv},\vec v\rangle}
\right)=
\\
&=\tfrac{\partial}{\partial v}(-\tfrac1{b}\cdot a_v)
-
\tfrac{\partial}{\partial u}(\tfrac1{a}\cdot b_u).
\end{align*}

\parbf{\ref{ex:conformal}.}
Примените \ref{lem:K(orthogonal)} при $a=b$, и упростите.

\parbf{\ref{ex:isom(geod)}.}
Убедитесь, что внутренняя изометрия переводит кратчайшие в кратчайшие, и примените \ref{prop:gamma''}.

\parbf{\ref{ex:K=0}.}
Применив \ref{prop:rauch}, убедитесь, что $\exp_p$ сохраняет длину, и воспользуйтесь этим.

\parbf{\ref{ex:K=1}.}
Переиначив доказательство \ref{prop:rauch}, получите, что $K\equiv 1$ означает, что $b(\theta,r)\z=\sin r$ для всех малых $r\ge 0$.

\parbf{\ref{ex:deformation}.}
Решив дифференциальное уравнение $x'(t)\z=\sqrt{1-|y'(t)|^2}$, можно найти $x(t)$.
Решение определено пока $y(t)>0$ и $|y'(t)|\z=|a\cdot \sin t|<1$.

Применив \ref{ex:principal-revolution:formula}, найдите гауссову кривизну поверхности $\Sigma_a$.

По \ref{ex:K=1}, малые диски $\Delta_a$ на $\Sigma_a$ внутренне изометричны диску на единичной сфере.
Чтобы показать, что $\Delta_a$ не конгруэнтен сферическому диску $\Delta_1$ при $a\ne 1$, покажите, что его главная кривизна не равна $1$ в какой-то точке.

\parbf{\ref{ex:line-cylinder}.}
Убедитесь в том, что $|f(p)-f(q)|\z\le |p-q|$ для любых двух точек $p$ и $q$ на $(u,v)$-плоскости.

Пусть $f(u,v)=(x(u,v),y(u,v),z(u,v))$;
можно предположить, что $z(0,v)=v$ для любого $v$.

Применив это наблюдение и неравенство треугольника к точкам $(0,-a)$, $(u,v_1)$, $(u,v_2)$ и $(0,a)$ при $a\to \infty$, покажите, что 
\begin{align*}
z(u,v_1)&=v_1,
&
z(u,v_2)&=v_2,
\\
x(u,v_1)&=x(u,v_2),
&
y(u,v_1)&=y(u,v_2).
\end{align*}
Сделайте вывод.


\stepcounter{chapter}
\setcounter{eqtn}{0}

\parbf{\ref{ex:wide-hinges}.}
Из неравенства треугольника $c_n\z\le a_n\z+b_n$ и условия $a_n,b_n>\epsilon$, выведите, что последовательность 
\[
s_n=\frac{a_n+b_n+c_n}{2\cdot a_n\cdot b_n}
\]
ограничена.
Убедитесь, что
\[
\frac{(a_n+b_n)^2-c_n^2}{2\cdot a_n\cdot b_n}=s_n\cdot (a_n+b_n-c_n)\to 0
\]
при $n\to\infty$.
Воспользовавшись этим и теоремой косинусов
\[
a_n^2+b_n^2-2\cdot a_n\cdot b_n\cdot\cos\tilde\theta_n -c_n^2=0,
\]
покажите, что $\cos\tilde\theta_n\to -1$, и, значит, $\tilde\theta_n\to \pi$ при $n\to\infty$.

\parbf{\ref{ex:thm:comp:cat:nsc}.}
Рассмотрите треугольник на гиперболоиде с вершинами $(1,0,0)$, $(-\tfrac12, \tfrac{\sqrt{3}}2,0)$ и $(-\tfrac12,-\tfrac{\sqrt{3}}2,0)$.
Убедитесь, что все его углы равны $\pi$, а все его модельные углы равны $\tfrac{\pi}3$.

\parbf{\ref{ex:diam-angle}.}
Убедитесь, что $\dist{p}{x}\Sigma\le \dist{p}{q}\Sigma$ и $\dist{q}{x}\Sigma\le \dist{p}{q}\Sigma$.
Выведите, что
\[\modangle xpq\ge \tfrac\pi3,\]
и примените \ref{thm:comp:toponogov}.

\parbf{\ref{ex:sum=<2pi}.}
Покажите, что 
$\measuredangle\hinge pxy+\measuredangle\hinge pyz+\measuredangle\hinge pzx\z\le2\cdot \pi$,
и примените \ref{thm:comp:toponogov}.

\parbf{\ref{ex:s-r}.}
Выберите такие векторы $\vec u$ и $\vec w$, что $|\vec w|\z=\dist{p}{x}{}$, $|\vec u|=1$, и $\measuredangle(\vec u,\vec w)=\measuredangle\hinge xpq$.
Рассмотрите функцию
$f\:t\mapsto t+|\vec w|-|\vec w-t\cdot \vec u|$.
Заметьте, что $f(0)\z=0$,
\begin{align*}
f(\dist{x}{y}{})&=\dist{x}{p}{}+\dist{x}{y}{}-\side\hinge xpy,
\\
f(\dist{x}{q}{})&=\dist{x}{p}{}+\dist{x}{q}{}-\side\hinge xpq.
\end{align*}
Покажите, что $f$ вогнута, и воспользуйтесь этим.

\parbf{\ref{ex:open-comparison}}; \ref{SHORT.ex:open-comparison:positive}.
Воспользуйтесь теоремой сравнения Рауха (\ref{prop:rauch}) и свойствами экспоненциального отображения в \ref{prop:exp}.

\parit{\ref{SHORT.ex:open-comparison:almost-min}.}
Рассуждая от противного,
предположите, что для любой точки $p\in\Sigma$ найдётся точка $q=q(p)\in \Sigma$ такая, что 
$\dist{q}{p}\Sigma<100\cdot R_p$
и
$R_q<(1-\tfrac1{100})\cdot R_p$.

Начав с любой точки $p_0$, постройте такую последовательность точек $p_n$, что $p_{n+1}\z\df q(p_n)$.
Покажите, что $p_n$ сходится, и $R_{p_n}\to 0$ при $n\to\infty$.
Придите к противоречию с \ref{SHORT.ex:open-comparison:positive}.

\parit{\ref{SHORT.ex:open-comparison:proof}.}
Повторите доказательство теоремы, предполагая, что для точки $p$ выполняется~\ref{SHORT.ex:open-comparison:almost-min}.

\parbf{\ref{ex:geod-convexity}.}
Примените \ref{angk>angk} и \ref{angk<angk}.

\parbf{\ref{ex:midpoints}.}
Используйте \ref{angk>angk} или \ref{angk<angk} дважды: 
сначала --- для треугольника $[pxy]$ и $\bar x\in [p,x]$, 
потом --- для треугольника $[p\bar xy]$ и $\bar y\in [p,y]$.
Затем примените монотонность угла (\ref{lem:angle-monotonicity}).

\parbf{\ref{ex:convex-dist}.}
Достаточно доказать неравенство Йенсена:
\begin{align*}
\dist{\gamma_1(\tfrac12)}{\gamma_2(\tfrac12)}{}
\le
\tfrac12\cdot&\dist{\gamma_1(0)}{\gamma_2(0)}{}
+
\\
+
\tfrac12\cdot&\dist{\gamma_1(1)}{\gamma_2(1)}{}.
\end{align*}

\begin{Figure}
\vskip-0mm
\centering
\includegraphics{mppics/pic-1650}
\vskip1mm
\end{Figure}

Пусть $\delta$ --- геодезический путь от $\gamma_2(0)$ до $\gamma_1(1)$.
Из \ref{ex:midpoints}
\begin{align*}
\dist{\gamma_1(\tfrac12)}{\delta(\tfrac12)}{}
&\le
\tfrac12\cdot\dist{\gamma_1(0)}{\delta(0)}{},
\\
\dist{\delta(\tfrac12)}{\gamma_2(\tfrac12)}{}
&\le
\tfrac12\cdot\dist{\delta(1)}{\gamma_2(1)}{}.
\end{align*}
Сложите их и примените неравенство треугольника.

\parit{Замечание.}
По модулю теоремы сравнения, случай евклидовой плоскости не проще.

\parbf{\ref{ex:disc+}.}
Пусть $B=\bar B[p,R]_\Sigma$.
Для данной точки $x\in B$ выберем геодезическую $[p,x]$ и обозначим через $\tilde x$ точку в $\T_p$, которая лежит на расстоянии $|p-x|_\Sigma$ от $p$ в направлении $[p,x]$.
Отображение $x\mapsto \tilde x$ отправляет $B$ в $R$-шар в $\T_p$;
по \mbox{\ref{mang>angk}}, это отображение не уменьшает расстояния,
и утверждение следует.

\parbf{\ref{ex:disc-}}.
Пусть $p\in\Sigma$ --- центр $\Delta$. 

\parit{\ref{SHORT.ex:disc-:kg}.}
Рассмотрим геодезический треугольник $[pxy]_\Sigma$ с $|p-x|_\Sigma=|p-y|_\Sigma=R$.
Пусть $\alpha=\measuredangle\hinge xpy$, 
$\beta=\measuredangle\hinge ypx$ и $\ell=|x-y|_\Sigma$.
Используя \mbox{\ref{mang<angk}}, покажите, что 
\[\pi-\alpha-\beta>\tfrac\ell R+o(\ell).\]

Пусть $\sigma$ --- дуга $\partial \Delta$ от $x$ до $y$.
Убедитесь, что $[x,y]_\Sigma$ лежит в области, ограниченной $[p,x]$, $[p,y]$ и $\sigma$.
Применив формулу Гаусса --- Бонне, покажите, что $\tgc\sigma\z>\pi-\alpha-\beta$.
Перейдите к пределу при $\ell\to0$, и сделайте вывод.

\parit{\ref{SHORT.ex:disc-:area}.}
Согласно \ref{mang<angk}, отображение $\exp_p\:\T_p\z\to\Sigma$ не уменьшает расстояния.
По \ref{prop:gamma''},  оно отображает $R$-шар в $\T_p$, с центром в нуле, на $\Delta$.
Отсюда следует утверждение.

\parbf{\ref{ex:moon-}.}
Повторите доказательство \ref{thm:moon-orginal}, используя \ref{ex:disc-:kg}.
Для последнего утверждения примените \ref{ex:disc-:area}.

\parit{Источник:}
Задача предложена Дмитрием Бураго.


\stepcounter{chapter}
\setcounter{eqtn}{0}

\parbf{\ref{ex:ell-infty}.}
Проверьте все условия определения метрики на странице \pageref{page:def:metric}.

\parbf{\ref{ex:B2inB1}}; \ref{SHORT.ex:B2inB1:a}.
Заметьте, что $\dist{p}{q}{\spc{X}}\le 1$.
Применив неравенство треугольника, покажите, что $\dist{p}{x}{\spc{X}}\le 2$ для любого $x\in B[q,1]$.
Сделайте вывод.

\parit{\ref{SHORT.ex:B2inB1:b}.}
Возьмите за $\spc{X}$ луч $[0,\infty)$ со стандартной метрикой; $p=0$ и $q=\tfrac45$.

\parbf{\ref{ex:shrt=>continuous}.}
Покажите, что условия в \ref{def:continuous} выполняются при $\delta=\epsilon$.

\parbf{\ref{ex:close-open}.}
Предположим, что дополнение $\Omega=\spc{X}\setminus Q$ открыто.
Тогда для каждой точки $p\in \Omega$ существует такое $\epsilon>0$, что $\dist{p}{q}{\spc{X}}>\epsilon$ для любого $q\in Q$.
Следовательно, $p$ \textit{не} может быть предельной точкой никакой последовательности из $Q$.
Значит, любой предел последовательности в $Q$ принадлежит $Q$,
и по определению, $Q$ замкнуто.

Теперь предположим, что $\Omega=\spc{X}\setminus Q$ не открыто.
Докажите, что существует точка $p\in \Omega$ и последовательность $q_n\in Q$ такая, что $\dist{p}{q_n}{\spc{X}}<\tfrac1n$ для любого~$n$.
Выведите отсюда, что $q_n\to p$ при $n\to \infty$, и, значит, $Q$ незамкнуто.

\spell{\end{multicols}}{}

}

\arxiv{\newgeometry{top=0.9in, bottom=0.9in,inner=0.9in, outer=0.7in}}{\newgeometry{top=1.025in, bottom=1.025in,inner=0.9in, outer=0.825in}}

\chapter{Послесловие}

Если после этой книжки вы перейдёте к изучению римановой геометрии, то половина материала в этой области окажется вам знакомой.
Но прежде придётся разобраться в тензорном исчислении; например, по книге Ричарда Бишопа и Самуэля Гольдберга \cite{bishop-goldberg}.

Далее мы приводим список введений в риманову геометрию, которые мы знаем и любим
от простых и подробных до сжатых и сложных:
\begin{itemize}
\item «Риманова геометрия»  Манфредо до Кармо \cite{carmo1992riemannian}.
\item «Введение в риманову геометрию» Юрия Бураго и Виктора Залгаллера \cite{burago-zalgaller}.
\item «Риманова геометрия в целом»  Дэтлефа Громолла, Вильгельма Клингенберга и Вернера Мейера \cite{gromoll-klingenberg-meyer,gromoll-klingenberg-meyer-ru}.
\item «Comparison geometry»  Джеффа Чигера и Дэвида Эбина \cite{cheeger-ebin}.
\item «Знак и геометрический смысл кривизны» Михаила Громова \cite{gromov-1991}.
\end{itemize}
Удачи.

\begin{flushright}
Антон Петрунин и Серхио Замора Баррера,\\
Стейт-Колледж, Пенсильвания, 11 мая 2021 года.
\end{flushright}

\arxiv{\newgeometry{top=0.83in, bottom=0.83in,inner=0.52in, outer=0.52in}}{\newgeometry{top=0.955in, bottom=0.955in,inner=0.52in, outer=0.645in}}

{
\newpage
\phantomsection
\footnotesize\sloppy
\input{curves-and-surfaces-ru.ind}

}

{

\sloppy

\def\emph{\textit}
\printbibliography[heading=bibintoc]

\fussy

}

\end{document}